\documentclass[12pt,a4paper]{article}
\usepackage[english]{babel}
\usepackage{newtxtext}
\usepackage[varvw]{newtxmath}
\usepackage[scaled=0.90]{helvet} % ss
\usepackage{courier} % tt
\normalfont
\usepackage[T1]{fontenc}
\usepackage{latexsym,amsmath}%amsthm}%,epsfig,pstricks}%,pst-node}
\usepackage{makeidx}
\usepackage{a4}
\usepackage{dsfont}
\usepackage[usenames,dvipsnames]{color}
\usepackage{comment}
\usepackage{pdfsync}
\newtheorem{Thm}{Theorem}
\newtheorem{Defi}[Thm]{Definition}
\newtheorem{Cor}[Thm]{Corollary}
\newtheorem{Lemma}[Thm]{Lemma}
\newtheorem{Prop}[Thm]{Proposition}
\newtheorem{Rem}[Thm]{Remark}
\newtheorem{Conj}[Thm]{Conjecture}
\newtheorem{Prelim}[Thm]{Preliminary}

\newenvironment{thm}[0]{\begin{Thm}\noindent}%
{\end{Thm}}
\newenvironment{defi}[0]{\begin{Defi}\noindent\rm}%
{\end{Defi}}
\newenvironment{cor}[0]{\begin{Cor}\noindent}%
{\end{Cor}}
\newenvironment{lemma}[0]{\begin{Lemma}\noindent}%
{\end{Lemma}}
\newenvironment{prop}[0]{\begin{Prop}\noindent}%
{\end{Prop}}
\newenvironment{rem}[0]{\begin{Rem}\noindent\rm}%
{\end{Rem}}
{\end{Conj}}
{\end{Prelim}}

\def\proof{\par\noindent{\it Proof.}{\ }{\ }}
\def\qed{~\hfill\hbox{$\Box$}\medbreak}

\def\naam#1{\label{#1}}

\def\medno{\medbreak\noindent}
\def\text#1{\;\;\;\;{\rm \hbox{#1}}\;\;\;\;}
\def\qquad{\quad\quad}

\def\itema{\vspace{-1mm}\item[{\rm (a)}]}
\def\itemb{\item[{\rm (b)}]}
\def\itemc{\item[{\rm (c)}]}

\def\msy#1{{\mathbb #1}}
\def\C{{\msy C}}
\def\N{{\msy N}}
\def\Z{{\msy Z}}
\def\R{{\msy R}}

\def\ga{\alpha}
\def\gb{\beta}
\def\gd{\delta}
\def\ge{\varepsilon}
\def\gf{\varphi}
\def\gg{\gamma}

\def\gl{\lambda}

\def\gs{\sigma}

\def\gz{\zeta}

\def\gD{\Delta}

\def\gS{\Sigma}

\def\got#1{\mathfrak #1}
\def\fa{{\got a}}
\def\fb{{\got b}}
\def\fg{{\got g}}
\def\fh{{\got h}}

\def\fk{{\got k}}
\def\fl{{\got l}}
\def\fm{{\got m}}
\def\fn{{\got n}}

\def\fp{{\got p}}

\def\ft{{\got t}}

\def\implies{\Rightarrow}
\def\to{\rightarrow}
\def\Re{{\rm Re}\,}
\def\Im{{\rm Im}\,}
\def\inp#1#2{\langle#1\,,\,#2\rangle}

\def\Ad{{\rm Ad}}
\def\End{{\rm End}}
\def\Hom{{\rm Hom}}

\def\ad{{\rm ad}}

\def\after{\,{\scriptstyle\circ}\,}

\def\pr{{\rm pr}}

\def\implies{\Leftarrow}

\def\tr{{\rm tr}\,}

\def\iC{{\scriptscriptstyle \C}}
\def\iR{{\scriptscriptstyle \R}}

\def\cA{{\mathcal A}}
\def\cB{{\mathcal B}}
\def\cC{{\mathcal C}}
\def\cD{{\mathcal D}}
\def\cE{{\mathcal E}}
\def\cF{{\mathcal F}}
\def\cH{{\mathcal H}}
\def\cI{{\mathcal I}}
\def\cK{{\mathcal K}}

\def\cM{{\mathcal M}}

\def\cO{{\mathcal O}}
\def\cP{{\mathcal P}}
\def\cR{{\mathcal R}}
\def\cS{{\mathcal S}}

\def\cY{{\mathcal Y}}

\def\bs{{\backslash}}

\numberwithin{Thm}{section}
\numberwithin{equation}{section}
\def\bfn{\bar\fn}
\def\bN{N}
\def\bn{n}
\def\Ind{{\rm Ind}}
\def\supp{{\rm supp}}
\def\dotvar{\, \cdot\,}
\def\Vtau{V_\tau}
\def\fZ{{\mathfrak Z}}

\def\implies{\Rightarrow}
\def\Wh{{\rm Wh}}
\def\ev{{\rm ev}}
\def\fadc{\fa_\C^*}

\def\embeds{\hookrightarrow}
\def\faPdc{{\fa_{P\iC}^*}}
\def\faPd{{\fa_P^*}}

\def\gL{\Lambda}

\def\fhdc{\fh_\iC^*}
\def\Cartan{\theta}

\numberwithin{equation}{section}

\def\fad{{\fa^*}}
\def\cl{{\rm cl}}
\def\fZ{{\mathfrak Z}}
\def\gL{\Lambda}
\def\ft{\mathfrak t}

\def\bfn{{\bar \fn}}
\def\fhdc{\fh_\iC^*}
\def\wt{{\rm wt}}
\def\bp{\!{\,}^\backprime}
\def\spec{{\rm spec}}

\def\faPdc{\fa_{P\iC}^*}
\def\faQdc{\fa_{Q\iC}^*}
\def\fh{{\mathfrak h}}
\def\Cartan{\theta}
\def\fds{{\rm fs}}
\def\fX{{\mathfrak X}}
\def\faP{\fa_P}
\def\cchidiff{\xi}
\def\cchiP{\cchi_{{\!}_P} }
\def\half{{\textstyle \frac12}}
\def\ulrmZ{\,\underline{\rm Z}_\mu}
\def\ulZ{\, \underline{\rm Z}}
\def\snorm{n}
\def\fv{{\mathfrak v}}
\def\IPgs{I_{P,\gs}^\infty}
\def\cAtwoP{\cA_{2,P}}
\def\circG{{}^\circ G}
\def\ord{{\rm ord}}
\def\Rinp#1#2{(#1,#2)}

\def\topping#1#2{\genfrac{}{}{0pt}{}{#1}{#2}}
\def\wh{{\rm wh}}
\def\col{\!:\!}

\def\jv{{}_vj}
\def\umu{\underline{\mu}}
\def\ideal{\vartriangleleft }
\def\spec{{\rm spec}}
\def\cPst{\cP_{\rm st}}
\def\MPds{\widehat M_{P,\rm ds}}
\def\cAtwoPgs{\cA_{2,P,\gs}}
\def\cchi{\chi}
\def\xx{x}
\def\xxi{\xi}
\def\gsP{\gs}%
\def\veec{\vee}%
\def\bfhyp{hyperplanes}
\def\trw{{}^{\scriptscriptstyle\veec}}
\def\hypp{{\mathfrak H}}
\def\CPT{{\rm CPT}}
\title{Uniform temperedness of Whittaker integrals\\for a real reductive group}
\author{Erik P. van den Ban}
\date{April 20, 2023}
\begin{document}
\maketitle
\hfill\mbox{\em In memory of my son Mark}
\tableofcontents
\section*{Introduction}
In 1982, Harish-Chandra announced the Whittaker Plancherel theorem for 
real reductive groups in an invited lecture at the AMS summer conference in Toronto. Because of 
his failing health, the lecture, with the title `On the theory of the Whittaker integral',  was delivered 
 on his behalf by V.S. Varadarajan.
As a consequence of Harish-Chandra's untimely death in 1983, the details of the proof remained 
unpublished until they finally appeared in the posthumous 5th volume of his collected papers \cite[pp.~141-307]{HCwhit}. That volume also contains the text of the 1982 announcement,
see \cite[\S 1.2]{HCwhit},

The proof of the Whittaker Plancherel theorem given in \cite{HCwhit} seems to be incomplete, mainly since it does not develop a complete theory of Fourier transform for the Whittaker Schwartz space. 
In particular the required uniformly tempered estimates for the Whittaker integral are not addressed.
In the present paper we give a proof of these estimates. 

Independently, N. Wallach developed a completely different approach to the Whittaker Plancherel theory 
in his book \cite{Wrrg2}. However, the treatment was flawed because of an erroneous 
estimate, pointed out in \cite[Remark 7.5]{BKcf}.  Wallach has made several attempts
to circumvent the error, see \cite{Wwhitrepair}, but the final status of his results 
seems unclear
at this point.

Clearly, the present paper has been inspired by both \cite{HCwhit} and \cite{Wrrg2}.
My desire to investigate the details of all arguments has led to 
a somewhat different and rather self-contained treatment of the theory needed for  
the derivation of a new functional equation for Whitttaker vectors in the generalized
principal series, which lies at the basis of the mentioned uniform tempered estimates.

Now that the results of the present paper are available, it is natural to develop a theory of the constant term for tempered families of Whittaker coefficients, as well as a theory of wave packets of Whittaker integrals, in analogy with the Plancherel theory for groups or symmetric spaces. This will be addressed in a follow up article.

We will now describe the contents of our paper in some detail. Throughout the paper, we assume that $G$ is a real reductive Lie group 
of the Harish-Chandra class, that $G = KAN_0$ is an Iwasawa decomposition and that $\cchi$ is a unitary 
character of $N_0.$ The character $\cchi$ is assumed to satisfy the regularity condition that for each simple root $\ga$ of $\fa$ in $\fn_0,$ 
the derivative $\cchi_*:= d\cchi(e)$ has a non-zero restriction to the root space $\fg_\ga.$  

Let 
$
C(G/N_0\col \cchi)$ be the space of continuous functions $f: G \to \C$ 
 transforming according to the rule
$$
f(xn) = \cchi(n)^{-1} f(x), \qquad (x \in G, n \in N_0),
$$%
and let $C_c(G/N_0\col \cchi)$ 
be the subspace of such functions with compact support modulo $N_0.$ The latter 
space has a natural left $G$-invariant pre-Hilbert structure, for which the completion
is denoted by $L^2(G/N_0\col \cchi).$ The representation of $G$ in this completion induced by the left regular action is the unitarily induced representation
\begin{equation}
\label{e: ind rep}
\Ind_{N_0}^G(\cchi).
\end{equation}
The Whittaker Plancherel formula 
concerns the unitary direct integral decomposition of (\ref{e: ind rep}).  It should be built from pairs $(\pi, \gl)$ 
with $\pi$ an irreducible unitary representation in a Hilbert space, and 
$\gl$ a continuous linear functional 
on the associated space of smooth vectors $H_\pi^\infty,$ transforming according the rule
$$
\gl \after \pi(n) = \cchi(n) \gl,\qquad (n\in N_0).
$$%
The functionals of this type are called 
Whittaker functionals of $\pi,$ and the space of these is denoted by $\Wh_\cchi(H_\pi^\infty).$ 
An element $\gl $ in the latter space determines a $G$-equivariant (Whittaker) matrix coefficient map 
$\wh_\gl: H_\pi^\infty \to C^\infty(G/N_0\col \cchi)$ given by 
$$ 
\wh_\gl(v)(x) = \gl(\pi(x)^{-1} v),\quad (v \in H_\pi^\infty , x \in G).
$$% 
In Sections \ref{s: whittaker vectors} 
and \ref{s: moderate estimates}, these Whittaker coefficients are discussed in more detail. They  
have moderate growth behavior towards infinity. In Lemma \ref{l: nilpotent trick}
we formulate a  technique which shows the importance of the regularity condition 
on $\cchi.$ As a consequence each Whittaker coefficient of the above type has
faster than exponential decay towards infinity in any closed cone disjoint from 
$\overline{ A^+}\setminus \{0\},$ see Corollary \ref{c: rapid descent}. In Section \ref{s: moderate estimates}  several related estimates are proven that are needed 
in the later sections.

Section \ref{s: W S space} concerns aspects of the Whittaker Schwartz space $\cC(G/N_0\col \cchi)$ as introduced by \cite{HCwhit} and \cite{Wrrg2}.  

In Section \ref{s: sharp estimates of whit coeff}
we discuss sharp estimates for a Whittaker 
coefficient $\wh_\gl: H_\pi^\infty \to C^\infty(G/N_0\col \cchi)$  in terms of a functional 
$\gL_V \in \fa^*$ attached to $\pi.$ More precisely, $\gL_V$ is defined in terms of 
the $\fa$-weights of $V / \fn_0 V,$ with $V$ 
the Harish-Chandra module of $K$-finite vectors of $\pi,$ see \ref{e: defi gL V}. In  \cite{Wrrg2} these estimates
were obtained on the positive chamber $A^+.$ In view of the results of 
Section \ref{s: moderate estimates}  
the estimates turn out to be valid on the entire group $A.$ In \cite{Wrrg2} the estimates on the positive 
chamber are obtained by using the method of estimate improvement along maximal parabolic subgroups. We use the same method, cast in the form of Lemma \ref{l: estimate imp}. This prepares
for the lengthy argumentation
in Section \ref{s: uni temp estimates}, where the uniformly tempered estimates are obtained.
The proof of Lemma \ref{l: estimate imp} is deferred to Section \ref{s: proof of est imp}.
We end Section \ref{s: sharp estimates of whit coeff} with Cor. \ref{c: schwartz estimate for ds coeff} which is due to both \cite{HCwhit} (on the $K$-finite level)
and \cite{Wrrg2}. It asserts that if $G$ has compact center and $\pi$ 
belongs to the discrete series of $G,$ then for every Whittaker vector $\gl \in \Wh_\cchi(H_\pi^\infty)$ 
the associated Whittaker coefficient $\wh_\gl$ is a continuous linear map
from $H_\pi^\infty$ into the Schwartz space
$\cC(G/N_0\col \cchi).$ 

In Section \ref {s: par ind reps} we discuss the space of smooth vectors 
for parabolically (normally) induced representations of the form 
$\Ind_P^G(\xi),$ with $P = M_P A_P N_P$ a parabolic subgroup of $G$ 
and $\xi$ a continuous representation of $P$ in a Hilbert space $H_\xi.$ 
For technical reasons we need to deal with this
in the generality of a representation of the form  $\xi = \gs \otimes \pi,$ with $\gs$ an irreducible  unitary representation
of $M_P,$ extended to $P$ by triviality on $A_P N_P$ and with $(\pi, F)$ a continuous  
representation of $P$ in a finite dimensional Hilbert space.
The main result is the characterization of Theorem
 \ref{t: identification of smooth vectors}, which asserts that the space of smooth vectors 
of $\Ind_P^G(\xi)$ equals 
\begin{equation}
\label{e: smooth vectors}
C^\infty(G/P:\xi):= \{f\in C^\infty(G,  H_\xi) \mid f(x man) = a^{-\rho_P} \xi(man )^{-1} f(x)\}.
\end{equation}  
The left regular representation of $G$ in this space is denoted by $\pi_{P,\xi}^\infty.$ 

In the subsequent Section \ref{s: generalized vectors for ind reps}, the space of generalized 
vectors for $\Ind_P^G(\xi)$ is defined as a conjugate continuous linear dual by 
\begin{equation}
\label{e: generalized vectors}
C^{-\infty}(G/P:\xi) := \overline{C^\infty(G/P:\xi^*)'}.
\end{equation}
Here $\xi^*$ is the continuous representation of $P$ in $H_\xi$ defined by $\xi^*(p):= \xi(p^{-1})^*.$ 
This definition has the advantage that (\ref{e: smooth vectors}) can be viewed as a subspace of 
(\ref{e: generalized vectors}), via a $G$-equivariant sesquilinear pairing defined by the usual
integration 
over $K/K\cap M_P.$ 

In Section \ref{s: whit for ind reps} we turn to the induced representations $\Ind_{\bar P}^G(\xi)$,
with $P$ a standard parabolic subgroup, and with $\xi = \gs \otimes (-\bar \nu) \otimes 1,$
 where $\gs$ is an irreducible unitary representation
of $M_P$ and $\nu \in \faPdc.$ The Whittaker functionals for the space of smooth vectors $C^\infty(G/\bar P:\gs: -\bar \nu)$ can then be identified with (conjugates of) elements of the space
of Whittaker vectors 
\begin{equation}
\label{e: chi invt gen sections ps} 
C^{- \infty}(G/\bar P: \gs:\nu)_\cchi:= \{ j \in C^{-\infty}(G/\bar P: \gs:\nu) \mid L_n j = \cchi(n) j, \; (n \in N_0)\}.
\end{equation}
The $N_0$-equivariance of an element $j$ of this space makes that on the open set $N_P \bar P$ 
it can be represented by a continuous function with values in $ H_\gs^{- \infty }:= \overline{H_\gs^\infty}{}'.$ 
Subsequent evaluation of this function in the identity element $e$ defines a linear map
\begin{equation}
\label{e: defi ev at e} 
\ev_e: C^{- \infty}(G/\bar P: \gs:\nu)_\cchi \to (H_\gs^{- \infty })_{\cchi_P},
\end{equation}
where $\cchi_P:= \cchi|_{M_P \cap N_0}.$ At this point we invoke the fundamental result 
\cite[Thm.~1]{HCwhit}, on which Harish-Chandra's entire treatment of the Whittaker theory is founded,
see Theorem \ref{t: main thm HC}. It allows us to conclude that the map
(\ref{e: defi ev at e}) is injective, see Corollary \ref{c: ev is injective}. 
 
Conversely, if $\gs$ is a representation of the discrete series of $M_P$, and if $\Re \nu$ 
is $P$-dominant, we define for each $\eta \in (H_\gs^{-\infty})_{\cchi_P}$ 
a continuous $H_\gs^{-\infty}$-valued function on $N_P \bar P$ which represents an element 
$j(\bar P, \gs, \nu, \eta)$ 
of (\ref{e: chi invt gen sections ps}) with $\ev_e j(\bar P, \gs, \nu, \eta) = \eta,$ see
Propositions \ref{p: integrability of j} and \ref{p: holomorphy j in dominant region}.
The element $j(\bar P, \gs, \nu, \eta)$  depends holomorphically on $\nu$ 
in the region where $\Re \nu$ is $P$-dominant.

In Section \ref{s: whit int} we discuss the close relation between $j(P, \gs, \nu,\eta)$ 
and the Jacquet integral introduced in
 \cite{Wrrg2}, see (\ref{e: comparison j and J}). Moreover,
we discuss the definition of Harish-Chandra's
Whittaker integral $\Wh(P,\psi, \nu)$ which is the analogue of the Eisenstein integral for a group or a symmetric space.  We show that the Whittaker integral  is expressible as a sum of Whittaker matrix coefficients involving Whittaker vectors $j(\nu) = j(\bar P, \gs, \nu, \eta),$ see
Corollary \ref{c: Wh as matrix coeff with j}.
 
In strong analogy with the theory of symmetric spaces one needs to extend the map $\nu \mapsto j(\nu)$ 
meromorphically in order to reach imaginary $\nu,$ which correspond to the unitary principal series. 
In addition one needs to establish uniformly tempered estimates in regions of the form 
$|\Re \nu|< \ge,$ with $\ge > 0$ a suitable constant. In the theory of symmetric spaces, 
tools for this program where initially developed in \cite{Bps2} for minimal $\gs$-parabolic subgroups 
and then extended to arbitrary $\gs$-parabolic subgroups in \cite{CD}.
It turns out that these tools from the theory of symmetric spaces are ideally suited for the Whittaker setting.
This unfolds in the final Sections \ref{s: fin dim reps} - \ref{s: uniform temp W int}.

In Section \ref{s: fin dim reps} we prepare  by reviewing the characterization
of irreducible finite dimensional spherical representations of $G$ with an $M_P$-fixed highest weight vector.
Then, in Section \ref{s: proj inf char} we consider the action of the center $\fZ$ of the universal enveloping 
algebra of $\fg_\iC$ on a tensor product of the form 
\begin{equation}
\label{e: ind tensor pi mu}
\Ind_{\bar P}^G(\gs\otimes \nu \otimes 1)\otimes \pi_\mu,
\end{equation} 
with $\pi_\mu$ 
an irreducible finite dimensional spherical representation of strictly $P$-dominant 
highest weight $\mu,$ with 
$M_P$ acting trivially on the highest weight space.
Let $\gL$ denote 
the infinitesimal character of $\gs,$ and let 
$p_{\gL + \nu +\mu}$ denote the projection in the space of (\ref{e: ind tensor pi mu}) onto the generalized weight space
for the infinitesimal character $\gL + \nu + \mu.$ 
 Then the main result of the section is that there exists a non-zero polynomial function $q$ on $\faPdc$ 
 such that 
$q(\nu) p_{\gL_\gs + \nu + \mu}$ can  be realized by the action of an element $\underline Z_\mu(\nu) \in \fZ$ wich depends polynomially on $\nu,$ see Corollary \ref{c: q p is Z nu}.

In Section \ref{s: functional equation} the element $\underline{Z}_\mu(\nu) \in \fZ$ is used to define a suitable differential operator
$D_\mu(\gs, \nu): C^{-\infty}(\bar P:\gs: \nu) \to C^{-\infty}(\bar P:\gs: \nu + \mu)$ such that 
one has a Bernstein-Sato type functional equation for the Whittaker vector, 
$$ 
j(\bar P,\gs , \nu ) = D_\mu(\gs, \nu) \after j(\bar P,\gs, \nu + \mu) \after R_\mu(\gs, \nu),
$$% 
see (\ref{e: functional equation}). Here $R_\mu(\gs)$ is a rational function on $\faPdc$ 
with values in $\End(H_{\gs,\cchi_P}^{-\infty}).$ The main problem is to show that 
$R_\mu(\gs, \dotvar)$ is generically invertible. This is done in the rest of Section 
\ref{s: functional equation} and the next.

In Section \ref{s: holo and ume}  the functional equation is used to obtain the meromorphic continuation
of $\nu \mapsto j(\nu) = j(\bar P, \gs, \nu)$ with estimates in terms of continuous seminorms 
on $C^{-\infty}(P, \gs, \nu)$, see Theorem \ref{t: estimate by functional equation}. 
The functional equation implies the existence of a non-zero polynomial $p_R$ 
such that $\nu \mapsto p_R(\nu) j(\nu)$ is holomorphic in the range $\inp{\Re\nu}{\ga} > R$ 
for all $\ga\in\gD.$ Any singularity of $j$ in this range is contained in the
zero set $M= p_R^{-1}(0).$ If it is contained in the regular part of $M$ then by a local 
analysis transversal to $M$ it can be shown that the singularity produces
a non-zero element of $C^{-\infty}(G/\bar P, \gs, \nu)_\cchi$ which must be zero at
$e$ hence zero. Therefore, singularities can only occur at singular points of $p^{-1}(0).$ The appearence of singularities would thus violate
a form of Hartog's principle formulated and proven in the appendix in Section \ref{s: Hartog thm}. 
Therefore, $j$ cannot have singularities. This gives a new proof of Wallach's
result \cite[Thm.~15.4.1]{Wrrg2} on the holomorphy of the Jacquet integral, but with strong estimates, see Theorems \ref{t: j is holomorphic} and \ref{t: uniform estimate of j}. This in turn leads to uniformly moderate
estimates for the associated family of Whittaker coefficients in Theorem \ref{t: uniformly moderate estimate}. 

In Section \ref{s: uni temp estimates}  it is shown that the uniformly moderate estimates
for the family of Whittaker coefficients
$(\wh_\nu)$ produced by $j_\nu = j(\bar P,\gs,\nu,\eta)$ can be improved to the so-called 
uniformly tempered estimates, see Theorem \ref{t: uniformly tempered estimates}. 
This is done by using the differential equations satisfied by $(\wh_\nu)$ and 
the method of estimate improvement along maximal parabolic subgroups, as in 
Section \ref{s: sharp estimates of whit coeff}, but now with suitable uniformity in the parameter $\nu.$ 
See Lemma \ref{l: improvement step} for the crucial stepwise improvement. 

In the final section, 16, the results obtained in the previous sections are applied 
to the Whittaker integral $\Wh(P, \psi, \nu) \in C^\infty(\tau\!:\! G/N_0\col \cchi).$ 
Here `$\tau\!:$' indicates that left $\tau$-spherical functions are considered, 
$\nu \in \faPdc$ and $\psi$ is an element of the finite dimensional space $\cA_{2,P}$ of 
$\fZ(\fm_P)$-finite functions in $\cC(\tau_P\col M_P/M_P\cap N_0\!: \!\cchi_P).$ 
The uniformly tempered estimates of Theorem \ref{t: uniformly tempered for Wh} thus 
obtained allow us to define 
a Fourier transform $\cF_P$ in terms of the Whittaker integral, and to show 
that $\cF_P$  defines a continuous linear map from the Whittaker Schwartz space
$\cC(\tau\!:\! G/N_0\col \cchi)$ to the Euclidean Schwartz space $\cS(i\faPd , \cA_{2,P}),$  see 
Theorem \ref{t: Fourier continuous into Schwartz}. 

\section{Whittaker vectors and matrix coefficients}
\label{s: whittaker vectors}
We consider a real reductive group $G$ of the Harish-Chandra class
and fix an Iwasawa decomposition $G = KAN_0.$ 
We denote by $\gS= \gS(\fg, \fa)$ 
the root system of $\fa$ in $\fg$ and by $\gS^+$ the positive system 
consisting of the roots $\ga \in \gS$ with $\fg_\ga$ contained in the 
Lie algebra $\fn_0$ of $N_0.$ The associated
collection of simple roots in  $\gS^+$ is denoted by $\gD.$ 

Thus,
if $M= Z_K(\fa)$, then $P_0:= MA N_0$ is the standard minimal parabolic
subgroup associated with $\gS^+.$ 
Here and in the rest of the paper we adopt the convention to 
denote Lie groups by Roman capitals and the associated Lie algebras by the corresponding
fraktur lower cases.

In this article $\cchi$ will always be a  unitary character of $\bN_0$. 
Following \cite[p.~142]{HCwhit}, we say that $\cchi$ is 
{\em regular} if its derivative $\cchi_*:= d\cchi(e) \in i \fn_0^*$ is non-zero on
each of the simple root spaces $\fg_\ga$, for $\ga \in \gD.$ 
Unless otherwise specified, it will always be assumed that $\chi$ is regular.
Note that the notion of regularity as defined here coincides with the notion 
of genericity in \cite[p.~371]{Wrrg2}.

We consider the function space 
\begin{equation}
\label{e: continuous whittaker space} 
C(G/\bN_0\col \cchi) :=\{ f \in C(G) \mid f(xn) = \cchi(n)^{-1} f(x), \;\;(x \in G, n \in \bN_0)\}.
\end{equation}
The subspace of functions with compact support modulo $\bN_0$ 
is denoted $C_c(G/\bN_0\col \cchi)$ and the subspace of smooth functions by $C^\infty(G/\bN_0\col \cchi).$ Finally, the intersection of the latter two is denoted $C^\infty_c(G/\bN_0\col \cchi).$ 

We fix a choice of positive invariant Radon measure $d \dot{x}$ on $G/N_0$ and define a
pre-Hilbert structure on $C_c(G/\bN_0\col \cchi)$  by the formula 
$$ 
\inp{f}{g} := \int_{G/N_0} f(x)\overline{ g(x)} \; d \dot{x},\quad (f,g \in C_c(G/N_0\col \cchi)).
$$%
The associated completion is denoted by $L^2(G/\bN_0\col \cchi).$ 
The Whittaker Plancherel formula concerns the unitary decomposition for the left regular representation of $G$ in the latter space.
Here we note that 
 $L^2(G/\bN_0\col \cchi)$ is  the space for the unitarily 
induced representation $\Ind_{\bN_0}^G(\cchi).$ Our notation is slightly different 
Harish-Chandra's,  who uses the notation $L^2(G/\bN_0\col \cchi)$ for the space of 
the induced representation
$\Ind_{\bN_0}^G(\cchi^\veec),$ with $\cchi^\veec: n\mapsto \cchi(n)^{-1},$ 
see \cite[p.~143]{HCwhit}.
Finally, note that Wallach \cite[p.~365]{Wrrg2} works with another realization of this representation space, namely $L^2(\cchi\col N_0 \bs G)$ equipped with the right regular representation. 

In the following we will need a bit of background from representation theory that we will now 
explain.

If $V$ is a locally convex (Hausdorff) space then a representation $\pi$ of a Lie group $L$ in 
$V$ is called smooth if it is continuous and if $V = V^\infty.$ 

If $V$ is  Fr\'echet space, its strong dual $V'$ is a complete locally convex space.
Suppose that $\pi$ is a continuous representation of $L$ in $V,$ then the homomorphism 
$\pi^\veec: L \to {\rm GL}(V')$ defined by 
\begin{equation}
\label{e: contragredient}
\pi^\veec(x)\xi = \xi\after \pi(x)^{-1}, \qquad (x \in L, \;\xi \in V').
\end{equation}
need not be a continuous representation. However, if $\pi$ is smooth, then $\pi^\veec$ is a smooth representation of $L$ in $V,$ called the contragredient
of $\pi.$ For details we refer to  \cite{Warn1}, Proposition 4.4.1.9 and the definition of $V^\veec$ preceding 
Proposition 4.1.2.1. 

To prepare for the treatment of conjugate representations, we first briefly discuss the notion
of conjugate space. Let $V$ and $W$ be complex linear spaces. We denote by $\bar V$
the real linear space $V$ equipped with the conjugate complex multiplication 
$\C \times V \to V,$ $ (z, v)\mapsto \bar z v.$ A map $T: V \to W$ is said
to be conjugate linear if it is real linear and satisfies $T(\gl v) = \bar\gl T(v)$ for $v \in V$ and $\gl\in \C.$ The complex linear space of conjugate linear maps 
$V \to W$ (equipped with the pointwise operations of scalar multiplication and addition) equals the complex linear space $\Hom_\C(\bar V, W).$ Given $T \in \Hom(V,W)$ we denote by 
$\bar T$ the map $V \to W$ viewed as an element of $\Hom(\bar V, \bar W).$ We note
that the map 
$$ 
T \mapsto \bar T, \quad \Hom(V,W) \to \Hom(\bar V, \bar W)
$$ 
is not complex
linear, but conjugate linear. Hence,  $T \mapsto \bar T$ is 
a complex linear isomorphism $\overline{\Hom(V,W)} \to \Hom(\bar V, \bar W).$ 
Since the map $T\mapsto \bar T$ is the identity on the set $\Hom(V,W)$ we have the following
identity of complex linear spaces 
\begin{equation}
\label{e: identity hom and bar}
\overline{\Hom(V,W)} = \Hom(\bar V, \bar W).
\end{equation} 
In the sequel we will encounter the conjugate dual space  $\overline{V^*} = \overline{\Hom(V,\C)}$ 
and  the dual conjugate space $(\bar V)^* = \Hom(\bar V, \C).$ These spaces
are not equal, but complex linearly isomorphic under the map $\gl \mapsto \overline{{}^c\! \gl}$ given 
by 
\begin{equation} 
\label{e: defi c gl}
{}^c \!\gl = {\bf c} \after \gl, 
\end{equation}
where ${\bf c}: \C \to \C$ denotes the conjugation
map $z \mapsto \bar z.$ 

Indeed, ${\bf c} \in \Hom(\bar \C, \C)$ so in view of the equality (\ref{e: identity hom and bar}) it follows that ${}^c\!\gl \in \Hom(\bar V, \C)$ $(\gl \in \Hom(V, \C)).$ 
It is now readily verified that $\gl \mapsto {}^c\!\gl$ 
is a conjugate linear map $\Hom(V, \C) \to \Hom(\bar V, \C)$ hence a complex linear
isomorphism
$$ 
\overline{\Hom(V, \C)}\;\; {\buildrel \simeq \over \longrightarrow}\;\;  \Hom(\bar V, \C).
$$ 

Given a representation $\pi$ of $L$ in a locally convex space $V$,
we denote by $(\bar \pi, \bar V)$ the conjugate of $\pi.$ Here the conjugate complex
linear space $\bar V$ is equipped
with the locally convex topology of $V.$  
 Furthermore, for $x \in L,$ 
$\bar \pi(x)$ equals the complex linear map $\bar {\pi(x)}: \bar V \to \bar V.$ 
It is clear that $(\bar \pi, \bar V)$ is a representation of $L$ in a locally convex 
space again, which is continuous if and only if $\pi$ is continuous.

We note that the spaces $\overline{V'}$ and  $\bar V '$ are
topologically complex linear isomorphic under the map $\gl \mapsto {}^c\! \gl$ given 
by (\ref{e: defi c gl}).

If $\pi$ is a smooth Fr\'echet representation, then $\bar \pi^\veec$ is a smooth continuous representation, and therefore, so is the equivalent representation $\overline{\pi^\veec}.$ 

We use the notation $U(\fl)$ for the universal enveloping algebra of the complexification
$\fl_\iC$ of $\fl.$ The canonical anti-automorphism of $U(\fl)$ is denoted by $u \mapsto u^\veec.$
 It is readily verified that for the associated infinitesimal representations
$\pi: U(\fl) \to \End(V)$ and $\pi^\veec: U(\fl) \to \End(V')$ we have
$$ 
\pi^\veec (u) \xi = \xi \after \pi(u^\veec), \quad (u \in U(\fl), \xi \in V').
$$%

Given a continuous Hilbert space representation $(\pi, H)$ of $L$ 
it is known that the contragredient $(\pi^\veec,H')$
is continuous, see \cite[Cor. 4.1.2.3]{Warn1}. 
Therefore, so are $\bar\pi^\veec$ and $\overline{\pi^\veec}.$ Both duals 
$(\bar H)'$ and $(H')^\veec$ come into play through the Hermitian
 inner product $b$ viewed as a  bilinear map $H \times \bar H \to \C.$ 
Let $b_1: H \to (\bar H)'$ be the linear map defined by $b_1(v)= b(v, \dotvar).$ 
Let $b_2: \bar H \to H'$ be the linear map defined by $b_2(v) = b(\dotvar, v).$ 
Then $b_2$ can be viewed as a linear map $H \to \overline{H'}.$ 
As such $b_1$ and $b_2$ are topological linear isomorphisms
from $H$ onto $(\bar H)'$ and $\overline{H'}$ respectively. 
Thus, here the isomorphism  $\overline{H'} \to \bar H'$ is given by $\gb= b_1\after b_2^{-1}.$ 
Using the conjugate symmetry of $b$, it readily follows that $\gb$ coincides with the 
isomorphism ${}^c(\dotvar)$ defined by (\ref{e: defi c gl}). Since $\gl \mapsto {}^c\! \gl$ 
intertwines the representations $\overline{\pi^\veec }$ and $\bar \pi^\veec,$ 
it follows that $b_2$ and $b_1$, respectively, intertwine these
representations with the same continuous representation $\pi^*$ of $L$ in $H,$ given by 
$$ 
\pi^*(x) = \pi(x^{-1})^* \qquad (x \in L),
$$ 
where the star indicates that the Hilbert adjoint is taken.

For such a continuous Hilbert representation 
we denote the associated Fr\'echet representation in the space
of smooth vectors by $(\pi^\infty, H^\infty).$ As we mentioned above, 
the continuous linear dual $(H^\infty)'$ of $H^\infty,$ equipped with the strong 
dual topology, is complete. The associated  contragredient of $\pi^\infty$ in
$(H^\infty)'$, denoted by $(\pi^{\infty\veec}, (H^\infty)')$ is a smooth representation of
$L.$ 
In this setting, with $G$ in place of  $L,$ it is of interest to consider the space of Whittaker functionals
\begin{equation}
\label{e: defi space of Wh functionals}
\Wh_\cchi(H^\infty) := \{\gl \in (H^\infty)' \mid \forall \bn \in \bN_0:\;\;  \gl \after \pi^\infty(\bn) = \cchi(\bn)
\gl\},
\end{equation} 
see \cite[15.3.4, p.~378]{Wrrg2}.
Equivalently, $\Wh_\cchi(H^\infty)$ consists of the functionals 
$\gl \in (H^\infty)'$ such that 
$$ 
\pi^{\infty\veec}(n) \gl = \cchi(n)^{-1} \gl, \qquad (n \in N_0).
$$% 
For a given Whittaker functional $\gl \in \Wh_\cchi(H^\infty)$, 
the matrix coefficient map $\wh_\gl: H^\infty \to C^\infty(G/N_0 \col \cchi)$, 
given by 
$$ 
\wh_\gl(v)(x) = \gl(\pi(x)^{-1} v),
$$%
is readily seen to be continuous and $G$-equivariant. Moreover, we
have the following easy lemma. 

\begin{lemma}
The matrix coefficient map $\gl \mapsto \wh_\gl$ is a bijection 
\begin{equation}
\label{e: matrix coefficient Whittaker}
\Wh_\cchi(H^\infty)\;\;
 {\buildrel \simeq \over \longrightarrow} \;\;\Hom_G(H^\infty, \; C^\infty(G/\bN_0\col \cchi)),
\end{equation}
where $\Hom_G$ indicates the space of intertwining continuous linear maps. 
The inverse of (\ref{e: matrix coefficient Whittaker}) 
is given by $T \mapsto \ev_e \after T$, where
$\ev_e: C^\infty(G/\bN_0\col \cchi) \to \C$ denotes evaluation at the identity.
\end{lemma}

The following result, valid for any continuous character $\chi$ of $N_0,$ is stated and proven in 
\cite[Cor 15.4.4]{Wrrg2}. 

\begin{lemma}
\label{l: finite dimensionality whit vectors}
If $(\pi, H)$ is admissible and of finite length, then 
$$ 
\dim \Wh_\chi(H^\infty) < \infty.
$$ 
\end{lemma}

Given a continuous representation $\rho$ of a Lie group in a complete locally convex space
$V,$ we use the notation $\bar\rho$ for $\rho$ viewed as a representation in 
the conjugate space $\bar V.$ Clearly, $\bar \rho$ is continuous again and the identifications
$\overline{\rho^\infty} = {\bar \rho}^\infty$ and $\overline{V^\infty} = \bar V^\infty$ are obvious.

In this paper it will be desirable to view the matrix coefficient map 
$\wh_\gl,$ for $\gl \in \Wh_\cchi(H^\infty)$ as a matrix coefficient with 
a suitable generalized vector. This is possible in the following setting of duality.

Let $(\pi_j, H_j)$ be two continuous Hilbert representations of a Lie group $L,$ for $j=1,2.$ 
By a perfect sesquilinear pairing of $\pi_1$ and $\pi_2$ 
we mean an equivariant continuous sesquilinear pairing 
\begin{equation}
\label{e: sesquilinear pairing} 
H_1 \times H_2 \to \C
\end{equation}
such that the induced maps $\ga_1: H_1 \to (\bar H_2)'$ and 
$\ga_2: H_2 \to \overline{H_1'}$ are unitary isomorphisms. Note that these maps are intertwining.
In particular, the restriction map
$$
r_1: \xi \mapsto \xi|_{H_1^\infty}, \qquad \overline{H_1'} \to \overline{{H_1^\infty}'}
$$%
is a continuous linear injection, intertwining $\overline{\pi_1^\veec}$ with 
$\overline{\pi_1^{\infty\veec}}.$ 

We put 
\begin{equation}
\label{e: defi H2 minus}
H_2^{-\infty} := \overline{{H_1^\infty}'},
\end{equation}
and accordingly denote by $\pi_2^{-\infty} := \overline{\pi_1^{\infty\veec}}$
the natural continuous representation on (\ref{e: defi H2 minus}).
We consider the canonically associated sesquilinear pairing  
\begin{equation}
 \label{e: pairing smooth and minus}
H_1^\infty  \times H_2^{-\infty} \to \C,\;\; (v,j) \mapsto \inp{v}{j}.
\end{equation}
This pairing is equivariant for $\pi_1^\infty$ and $\pi_2^{-\infty}$ and induces the inverse of the continuous linear isomorphism 
(\ref{e: dual H1 infty to H2 minus infty}).
Put $\iota_2: = r_1 \after\ga_2.$ Then 
the map 
\begin{equation}
\label{e: embedding iota 2}
\iota_2:  H_1 \embeds H_2^{-\infty}
\end{equation} 
is a continuous linear injection, intertwining $\pi_1$ with $\pi_1^{-\infty}.$ 
We will use it to identify the first of these spaces as an invariant subspace of 
the second. 
This allows us to view the elements of (\ref{e: defi H2 minus})
as generalized vectors for $\pi_2.$ 
We now note that for $(v,w)\in H_1^\infty \times H_2,$ we have
$$ 
\inp{v}{\iota_2(w)} = [r_1(\ga_2(w))](v) = [\ga_2(w)](v) = \inp{v}{w},
$$ 
where the last mentioned pairing is  (\ref{e: sesquilinear pairing}). We thus see that 
the sesquilinear pairing (\ref{e: pairing smooth and minus}) is an extension 
of the pairing $H_1^\infty \times H_2 \to \C$ given by restricting (\ref{e: sesquilinear pairing}).

In the present context it is sometimes  convenient to also use the sesquilinear pairing
\begin{equation} 
\label{e: reversed order pairing}
H_2^{-\infty} \times H_1^\infty \to \C, \quad  (j, v) \mapsto \inp{j}{v}:= \overline{\inp{v}{j}}.
\end{equation}
Finally, we note that the above definitions imply that the pairing (\ref{e: pairing smooth and minus})
induces a topological linear isomorphism 
\begin{equation}
\label{e: dual H1 infty to H2 minus infty}
\overline{H_2^{-\infty}}  \;\; {\buildrel \simeq \over \longrightarrow}\;\; (H_1^\infty)'.
\end{equation}
This isomorphism intertwines $\overline{\pi_2^{-\infty}}$ with $\pi_1^{\infty\veec}.$
In the Whittaker setting, with $L = G$ and $\pi_1$ and $\pi_2$ of finite length,
this equivariance implies that (\ref{e: dual H1 infty to H2 minus infty}) restricts
to a bijective conjugate linear 
map 
\begin{equation}
\label{e: from j to wj}
H_{2,\chi}^{-\infty} \;\; {\buildrel \simeq \over \longrightarrow}\;\;   \Wh_\chi(H_1^\infty),\quad
j \mapsto \trw\! j,
\end{equation}
where
\begin{equation}
\label{e: whittaker vectors}
 H_{2,\chi}^{-\infty}:= \{j \in H_{2}^{-\infty}\mid \pi_2^{-\infty}(n) j = \chi(n) j\;  (n \in N_0)\}.
\end{equation} 
For obvious reasons, we agree to call (\ref{e: whittaker vectors})
the space of Whittaker vectors for $\pi_1.$ 
Given $j \in (H_2^{-\infty})_\cchi$ there is the associated Whittaker coefficient map
$\wh_j = \wh_{(\trw\! j)}: H_1^\infty \to C^\infty(G/N_0\col \cchi),$ given by 
$$ 
\wh_{j}(v)(x) = \trw\!j( \pi_1(x)^{-1} v) = \inp{ \pi_1(x)^{-1} v}{ j}, \;\;\;  (v \in H_1^\infty , \; x \in G).
$$%

\begin{rem}
\label{r: dual of pi}
If $(\pi, H)$ is a continuous Hilbert representation of  a Lie group $L$ we define the representation
$\pi^*$ of $L$ in $H$ by $\pi^*(x) = \pi(x^{-1})^*$ for $x \in L.$ It is readily verified
that the isometry ${\rm i}: H \to \bar H'$ induced by the Hilbert inner product 
intertwines $\pi^*$ with $\bar \pi^\veec.$ Thus, $\pi^*$ is continuous since 
$\bar\pi^\veec$ is continuous. Let now $H_\pi$ and $H_{\pi^*}$ denote
$H$ equipped with the representations $\pi$ and $\pi^*,$ respectively. 
Then the inner product of $H$ gives an equivariant perfect sesquilinear pairing
\begin{equation}
\label{e: pairing pi with pi star}
H_\pi \times H_{\pi^*}\to \C.
\end{equation}
Conversely, any equivariant perfect pairing of the form (\ref{e: sesquilinear pairing}) 
can be transfered to a pairing as (\ref{e: pairing pi with pi star}) by putting $H = H_1,  \pi = \pi_1$ 
and using the equivariant unitary isomorphism ${\rm i}_2^{-1} \after \ga_2: H_2 \to H_1.$
The pairing (\ref{e: pairing pi with pi star}) gives rise to an intertwining 
injective linear map $H_{\pi^*}  \to H_{\pi^*}^{-\infty} = \overline{{H}_{\pi}^\infty{}'}.$ The representation
$\pi$ is unitary if and only if $\pi = \pi^*.$ In that case we obtain the equality 
$H_\pi^{-\infty} = \overline{H_{\pi}^\infty{}'}$ which is compatible with an existing convention in the literature.
\end{rem}

\begin{rem} 
The point of view explained above Remark \ref{r: dual of pi} will be of particular importance in the setting
of parabolically induced representations of the form $\Ind_P^G(\xi),$ 
with $\xi$ a continuous Hilbert representation of a parabolic subgroup $P$ of $G.$ 
Let $L^2(G/P:\xi )$ be the Hilbert space in which $\Ind_P^G(\xi)$ is 
realized by the left regular action. Then, with the similar notation for $\xi^*,$  there exists a natural $G$-equivariant perfect sesquilinear pairing
$$ 
L^2(G/P:\xi ) \times L^2(G/P:\xi^*) \to \C.
$$% 
Applying the formalism introduced above, one obtains a compatible equivariant
sesquilinear pairing
$$
L^2(G/P:\xi)^\infty  \times L^2(G/P:\xi^*)^{-\infty} \to \C
$$% 
which induces an equivariant continuous linear isomorphism 
$$ 
L^2(G/P:\xi^*)^{-\infty} \;\;
{\buildrel  \simeq \over \longrightarrow}\;\; 
(\overline{L^2(G/P:\xi)^\infty{}'}
$$%
The associated space of Whittaker vectors,
$
(L^2(G/P:\xi^*)^{-\infty})_\cchi,
$
can thus be viewed as a space of generalized sections of a Hilbert bundle.
\end{rem}

%%%%%%%%%%%%%%%%%%%%%%%%%%%%%%%%%%%%%%%%%%
\section{Moderate estimates for Whittaker coefficients}\label{s: moderate estimates}
We fix a non-degenerate $\Ad(G)$-invariant symmetric bilinear form
\begin{equation} 
\label{e: intro of B}
B: \fg \times \fg \to \R
\end{equation}
which is negative definite on $\fk,$ positive definite on $\fp$ and which restricts to the Killing form on the semisimple 
part $[\fg, \fg].$ Let $\theta$ denote the Cartan involution on $\fg$ associated with $K.$ 
We define the positive definite $\Ad(K)$-invariant inner product 
$\inp{\dotvar}{\dotvar}$ on $\fg$ by 
\begin{equation}
\label{e: defi inp}
\inp{X}{Y}:= - B(X, \Cartan Y), \qquad (X,Y\in\fg).
\end{equation}
The restriction of this inner product to $\fa$ induces a dual inner product
on $\fa^*,$ The latter's extension to a complex biliinear form on $\fadc$ is also denoted 
$\inp{\dotvar}{\dotvar}.$  The associated norms on $\fa$ and $\fad$ are denoted by 
$|\dotvar|.$ Finally we extend the norm on $\fad$ to the norm $|\dotvar|$ on $\fadc$ associated with the 
Hermitian extension of the inner product on $\fad.$ Accordingly, 
\begin{equation}
\label{e: defi norm on fadc}
|\nu|^2 = |\Re \nu|^2 + |\Im \nu|^2, \qquad (\nu \in \fadc).
\end{equation}%
If $\fv \subset \fa$ is a linear subspace, we will use the inner product on $\fa$ 
to identify the real linear dual $\fv^*$ with a subspace of $\fa^*,$ unless
otherwise specified.

We define ${}^\circ G$ to be the intersection of the kernels $\ker \xi$ 
where $\xi$ ranges over the characters $G \to \R_{>0}.$ Let 
$\fa_\gD = \cap_{\ga \in \gD} \ker\ga$ and put $A_\gD = \exp (\fa_\gD).$ 
Then multiplication induces an isomorphism of Lie groups
$$ 
G \simeq {}^\circ G \times A_\gD.
$$ 
It follows that $G = {}^\circ G$ if and only if $G$ has compact center.

We define ${}^\circ \|\dotvar\|: {}^\circ G \to ]0,\infty[$ by 
${}^\circ\|x\| = \|\Ad(x)\|_{\rm op}$ $(x \in G),$ 
where the subscript `${\rm op}$' indicates that the operator norm with 
respect to the inner product (1.2) has been taken.

We put ${}^*A_\gD := G^0\cap A.$ Then via the direct sum $\fa = {}^*\fa_\gD \oplus \fa_\gD$ 
we identify the elements of the real duals ${}^*\fa_\gD^*$ and $\fa_\gD^*$  with elements of $\fa^*.$ 
We select a basis $\cB$ of $\fa_\gD^*$ and define $\|\dotvar \|_\gD: A_\gD \to [1, \infty[$ 
by  
$$ 
\|b\|_\gD = \max_{\beta \in \pm \cB}  b^{\gb}, \qquad (b \in A_\gD).
$$ 
Finally we define $\|\cdot \|: G \to [1,\infty[$ by 
$$ 
\| x b\|  := \max( {}^\circ \|x\|, \|b\|_\gD),\qquad (x \in {}^\circ G , \; b \in A_\gD).
$$ 
Put $\gS_e:= \gS \cup \cB \cup (-\cB).$ Then it is easily verified that for $k_1, k_2 \in K$ and $a \in A,$ 
\begin{equation}
\label{e: norm of kak}
\|k_1 a k_2\| = \max_{\ga \in \Sigma_e} a^\ga.
\end{equation}
See also \cite[Lemma 2.1]{BSbv}, where the definition is given for $G$ with compact center. 

From the above definitions and (\ref{e: norm of kak}) it readily follows that  $\|\dotvar \|$ is a norm on $G$ in the sense 
of \cite[Lemma 2.A.2.1]{Wrrg1}.
 
\begin{lemma}
If $(\pi, H)$ is a (continuous) representation of $G$
in a Hilbert space, then there exist constants $r(\pi)  \geq 0$ and $C > 0$ 
such that
\begin{equation}
\label{e: norm pi g op}
\|\pi(g)\|_{\rm op} \leq  C \|g\|^{r(\pi)}, \qquad (g \in G),
\end{equation} 
where $\|\dotvar\|_{\rm op}$ indicates the operator norm. 
\end{lemma}
 \proof See \cite[Lemma 2.A.2.2]{Wrrg1}.
\qed

\begin{lemma}
Let $(\pi, H)$ be of finite length, and $\gl \in \Wh_\cchi(H^\infty).$ There exists
a constant $r > 0$ and a continuous seminorm ${\rm n}$ on $H^\infty$ 
such that the Whittaker coefficient $\wh_\gl$ satisfies 
\begin{equation}
\label{e: first estimate Whittaker coeff}
|\wh_\gl(v)(a)|\leq e^{r|\log a|} {\rm n}(v),
\end{equation} 
for all $v \in H^\infty $ and $a \in A.$ 
\end{lemma}

\proof 
Let $r(\pi)$ and $C$ be as in (\ref{e: norm pi g op}). Put $m = \max_{\ga \in \gS_e} |\ga|.$
Then it follows that for $a \in A$ we have
$$ 
\|a\|^{r(\pi)}  \leq e^{m r(\pi)| \log a|}.
$$ 

By continuity of $\gl,$ there exist a finite subset $S \subset U(\fg)$ 
such that 
$|\gl(v)| \leq \sum_{u \in S} \|\pi^\infty(u) v\|, $ for all $v \in V^\infty.$
By decomposing each element of $S$ as a sum of weight vectors for 
$\ad(\fa)$ it is readily seen that we may assume $S$ to consist of weight vectors
from the start. Let $\xi_u$ denote the weight by which $\ad(\fa)$ acts on $u.$ 
Then it follows that, for all $v \in H^\infty$ and $a \in A,$ 
\begin{eqnarray*}
| \gl(\pi(a) v)| & \leq & \sum_{u\in S} \|\pi^\infty(u) \pi(a)v \| \\
&\leq&
C\, \|a\|^{r(\pi)}  \sum_{u \in S} a^{-\xi_u} \|\pi^\infty(u) v\| \\
&\leq & 
C\, e^{(m r(\pi)  + s) |\log a| }  \sum_{u \in S} \| \pi^\infty(u) v\|,
\end{eqnarray*}
where $s = \max_{u \in S} |\xi_u|.$ 
The result follows with $r = m r(\pi) + s.$ 
\qed

The above proof does not use the assumption that $\cchi \in \widehat N_0$ 
is regular.
If  $\cchi$ is regular, then the above estimate
gives rise to remarkable new exponential estimates. The argumentation for this
is suggested by the following lemma.

\begin{lemma}
\label{l: nilpotent trick}
Let $u \in U(\fn_0)$ have weight $\eta \in \fa^*$ for the adjoint action
of $\fa.$ Then for all $f \in C^\infty(G/N_0\col \cchi)$ and $a \in A,$ 
$$ 
L_u f(a) = a^{-\eta} \cchi_*(u) f(a).
$$% 
\end{lemma}

\proof 
Note that 
$$
L_u f(a) = [R_{\Ad(a)^{-1} (u^\veec)} f](a)  = a^{-\eta} R_{u^\veec} f(a) = a^{-\eta} \cchi_*(u) f(a).
$$ 
\vspace{-40pt} {\ }

\qed
\medbreak
The regularity of the character implies that Whittaker coefficients have fast
decay outside the closed positive Weyl chamber $\cl (A^+).$ 

\begin{cor}
\label{c: rapid descent}
Let $(\pi, H)$ be an admissible Hilbert representation of $G$ of finite length
and let $\gl \in \Wh_\chi(H^\infty).$ 
Let $\Gamma$ be a closed cone in $\fa$ which is disjoint from $\cl(\fa^+)\setminus \{0\}.$ 

Then for every $s > 0$ there exists a continuous seminorm $n$ of $H^\infty$ such 
that, for all $v \in H^\infty,$ $k \in K$ and $a \in \exp (\Gamma),$ 
\begin{equation}
\label{e: rapid decay wh} 
|\wh_\gl(v)(ka)| \leq  e^{-s|\log a\, |} n(v).
\end{equation}  
\end{cor}

\proof 
Let $r > 0$ and $n$ be as in (\ref{e: first estimate Whittaker coeff}). 
Let $S$ be the unit sphere in $\fa$. Then by compactness of $S \cap \Gamma$ 
it suffices to show that for every $H_0 \in S\cap \Gamma$ there exists a closed neighborhood
$\omega \ni H_0$ in $S$ such that (\ref{e: rapid decay wh}) holds  with suitable $s' $ and $n'$ 
in place of $s$ and $n$, for all $v \in H^\infty,$ all $k \in K$ and all $a \in \exp(\R_{\geq 0} \omega).$ 

Let $H_0 \in S \cap \Gamma$ be given. Then $H_0\notin \cl(\fa^+)$ and it follows that there exists a simple root $\ga \in \gD$ 
such that $\ga(H_0) < 0.$ We may fix $p \in \N$ such that for $H = H_0$ 
we have
\begin{equation}
\label{e: estimate r + p H}
r|H| + p \ga(H) < - s|H|.
\end{equation} 
We may now fix a closed neighborhood $\omega$ in $S$ such that this estimate
holds for $H  \in \omega.$ By positive homogeneity (\ref{e: estimate r + p H})  holds for 
$H \in \R_+ \omega.$  Let now $X \in \fg_\ga$ be such that $\chi_*(X) = 1.$ 
Put $u = X^k \in U(\fg).$ Then by application of Lemma \ref{l: nilpotent trick} it follows that,
for all $v \in H^\infty,$ $k \in K$ and $a \in \exp (\R_{\geq\, 0} \omega),$ 
\begin{eqnarray*}
|\wh_\gl(v)(k a) |  & = & a^{k\ga} |L_u(\wh_\gl(\pi(k)^{-1}v))(a)| 
\\
&=& e^{k\ga(\log a)} |\wh_\gl(\pi(u)\pi(k)^{-1}v)(a)|\\
& \leq &  e^{k\ga(\log a) + r|\log a|} n (\pi(u) \pi(k)^{-1}v) \leq e^{-s|\log a|} n'(v),
\end{eqnarray*}
where $n' (v) = \sup_{k\in K} n ( \pi(u)\pi(k)^{-1} v).$ 
\qed

\begin{lemma}
\label{l: basic comparison of estimates}
Let $\vartheta: \fa \to \R$ be either linear, or of the form $\vartheta = r|\dotvar|,$  with $r> 0.$ 
Let $\xi \in \fad$ and assume that $\xi \geq  \vartheta$ on $\fa^+.$ 
Then there exists a finite subset $\Theta \subset U(\fn)$ 
such that for all $f \in C^\infty(G/N_0\col \cchi)$ and all $a\in A$ we have the estimate
\begin{equation}
\label{e: basic estimate f from exponent}
 a^{-\xi} |f(a)| \leq e^{-\vartheta(\log a)} \max_{u \in \Theta} |L_uf(a)|.
\end{equation}
\end{lemma}

Before we start with the proof, we need to introduce suitable notation. As usual, for $\Phi \subset \gD$ we define
$\fa_\Phi$ to be the intersection of the spaces $\ker \ga,$ for $\ga \in \Phi.$ 
In particular, $\fa_\gD$ equals the centralizer of $\fg$ in $\fa.$ We agree to 
write ${}^*\fa_\Phi$ for the orthocomplement of $\fa_\Phi $ in $\fa.$ Then 
$\fa = {}^*\fa_\Phi \oplus \fa_\Phi.$ 

The collection of restrictions $\ga|_{{}^*\fa_\gD},$ 
for $\ga \in \gD,$ is a basis of ${}^*\fa_\gD^*.$ The associated dual basis 
of ${}^*\fa_\gD$ is denoted by $\{h_\ga\mid \ga \in \gD\}.$ We define
$
\bp \fa_\Phi : ={\rm span}_\R\{h_\ga \mid \ga \in \Phi\}.
$
Then we have the following direct sum decomposition
\begin{equation}
\label{e: bp deco of fa}
\fa = \bp\fa_\Phi \oplus \fa_\Phi.
\end{equation} 
This decomposition will be important in the proof of Lemma 
\ref{l: basic comparison of estimates}.

We denote by $\bp \fa_\Phi^-$ the interior in $\bp\fa_\Phi$ of the closed
cone spanned by the elements $- h_\ga$ for $\ga \in \Phi.$ Then 
$$ 
\bp\fa_\Phi^- =  \{ H \in \bp\fa_\Phi \mid (\forall \ga\in \Phi): \ga(H) < 0\}.
$$% 
In addition, we define $\fa(\Phi):= \bp \fa_\Phi^- + \cl(\fa_\Phi^+).$
\begin{lemma}
The set $\fa$ is the disjoint union of the sets $\fa(\Phi),\,$ for $\Phi \subset \gD.$ 
\end{lemma}
\proof
We write $\R^\gD$ for the real linear space of functions $\gD \to \R.$ 
For $\Phi\subset \gD$ we define  $\R^\gD(\Phi)$ to be the subset of 
$\R^\gD$ consisting of $x \in \R^\gD$ with $x_\ga<0$ for $\ga\in \Phi$ 
and $x_\gb \geq 0$ for $\gb\in \gD\setminus \Phi.$ It is clear that
$\R^\gD$ is the disjoint union of the sets $\R^\gD(\Phi).$ 
Consider the linear map $p: \fa \to \R^\gD$ defined by $p(H)_\ga = \ga(H)$ 
for $\ga \in \gD.$ Then $p$ is a surjective linear map, hence
$\fa$ is the disjoint union of the sets $p^{-1}(\R^\gD(\Phi)),$ for $\Phi \subset \gD.$ 

We  will finish the proof by showing that $p^{-1}(\R^\gD(\Phi)) = \fa(\Phi).$ 
For this, suppose $H \in \fa$ and consider the decomposition
$H = \bp H + H_\Phi,$ according to (\ref{e: bp deco of fa}). Then  
$H \in p^{-1}(\R^\gD(\Phi)) $ is equivalent to the assertion that
$\ga(H) < 0$ and $\gb(H) >0$ 
for all $\ga \in \Phi$ and $\gb \in \gD\setminus \Phi.$  
This in turn is equivalent to the assertion that $\ga(\bp H) < 0$ and $\gb(H_\Phi) \geq 0$ 
for all $\ga \in \Phi$ and $\gb \in \gD\setminus \Phi,$ hence to 
$\bp H \in \bp \fa_\Phi^-$ and $H_\Phi \in \cl(\fa_\Phi^+).$ 
By definition, the latter is equivalent to $H \in \fa(\Phi).$ 
\qed

\medno
{\em Proof of Lemma \ref{l: basic comparison of estimates}.\ }
We may fix $k \in \N$ sufficiently large, such that for every $\ga \in 
\gD$ we have $k \ga +  \vartheta \leq \xi$ on $- h_\ga.$  For $\Phi \subset \gD$ we put 
$$ 
\gs_\Phi:= \sum_{\ga \in \Phi} \ga.
$$% 
Then it is readily verified that $k\gs_\Phi + \vartheta \leq \xi$ on 
$\bp\fa_{\Phi}^{-}$. Since $\gs_\Phi$ vanishes on $\fa_\Phi^+,$ whereas $\fa_\Phi^+\subset \cl(\fa^+),$
it follows from the hypothesis that the same estimate is valid on $\fa_\Phi^+.$ Using the subadditivity 
of $\vartheta$ and the linearity of $\gs_\Phi$ and $\xi$ we now find that 
\begin{equation}
\label{e: estimate r versus xi}
 k \gs_\Phi  + \vartheta \leq \xi \qquad{\rm on}\quad \fa(\Phi).
\end{equation}
The idea is now to derive  suitable estimates on the set $A(\Phi): = \exp (\fa(\Phi))$ 
for $\Phi \subset \gD$ fixed, by using Lemma \ref{l: nilpotent trick}.
Put $u=u_\Phi= \prod_{\ga \in \Phi} X_\ga^k$, where an arbitrary fixed ordering in the product
may be taken, and where $X_\ga \in \fg_\ga$ are such that $\cchi_*(X_\ga) =1.$ 

Let  $f\in C^\infty(G/N_0\col \cchi)$ and $a  \in A(\Phi),$ then it follows that 
$$ 
|f(a)| = a^{k\gs_\Phi} |\cchi_*(u)|^{-1} |L_u f(a)|= a^{k\gs_\Phi} |L_u f(a)|.
$$% 
Therefore, 
$$ 
a^{-\xi} |f(a)| = a^{-\xi + k\gs_\Phi } |L_u f(a)| \leq 
e^{-\vartheta(\log a)} |L_{u} f(a)|.
$$%
As the sets $\fa(\Phi)$ cover $\fa$, we find 
the desired estimate with $\Theta= \{ u_\Phi\mid \Phi\subset \gD\}.$
\qed

\begin{cor}
\label{c: estimate whittaker coeff with xi}
Assume that $G$ has compact center and let $(\pi, H)$ be an admissible Hilbert $G$-representation of finite length.  Let $\gl \in \Wh_\cchi(H^\infty).$ Then there
exists a $\xi \in \fad$ and a continuous seminorm ${\rm n}$ on $H^\infty$ such that 
for all $v \in H^\infty$ and $a \in A$ we have 
$$ 
|\wh_\gl(v)(a)| \leq a^\xi {\rm n}(v).
$$% 
\end{cor}

\proof 
It follows from (\ref{e: first estimate Whittaker coeff}) that there exists a continuous
seminorm ${\rm n_0}$ on $H^\infty$ 
such that for all $v \in H^\infty$ and $a \in A,$ 
$$ 
|\wh_\gl(v)(a)| \leq e^{r|\log a|} {\rm n}_0(v) . 
$$% 
Since $G$ has compact center, $\fa_\gD = 0,$ so that $\cl(\fa^+)$ 
is a proper closed cone in $\fa.$ Hence, there
exists a linear functional $\xi \in \fad$ such that $\xi > 0$ on $\cl(\fa^+)\setminus \{0\}.$ 
Let $S$ be the unit sphere in $\fad,$ then by compactness of $S \cap \cl(\fa^+)$ we may 
multiply $\xi$ by a positive scalar to arrange that $\xi > r$ on $S \cap \cl(\fa^+).$ 
This implies that  $r |\dotvar| \leq \xi$ 
on $\cl(\fa^+)$. Let now $\Theta \subset U(\fn_0)$ 
be a finite subset as in Lemma \ref{l: basic comparison of estimates}.
Then $L_u w (v) = w( \pi(u) v),$ so that 
$$
a^{-\xi} |\wh_\gl (v)(a)|  \leq \max_{u\in \Theta} e^{- r |\log a|} |\wh_\gl(\pi (u)v)(a)|
\leq 
\max_{u\in \Theta} n_0(\pi (u)v).
$$% 
The required estimate now follows with the continuous seminorm defined by 
$n(v):= \max_{u\in \Theta}n_0(\pi(u)v).$ 
\qed

At a later stage we will also need the following result. We retain the assumption 
that $G$ has compact center.

\begin{lemma}
\label{l: moderate estimate with negative exp}
Let $\mu \in \fad$ be 
such that $\mu(h_\ga)  > 0$ for all $\ga \in \gD.$ Then there exists
a constant $s > 0$ and a finite subset $\Theta \subset U(\fn_0)$ such that
for all $f \in C^\infty(G/N_0\col \cchi)$ and all $a \in A$ we have the estimate
\begin{equation}
\label{e: estimate f from exponent on aplus}
a^{-\mu}  | f(a)| \leq  e^{-s|\log a|} \,\max_{u\in \Theta}  |L_u f(a)|.
\end{equation}
\end{lemma}

\proof Since $\cl(\fa^+)$ is the cone spanned by the elements $h_\ga,$ for $\ga \in
\gD,$ it follows that there exists $s > 0$ such that 
$\mu \geq s$ on $\cl(\fa^+) \cap S,$ where $S$ is the unit sphere
in $\fa.$ This implies that $\mu(H) \geq s|H|$ for all $H \in \fa^+.$ 
The result now follows by application of Lemma \ref{l: basic comparison of estimates}.
\qed

\section{The Whittaker Schwartz space}\label{s: W S space}
We denote the map $G \to \fa$ associated with the Iwasawa decomposition 
$G= K A N_0$ by $H.$ Thus, for $k\in K, a\in A$ and $n_0 \in N_0,$
\begin{equation}
\label{e: Iwasawa H}
H(ka n_0) = \log a.
\end{equation}
Let $\rho \in \fa^*$ be defined by $\rho(H) := \frac12 \tr [\ad(H)|_{\fn_0}].$ 

Following Harish-Chandra \cite[\S 1.3]{HCwhit} and Wallach \cite[\S 15.3.1]{Wrrg2}
we define the Whittaker Schwartz space 
$\cC(G/N_0\col \cchi)$ to be the space of functions $f \in C^\infty(G/N_0\col \cchi)$ 
such that for all $u \in U(\fg)$ and $N > 0$
$$ 
n_{u,N}(f):= \sup_{x \in G} (1 + |H(x)|)^N e^{\rho H(x)}\cdot  | L_u f(x)| < \infty.
$$% 
The indicated seminorms $n_{u,N}$ induce a Fr\'echet topology on $\cC(G/N_0\col \cchi).$ 
It is readily verified that
$
C^\infty_c (G/N_0\col \cchi)\subset \cC(G/N_0\col \cchi) \subset C^\infty(G/N_0\col \cchi),
$
with continuous inclusion maps.

\begin{lemma}
The space $C_c^\infty(G/N_0\col \cchi)$ is dense in $\cC(G/N_0\col \cchi).$ 
\end{lemma}

\proof
For $t >0$ we define 
$$ 
B_t:= \{x \in G\mid |H(x)| \leq t\}
$$% 
Then $B_t$ is right $N_0$-invariant, with compact image in $G/N_0.$ 
Adapting the argument given in \cite[p.~343, Lemma 1]{VarHA} in an obvious fashion,
we infer that there exist left $K$-invariant
functions $\psi_t \in C_c^\infty(G/N_0),$ for $t> 0$ such that 
$0\leq \psi_t \leq 1,$ $\psi_t = 1$ on $B_t,$ 
$\supp \,\psi_t \subset B_{t+1}$ for all $t > 0$ and such that, in addition, 
for every $u \in U(\fg)$ there exists a constant $C_u > 0$ such that 
$$ 
|L_u(\psi_t)(x)| \leq C_u  \text{for {\ }all}\!\! t  > 0,\; x \in G. 
$$% 
Adapting the argument of \cite[p.~343, Thm.\ 2]{VarHA}, again in an obvious way, 
we deduce that
for every $u \in U(\fg),$ $N > 0$ there exists a finite subset 
$V \subset U(\fg)$ such that for all $t \geq 1,$ 
$$ 
n_{u,N}(f - \psi_t f) \leq \sum_{v \in V} (1 + t)^{-1} n_{v, N+1} (f).
$$% 
From this it follows that $\psi_t f \to f$ in $\cC(G/N_0\col \cchi) $ as $t \to 
\infty$. 
\qed

\begin{lemma}
The space $\cC(G/N_0\col \cchi)$ is invariant under left translation by elements
of $G.$ The associated left regular representation $L$ of $G$ on it is continuous.
\end{lemma}

\proof
We start with the observation that for $x\in G/N_0$ and $g \in G$ one has
$$ 
H(g x) = H(g k(x)) + H(x),
$$% 
where $k(x)$ is determined by $x \in k(x)AN_0.$ 
It follows from this that for every compact subset $S \subset G$ and every $N \in \N$ 
there exists a 
constant $C_{S,N} > 0$ such that 
$$ 
e^{\rho H(gx)} (1 + |H(gx)|)^N \leq C_{S, N} e^{\rho H(x)}(1 + |H(x)|)^N.
$$% 
This implies that, for $g \in S$ and $f \in \cC(G/N_0\col \cchi),$
$$
n_{1, N} (L_g f) \leq C_{S,N} n_{1, N+1}(f).
$$% 
Noting that $L_u (L_g f) = L_{\Ad(g^{-1}) u} f$ and observing that $\Ad(S^{-1}) u$ 
is a bounded subset of a finite dimensional subspace of $U(\fg),$ 
we deduce the existence of a finite subset $V \subset U(\fg)$ such that for all
 $g \in S$ and $f \in \cC(G/N_0\col \cchi)$ we have
 $$ 
 n_{u, N} (L_g f) \leq \sum_{v \in V} n_{v, N}(f).
 $$% 
 This implies that $\cC(G/N_0\col \cchi)$ is invariant for the left regular representation 
 and that the set of linear maps $L_g,$ for $g\in S,$ is equicontinuous
  in $\End(\cC(G/N_0\col \cchi)).$ If $f_0 \in C_c^\infty(G/N_0\col \cchi)$ 
 then for $g \to e,$ the function $L_g f_0$ tends to $f_0$ in $C_c^\infty(G/N_0\col \cchi)$
 hence in $\cC(G/N_0\col \cchi).$ Using the density of $C_c^\infty(G/N_0\col \cchi)$ 
 in $\cC(G/N_0\col \cchi)$ it follows by a standard argument that for all $f \in \cC(G/N_0\col \cchi)$ we have 
 $$
 \lim_{g \to e} L_g f = f \text{in} \cC(G/N_0\col \cchi).
 $$% 
Invoking the equicontinuity mentioned above, it now follows by a 
standard argument that the map $(g,f)\mapsto L_g f$ is continuous
$G \times \cC(G/N_0\col \cchi) \to \cC(G/N_0\col \cchi).$
\qed
\begin{lemma}
\label{l: convergence integral over G/N_0}
If $\ell > \dim A$ then
$$ 
\int_{G/N_0} (1 + |H(x)|)^{-\ell} e^{- 2 \rho H(x)} \, dx < \infty.
$$ 
\end{lemma}

\proof By substitution of variables a measurable function $\gf: G/N_0 \to \C$ 
is absolutely integrable if and only if the function $(k, a)\mapsto  \gf(ka) a^{2\rho}$ 
is absolutely integrable over $K \times A.$ If so, the integrals 
$\int_{G/N_0} \gf(x) \, d\dot{x}$ and $\int_{K\times A} \gf(ka) a^{2\rho} \; dk da$ 
are equal, provided the invariant measures are suitably normalized. From this, the proof 
is immediate.
\qed

\begin{cor} $\cC(G/N_0:\chi) \subset L^2(G/N_0:\chi),$ with continuous linear inclusion map.
\end{cor} 

We end this section with a result that will be applied in
the next section.
It is assumed that $G$ has compact center. Then $\gD$ is a linear basis of 
$\fad;$ the associated dual basis of $\fa$ is denoted by $\{h_\ga\mid \ga \in \gD\}.$

\begin{lemma}
\label{l: exp condition for coeff Schwartz}
Suppose that $G$ has compact center and let
$\xi \in \fad$ be such that $\xi(h_\ga ) < -\rho(h_\ga)$ for all $\ga \in \Delta.$ 
Let $(\pi, H)$ be an admissible continuous representation of finite length
of $G$ in a Hilbert space and let $\gl \in \Wh_\cchi(H^\infty).$ Assume there exist
 a continuous
seminorm ${\rm n}$ on $H^\infty$ and a constant $d \in \N$ such that,
for all $v \in H^\infty$ and $a \in A,$ 
$$ 
|\wh_\gl(v)(a) | \leq a^\xi (1 + |\log a|)^d \; {\rm n} (v).
$$%  
Then the Whittaker coefficient map $\wh_\gl$ is continuous 
$H^\infty \to \cC(G/N_0\col \cchi).$ 
\end{lemma}

\proof
We put $\mu = -\rho - \xi.$ Then it follows that $\mu(h_\ga) > 0$ 
for all $\ga \in \gD.$ By Lemma \ref{l: moderate estimate with negative exp} there exists a finite set $\Theta\subset U(\fn_0)$ and  a 
constant $s > 0$ such that for all $f \in C^\infty(G/N_0\col \cchi)$ and all $a \in A$ we have the estimate
$$ 
a^{-\mu} |f(a)| \leq e^{-s|\log a|} \max_{u \in \Theta} |L_u f(a)|.
$$% 
This implies that 
\begin{eqnarray*}
a^\rho |f(a)| 
&= &  a^{-\xi}  a^{-\mu} |f(a)| \\
&\leq & \max_{u \in \Theta}  a^{-\xi} e^{-s|\log a|} |L_u f (a)|\\
&\leq & e^{-s| \log a|} (1 +|\log a|)^{d} \max_{u \in \Theta}  a^{-\xi} (1+|\log a|)^{-d} |L_u f(a)|.\end{eqnarray*}
Using the above estimate for $f = \wh_\gl(v),$ with $v \in H^\infty,$ we find that
\begin{equation}
\label{e: ds estimate 1}
a^\rho |\wh_\gl(v)(a)| \leq a^{-s|\log a|}(1 +|\log a|)^{d} \max_{u\in \Theta} {\rm n}(\pi(u) v).
\end{equation}
For $N \in \N$ we define the positive number
\begin{equation}
\label{e: CN as sup over t} 
C_N = \sup_{t \geq 0} e^{-st} (1 + t)^{N +d}.
\end{equation} 
It follows from (\ref{e: ds estimate 1}) and (\ref{e: CN as sup over t}) that, for all $v \in H^\infty$ and $a\in A,$ 
$$ 
a^\rho (1 + | \log a|)^N |\wh_\gl(v)(a)| \leq C_N \, {\rm n_1}(v),
$$% 
where ${\rm n}_1$ is the continuous seminorm on $H^\infty$ given by 
${\rm n}_1 (v) = \max_{u\in \Theta} (\pi(u) v).$ 
Finally, the last displayed estimate implies that,
for $u \in U(\fg)$ and for all 
$v \in H^\infty$ and $x \in G,$
\begin{eqnarray*}
e^{\rho H(x)} ( 1 + | H(x)|)^N |L_u (\wh_\gl(v))( x)| &=& e^{\rho H(x)} ( 1 + | H(x)|)^N |(\wh_\gl(\pi(u) v))( x)|\\
&\leq & 
 {\rm n}_u (v),
\end{eqnarray*}
 with ${\rm n}_u$ the continuous seminorm on $H^\infty$ 
given by ${\rm n_u}(v) = C_2 \sup_{k \in K} {\rm n}_1( \pi(k)\pi(u)v)).$ 
\qed

\section{Sharp estimates for Whittaker coefficients}
\label{s: sharp estimates of whit coeff}
In this section we assume that $G$ has compact center. 
We shall derive certain growth properties of Whittaker 
coefficients, building on results and ideas of Wallach \cite{Wrrg2} and Harish-Chandra
\cite{HCwhit}. 

We assume that $(H, \pi)$ is an admissible continuous representation 
of finite length of $G$ in a Hilbert space. Let $V = H_K$ be the associated Harish-Chandra module, and let $V^\sim$ denote  the associated dual Harish-Chandra module.
Then it is well known that the natural map $(H^\infty)'  \to V^*$ induces an isomorphism 
$(H^\infty)'_K \simeq V^\sim.$ 

We agree to write $E(P_0, V )$ for the set of generalized weights of the 
finite dimensional $\fa$-module $V/\fn_0 V.$ 

Let $\gD$ be the collection of simple roots in $\gS(\fn_0, \fa)$ and let  
$\{h_\ga \mid \ga\in \gD\}$ be the associated dual basis of $\fa.$ 
We define $\gL_V \in \fad$ by 
\begin{equation}
\label{e: defi gL V}
\gL_V(h_\ga) := \max\{- \Re \mu (h_\ga) \mid \mu \in E(P_0, V)\}.
\end{equation} 

\begin{rem}
\label{r: ds and gL}
At a later stage it will be of crucial importance to us that for  $\pi$ irreducible unitary, the following two conditions are equivalent
\begin{enumerate}
\itema $\pi$ is equivalent to a direct sum of representations from the discrete series of $G;$ 
\itemb $\gL_V(h_\ga) < -\rho(h_\ga)$ for all $\ga \in \gD.$ 
\end{enumerate}
For a proof of this well known result, we refer to \cite{Wrrg1} as follows.
The element $\gL_V$ corresponds to $\gL_{V^\sim}$ 
as defined in \cite[\S 4.3.5]{Wrrg1}. In the terminology of \cite[\S 5.1.1]{Wrrg1} assertion 
(b) means that $V^\sim$ is rapidly decreasing. According to \cite[Prop.\ 5.1.3 \&\ Thm. 5.5.4]{Wrrg1} 
this is equivalent to the assertion that $(\pi^\veec, H')$ is a direct sum of square integrable
representations, which in turn is equivalent to assertion (a).
\end{rem}

The following result is due to  Wallach \cite[15.2.2]{Wrrg2}. The regularity
of $\chi$ is not required.

\begin{prop}
\label{p: first estimate whittaker coefficient}
Put $\gL:= \gL_V$ and let $\gl \in \Wh_\cchi(H^\infty).$ 
Then there exists a constant $d \geq 0$ 
and for every $v \in H_K$ a constant $C_v > 0$ such that 
$$ 
|\gl(\pi(a^{-1}) v) |\leq C_v ( 1+ |\log a|)^d \,a^\gL.
$$% 
for all $a \in \cl(A^+).$ 
\end{prop}

Wallach's proof of this result follows the lines of the proof of an earlier result,
stated in  \cite[Thm. 4.3.5]{Wrrg1}. That result, applied to the contragredient
 representation
$\pi^\veec,$  asserts that for a given $v \in ((H')^\infty)'_K\simeq H_K$ there
exists a constant $d > 0,$ a continuous seminorm $\gs_v $ on $(H')^\infty$ such that for all 
$\gl \in (H')^\infty$ one has 
$$ 
|v(\pi^\veec(a)\gl)| \leq ( 1 + |\log a|)^d a^\gL  \gs_v(\gl).
$$ 
As $v(\pi^\veec(a)\gl) = \gl(\pi(a^{-1}) v),$ the result \cite[Thm. 4.3.5]{Wrrg1} implies
that Proposition \ref{p: first estimate whittaker coefficient} is valid for $\gl$ 
a smooth vector in $H'.$ 

The proof of \cite[Thm. 4.3.5]{Wrrg1} makes use of initial estimates and 
of estimate improvement
through asymptotic behavior along maximal standard parabolic subgroups $P_\Phi,$ 
with $\Phi = \gD \setminus \{\ga\}.$ It exploits a system of differential
equations coming from the observation that $V /\fn_\Phi V$ is an admissible
$(\fm_{1\Phi}, K_\Phi)$-module  ($\fm_{1\Phi} = \fm_\Phi +\fa_\Phi$), so that $\fa_\Phi$ acts finitely on it. 
From the proof one sees that only estimates of 
$[\pi^\veec(U)\gl](\pi(a) w) $ for $U \in U(\fn_0),$ $w \in H_K$ and $a \in A^+$ 
are needed to make the approach work. This is the condition of $(P_0, A)$-tameness of \cite[\S 15.2.1]{Wrrg2}.
 If $\gl \in \Wh_\chi(H^\infty)$ then 
$\pi^\veec(U) \gl = \chi_*(U^\veec)\gl$ so that the needed tameness is trivially guaranteed.
Therefore, essentially the same approach gives the validity 
of Proposition \ref{p: first estimate whittaker coefficient}. For further details, see the proof of 
\cite[Thm.~15.2.2]{Wrrg2}, assertion (1). 

For regular $\chi$ the following refinement  of Proposition 
\ref{p: first estimate whittaker coefficient} will be of crucial importance to us, since it gives an estimate for the Whittaker coefficient $\wh_\gl(v)$ for every smooth $v\in H^\infty$ on all Weyl chambers of $A.$ 
The key idea is to now focus on the
$\fa_\Phi$-actions on the modules $U(\fg)\gl/ \bar \fn_\Phi^k U(\fg)\gl,$ for 
$k\geq 1,$ making use of the information provided by 
Proposition \ref{p: first estimate whittaker coefficient} to exclude the contribution 
of $\fa_\Phi$-weights that are not dominated by $\gL|_{\fa_\Phi}.$  

\begin{thm}
\label{t: second estimate whittaker coefficient}
Suppose that $G$ has compact center and let $\gL:= \gL_V.$ 
Assume that $\cchi$ is a regular unitary character of $N_0$ and that 
$\gl \in \Wh_\cchi(H^\infty).$  
Then there exists a constant $d \geq 0$ and a continuous seminorm ${\rm n}$ on 
$H^\infty$ such that for every $v \in H^\infty$, 
$$ 
|\gl(\pi(a^{-1}) v) |\leq ( 1+ |\log a|)^d \,a^\gL \, {\rm n}(v),
$$% 
for all $a \in A.$ 
\end{thm}

\begin{rem} \cite[Thm.~15.2.5]{Wrrg2} gives this estimate for  $a \in \cl(A^+).$ 
\end{rem}

The rest of this section is devoted to the proof of Theorem
\ref{t: second estimate whittaker coefficient}. We start with the fact that 
there exists a $\xi \in \fa^*$,  a constant $d \geq 0$ 
and a continuous seminorm ${\rm n}$ on $H^\infty$ 
such that 
\begin{equation} 
\label{e: dominating estimate with xi}
|\wh_\gl(v)(a)| \leq a^\xi (1 + |\log a |)^d \, {\rm n}(v)\qquad (v \in H^\infty, \;a \in A).
\end{equation}
Indeed, according to  Corollary \ref{c: estimate whittaker coeff with xi}
this estimate holds  with $d = 0$ for a suitable choice of $\xi.$ 

In the following we shall say that $\xi \in \fad$ dominates the Whittaker coefficient 
$\wh_{\gl}$ if there exist a $d \in \N$ and a continuous seminorm ${\rm n}$ on 
$H^\infty$ such that the estimate (\ref{e: dominating estimate with xi}) is valid 
for all $v \in H^\infty$ and $a \in A.$ The following result concerning domination
will be useful.

\begin{lemma}
\label{l: domination on fap}
Let $\vartheta, \xi  \in \fad$ be such that $\vartheta \leq \xi $ on $\fa^+.$ 
If $\vartheta $ dominates the Whittaker coefficient $\wh_\gl,$ then so does $\xi.$ 
\end{lemma}

\proof 
Let $\Theta$ be a finite subset of $U(\fn_0)$ with the properties guaranteed by 
Lemma \ref{l: basic comparison of estimates}.  
Assume $\vartheta$ is dominating. Then there exist a constant $d\geq 0$ and a continuous seminorm $\rm n$ on $H^\infty$ 
such that for all $v \in H^\infty$ and all $a\in A$ we have
$$ 
a^{-\vartheta} |\wh_\gl(v)(a)| \leq (1 + |\log a|)^d {\rm n}(v).
$$% 
By applying Lemma \ref{l: basic comparison of estimates}  
with $f = \wh_\gl(v),$ and using that 
$L_u(\wh_\gl(v)) = \wh_\gl(\pi_*(u) v)$ we infer that, for all $v \in H^\infty$ and
$a \in A,$  
\begin{eqnarray*}
a^{-\xi} |\wh_\gl(v)(a)|&  \leq & \max_{u \in \Theta} a^{-\vartheta} |\wh_\gl(\pi_*(u)v)(a)|\\
&\leq &
(1+ |\log a|)^d \max_{u \in \Theta}  {\rm n}(\pi_*(u) v).
\end{eqnarray*}
Since ${\rm n}': v \mapsto \max_{u \in \Theta} {\rm n}(\pi_*(u) v)$ is a continuous
seminorm on $H^\infty$ it follows that $\xi$ dominates $\wh_\gl.$ 
\qed

The idea is now to show that a dominating $\xi$ can be improved by using asymptotic 
expansions along maximal standard parabolic subgroups derived
from suitable differential equations. This method, inspired by \cite[\S 4.4]{Wrrg1} and \cite[\S 15.2]{Wrrg2}, leads to the following lemma, which 
is the main step in our proof of Theorem \ref{t: second estimate whittaker coefficient}.
We write $\gL= \gL_{H_K}.$ 

\begin{lemma}{\rm (estimate improvement)\ }
\label{l: estimate imp}
Assume that the Whittaker coefficient $\wh_\gl$ is dominated by $\xi \in \fad.$
Let $\ga \in \gD.$  
Then $\wh_\gl$ is also dominated by $\xi',$ where $\xi'\in \fad$ is defined 
by 
\begin{enumerate}
\itema $\xi' =\xi$ on $\ker \ga;$ 
\itemb  $\xi'(h_\ga) = \min(\xi(h_\ga), \gL(h_\ga)).$
\end{enumerate}
\end{lemma}
The proof of this lemma will be given in the next section.
Here we note that  
using the lemma successively for all simple roots $\ga \in \gD$ we obtain 
the following corollary. 

\begin{cor}{\rm (estimate improvement)\ } \label{c: est improvement}
Assume that the Whittaker coefficient $\wh_\gl$ is dominated by $\xi \in \fad.$
Then $\wh_\gl$ is also dominated by $\xi'',$ where $\xi'' \in \fad$ is defined by 
$\xi''(h_\ga) = \min (\xi(h_\ga), \gL(h_\ga))$ for every $ \ga \in\gD.$ In particular,
$\xi''(h_\ga) \leq \gL(h_\ga)$ for all $\ga \in \gD.$ 
\end{cor}
\medno
The above results allow us to establish the main result of this section. 
\medno
{\em Proof of Theorem \ref{t: second estimate whittaker coefficient}.\ }
According to Corollary \ref{c: estimate whittaker coeff with xi} there exists
a $\xi\in \fad$ which dominates $\wh_\gl.$ By Corollary \ref{c: est improvement}
the coefficient $\wh_\gl$ is also dominated by $\xi''.$ From the definition of $\xi''$ one 
sees that $\xi''\leq \gL$ on $\fa^+.$ By application of Lemma \ref{l: domination on fap}
it now follows that $\gL$ dominates $\wh_\gl.$ This
establishes the assertion of Theorem \ref{t: second estimate whittaker coefficient}.
\qed

\begin{cor}
\label{c: schwartz estimate for ds coeff}
Suppose that $(\pi, H)$ belongs to the discrete series of $G$ and let 
$\gl \in \Wh_\cchi(H^\infty).$ Then the associated Whittaker coefficient
$\wh_\gl$ defines a continuous linear map $H^\infty \to \cC(G/N_0\col \cchi).$ 
\end{cor}

\proof 
From Remark \ref{r: ds and gL} it follows that $\gL(h_\ga) < -\rho(h_\ga)$ 
for $\ga \in \gD.$ The result now follows by combining Theorem 
\ref{t: second estimate whittaker coefficient}
with Lemma \ref{l: exp condition for coeff Schwartz}.
\qed

\section{Proof of Lemma \ref{l: estimate imp}: improvement of estimates}
\label{s: proof of est imp}
As in the previous section, we assume that $G$ has compact center.
Furthermore, $(\pi, H)$ is an admissible continuous representation of finite length of $G$ in 
a Hilbert space, and $V = H_K$ is the associated Harish-Chandra module.
We define $\gL = \gL_V$ as in (\ref{e: defi gL V}) and assume that 
$\gl \in \Wh(H^\infty).$ Let  $\wh_\gl: H^\infty \to C^\infty(G/N_0\col\cchi)$ be 
the associated Whittaker coefficient. The purpose of this section is to
prove Lemma  \ref{l: estimate imp}.

Our starting assumption is that $\wh_\gl$ is dominated by $\xi \in \fad.$ This means
that there exists a constant $d \in \N$ and a continuous seminorm ${\rm n}$ 
on $H^\infty$ such that for all $v \in H^\infty$ and all $a\in A$ we have
\begin{equation}
\label{e: initial domination}
|m_\gl(v)(a)| \leq ( 1 + |\log a|)^d a^\xi {\rm n}(v).
\end{equation}

We will improve upon this estimate by using a system of differential equations
satisfied by $\wh_\gl.$ Our first goal is to set up this system.

Given a finite dimensional real linear space $\fv$ we denote by $S(\fv)$ the symmetric algebra
of its complexification $\fv_{\iC}$ and by $P(\fv)$ the algebra of polynomial
functions $\fv_\iC \to \C.$

Given a real Lie algebra $\fl$ we denote by $U(\fl)$ the universal enveloping algebra of
its complexification $\fl_\iC$, and by $\fZ(\fl)$ the center of $U(\fl).$ Furthermore, 
$U(\fl)$ is equipped with the standard filtration by order, $(U(\fl)_n)_{n \geq 0}.$ 
The center $\fZ(\fl)$ is equipped with the induced filtration. 

For a parabolic subgroup $P$ of $G$ we denote its Langlands decomposition
by $P = M_P A_P N_P$ and we  write $M_{1P}:= M_P A_P.$  
We agree to use the abbreviated notation $\fZ = \fZ(\fg),$ $\fZ_{1P}= \fZ(\fm_{1P})$ 
and $\fZ_P:= \fZ( \fm_P) .$

Given a subset $\Phi \subset \gD$ 
we denote by $P_\Phi$ the associated standard parabolic subgroup of $G.$ 
Its Langlands components  are denoted by $M_\Phi, A_\Phi, N_\Phi.$ Furthermore, $M_{1\Phi} := M_\Phi A_\Phi,$ $\fZ_{1\Phi}:= \fZ(\fm_{1\Phi})$ and $\fZ_\Phi:= \fZ(\fm_\Phi).$ 

We consider 
the $U(\fg)$-submodule 
$$ 
\cY:= U(\fg)\gl
$$% 
of $(H^\infty)'.$ 
Since $H$ is admissible and of finite length, $\fZ \gl$ is a finite
dimensional subspace of $\cY.$ 

We fix $\Phi \subset \gD.$  
By the PBW theorem, $U(\fg) = U(\fm_{1\Phi}) \oplus (\bar \fn_\Phi U(\fg) + U(\fg) \fn_\Phi).$ 
The associated projection $U(\fg) \to U(\fm_{1\Phi}),$ restricted to $\fZ,$ defines an 
algebra homomorphism
\begin{equation}
 \label{e: defi homomorphism p}
p: \fZ \to \fZ_{1\Phi}.
\end{equation} 
It is well known that $p$ is injective and preserves the filtrations
induced by the standard filtration $(U(\fg)_n)_{n \in \N}$ by order on $U(\fg).$ 

We fix a maximal torus $\ft \subset \fm;$ then $\fh:= \ft \oplus \fa$ is a Cartan subalgebra
of $\fg. $ Let $W(\fh)$ denote the Weyl group of the root system $R(\fg, \fh)$ 
of $\fh_\iC$ in $\fg_\iC.$ 
Furthermore, let $W_\Phi(\fh)$ denote the Weyl group of $R(\fm_{1\Phi}, \fh).$ 
Then $W_\Phi(\fh)$ equals the centralizer of $\fa_\Phi$ in $W(\fh).$ 

We denote by $S(\fh)^{W(\fh)}$ and $S(\fh)^{W_\Phi(\fh)}$ the associated
subalgebras of Weyl group invariants in $S(\fh).$ Then it is well known that
$S(\fh)^{W_\Phi(\fh)}$ 
is a free $S(\fh)^{W(\fh)}$-module of rank $\ell = [W(\fh): W_\Phi(\fh)]$ 
with free homogeneous generators, $q_1 =1, q_2, \ldots , q_\ell.,$ 
see \cite[Thm 2.1.3.6]{Warn1}. 

Let $\gg: \fZ \to S(\fh)^{W(\fh)}$ and $\gg_\Phi: \fZ_{1\Phi} \to S(\fh)^{W_\Phi(\fh)}$ 
denote the associated Harish-Chandra isomorphisms. These are known to
be isomorphisms of filtered algebras. Furthermore,
$$ 
p = T_{-\rho_\Phi} \after \gg_\Phi^{-1}\after  \gg,
$$ 
where $T_{-\rho_\Phi}$ is the automorphism $T$ of $\fZ_{1\Phi}\simeq \fZ_\Phi \otimes S(\fa_\Phi)$ determined by $T = I $ on $\fZ_\Phi$ and by $T(X)= X - \rho_\Phi(X),$ for $X \in \fa_\Phi.$ 
Here $\rho_\Phi \in \fa^*$ is given by $\rho_\Phi(X) = \frac12 \tr (\ad(X)|_{\fn_\Phi}).$ 
For $1\leq i \leq \ell,$ put  
$$
u_i := T_{-\rho_\Phi}\gg_\Phi^{-1}(q_i) \in \fZ_{1\Phi}.
$$
Then we see that $\fZ_{1\Phi}$ is a free $p(\fZ)$-module with basis $u_1=1 , u_2, \ldots, u_\ell.$
Furthermore, since the $q_i$ are homogeneous and since $\gg$, $\gg_\Phi$ and $T_{-\rho_\Phi}$ are isomorphisms
of filtered algebras, the following is valid.
\medno\sl
For every $u \in \fZ_{1\Phi} \cap U(\fg)_n,$ $(n \in \N),$ let $Z_1, \ldots , Z_\ell\in \fZ$ 
be the unique elements such that $u = \sum_{i=1}^\ell u_i p(Z_i).$ Then 
for all $1\leq i \leq \ell,$
\begin{equation}
\label{e: estimation orders}
\ord(u_i) + \ord(Z_i) \leq n 
\end{equation}\rm
 
Let $E_\Phi $ be the complex linear
span of the elements $\{u_i\mid 1\leq i \leq \ell\}$ in $\fZ_{1\Phi}.$ Then 
the map $(u, Z)\mapsto  u p(Z)$ induces a linear isomorphism
$E_\Phi\otimes \fZ \simeq \fZ_{1\Phi}. $
\begin{lemma}
With $E_\Phi \subset \fZ_{1\Phi}$ as above, we have   
\begin{equation}
\label{e: deco of Ug with Z and E}
U(\fg) = U(\bfn_\Phi) U(\fm_\Phi) E_\Phi \fZ U(\fn_\Phi).
\end{equation} 
\end{lemma}

\proof 
By induction on $n \in \N$ we will show that $U(\fg)_n$ is contained in 
the space on the right-hand side of (\ref{e: deco of Ug with Z and E}). 
For $n=0$ the inclusion is obvious. Thus, let $n \geq 1$ and assume the inclusion
has been established for strictly smaller values of $n.$   
We observe that by the PBW theorem, every element of $U(\fg)_n$ may be written
as a sum of an element of $U(\fg)_{n-1}$ and a finite sum of products 
$w v u y $ with $w \in U(\bfn_\Phi)$, $v \in U(\fm_\Phi),$ $ u \in  U(\fa_\Phi)$ 
and $y \in U(\fn_\Phi)$ such that 
\begin{equation}
\label{e: sum orders at most n}
\ord (w) + \ord (v) + \ord (u) + \ord (y) \leq n.
\end{equation}
In view of the induction hypothesis, it suffices to 
show that each such product $wvuy$ with (\ref{e: sum orders at most n})
belongs to the space on the right of 
(\ref{e: deco of Ug with Z and E}). 

Since  $U(\fa_{\Phi})\subset \fZ_{1\Phi},$ the element $u$ may be
expressed as a sum of elements $u_i p(Z_i)$ with $\ord(u_i) + \ord(Z_i) \leq
\ord (u).$ Now $p(Z_i) - Z_i \in U(\fg)_{n_i -1}\fn_{\Phi},$ where 
$n_i: = \deg Z_i \leq n.$ It follows that for every $1 \leq i \leq \ell,$ 
$$
wvu_i (p(Z_i) - Z_i )y \in U(\fg)_{n-1}\fn_\Phi.
$$% 
Summing over $i$ and applying the induction hypothesis, we find
$$ 
wvu y \in w v \sum_{i} u_i Z_i y + U(\fg)_{n -1} \fn_\Phi \subset U(\bfn_\Phi) U(\fm_\Phi) E_\Phi \fZ U(\fn_\Phi). 
$$% 
\qed

For $k \geq 1,$ the quotient $\cM_k:= \cY/\bfn_\Phi^k \cY$ is a 
left $U(\fm_{1\Phi})$-module. 

\begin{lemma}
The subspace $E_\Phi \fZ \gl$ of $\cY$ has a finite dimensional 
$\fa_\Phi$-invariant image in $\cM_1,$ which generates the $U(\fm_{\Phi})$-module $\cM_1.$ 
\end{lemma}

\proof
Since $E_\Phi$ and $\fZ \gl$ are finite dimensional, the mentioned image
$F_1$ of $E_\Phi \fZ \gl$ in $\cM_1$ is finite dimensional. From 
$$ 
U(\fa_\Phi)E_\Phi  \subset \fZ_{1\Phi } \subset  E_\Phi p(\fZ)  \subset E_\Phi \fZ + \bar \fn_\Phi U(\fg)
$$% 
it follows that  $F_1$ is finite dimensional and $\fa_\Phi$-invariant. 
\qed

In particular, it follows that the set $S_\Phi$ of generalized $\fa_\Phi$-weights
of $\cM_1$ is finite and that $\cM_1$ is the direct sum of the associated
generalized $\fa_\Phi$-weight spaces.

\begin{lemma}
For $k\geq 1$ the set $\wt(\cM_k)$ of generalized $\fa_\Phi$-weights of 
$\cM_k$ is finite. Each of its elements is of the form  
$\gs - (\ga_1 + \cdots + \ga_l)$, with $\gs \in S_\Phi$, $0 \leq l < k$ 
and $\ga_j \in \gS(\fn_\Phi, \fa_\Phi)$ for all $1 \leq j \leq l.$ 
\end{lemma}

\proof 
For $k=1$ the result has been established above. For  $k\geq 2$ we notice
that the natural map $q_k: \cM_{k} \to \cM_{k-1}$ is a surjective morphism of 
$\fa_\Phi$-modules. Furthermore, its kernel is isomorphic as an $\fa_\Phi$-module 
to the quotient
$$ 
Q_k := \bfn^{k-1}_\Phi U(\fg) \gl / \bfn_\Phi^{k}U(\fg) \gl. 
$$% 
Define $Q_1: =  \cM_1,$ then $\wt(Q_1) = S_\Phi$ and $Q_1$ is the associated
direct sum of generalized $\fa_\Phi$-weight spaces. 
The natural map $\bar\fn_\Phi \otimes Q_{k-1} \to Q_k$ is surjective, for $k \geq 1.$ 
Thus, if the $\fa_\Phi$-module $Q_{k-1}$ is the direct sum of its generalized weight spaces, then so is $Q_k,$ and $\wt(Q_k) \subset \wt(\bfn_\Phi) + \wt(Q_{k-1}).$ 
It follows by induction on $k$ that each $Q_k$ is the direct sum of its
generalized weight spaces, and that for $k \geq 1$ each weight of $Q_k$ is 
of the form $\gs - (\ga_1+  \cdots + \ga_\ell),$ where $\gs \in S_\Phi,$ $\ell < k$ and $\ga_j \in \gS(\fn_\Phi, \fa_\Phi),$ for $1 \leq j \leq \ell.$ 

For $k \geq 2$ we now have the short exact sequence of $\fa_\Phi$-modules
$$ 
0 \to Q_k \to \cM_k \to \cM_{k-1} \to 0.
$$%
Here $Q_k$ is the direct sum of finitely many generalized $\fa_\Phi$-weight spaces. 
If $\cM_{k-1}$ is the direct sum of finitely many weight spaces, then it follows
that the module $\cM_k$ is the direct sum of finitely many weight spaces as well,
while $\wt(\cM_k) \subset \wt (Q_k ) \cup \wt(\cM_{k-1}).$ The asserted result
now follows by induction on $k.$ 
\qed

After this preparation, we proceed with setting up the system of differential equations. 
Fix $k \geq 1.$ At a later stage we shall impose a condition on the magnitude of $k.$ 
The space $U(\fa_\Phi)\gl$ maps to a finite dimensonal
subspace $F$ of $\cM_k= \cY/\bfn_\Phi^k \cY.$ We fix elements $u_1=1, u_2, \ldots , u_p$ 
of $U(\fa_\Phi)$ such that the images $[u_j \gl]$ in $\cM_k$ form a basis 
of $F.$

For $H \in \fa_\Phi$ we denote by $B(H)$ the transposed of the 
matrix of the action of $H$ 
on $F$ relative to this basis. Then there exist linear maps
$y_j: \fa_\Phi \to  \bfn_\Phi^k U(\bfn_0 +\fm_{1}),$ 
($\fm_1 := \fm + \fa$) such that for every $H \in \fa_\Phi$ we have
the following identities in $\cY = U(\fg) \gl,$ 
$$ 
H u_j \gl = \sum_{i=1}^p B(H)_{ji} u_i \gl + y_j(H)\gl, \qquad  (1 \leq j \leq p).
$$% 
We now define the functions $F: H^\infty \times A \to \C^p$ and $R: V \times 
A \times \fa_{\Phi} \to \C^p$ 
by 
$$
F_j(\gf, a) = u_j \gl( \pi(a)^{-1}v  ), \quad R_j(v, a, H) = y_j(H)\gl (\pi(a)^{-1}v), 
\qquad (1 \leq j \leq p).
$$% 
We use the decomposition (\ref{e: bp deco of fa}) and put $\bp A_\Phi :=\exp (\bp \fa_\Phi),$ 
so that $A \simeq \bp A_\Phi A_\Phi.$ 
Then,
$$ 
\frac{d}{dt} F(v, \bp a \exp ( tH)) =  B(H) F(v, \bp a \exp tH) + 
R(v, \bp a \exp tH, H),
$$% 
for all $\bp a \in \bp A_{\Phi},$ $H \in \fa_\Phi$ and $t \in \R.$ This equation in turn leads
to 
$$ 
\frac{d}{dt} e^{-tB(H)} F(v, \bp a \exp (tH) ) =  e^{-tB(H)} R(v, \bp a \exp tH, H).
$$% 
Finally, by integrating with respect to $t$ we obtain the equivalent equation
\begin{equation}
\label{e: integral expression F}
F(v, \bp a \exp (tH))  = e^{tB(H)} F(v, \bp a) + e^{tB(H)} \int_0^t e^{-\tau B(H)} R(v, \bp a \exp (\tau H) , H) d \tau. 
\end{equation}
The domination estimate (\ref{e: initial domination}) leads
to estimates of $F$ and $R.$ 

\begin{lemma}
There exists a constant $d \in \N$ and a continuous seminorm ${\rm n}$ on $H^\infty$ such that for all $v \in H^\infty$ and
all $a \in A,$ 
\begin{equation}
\label{e: initial estimate F}
\|F(v, a)\| \leq a^\xi (1 + |\log a|)^d \; { \rm n}(v).
\end{equation}
\end{lemma}
\proof 
Using that $F_j(v,a ) = u_j \gl (\pi(a)^{-1} v) = \wh_\gl(\pi(u_j^\veec) v)(a)$ 
in combination with the estimate (\ref{e: initial domination}) we find that
the estimate (\ref{e: initial estimate F}) is valid with a seminorm of the 
form ${\rm n}(v) = C \max_j {\rm n}_0(\pi(u_j^\veec)v),$ where 
 ${\rm n}_0$ denotes the seminorm of (\ref{e: initial domination}).
\qed
In the following we will use the abbreviated notation
$$ 
|(\bp a, H)| := (1 + |\log \bp a|)(1 + |H|) , \qquad (\bp a \in \bp A_\Phi, H \in \fa_\Phi).
$$%
\begin{lemma}
\label{l: initial estimate R}
Let $\xi \in \fa^*$ be as in (\ref{e: initial domination}).
Then
there exist $N \in \N$ and a continuous seminorm ${\rm n}$ on $H^\infty$ 
such that for all $v \in H^\infty,$ $\bp a \in \bp A_\Phi$ and $H \in \fa_\Phi^+,$ $t \geq 0,$ 
\begin{equation}
\label{e: first estimation R}
\|R(v, \bp a \exp tH, H)\| \leq  |(\bp a, tH)|^N (1 + |H|)(\bp a)^\xi e^{(\xi -k\gb_\Phi)(tH)}
\, {\rm n }(v).
\end{equation}
Here $\gb_\Phi: \fa_\Phi^+ \to [0,\infty[$ is defined by 
$
\gb_\Phi(H) = \min_{\ga \in \gD\setminus \Phi} \ga(H).
$ 
\end{lemma}

\proof
Let $y \in \bar \fn_\Phi^k U(\bar \fn_0 + \fm_1).$ Then we may express $y$ as a finite 
sum of weight vectors $U$ for the adjoint action of $\fa.$ The weight associated with
such a term $U$ is of the form $-\mu = - \sum_{\ga \in \gD} \mu_\ga\, \ga$ with  $\mu_\ga \in \N$ 
and with $\sum_{\ga \in \gD\setminus \Phi} \mu_\ga \geq k.$ It follows
that 
$$
U \gl (\pi(\bp a \exp tH)^{-1} v) = (\bp a\exp tH)^{-\mu} \gl (\pi(\bp a \exp tH)^{-1} U^\veec v).
$$% 
Now there exists a finite sequence $\ga_1, \dots, \ga_q$ of simple roots
from $\Phi$ such that $\mu = \ga_1 + \cdots + \ga_q$ on $\bp \fa_\Phi.$ 
Let $X_{\ga_j} \in \fg_{\ga_j}$ be such that $\cchi_*(X_{\ga_j}) = 1$ 
for each $j$ and put $X = X_{\ga_1}\cdots X_{\ga_q},$ then $\cchi_*(X) = 1.$ 
Hence,  $X \gl = \gl$ and  it follows that
\begin{eqnarray*}
U \gl (\pi(\bp a \exp tH)^{-1} v) & = &  U X \gl (\pi(\bp a \exp tH)^{-1} v)\\
&=&  
(\exp tH)^{-\mu} \gl (\pi(\bp a \exp tH)^{-1} X^\veec U^\veec v).
\end{eqnarray*} 
Now  
$$
-\mu(H)  \leq -k \gb_\Phi(H), \qquad (H \in \fa_\Phi^+),
$$% 
and we see, by using the initial estimate (\ref{e: initial domination}) that there exists a continuous seminorm $n_0$ on $H^\infty$ 
such that for all $v \in H^\infty,$ $\bp a\in \bp A_\Phi,$ $H \in \fa_\Phi^+, t\geq 0$ 
\begin{eqnarray*}
|U \gl (\pi(\bp a \exp tH)^{-1} v)| 
& \leq & e^{-k\gb_\Phi( t H)} |\wh_\gl(\pi(\bp a \exp tH) X^\veec U^\veec v )|\\
& \leq & 
e^{-k\gb_\Phi( t H)} (\bp a \exp tH)^\xi |(\bp a, tH)|^d n_0(v).
\end{eqnarray*}
It follows that a similar estimate holds for each $y \in \bfn_\Phi^kU(\bfn_0 + \fm_1).$ 
Using the linear dependence of  $y_i(H)$ on $H$ for $1 \leq i \leq p$ 
and using the above estimates with $y_i(H)$ in place of $y,$ 
we infer the existence of a continuous seminorm ${\rm n}$ on $H^\infty$ such that
the asserted estimation (\ref{e: first estimation R}) is valid.
\qed

We now assume that  $\ga \in \gD$ and $\Phi = \gD\setminus \{\ga\},$ 
so that $P_\Phi$ is a maximal standard parabolic subgroup. Clearly,
$\fa_\Phi = \R h_\ga,$ where $h_\ga$ is defined as in the beginning 
of Section \ref{s: sharp estimates of whit coeff}.

We may and will assume $\gL(h_\ga) < \xi(h_\ga)$ 
since otherwise $\xi' = \xi$ and there would be nothing to prove.
Accordingly,
\begin{equation}\label{e: xi prime on h ga}
\xi'(h_\ga) = \gL(h_\ga),\qquad  \xi'|_{\bp \fa_\Phi}  = \xi|_{\bp \fa_\Phi}.
\end{equation} 
We now observe that $\gb_\Phi(h_\ga) = 1$ and 
impose the mentioned condition on $k \in \N$  that 
\begin{equation}
\label{e: condition on k}
(\xi  - k \gb_\Phi)(h_\ga) = \xi(h_\ga) - k  < \gL(h_\ga).
\end{equation}
 
The spectrum of $B(h_\ga)$ is contained in the set
$$
X : = \{ [\gs - (\ga_1 +  \ldots + \ga_\ell)](h_\ga)\mid 0\leq \ell < k, \;\forall j: \ga_j \in \gS(\fn_\Phi, \fa_\Phi) \}.
$$% 
For $x \in X$ we denote by $P_x: \C^p \to \C^p$ the projection 
onto the associated generalized weight space of $B(h_\ga),$ along the generalized
weight spaces for the eigenvalues different from $x.$ 
If $x$ is not an eigenvalue for $B(h_\ga),$ 
then $P_x = 0.$ Then it is clear from (\ref{e: initial estimate F}), 
possibly after adaptation of $d$ and ${\rm n}$ 
that, 
\begin{equation}
\label{e: estimate Px F}
\|P_{\!x} F(v, a)\| \leq a^\xi ( 1+ |\log a|)^d \;{\rm n}(v),
\end{equation} 
for all $v \in H^\infty$ and $a \in A$. 
Likewise, for the components $P_{\! x} R$ we have estimates 
of the form (\ref{e: first estimation R}), with adapted $N'$ and $n'$ if necessary. 

In the following we agree to write
$$ 
a = \bp a a_t, \;\;\;\bp a \in \bp A_\Phi\quad  {\rm and} \quad \; a_t = \exp (th_\ga), \;\;\;t \in \R.
$$%

\begin{rem}
\label{r: restricted domain of estimate}
Clearly, to finish the proof of Lemma \ref{l: estimate imp} it suffices to prove the estimate (\ref{e: estimate Px F}) with $\xi$ replaced by $\xi'.$  Now $\xi'(h_\ga) \leq \xi(h_\ga)$ implies that for $\bp a \in \bp A_\Phi$ and $t \leq 0$ we have
$$ 
a^\xi = (\bp a)^\xi (a_t)^{\xi} \leq (\bp a)^\xi (a_t)^{\xi'}= a^{\xi'},
$$% 
so that the required estimate (\ref{e: estimate Px F}) with $\xi'$ in place of $\xi$ 
is automatically fulfilled for $a$ outside $\bp A_\Phi  A_\Phi^+.$ It therefore
suffices to prove the estimate for $a = \bp a a_t,$ $t > 0.$ 
\end{rem}

It is well known that there exist unique polynomial maps 
$Q_x: \R \to \End (\C^p),$ for $x \in X,$ such that $Q_x(t)$ commutes with 
$B(h_\ga)$ for all $t \in \R,$ hence with all the weight space projections
for $B(h_\ga),$ and satisfies $P_x Q_x(t) = Q_x(t) = Q_x(t) P_x$ and 
\begin{equation}
e^{t B(h_\ga)} P_x = e^{tx} Q_x(t), \qquad (t \in \R).
\end{equation}
Furthermore, the polynomial degrees of $Q_x$ are at  most $p,$ so that
there exists a suitable constant $C_0 > 0$ such that, for all $x \in  X ,$ 
\begin{equation}
\label{e: estimate Q x}
\|Q_x(t)\|_{\rm op} \leq C_0 (1 + |t|)^p , \qquad (t \in \R). 
\end{equation}

We recall that  $\gb_\Phi(h_\ga) = 1$ and that   
Lemma \ref{l: initial estimate R} with $H = h_\ga$ implies the 
existence of a continuous seminorm 
${\rm n}$ on $H^\infty$ such that for all  $v \in H^\infty,$ $\bp a\in \bp A_\Phi$ and $t \geq 0$ we have 
\begin{equation}
\label{e: estimate R along max par}
\|R(v,\bp aa_t , h_\ga) \| \leq |(\bp a , t h_\ga)|^N (\bp a)^\xi e^{t[\xi(h_\ga) - k]} {\rm n}(v).
\end{equation}
Writing 
$$ 
R(v, \bp a\exp t h_\ga, h_\ga) = e^{t [\xi (h_\ga) - k]} R_0(v, \bp a \exp(t h_\ga))
$$% 
we find that the estimate (\ref{e: estimate R along max par}) becomes
\begin{equation}
\label{e: estimate Px of R0}
\|P_x R_0(v, \bp a \exp(t h_\ga))\| \leq |(\bp a, th_\ga)|^N (\bp a)^\xi {\rm n} (v).
\end{equation} 
Finally, writing $F_x = P_x F$, 
formula (\ref{e: integral expression F}) leads to 
\begin{eqnarray}
\label{e: F x and integral}
\lefteqn{F_x(v, \bp a\exp (th_\ga))} \\ \nonumber
&=& e^{tx} Q_x(t) F(v,\bp a) +
e^{t x} Q_x(t) \int_0^t Q_x(\tau) e^{\tau(- x + [\xi(h_\ga) - k])} 
R_{0} (v, \bp a \exp \tau h_\ga)\; d\tau.
\end{eqnarray}

We will need to distinguish two cases
depending on the position of the real number $x$,
namely: 
\begin{enumerate}
\itema $-x + [\xi(h_\ga)  - k]  < 0,$ 
\itemb $-x + [\xi(h_\ga) - k] \geq 0.$ 
\end{enumerate}
In case (a) we will need to distinguish the subcases 
(a.1): $x \leq \gL(h_\ga)$ and 
(a.2): $x >  \gL(h_\ga). $ 
In case  (b) we automatically have $ x < \gL(h_\ga),$ in view of 
(\ref{e: condition on k}).
\medno
{\bf Case (a)\ \ } In this case the integrand of (\ref{e: F x and integral}) is integrable over $[0, \infty[$ 
and we find that the expression on the right-hand side of (\ref{e: F x and integral}) 
becomes 
\begin{equation}
\label{e: deco in three exponential terms}
e^{tx} Q_x(t) F(v,\bp a) +  e^{tx} Q_x(t)   I_0(v, \bp a)  - 
e^{t x}Q_x(t) I_t (v,a ').
\end{equation}
where, for $t \in [0,\infty[,$ 
$$ 
I_t( v, \bp a) := \int_t^\infty  Q_x(\tau) e^{\tau(- x + [\xi(h_\ga) - k])} 
R_{0} (v, \bp a \exp \tau h_\ga)\; d\tau\,.
$$% 
From (\ref{e: estimate Px F}) and (\ref{e: estimate Q x}) we infer that  
$$ 
\|Q_x(t) F(v,\bp a)\| \leq C_0 (1 + |t|)^p (1 + |\log \bp a|)^d (\bp a)^\xi.
$$% 
We select $\ge > 0$ such that 
$$ 
-x + [\xi(h_\ga) - k] +\ge < 0 \quad {\rm and}\quad [ \xi(h_\ga)  - k] +\ge < \gL(h_\ga).
\quad 
$$% 
Then
\begin{eqnarray}
\|I_t( v, \bp a)\| & \leq  & e^{(-x +[ \xi(h_\ga) - k]+\ge)t } \int_t^\infty \|Q_x(\tau)\| e^{-\ge \tau}\; |R_{0} (v, \bp a \exp \tau h_\ga)| d\tau \nonumber \\
& \leq & C_1 e^{(-x +[ \xi (h_\ga) - k]+\ge)t }(1 + |\log \bp a|)^N (\bp a)^\xi {\rm n } ( v)\nonumber \\ \nonumber
&\leq & C_1 e^{(-x + \gL(h_\ga) )t }(1 + |\log \bp a|)^N (\bp a)^\xi {\rm n } ( v)
\end{eqnarray}
with 
$$ 
C_1 = \int_0^\infty \|Q_x(\tau)\| e^{-\ge \tau}\; ( 1 + |\tau|)^N d\tau < \infty. 
$$% 
For $t=0$ this leads to 
$$
\|I_0(v,\bp a)\|\leq C_1 (1 + |\log \bp a|)^N (\bp a)^\xi {\rm n } ( v).
$$% 
For general $t\geq 0$ this leads to 
\begin{equation} 
\label{e: estimate I t v bp a}
e^{xt} \|I_t(v,\bp a)\| \leq C_1 e^{t \gL(h_\ga)}(1 + |\log \bp a|)^N (\bp a)^\xi {\rm n } ( v)
\end{equation} 
We will now consider the two subcases (a.1) and (a.2).
\medno
{\bf Subcase (a.1)\ \ } In this subcase, $x \leq  \gL(h_\ga).$ Then from the above estimates it follows that the norm of each term of (\ref{e: deco in three exponential terms}) can be estimated by
$C_2:= \max (C_0, C_1)$ times 
$$ 
e^{t \gL(h_\ga)}(1 + |t|)^p (1 + |\log \bp a|)^N (\bp a)^\xi {\rm n } ( v)
$$% 
There exists a constant $C_3 > 0$ such that the above expression
can be estimated by 
$$ 
C_3 e^{t \gL(h_\ga)} (\bp a)^\xi (1+ |\log (\bp a)+ th_\ga|)^{p+N} {\rm n}(v)
=
C_3 (\bp a a_t)^\xi \, (1 +|\log(\bp a a_t)|)^{p+N} {\rm n}(v).
$$% 
for all $\bp a \in \bp A_\Phi$ and $t \geq 0.$ It follows that $F_x$ satisfies 
the required estimate for all $v \in H^\infty$ and $a = \bp aa_t$ with $t \geq 0.$ 
In view of Remark \ref{r: restricted domain of estimate} this establishes the required estimate in the present subcase.

\medno
{\bf Subcase (a.2)\ \ } In this case we have $x > \gL(h_\ga).$ If follows from 
(\ref{e: F x and integral}) that for $\bp a \in \bp A_\Phi$ and $t \geq 0$ we have 
\begin{equation}
\label{e: different growth left right}
F_x(v, \bp a a_t) -  e^{tx} Q_x(t) I_t (v,a ') = e^{tx} Q_x(t) [F(v,\bp a) +    I_0(v, \bp a)] 
\end{equation}
In view of Proposition \ref{p: first estimate whittaker coefficient} there exists a constant $d' \in \N$ and for every $v \in H_K$ a constant $C_v >0$ such that
for $ \bp aa_t \in \cl(A^+)$ we have 
\begin{equation}
\label{e: second estimate F x v}
\|F_x(v, \bp a a_t)\| \leq C_v  (1 +| \log \bp a|)^{d'} (\bp a)^{\gL} (1 + |t|)^{d'} e^{t \gL(h_\ga)}.
\end{equation} 
Write $\bp \fa_{\Phi}^+$ for the set of elements $H \in \bp\fa_\Phi$ such that 
$\gb(H) >  0$ for all $\gb \in \Phi.$ Then this set has non-empty interior
in $\bp \fa_{\Phi}.$ Let $\bp A_\Phi^+: = \exp (\bp\fa_\Phi^+).$ Then it follows that  
$\bp A_\Phi^+ a_t \subset A^+$ for all $t> 0.$ 

Combining (\ref{e: second estimate F x v}) with the estimates (\ref{e: estimate I t v bp a}) and (\ref{e: estimate Q x}) we infer that the norm of the sum on the left-hand side of (\ref{e: different growth left right}) allows for every $v \in H_K $ 
and $\bp a \in \cl(\bp A_\Phi^+)$ an estimation 
by a constant times $(1 + |t|)^{N'} e^{t\gL(h_\ga)},$ for all $t \geq 0.$
On the other hand, the expression on the right-hand side is of exponential 
polynomial type with exponent $x > \gL(h_\ga).$  
By uniqueness of asymptotics this implies that for $v \in H_K$, 
$\bp a \in \bp A_\Phi^+$ and $t >0$ the expression on the right-hand side
in (\ref{e: different growth left right}) is zero. 
Hence,
\begin{equation}
\label{e: F x equals I t}
F_x(v, \bp a a_t) =  e^{tx} Q_x(t) I_t (v,\bp a) 
\end{equation}  
for $v\in H_K,$ $\bp a \in \bp A_\Phi^+$ and $t > 0. $

From the definitions of $F_x$ and $I_t$ it easily follows that for $v \in H^\infty,$ 
$\bp a\in \bp A_\Phi$ and $t > 0$,
\begin{equation}
\label{e: bp a in action}
F_x(v, \bp a a_t) = F_x(\pi(\bp a)^{-1}v, a_t) \text{and} I_t(v, a') = I_t(\pi(\bp a)^{-1}v, e).
\end{equation}
Furthermore, from the estimates (\ref{e: estimate Px F}) and (\ref{e: estimate I t v bp a}) it follows that for each $t > 0$ the maps 
$v \mapsto F_x(v,a_t)$ and $v \mapsto I_t(v, e)$ are continuous linear $V^\infty \to \C^p.$ If $v$ is $K$-finite in $H$ then $v$ is an analytic vector for $H^\infty,$ 
so that the map $\bp a \mapsto \pi(\bp a)^{-1} v$ is analytic $\bp A_\Phi \to H^\infty.$ 
We may now conclude that the maps $\bp a \mapsto F_x(\pi(\bp a)^{-1}v, a_t)$ 
and $\bp a \mapsto I_t(\pi(\bp a)^{-1}v, e)$ are analytic $\bp A_\Phi \to \C^p.$ 
In view of (\ref{e: bp a in action}) it now follows by analytic continuation in the variable $\bp a$ 
that the validity of the identity (\ref{e: F x equals I t}) extends to all $v \in H_K,$ $\bp a \in \bp A_\Phi$ and $t > 0.$ Using that both members of (\ref{e: F x equals I t}) depend continuous linearly on $v \in H^\infty,$ whereas $H_K$ is dense in the latter space, we conclude that (\ref{e: F x equals I t}) is actually valid
for all $v \in H^\infty,$ $\bp a \in \bp A_\Phi$ and $t >0.$ 

Using (\ref{e: estimate Q x}) and (\ref{e: estimate I t v bp a}) to estimate  
the expression on the right-hand side of (\ref{e: F x equals I t}), we obtain,
for $v \in H^\infty$, $\bp a \in \bp A_\Phi$ and $t > 0,$ 
\begin{eqnarray*}
\|F_x(v, \bp a a_t)\| & \leq &  C_0 C_1 (1 + |t|)^p e^{t \gL(h_\ga)}(1 + |\log \bp a|)^N (\bp a)^\xi {\rm n } ( v)\\
&\leq& 
C_3  (1 + |\log (\bp a a_t)|)^{2(N + p)} (\bp a a_t)^{\xi '} {\rm n}(v),
\end{eqnarray*}
with $C_3>0$ uniform with respect to $v \in H^\infty$, $\bp a\in \bp A_\Phi$ 
and $t> 0.$ This gives the required estimate for $F_x$ on $H^\infty \times \bp A_\Phi A_\Phi^+.$ In view of Remark \ref{r: restricted domain of estimate} this  completes
the discussion in the present subcase.

\medno
{\bf Case (b)\ \ } In this case, $- x + [\xi(h_\ga) - k] \geq 0,$ and
$x < \gL(h_\ga).$ 
The identity  (\ref{e: F x and integral}) can be rewritten as
\begin{equation}
\label{e: Fx in terms of J}
F_x(v, \bp a a_t) = e^{tx} Q_x(t) [ F(v,\bp a) + J(t, v,\bp a)].
\end{equation} 
with 
$$ 
J(t, v,\bp a) := \int_0^t Q_x(\tau) e^{\tau(- x + [\xi(h_\ga) - k])} 
R_{0} (v, \bp a \exp \tau h_\ga)\; d\tau.
$$% 
By a straightforward estimation of the integral defining $J,$ we find, 
for $v \in H^\infty, $ $\bp a\in \bp A_\Phi$ and $t \geq 0,$ 
using (\ref{e: estimate Q x}) and (\ref{e: estimate Px of R0})  that 
\begin{eqnarray*}
\|J(t, v, \bp a)\|& \leq & C_0 e^{t(- x + [\xi(h_\ga) - k])}  \int_0^t (1 + |\tau|)^p \|R_0(v, \bp a\exp \tau h_\ga)\| d\tau\\
&\leq & C_0 C^N e^{t(- x + [\xi(h_\ga) - k])} (1+|t|)^{N+p+1} (1 + |\log \bp a|)^N (\bp a)^\xi{\rm n}(v),
\end{eqnarray*}
with $C = \sup ( 1 + |t||h_\ga|)(1 +|t|)^{-1}.$ This implies that
\begin{eqnarray*}
\|e^{tx} Q_x(t) J(t,v,\bp a)\| & \leq &  C_0^2 C^N e^{t([\xi(h_\ga) - k])} (1+|t|)^{2N+p+1} 
(1 + |\log \bp a|)^N (\bp a)^\xi{\rm n}(v)\\
&\leq & 
C_0^2 C^N e^{t \gL(h_\ga)} (1 + |t|)^{2N + p+1} (1 + |\log \bp a|)^N (\bp a)^\xi{\rm n}(v)
\\
&\leq & C_2 (\bp a a_t)^{\xi'} (1 + |\log \bp a a_t|)^{6N + 2 p +2} {\rm n}(v),
\end{eqnarray*}
see also (\ref{e: condition on k}), with $C_2 >0$ uniform with respect to 
$v \in H^\infty,$ $\bp a \in \bp A_\Phi$ and $ t > 0.$

On the other hand, since $x < \gL(h_\ga) =\xi'(h_\ga),$ see (\ref{e: xi prime on h ga}),
\begin{eqnarray*}
\| e^{tx} Q_x(t) F(v,\bp a)\|  &   \leq  & C_0 e^{\gL(h_\ga)} (1 +|t|)^p (\bp a)^\xi(1 +| \log \bp a|)^d 
{\rm n}(v)\\
& \leq & 
C_3 (\bp a a_t)^{\xi'} (1 + |\log \bp aa_t|)^{2(p+d)} {\rm n}(v),
\end{eqnarray*}
with $C_3 > 0$ uniform with respect to $v \in H^\infty,$ $\bp a \in \bp A_\Phi$ and $t\geq 0.$ 
Combining the two latter estimates we find that there exist a constant $N'>0$ 
and a continuous seminorm ${\rm n}'$ on $H^\infty$ such that 
$$ 
\|P_x F(v,\bp a a_t) \|\leq (\bp aa_t)^{\xi'} ( 1+ |\log (\bp aa_t)|)^{N'} {\rm n}'(v)
$$% 
for all $v \in H^\infty, $ $\bp a \in \bp A_\Phi$ and $t \geq 0.$ 
In view of Remark \ref{r: restricted domain of estimate} this gives the required estimates for $F_x$ in case (b) and completes the proof of Lemma \ref{l: estimate imp}.
\qed

\section{Parabolically induced representations}
\label{s: par ind reps}
In this section we will describe the space of smooth vectors for  
parabolically induced representations of the form
\begin{equation}
\label{e: parabolically induced rep}
\Ind_{P}^G(\xi),
\end{equation}
where 
 $P = M_P A_P N_P$ is a parabolic subgroup of $G$ with the indicated Langlands decomposition, 
 and $\xi$ a continuous representation of $P$ in a Hilbert space $H_\xi.$ 
 
 \begin{rem}
 \label{r: special form xi}
 In particular we will be interested in the situation that
 $\xi = \gs \otimes \nu \otimes 1,$ with $\gs$ a unitary representation in $H_\gs$ 
 and $\nu \in \faPdc.$ 
 The representation $\xi$ is now given by $\xi(man)v = a^\nu \gs(m)v$ for all
 $v \in H_\gs$ and $(m,a,n) \in M_P \times A_P \times N_P.$  For technical 
 reasons we wish to have the possibility to tensor representations 
 of this particular form with finite dimensional representations of $P,$ whence the greater generality.
 \end{rem}
 
We put $K_P:= K\cap M_P = K \cap P.$ 
By averaging over $K_P$ we may replace the inner product on $H_\xi$ 
with a $K_P$-invariant inner product for which the associated norm is equivalent to  the original norm.
Accordingly, we may and will assume that $\xi|_{K_P}$ is unitary.

Let $\gd_P: P \to [0,\infty[$ be the character of  $P$ defined by 
$$ 
\gd_\gd(p) = |\det [\Ad(p)|_{{\rm Lie}(P)}]|^{1/2},\qquad (p\in P).
$$% 
Then 
$$ 
\gd_P(man) = a^{\rho_P}, \qquad ((m,a,n) \in M_P \times A_P \times N_P),
$$% 
where $\rho_P\in \fad$ is defined by $\rho_P(X) =\frac12 \tr (\ad (X)|_{\fn_P}),$ 
for $X \in \fa.$ 

We denote by  $C(G/P: \xi)$ the Fr\'echet space of continuous functions
$f: G \to H_\xi$ transforming according to the rule
$$ 
f(g p) =   \gd_P(p)^{-1} \xi(p)^{-1} f(x), \qquad (x \in G, p \in P).
$$% 
This space is equipped with the left regular representation $L$ of $G.$ 
Let $\bp H_\xi$ denote the Hilbert space $H_\xi$ equipped with the 
representation $\bp \xi$ of $P$ given by 
$
\bp \xi(p) := \gd_P(p) \xi(p). 
$
The group  $P$ naturally acts on $C(G, \bp H_\xi)$ by the formula
\begin{equation}
\label{e: action of P on function space}
[p\gf] (g) = \bp\xi(p) [ \gf(gp)], \qquad (g \in G, p\in P).
\end{equation}
Accordingly, 
$$ 
C(G/P:\xi) = C(G, \bp H_\xi)^P.
$$% 

\begin{rem}
In the particular case  $\xi = \gs \otimes \nu \otimes 1,$ we write 
$$ 
C(G/P:\xi) = C(G/P: \gs :\nu).
$$% 
\end{rem}
We write 
$C(K/K_P: \xi)$ for the Fr\'echet space of continuous functions $\gf: K \to H_\xi$ such that
\begin{equation}
\label{e: rule compact picture}
\gf(km) = \xi(m)^{-1} \gf(k), \qquad (k\in K, m \in K_P).
\end{equation}
Using the decomposition $G \simeq K \times_{K_P} M_P \times A_P \times N_P,$ 
given by the multiplication map, we readily see that restriction 
to $K$ induces a $K$-equivariant topological linear isomorphism
\begin{equation} 
\label{e: compact picture iso}
r: C(G/P:\xi) \;\;{\buildrel \simeq\over \longrightarrow } \;\; C(K/K_P:\xi), \quad f\mapsto f|_K.
\end{equation}
Via this isomorphism we may transfer the representation $L$ of $G$ on 
$C(G/P:\xi)$ to a continuous representation of $G$ in $C(K/K_P: \xi),$ 
denoted $\pi_{P,\xi}.$ 
\begin{rem}
In case $\xi = \gs \otimes \nu \otimes 1$ we will use the notation $\pi_{P,\gs, \nu} = \pi_{P,\xi}.$ 
\end{rem}
Let $d\dot{k}$ be the choice of a $K$-invariant Radon measure on $K/K_P$ normalized by 
$\int_{K/K_P} d\dot{k} = 1.$ 
Then it follows from the isomorphism (\ref{e: compact picture iso}) that the 
sesquilinear pairing
\begin{equation}
\label{e: pairing induced reps}
C(G/P:\xi) \times C(G/ P:\xi) \to \C
\end{equation} 
given by 
\begin{equation}
\label{e: pre Hilbert structure}
\inp{f}{g}:= \int_{K/K_P} \inp{f(k)}{g(k)}_\xi \; d\dot{k}
\end{equation} 
is a $K$-equivariant pre-Hilbert structure. We denote the associated norm 
by $\|\dotvar \|_2. $ The associated 
Hilbert completion is denoted by $L^2(G/P: \xi).$ The restriction map 
(\ref{e: compact picture iso}) induces an isometric isomorphism from this 
completion onto $L^2(K/K_P:\xi ),$ the completion of $C(K/K_P:\xi)$ with respect to 
the pre-Hilbert structure given by (\ref{e: pre Hilbert structure}).
We note that a different but equivalent choice of $K$-invariant inner product on $H_\xi$ 
gives rise to the same completed spaces (as topological linear spaces), with equivalence of the 
Hilbert inner products. 

Out next objective is to show that the representation $(L, C(G/P:\xi))$ has a 
unique extension to a continuous representation of $G$ in 
the Hilbert space $L^2(G/P:\xi).$ To prepare for the proof we start by recalling a  well known result
involving the representation in $\C$ given by the character $\gd_P$ of $P$.
The associated space $C(G/P: \gd_P)$ consists
of the continuous functions $\gf: G \to \C$ such that  
$$ 
\gf(x man) = a^{- 2\rho_P} \gf(x), \qquad (x\in G, (m,a,n) \in M_P\times A_P \times N_P).
$$% 
\begin{lemma}
\label{l: around gd P}
For all $\gf \in C(G/P:\gd_P)$ and $g \in G,$ 
$$ 
 \int_{K/K_P} \gf(g k)\; d\dot{k}  = \int_{K/K_P} \gf(k)\; d\dot{k} .
$$%
There exists a choice of Haar measure $d \bar n$ on $\bar N_P$ such that 
for all $\gf \in C(G/P:\gd_P),$ 
$$
\int_{K/K_P} \gf(k)\; d\dot{k} = \int_{\bar N_P} \gf(\bar n)\; d\bar n.
$$% 
\end{lemma}

The following corollary will be useful for the ongoing discussion.
We define $\kappa_P: G\to K$ $\mu_P: G \to \exp(\fm_P \cap \fp),$ $h_P: G\to A_P$ and
$n_P: G \to N_P$  to be the unique maps determined by 
$$ 
x = \kappa_P(x) \mu_P(x) h_P(x) n_P(x), \qquad (x \in G).
$$% 
These maps  are all analytic maps between the indicated analytic manifolds. 

\begin{cor}
\label{c: L one boundedness}
Let $\omega$ be a bounded subset of $G.$ Then there exists a constant $C_\omega > 0$ such that
for all $\psi \in C(K/K_P)$ and $g \in \omega,$ 
$$ 
\int_{K/K_P} |\psi (\kappa_P (g k))| \; d\dot k \leq C_\omega \int_{K/K_P} |\psi(k)|\;d\dot{k}.
$$% 
\end{cor} 

\proof
Consider the function $\gf: G \to H_\xi$ defined by $\gf(kp) = \gd_P(p)^{-2}\psi(k).$ 
Then $\gf \in C(G/P: \gd_P).$ 
Put $C_\omega := \sup_{x \in \omega  K}  \gd_P(x)^{-2}.$
Then it follows by application of Lemma \ref{l: around gd P}   that, for $g \in \omega,$ 
\begin{eqnarray*}
\int_{K/K_P} |\psi (\kappa_P (g k))| \; d\dot{k} & = & 
\int_{K/K_P} \gd_P(k)^{-2} |\gf (g k))| \;d\dot{k}\\
&\leq&
C_\omega \int_{K/K_P} |\gf (g k)| \;d\dot{k}
= C_\omega \int_{K/K_P} |\gf(k)|\;d\dot k\\
&=& C_\omega \int_{K/K_P} |\psi(k)|\;d\dot k.
\end{eqnarray*}
 \vspace{-30pt}
 
 \qed

The following result as well as its proof are contained in \cite[III.7]{BWcc}.

\begin{prop}
\label{p: induced rep is cts Hilbert}
The left regular representation $L$ of $G$ in $C(G/P:\xi)$ 
has a unique extension to a continuous representation of $G$ in the Hilbert space
$L^2(G/P:\xi).$ 
\end{prop}
\proof
By the principle of uniform boundedness, the operator norm $\|\bp \xi(p)\|_{\rm op}$ of $\bp\xi(p)= 
\gd_P(p) \xi(p) \in \End(H_\xi)$ is locally bounded
as a function of $p \in P.$ 
For the purpose of this proof, we define the analytic map $p_P: G \to P$ by $p_P = \mu_P h_P n_P.$ 
Then $\|\bp \xi (p_P(g)\|_{\rm op} $ is locally bounded as a function of $g \in G.$ For a bounded subset $S \subset G$
let $\bp C_S >0$ be the supremum of the values  $\|\bp \xi(p_P(x^{-1} k)^{-1})\|_{\rm op}$ for $x \in S$ 
and $k\in K.$ 

Let $f \in C(G/P:\xi)$ and $x \in S.$ Since $\xi|_{K_P}$ is unitary, it
follows that the function $\psi: k \mapsto \|f(k)\|_\xi^2 $ belongs to $C(K/K_P).$ 
Hence, by application of Corollary \ref{c: L one boundedness} with $\omega = S^{-1},$ 
\begin{eqnarray*}
\int_{K/K_P} \|L_x f(k)\|_\xi^2\; d\dot{k}&=& 
\int_{K/K_P} \| \bp \xi(p_P(x^{-1} k)^{-1}) f(\kappa_P(x^{-1}k)\|_\xi^2\; d\dot{k} 
\\
&\leq &
\bp C_S^2 \int _{K/K_P} |\psi(\kappa_P(x^{-1}k)|\; d\dot{k} \\
&\leq & \bp C_S^2 C_{\omega} \int_{K/K_P} |\psi(k)| \;d\dot{k}
=
 \bp C_S^2 C_{\omega}  \int_{K/K_P} \|f(k)\|_\xi^2 \; d\dot{k}.
\end{eqnarray*}

Let $\|\dotvar \|_2$ denote the norm associated with the Hilbert structure on $L^2(G/P:\xi).$  
It follows from the estimate above that the map $L_x$ is continuous with respect to the norm $ \|\dotvar\|_2$ 
on $C(G/P:\xi)$  with operator norm that is locally bounded relative to $x \in G.$ 
This implies that $L_x$ has a unique continuous linear extension to a bounded 
map of the Hilbert completion $L^2(G/P:\xi),$ with locally bounded operator norm.
It remains to be shown that the 
associated action map $G \times L^2(G/P:\xi) \to L^2(G/P:\xi)$ is continuous.

Let $f \in C(G/P:\xi).$ Then it is readily checked that $\sup_{k \in K} \|L_x f(k) - f(k)\|_\xi \to 0$ 
for $G\ni x \to e.$ This implies that $L_x f\to f$ in $L^2(G/P: \xi).$ 
If $f \in L^2(G/P:\xi)$ and $f_0 \in C(G/P: \xi)$ then 
\begin{eqnarray*}
\|L_x f - f\|_2 & \leq& \|L_x (f - f_0)\|_2 + \|L_x f_0 - f_0 \|_2 + \|f_0 - f\|_2 \\
& \leq & (\|L_x\|_{\rm op} +1 )\, \|f -f_0\|_2 + \|L_x f_0 - f_0\|_2.
\end{eqnarray*}
Using density of $C(G/P:\xi)$ in $L^2(G/P:\xi)$ we infer from the results obtained in the first part of this proof  that $\|L_x f -f \|_2 \to 0$ for $x\to e$ in $G.$ 

Finally, let  $f_0 \in L^2(G/P:\xi).$ 
Then 
\begin{eqnarray*} 
\|L_x f -  f_0 \|_2 & \leq &\|L_x f - L_x f_0 \|_2 + \|L_x f_0 - f_0\|_2 \\
& \leq & \|L_x\|_{\rm op} \|f -f_0\|_2+ \|L_x  f_0 - f_0\|_2.
\end{eqnarray*}
By what we have established above it follows that both terms in the latter sum 
tend to zero as $(x, f) \to (e, f_0)$ in $G \times L^2(G/P:\xi).$ 
Thus the action map $G \times L^2(G/P:\xi)  \to L^2(G/P:\xi)$ is continuous at every point of $\{e\} \times L^2(G/P:\xi).$ 

Since the operator norm of $L_x$ is locally bounded relative to $x\in G$, the continuity
of the action map at any point of $G \times L^2(G/P:\xi)$ follows.
\qed
 
The representation (\ref{e: parabolically induced rep}) is defined to be the unique extended representation of Proposition \ref{p: induced rep is cts Hilbert}. Under the (unitary) restriction map $\gf \mapsto \gf|_{K}$, this 
representation is transfered to a continuous representation of $G$ in  $L^2(K/K_P:\xi)$ 
which extends $\pi_{P,\xi}$ and is denoted by the same symbol. The latter representation  is called
the compact picture of (\ref{e: parabolically induced rep}).
\medbreak

We are now prepared to determine the space of smooth vectors of the representation
(\ref{e: parabolically induced rep}) . 
We define $C^\infty(G/P:\xi): = C^\infty(G, \bp H_\xi)  \cap C(G/P, \xi).$
Equivalently, this can be expressed in terms of the  action of $P$ on $C^\infty(G, \bp H_\xi)$ given by formula (\ref{e: action of P on function space}): 
\begin{equation}
\label{e: smooth vectors as P invariants}
C^\infty(G/P:\xi) = C^\infty(G, \bp H_\xi)^P
\end{equation}
This is a closed subspace of the Fr\'echet space $C^\infty(G, \bp H_\xi),$ 
hence a Fr\'echet space of its own right. The left regular representation
of $G$ in the first space in (\ref{e: smooth vectors as P invariants}) is smooth, hence restricts to a smooth representation of $G$ in the second space.

\begin{thm}
\label{t: identification of smooth vectors}
The space of smooth vectors in $(L, L^2(G/P:\xi)),$ equipped with its usual 
Fr\'echet topology is given by 
\begin{equation}
\label{e: space of smooth vectors} 
L^2(G/P:\xi)^\infty = C^\infty(G/P:\xi).
\end{equation}
\end{thm}

\proof  
By Fubini's theorem and compactness of $G/P$ it follows that
\begin{equation}
\label{e: P invariants in Ltwoloc}
L^2(G/P:\xi)  =  L^2_{\rm loc}(G, \bp H_\xi)^{P},
\end{equation}
with equality of the usual locally convex topologies; here superscript $P$ indicates
the space of invariants for the obvious action of $P$ on $L^2_{\rm loc}(G, H_\xi),$ 
described by the formula given in (\ref{e: action of P on function space}). This $P$-action is by continuous linear maps which commute with the left regular action of $G.$  Hence, the space on the right of (\ref{e: P invariants in Ltwoloc}) is a closed $G$-invariant subspace of 
$L^2_{\rm loc}(G, H_\xi).$ From this it readily follows that
\begin{equation}
\label{e: smooth vectors of Ltwo}
L^2(G/P: \xi)^\infty = L^2_{\rm loc}(G, \bp H_\xi)^{P}\cap L^2_{\rm loc}(G, \bp H_\xi)^\infty.
\end{equation}
By the Sobolev embedding theorem we have that 
\begin{equation}
\label{e: smooth vectors of Ltwoloc}
L^2_{\rm loc}(G, \bp H_\xi)^\infty = C^\infty(G, \bp H_\xi),
\end{equation}
with equality of the usual locally convex topologies. 
Combining (\ref{e: smooth vectors of Ltwo}) and (\ref{e: smooth vectors of Ltwoloc})
we find that $L^2(G/P:\xi)^\infty = C^\infty(G, \bp H_\xi)^P.$ In view of 
(\ref{e: smooth vectors as P invariants}) this completes the proof.
\qed

\begin{cor}
The left regular representation $L$ of $G$ in $C^\infty(G/P:\xi)$ is smooth.
\end{cor}

\begin{rem}
For $\xi = \gs \otimes \nu \otimes 1$ with $\gs$ a unitary representation of 
$M_P$ and $\nu \in \faPdc$ the  characterisation (\ref{e: space of smooth vectors})
was used in \cite{BScfunc},
with a reference to \cite[III.7.9]{BWcc}. However, the characterization of 
$L^2(G/P:\xi)^\infty$ in \cite{BWcc} was different. It is presented in the lemma below.
\end{rem}
\medbreak

\begin{lemma} 
\label{l: equality of spaces of smooth vectors}
The following equality of locally convex spaces is valid:
\begin{equation}
\label{e: equality of spaces of smooth vectors}
C^\infty(G, \bp H_\xi)^P = C^\infty(G, \bp H_\xi^\infty)^{P}.
\end{equation} 
\end{lemma}

\proof
We consider the Fr\'echet space $C^\infty(P, \bp H_\xi)$ equipped with the
$P$-action given by $p\psi = \bp\xi(p) L_p \psi.$ The action is by continuous linear maps
so that the subspace of $P$-invariants is closed, hence Fr\'echet.

The proof is motivated by the observation that
\begin{equation} 
\label{e: motivating iso}
\bp H_\xi^\infty \simeq C^\infty(P, \bp H_\xi)^P
\end{equation} 
as topological linear spaces. The isomorphism from left to
right is given by $\ga: v \mapsto \psi_v,$ where $\psi_v: P \to \bp H_\xi$ 
is given by $\psi_v(p) = \bp \xi(p) v.$ The inverse $\ge$ is given 
by evaluation at the identity.
The isomorphism intertwines the $P$-action given by $\bp \xi$ on the first space
with the $P$-action obtained from restriction of the right regular $P$-action on 
$C^\infty(P, \bp H_\xi).$ 

The idea of the proof is to establish the following sequence of 
topological linear isomorpisms:
\begin{eqnarray}
\label{e: first iso in array}
C^\infty(G, \bp H_\xi)^P & \simeq & C^\infty(G\times P, \bp H_\xi)^{P\times P} \\
\label{e: second iso in array}
& \simeq & C^\infty(G, C^\infty(P, \bp H_\xi)^P)^P\\
\label{e: third iso in array}
& \simeq &  C^\infty(G, \bp H_\xi^\infty)^{P}
\end{eqnarray}
The isomorphism (\ref{e: first iso in array}) is obtained by restriction of 
the map 
$$
S: C^\infty(G, \bp H_\xi) \to C^\infty(G\times P, \bp H_\xi)
$$% 
given by  $S(\gf)(g, q ) = \gf (g q ),$ for $\gf \in  {\rm dom}(S)$ and $(g,q)\in G \times P.$ 
The image of $S$ consists of the space of invariants
for the $P$ action on its codomain given by $[p_1 \cdot \psi](g, q) = \psi(gp_1, p^{-1}_1 q),$ 
for $\psi \in C^\infty(G\times P, \bp H_\xi), $ $ g \in G$ and $p_1,q \in P.$ As the action takes
place by continuous linear maps, the image of $S$ is closed hence Fr\'echet 
and it is readily verified that $S$ is a topological isomorphism onto its image.

There is a second action of $P$ on the codomain of $S,$ given by 
$p_2 \cdot \psi(g,q) = \bp\xi(p_2)\psi(g, q p_2),$ for $g \in G$ and $q, p_2 \in P.$ 
This second action commutes with the first one, hence leaves ${\rm im}(S)$ invariant.
The map $S$ intertwines the usual $P$-action on its domain with 
the second $P$-action on its codomain. The associated spaces of invariants
are closed, hence Fr\'echet, and it follows that $S$ induces a topological linear
isomorphism between these spaces. This is the isomorphism (\ref{e: first iso in array}).

To understand the map (\ref{e: second iso in array}) we consider the map 
$$
T: C^\infty(G\times P, \bp H_\xi)\;\; {\buildrel \simeq \over \longrightarrow}\;\; 
C^\infty(G, C^\infty(P, \bp H_\xi)),
$$% 
given by $[T\psi](g)(p) = \gf(g,p),$  It is well known that this map is a 
topological linear isomorphism of Fr\'echet spaces.
It is readily verified that $T$ intertwines the given action of $P \times P$ 
on its domain with the action on its codomain given by 
$$ 
(p_1, p_2)\cdot \vartheta(g)(q) := \bp \xi(p_2) \vartheta(gp_1)(p_1^{-1} q p_2),
$$% 
for $\vartheta \in C^\infty(G, C^\infty(P, \bp H_\xi)),$  $g \in G$ and $q,p_1,p_2 \in P.$ 
The associated space of invariants for $\{e\}\times P$ equals
$C^\infty(G, C^\infty(P, \bp H_\xi)^P).$ Taking the remaining action of $P\times \{e\}$ 
into account we infer that $T$ induces a topological linear isomorphism
(\ref{e: second iso in array}). 
Finally, the isomorphism (\ref{e: motivating iso}) induces the 
topological linear isomorphism  (\ref{e: third iso in array}).
\qed
\begin{rem}
\label{r: generality of result}
From the proof it is clear that  Lemma \ref{l: equality of spaces of smooth vectors} is valid for any 
triple $(G, P, \bp \xi),$ with $G$ a Lie group, $P$ a closed subgroup and 
$(\bp \xi, \bp H_\xi)$ a continuous Hilbert representation of $P.$
\end{rem}

\begin{rem}
If $\xi = \gs \otimes \nu \otimes 1$ as in Remark \ref{r: special form xi},
then the space of smooth vectors for $\gs$  in $H_\gs$ equals the space
of smooth vectors  for $\xi$ in $H_\xi = H_\gs.$ 
Hence, in this setting the equality (\ref{e: equality of spaces of smooth vectors})
becomes
$$ 
C^\infty(G/P: \gs :\nu) =  C^\infty(G, H_{\gs,\nu}^\infty)^P
$$% 
where $H_{\gs, \nu}^\infty$ denotes  $H_\gs^\infty$ 
equipped with the representation of $P = M_P A_P N_P$ given by 
$(man,v) \mapsto a^{\nu + \rho_P}\gs(m) v.$ 
\end{rem}

In the sequel we will need the description of the space of smooth vectors for 
parabolically induced representations in the compact picture. To accomodate this we agree to
write $\cR(P)$ for the set of (equivalence classes of) 
continuous Hilbert space  representations $(\xi, H_\xi)$ of $P$ such that
the space of smooth vectors for $P$ in $H_\xi$ equals the space of smooth vectors
for its compact subgroup $K_P.$

\begin{lemma}
\begin{enumerate}
\itema
Any $\xi = \gs \otimes \nu \otimes 1$ with $\gs$ a unitary representation of $M_P$ and
$\nu \in \faPdc$ belongs to $\cR(P).$
\itemb
If $\xi \in \cR(P)$ and if $(\pi, F)$ is a finite dimensional continuous representation of 
$P,$ then $\xi \otimes \pi \in \cR(P).$
\end{enumerate} 
\end{lemma}

\proof 
Assertion (a) follows from \cite[p. 3]{Wrrg2}. For (b) we first note that
if $L$ is any Lie group, $\xi$ a continuous representation of $L$ in a Hilbert 
(or more generally quasi-complete locally convex) space
$E$ and $(\pi, F)$ a finite dimensional continuous representation of $L,$ 
then 
\begin{equation}
\label{e: smooth of tensor}
(E \otimes F)^\infty = E^\infty \otimes F. 
\end{equation}
To prove this, we first note that $\pi: L \to \End(F)$ is smooth. By finite dimensionality 
of $F$ it now follows that $(1 \otimes \pi): L \to \End(E)\otimes \End(F) \simeq \End(E \otimes F)$ 
is smooth. Let $T \in E \otimes F$ be a smooth vector.
Then it follows that 
$$ 
x \mapsto (\xi(x) \otimes 1)T =   (1 \otimes \pi(x^{-1}))(\xi(x) \otimes \pi(x)) T
$$% 
is smooth from $L$ to $E \otimes F.$ 
By finite dimensionality of $F$ this implies that $T \in E^\infty \otimes F.$ This shows that 
the space on the left in (\ref{e: smooth of tensor}) is contained in the space on the right.

For the converse inclusion, let $e \in E^\infty$ and $f \in F,$ then $x \mapsto \xi(x) e$,  $L \to F$ and 
$x \mapsto \pi(x)f,$ $  L \to F$ are smooth. 
By finite dimensionality of $F$ it now follows that $e \otimes f$ is smooth.
The claim follows.

Returning to (b), and letting $\infty(P)$ indicate the smooth vectors for $P$ 
and $\infty(K_P)$ those for $K_P,$ we see that
$$ 
(H_\xi \otimes F)^{\infty(P)}= H_\xi^{\infty(P)} \otimes F = H_\xi^{\infty(K_P)} \otimes F
=
(H_\xi \otimes F)^{\infty(K_P)}.
$$% 
\qed 

We denote by $C^\infty(K/K_P:\xi)$ the Fr\'echet space of smooth functions $\gf: K \to H_\xi$ 
transforming according to the rule (\ref{e: rule compact picture}). 
Then clearly the isomorphism (\ref{e: compact picture iso}), induced by restriction to $K$, 
restricts to an injective continuous linear map from $C^\infty(G/P:\xi)$ into $C^\infty(K/K_P:\xi).$

\begin{lemma}
\label{l: restriction of smooth vectors and the compact picture}
If $\xi \in \cR(P),$  then
the restriction map $ f \mapsto f|_K$ defines a $K$-equivariant 
topological linear isomorphism
\begin{equation}
\label{e: compact picture smooth vectors} 
r: C^\infty (G/ P:\xi) \;\; {\buildrel \simeq \over \to} \;\; C^\infty(K/K_P:\xi).
\end{equation}
\end{lemma} 
See also \cite[\S 10.1.1]{Wrrg2}, or \cite[Cor. III.7.9]{BWcc}).
\medno
Before proceeding with the proof we notice that there exists a unique representation $\pi_{P,\xi}^\infty$ 
of $G$ in $C^\infty(K/K_P:\xi)$ which makes the map (\ref{e: compact picture smooth vectors})
$G$-equivariant.
This representation is called the compact picture of the induced representation on the level
of smoooth vectors.
\medbreak

\proof
Suppose $\xi \in \cR(P).$ 
In view of the closed graph theorem for Fr\'echet spaces, it suffices
to prove the surjectivity of the map above. 
Let $V$ denote the space of $v \in H_\xi$ which are smooth for the restricted
unitary representation $\xi|_{K_P},$ equipped with its natural topology,
and let $H_\xi^\infty$ denote the space of smooth vectors for $\xi.$ 
Then $H_\xi^\infty$ is contained in $V$ with continuous inclusion 
map. 

By assumption on $\xi$ all elements of $V$ are smooth for $P,$ so that
$V = H_\xi^\infty$ as sets.
By application of the closed graph theorem for Fr\'echet spaces, it now follows
that $V = H_\xi^\infty$ as Fr\'echet spaces. By application
of Lemma \ref{l: equality of spaces of smooth vectors} and Remark \ref{r: generality of result}
we have 
$$ 
C^\infty(K/K_P:\xi) = C^\infty(K , H_\xi^\infty)^{K_P}.
$$% 
Hence, if  $\gf \in C^\infty(K/K_P:\xi),$ we infer that 
the function $F: K \times P \to H_\xi $ 
given by 
$$ 
F(k, p) = \bp\xi(p^{-1}) [\gf(k)] 
$$% 
is smooth. It factors through a smooth function 
$$ 
f: G \simeq K\times_{K_P} P \to H_\gs.
$$% 
Hence, $f: G \to H_\xi$ is smooth and belongs to $C^\infty(G/P:\xi).$ 
Moreover, $f|_K = \gf.$ 
\qed

Finally, we will need a few results related the nilpotent picture of $\Ind_P^G(\xi).$ 
Given a compact subset $S \subset \bar N_P$, we denote
by $C_S^\infty(\bar N_P, H_\xi^\infty)$ the subspace of functions in  $C^\infty(\bar N_P, H_\xi^\infty)$ 
with support contained in $S.$ Let  
$$
C_S^\infty(G/P : \xi): = \{ f\in C^\infty(G/ P : \xi)\mid \supp(f) \subset S P\}.
$$ 
This space is a closed subspace of $C^\infty(G/ P : \xi) = C^\infty(G, \bp H_\xi^\infty)^P.$ 
It follows that restriction from $G$ to $\bar N_P$ induces a continuous linear map
\begin{equation}
\label{e: defi nr}
{}^n r: C_S^\infty(G/P : \xi) \to C_S^\infty(\bar N_P, H_\xi^\infty)
\end{equation}
Since the multiplication map $\bar N_P \times P \to G$ is a diffeomorphism
onto the dense open subset $\bar N_P P$ of $G,$ 
it is readily seen that  the map (\ref{e: defi nr}) is injective. 

\begin{prop}
The map (\ref{e: defi nr}) is a topological linear isomorphism of Fr\'echet spaces.
\end{prop}

\proof 
By the closed graph theorem for Fr\'echet spaces, it is sufficient to show that
${}^n r$ is surjective. Let $\psi \in C^{\infty}_S(N_P, H_\xi^\infty).$ 
The map $\iota_\xi: H_\xi^\infty \to C^\infty(P, H_\xi)$ defined by 
$$
\iota_{\bp\xi}(v)(p) = \bp\xi(p)^{-1}v , \qquad (v \in H_\xi^\infty, p \in P),
$$ 
is a continuous embedding onto the closed subspace $C^\infty(P, \bp H_\xi)^P.$  It follows that 
$$ 
\iota_\xi \after \psi: \bar N_P \to C^\infty(P, H_\xi)
$$ 
is a smooth map. This implies that the function $\bar N_P \times P \to H_\xi,$ 
$$
(\bar n, p)\mapsto  \iota_{\bp\xi} (\psi(\bar n))(p) = \bp\xi(p)^{-1}\gf(\bar n)
$$
is smooth. This function has support contained in $S \times P.$ 
It follows from this that the function $\gf: G \to H_\xi$ defined by $\gf = 0$ on $G\setminus SP$
and by 
$$ 
\gf(\bar n p) = \bp\xi(p)^{-1}\gf(\bar n), \qquad ((\bar n, p) \in S \times P),
$$ 
belongs to $C^\infty_{SP} (G, H_\xi).$ It is now readily seen that $\gf \in C^\infty_S(G:\xi)$ 
and that ${}^n r(\gf) = \psi.$
\qed

The inverse of the above map ${}^n r$ will be denoted by 
$i_{\xi, S}:  C_S^\infty(\bar N_P, H_\gs^\infty) \to C_S^\infty(G/P : \gs: \nu).$ 
It is a continuous linear isomorphism of Fr\'echet spaces. We define
$C_c^\infty(\bar N_P, H_\xi^\infty)$ as the locally convex direct limit
of the Fr\'echet spaces $C_S^\infty(\bar N_P, H_\xi^\infty).$ 
Then it follows that the maps $i_{\xi, S}$, for all $S\subset \bar N_P$ compact,
are the restrictions of a single continuous linear map 
\begin{equation} 
\label{e: i nu}
i_{\xi}:  C_c^\infty(\bar N_P, H_\gs^\infty) \to C^\infty(G/P : \xi).
\end{equation}  
For every compact set $S \subset \bar N_P$ 
the natural bilinear map $C^\infty_S(\bar N_P ) \times H_\gs^\infty \to C^\infty_S(\bar N_P, H_\xi),$ $(\psi , v)\mapsto \psi\otimes v$ is jointly continuous. However, we warn the reader that 
this need not be true for the similar bilinear map $C_c^\infty(\bar N_P)\times H_\xi^\infty \to C_c^\infty(\bar N_P, H_\xi^\infty),$ see \cite[Cor. 4.18]{Osb}.

\section{Generalized vectors for induced representations}
\label{s: generalized vectors for ind reps}
We retain the notation of the previous section. In particular,
$P= M_P A_P N_P$ is a parabolic subgroup of $G$ containing $A$ and 
$(\xi, H_\xi)$ is a continuous Hilbert space representation of $P.$ 
Without loss of generality we may assume that $\xi|_{K_P}$ is unitary, see
the text subsequent to (\ref{e: pre Hilbert structure}).

We will now introduce a $G$-equivariant pairing between induced 
representations which is well known for the particular case 
$\xi = \gs \otimes \nu \otimes 1.$ 

We denote by $\xi^*$ the conjugate to $\xi$ in $H_\xi,$  which is the  representation 
of $P$ in $H_\xi$ defined by 
$$ 
\inp{\xi^*(p) v}{w} = \inp{v}{\xi(p^{-1})w} \qquad (p \in P, \;v,w \in H_\xi).
$$%  
Thus, $\xi^*(p)$ equals the Hilbert adjoint  of $\xi(p^{-1}).$ Clearly,  $\xi$ is unitary if and only if $\xi^* = \xi.$ In Remark \ref{r: dual of pi} we mentioned that $\xi^*$ is a continuous representation 
of $G$ in $H_\xi.$ 
  
If $f \in C(G/P: \xi)$ and $g \in C(G/P: \xi^*),$ we define
the function  $\inp{f}{g}_\xi: G \to \C$ by 
$$ 
\inp{f}{g}_\xi(x)= \inp{f(x)}{g(x)}_{\xi}, \qquad (x \in G),
$$% 
where the expression $\inp{\dotvar}{\dotvar}_\xi$ on the right denotes the inner product of $H_\xi.$ 
Since the restriction of $\xi$ to $K_P$ is unitary, $\xi^*|_{K_P}  =\xi|_{K_P},$ and we see
that restriction of the function $\inp{f}{g}_\xi$ to $K$ belongs to $C(K/K_P).$ 
This allows us to define the sesquilinear pairing 
\begin{equation}
\label{e: sesquilineair pairing with contragredient} 
\inp{\dotvar}{\dotvar}:  C(G/P: \xi) \times C(G/P: \xi^*) \to \C,
\end{equation} 
by the formula
\begin{equation}
\label{e: pairing by int over K} 
 \inp{f}{g}:=  \int_{K/K_P} \inp{f(k)}{g(k)}_\xi \; d\dot{k}.
\end{equation}

\begin{lemma}
\label{l: G equivariance of sesquilinear pairing}
The sesquilinear pairing (\ref{e: sesquilineair pairing with contragredient}) is $G$-equivariant.
\end{lemma}

\proof 
Let $f ,g \in C(G/P: \xi) \times C(G/P: \xi^*)$ and define $\gf: G \to \C$ by 
$
\gf(y) = \inp{f(y)}{g(y)}_\xi.
$ 
Then one readily verifies that $\gf \in C(G/P:\gd_P).$ 
Using Lemma \ref{l: around gd P} we infer, for $x \in G,$ that
$$ 
\inp{L_x f}{L_x g} = \int_{K/K_P} \gf(x^{-1} k) \; d \dot{k} = \int_{K/K_P} \gf(k)\;d\dot{k} = \inp{f}{g}.
$$% 
\vspace{-36pt}

\qed

\vspace{24pt}

The pairing (\ref{e: sesquilineair pairing with contragredient}) 
obviously extends to a continuous sesquilinear pairing 
\begin{equation}
\label{e: sesquilinear pairing L2 induced}
L^2(G/P: \xi) \times L^2(G/P: \xi^*) \to \C,
\end{equation}
given by the same formula. By density and continuity, the extended 
pairing is $G$-equivariant. In particular, we see again that if 
$\xi$ is unitary then $\xi^* = \xi$ and the representation
$(L, L^2(G/P:\xi))$ is unitary. In general, without the requirement that $\xi$ be unitary, the following result is valid.

\begin{lemma}
The Hermitian pairing (\ref{e: sesquilinear pairing L2 induced}) is a perfect pairing of Hilbert spaces,
realizing each of them $G$-equivariantly as the conjugate dual of the other one.
\end{lemma}

\proof
Since $\xi|_{K_P}$ is unitary, $\xi^*$ and $\xi$ are equal on $K_P.$ 
Accordingly, restriction to $K$ induces isometric isomorphisms
$L^2(G/P:\xi) \simeq L^2(K/K_P:\xi|_{K_P})$ and $L^2(G/P:\xi^*) \simeq L^2(K/K_P:\xi|_{K_P}).$
Via these isomorphisms, the pairing ((\ref{e: sesquilinear pairing L2 induced}) becomes
the Hermitian pairing of $L^2(K/K_P:\xi|_{K_P})$ with itself, given by (6.5). As that pairing 
is perfect, so is (\ref{e: sesquilinear pairing L2 induced}). The final assertion follows
from the $G$-equivariance of (\ref{e: sesquilinear pairing L2 induced}).
\qed

At various points in this article we will need the following description of  the equivariant pairing in terms of the nilpotent group $\bar N_P.$  

\begin{lemma}
If  $f\in C(G / P : \xi) $ and $g \in C(G / P :  \xi^*)$ 
then   
$$ 
\inp{f}{g} = \int_{\bar N_P} \inp{f(\bar n)}{g(\bar n)}_{\xi} \; d\bar n.
$$% 
\end{lemma}

\proof
Since $\gf:= \inp{f}{g}_\xi$ belongs to $C(G/P:\gd_P),$ see also the proof of 
Lemma \ref{l: G equivariance of sesquilinear pairing}, we have,
by application of Lemma \ref{l: around gd P}, that 
$$ 
\inp{f}{g} = \int_{K/K_P} \gf(k)\;dk = \int_{\bar N_P} \gf(\bar n)\; d\bar n,
$$% 
whence the required identity.
\qed

Being perfect, the Hermitian pairing (\ref{e: sesquilinear pairing L2 induced}) induces 
a $G$-equivariant topological linear isomorphism 
\begin{equation}
\label{e: Ltwo onto conjugate dual of Ltwo}
L^2(G/P: \xi)\; {\buildrel \simeq \over \longrightarrow}  \; \overline{L^2(G/P: \xi^*){\,}'}.
\end{equation} 
Since $C^\infty(G/P:\xi^*)$ is a dense subspace of $L^2(G/P: \xi^*),$ the transpose
of the associated inclusion map induces an injective continuous linear map 
\begin{equation}
\label{e: Ltwo dual into Cinfty dual}
\overline{L^2(G/P: \xi^*){\,}'} \;\embeds \;\overline{C^\infty(G/P:\xi^*){\, }'}.
\end{equation} 
Here the second space is equipped with the strong dual (locally convex) topology. The map 
is given by restriction to $C^\infty(G/P:\xi^*).$

\begin{defi}
\label{d: defi generalized vectors induced}
Let $(\xi, H_\xi)$ be a continuous representation of $P$ in a Hilbert space. 
Then by 
$C^{-\infty}( G/P: \xi)$ we denote the conjugate continuous linear dual  
$C^{\infty} (G/P: \xi^*){}',$ given as the second space in (\ref{e: Ltwo dual into Cinfty dual}),
equipped with the strong dual  topology. 
\end{defi}

\begin{rem}
Being the strong dual of a Fr\'echet space, $C^{-\infty}(G/P:\xi)$ 
is a complete locally convex Hausdorff space. Since the induced representation 
$\pi_\xi = L$ of $G$ in $C^\infty(G/P:\xi)$ is smooth, it follows that 
the natural representation  $\pi^{-\infty}_\xi$ 
of $G$ in $C^{-\infty}(G/P:\xi)$ is continuous and even smooth. 
See the text around (\ref{e: contragredient}).

Furthermore, the associated derived representation of $U(\fg)$, 
also denoted $\pi^{-\infty}_\xi,$ is given by 
$$ 
\pi^{-\infty}_\xi(u): \gf\mapsto \gf \after \pi^{\infty}_\xi(\bar u^\veec)
$$% 
for $u \in U(\fg).$ Here  $u \mapsto \bar u^\veec$ is the conjugate linear automorphism
of $U(\fg)$ that for (real) $X \in \fg$ is given by $\bar X^\veec = - X.$ 
\end{rem}

 The composition of the two maps (\ref{e: Ltwo onto conjugate dual of Ltwo}) 
 and (\ref{e: Ltwo dual into Cinfty dual}) leads to 
 the $G$-equivariant continuous linear embedding 
\begin{equation}
\label{e: Ltwo embedded into generalized sections}
L^2(G/ P : \xi) \;\; \embeds \;\;  C^{-\infty}(G/ P : \xi).
\end{equation} 
given by $f \mapsto \inp{f}{\dotvar}.$ The elements
of the latter space will be called generalized vectors for the induced representation $\Ind_P^G(\xi).$ 
In the sequel, we will use the map (\ref{e: Ltwo embedded into generalized sections})
to identify $L^2(G/P:\xi)$ with a subspace of $C^{-\infty}(G/P:\xi).$ 

\begin{rem}
\label{r: dual of gs otimes nu} 
In particular, if $\xi = \gs \otimes \nu \otimes 1,$ with $(\gs, H_\gs)$ a unitary representation of $M_P$ and $\nu \in \faPdc$ we obtain that 
$\xi^* = \gs \otimes -\bar\nu \otimes  1$ so that $C^{-\infty}(G/P: \gs : \nu)$ is defined as the conjugate continuous 
linear dual of $C^\infty(G/P: \gs :  -\bar\nu).$ 
\end{rem}

If $(\xi_j, H_j)$ are continuous representations of $P$ in Hilbert spaces,
 for $j=1,2,$ 
then any continuous linear intertwining operator $\gf: H_1 \to H_2$ induces the $G$-equivariant continuous linear map
$$ 
\Ind_P^G(\gf): L^2(G/P:\xi_1) \to L^2(G/P:\xi_2), \;\;f \mapsto \gf \after f.
$$% 
This map restricts to a $G$-equivariant continuous linear map $C^\infty\Ind_P^G(\gf)$ from the space 
$C^\infty(G/P:\xi_1)$ to the space 
$C^\infty(G/P:\xi_2).$

\begin{lemma}
\label{l: functoriality C min infty}
The map $\Ind_P^G(\gf)$ has a unique continuous linear extension to a map 
$$
C^{-\infty}\Ind_P^G(\gf) : C^{-\infty}(G/P:\xi_1) \to C^{-\infty}(G/P:\xi_2).
$$% 
The extension is $G$-equivariant.
\end{lemma} 

\proof 
The conjugate map $\gf^*: H_2 \to H_1$ is continuous linear and intertwines $\xi_2^*$ with $\xi_1^*.$ 
We consider the transpose $T'$ of the map $T:= C^\infty\Ind_P^G(\gf^*).$ This map, given by the formula
$\vartheta \mapsto \vartheta \after T,$  is $G$-equivariant and continuous linear 
$C^{-\infty}(G/P:\xi_1) \to C^{-\infty}(G/P: \xi_2)$ (use that $\xi_j^{**} = \xi_j$). We will proceed by establishing
the claim that this map restricts to $\Ind_P^G(\gf).$ Indeed, let $f\in L^2(G/P:\xi_1).$ Then for 
$g \in C^\infty(G/P:\xi_2)$ we have
\begin{eqnarray*}
\inp{T'(f)}{g} & = & \inp{f}{Tg} = \int_{K/K_P} \inp{f(k)}{\gf^*\after g(k)}_1\; d\dot{k}\\
&= & 
\int_{K/K_P} \inp{\gf\after f(k)}{g(k)}_2\; d\dot{k} =
\inp{\Ind_P^G(\gf)(f)}{g}.
\end{eqnarray*}
Hence $T' f= \Ind_P^G (\gf) (f)$ for  $f \in L^2(G/P: \xi_1),$ establishing the claim. 
This settles the existence of the continuous linear extension of $\Ind_P^G(\xi).$
The uniqueness and $G$-equivariance follow from the density of 
$C^\infty(G/P: \xi_1)$ in $C^{-\infty}(G/P: \xi_1).$ 
\qed

We agree to identify the open right $P$-invariant subsets of $G$ with the open subsets of $G/P$ 
via the canonical projection $G \to G/P.$ Likewise, the closed right $P$-invariant 
subsets of $G$ are identified with the closed subsets of $G/P.$ 
Accordingly, if $S\subset G/P$ is closed, we 
denote by 
$C^\infty_S(G/P:\xi^*)$ the closed subspace of $f \in C^\infty(G/P:\xi^*)$ such that  
$\supp \, f \subset S.$ For a given open subset $\Omega \subset G/P,$ we write
${\rm CPT}(\Omega)$ for the collection right $P$-invariant subsets of $\Omega$
which are closed, hence compact, as subsets of $G/P,$ 
and put 
$$ 
C_c^\infty(\Omega :\xi^*):= \cup_{S \in \CPT(\Omega)}  \;C_S^\infty (G/P:\xi^*).
$$ 
Accordingly, we equip $C_c^\infty(\Omega :\xi^*)$ with the direct limit locally convex topology.
Thus, a seminorm $\gs$ on $C_c^\infty(\Omega :\xi^*)$ is continuous if and only if 
for every $S\in \CPT(\Omega)$ the restriction of $\gs$ to $C_S^\infty(G/P:\xi^*)$
is continuous. 
For $\Omega \subset G/P$ open we define  
\begin{equation}
\label{e: defi gen fun on Omega}
C^{-\infty}(\Omega:\xi):= \overline{ C^\infty_c(\Omega:\xi^*)'},
\end{equation}
equipped with the strong dual topology. We will view $C(\Omega: \xi)$ as a linear subspace of  $ C^{-\infty}(\Omega: \xi)$ 
via the map $g\mapsto \inp{\dotvar}{g}$ given by the sesquilinear pairing
$$ 
C^\infty_c(\Omega:\xi^*) \times C(\Omega:\xi) \to \C, \quad (f,g)\mapsto \inp{f}{g} = \int_{K/K_P} \inp{f(k)}{g(k)}_\xi \; d\dot{k}.
$$
Accordingly, the natural sesquilinear pairing associated with (\ref{e: defi gen fun on Omega}) will be denoted by 
$$ 
\inp{\dotvar}{\dotvar}: \;\; C_c^{\infty}(\Omega:\xi^*) \times C^{-\infty}(\Omega:\xi) \to \C
$$ 
We denote by $C^\infty(\Omega)$ the space of smooth right $P$-invariant functions
$\Omega \to \C.$ 
For each $\gf$ in this space the multiplication map $g\mapsto \gf g,$ $C(\Omega:\xi) \to C(\Omega: \xi)$ has a unique continuous linear extension to a map $C^{-\infty}(\Omega:\xi) \to C^{-\infty }(\Omega: \xi).$ This map, denoted $u \mapsto \gf u,$  is given by 
$$ 
\inp{g}{\gf u} = \inp{\bar \gf  g}{u} , \quad (u \in C^{-\infty}(\Omega:\xi), \;g \in C^\infty(G/P:\xi^*)).
$$ 
If $\Omega_1\subset \Omega_2$ are open 
subsets of $G/P,$ the transpose of the inclusion map $C^\infty_c(\Omega_1:\xi) \embeds C^\infty_c(\Omega_2:\xi)$ gives us the continuous linear restriction map
$$ 
\rho^{\Omega_2}_{\Omega_1}: u \mapsto u|_{\Omega_1},\;\; C^{-\infty}(\Omega_2:\xi)\to C^{-\infty}(\Omega_1:\xi).
$$ 
Together with these restriction maps, the assignment $\Omega \mapsto C^{-\infty}(\Omega:\xi)$ 
defines a presheaf of  $C^\infty(G/P)$-modules on $G/P.$ 
By using smooth partitions of $1$ over $G/P,$ is readily seen that this presheaf is in fact a sheaf,
as it has the following required restriction and glueing properties, for any open covering
$\{\Omega_i\mid i \in I\}$ of $G/P.$ 
\medno
{\bf Restriction poperty.\ } If $u \in C^{-\infty}(G/P:\xi)$ and $u|_{\Omega_i} = 0$ for all 
$i\in I,$ then $u = 0.$ 
\medno
{\bf Glueing property.\ } If $u_i \in C^{-\infty}(\Omega_i:\xi)$ for $i \in I$ and
$u_i|_{\Omega_i\cap \Omega_j} = u_j|_{\Omega_i\cap \Omega_j}$ for all $i,j \in I,$ 
then there exists a $u \in C^{-\infty}(G/P:\xi)$ such that $u|_{\Omega_i}= u_i$ 
for all $i \in I.$ 
\medbreak
We will finish this section by introducing a certain direct limit topology
 on the spaces $C^{-\infty}(G/P:\xi),$ assuming that 
both $\xi$ and $\xi^*$ belong to $\cR(P).$                                                                                                                                                                 By Lemma \ref{l: restriction of smooth vectors and the compact picture}
restriction to $K$ induces a topological linear isomorphism $r: C^\infty(G/P:\xi^*) \to 
C^\infty(K/K_P: \xi^*).$ Denote the conjugate continuous linear dual of the latter
space by $C^{-\infty}(K/K_P: \xi).$ Since $\xi|_{K}$ is unitary,
it follows that $\xi^*$ and $\xi$ are equal on $K,$ so that the latter two spaces do not 
change if $\xi$ is  replaced by $\xi^*.$
By transposition we obtain a topological linear isomorphism 
\begin{equation}
\label{e: generalized compact picture}
r^*: C^{-\infty} (K/K_P:\xi)\;\;{\buildrel \simeq \over \longrightarrow} \;\;C^{-\infty} (G/P:\xi).
\end{equation}% 
The map $r^*$ is equivariant for a uniquel representation of $G$ in the space $C^{-\infty}(K/K_P:\xi)$ 
which we denote by $\pi_{P, \xi}^{-\infty}.$ This representation is called the compact picture of the induced 
representation $\Ind_P^G(\xi)$ on the level of generalized vectors.

In the sequel it will be important to use a specific set of continuous norms
on $C^{\infty}(K/K_P\col \xi)$, hence on $C^{\infty} (G/P\col \xi).$ These are introduced as follows.
We fix a basis $X_1, \ldots, X_m$ of $\fk$ and use the notation $X^\mu := X_1^{\mu_1} \cdots X_n^{\mu_n} \in U(\fk),$ 
for $\mu \in \N^n$ a multi-index. For $s \in \N$ the space $C^s(K, H_\xi)$ of $C^s$-functions 
$K \to H_\xi$ is a Banach space for the norm
\begin{equation}
\label{e:  intro norm f s} 
f\mapsto \|f\|_s:= \sum_{|\mu| \leq s} \sup_{k \in K} \|R_{X^\mu} f(k)\|_{H_\xi}.
\end{equation}
The space $C^s(K/K_P\col \xi):= C^s(K, H_\xi)\cap C(K/K_P\col \xi)$ is a closed subspace. 
Hence, equipped with the norm $\|\dotvar \|_s$ it is a Banach space of its own right.
Clearly, the Fr\'echet topology on  $C^\infty(K/K_P:\xi^*)$ is induced by the restrictions
of the norms $\|\dotvar\|_s,$ for $s \in \N.$ 

For each $s \in \N,$ the conjugate continuous linear dual of $C^s(K/K_P:\xi^*)$ will
be denoted by $C^{-s}(K/K_P:\xi).$ Equipped with the dual norm $\|\dotvar\|_{-s}$ , it is 
a Banach space of its own right. The transpose of the inclusion $C^{\infty}(K/K_P : \xi^*)
\to C^{s}(K/K_P:\xi^*)$  is an injective continuous linear map 
\begin{equation}
\label{e: limit topology}
C^{-s}(K/K_P:\xi) \embeds C^{-\infty}(K/K_P:\xi)
\end{equation}% 
via which we shall identify elements. As the norms $\|\dotvar\|_s, (s \in \N),$ 
induce the topology of $C^\infty(K/K_P:\xi^*),$ it follows that $C^{-\infty}(K/K_P:\xi)$ 
is the union of the spaces $C^{-s}(K/K_P:\xi),$ for $s \in \N.$
The latter spaces increase with 
$s,$ and constitute the so-called filtration by order.
The associated inclusion maps
$C^{-s} \to C^{-s-t},$ for $s ,t \in \N,$ are continuous, so that the spaces form a directed
family of locally convex spaces.
We now observe that, as a linear space,
$C^{-\infty}(K/K_P:\xi)$ is the direct limit of the directed family consisting of the spaces 
$C^{-s}(K/K_P:\xi),$ for $s \in \N.$ We may therefore equip $C^{-\infty}(K/K_P:\xi)$ 
with the associated direct limit locally convex topology.  Since the  natural maps
$C^{-s}(K/K_P:\xi) \to C^{-\infty}(K/K_P: \xi)$ are continuous for the strong
dual topologies, it follows that the direct limit locally convex topology on 
$C^{-\infty}(K/K_P: \xi)$ is finer than (or equal to) the strong dual topology. 

\section{Whittaker vectors for induced representations}
\label{s: whit for ind reps}
In this section we will initiate our study of Whittaker vectors
for induced representations of the form $\Ind_{\bar P}^G (\gs \otimes \nu \otimes 1),$ 
with $P$ a standard parabolic subgroup of $G,$ i.e., $P \supset P_0 = M A N_0.$ 
Here $(\gs, H_\gs)$ is an irreducible unitary representation of $M_P$ 
and $\nu \in \faPdc.$ 
At a later stage, 
$\gs$ will be assumed to belong to the 
discrete series $\MPds$ of $M_P,$ i.e. the set of equivalence classes of irreducible square
integrable representations of $M_P.$ Implicitly it is then assumed that $P$ is cuspidal.
From Remark \ref{r: dual of gs otimes nu} we recall that $(\gs \otimes \nu \otimes 1)^* = 
\gs \otimes (-\bar \nu) \otimes 1.$ 

We note that $\fa_P\subset \fa$, $M_P \supset M$ and $N_P \subset N_0.$ 
Accordingly, the set $N_0 \bar P$ is open (and dense) in $G.$ 
The space of Whittaker functionals for the induced representation 
$\Ind_{\bar P}^G (\gs \otimes - \bar \nu \otimes 1)$ is denoted 
\begin{equation}
\label{e: wh of induced}
\Wh_\cchi (L^2(G/ \bar P:\gs: -\bar \nu)^\infty),
\end{equation}
cf.\  (\ref{e: defi space of Wh functionals}). In view of
and (\ref{e: defi space of Wh functionals}) and Theorem \ref{t: identification of smooth vectors}
the space (\ref{e: wh of induced}) consists of the continuous linear 
functionals $\eta \in  C^\infty(G/\bar P:\gs: -\bar\nu)'$ 
such that
$$
\eta \after L_n = \cchi(n) \eta,\qquad (n \in N_0).
$$% 
In view of (\ref{e: whittaker vectors}) and Definition \ref{d: defi generalized vectors induced} with the subsequent remark, 
the space of 
Whittaker vectors of  $\Ind_{\bar P}^G (\gs \otimes -\bar \nu \otimes 1)$ equals 
the space 
\begin{equation}
\label{e: conjugate whittaker for induced rep} 
C^{-\infty}(G/\bar P: \gs: \nu)_\cchi,
\end{equation}
consisting of $j \in  C^{-\infty}(G/\bar P: \gs: \nu)$ transforming according to the rule
$$ 
L_n j  = \cchi(n) j ,\qquad (n \in N_0),
$$% 
in the sense of generalized functions,
see also (\ref{e: whittaker vectors}).

The following result of Harish-Chandra \cite[Thm.~1, p.~143]{HCwhit}, see also \cite{KVwhit} for a proof,
will be crucial for the determination of the space (\ref{e: conjugate whittaker for induced rep}). In fact, it is 
crucial for the entire Whittaker theory.  
If $Q$ is a parabolic subgroup of $G$ we write 
$$ 
M_{1Q} = M_Q A_Q
$$% 
for its $\theta$-stable Levi component.

\begin{thm}{\rm (Harish-Chandra \cite{HCwhit})\ }
\label{t: main thm HC}
Let $Q$ be a standard parabolic subgroup. Then $N_0 M_{1Q} \bar N_Q$ is open 
in $G.$ Let $\cchi$ be regular and suppose that $u$ is a distribution on $G$ such that 
$$ 
R_{\bar n} u = u \;\;\mbox {\rm and}\;\; L_{n_0}u = \cchi(n_0) u,\qquad (\bar n \in \bar N_Q, n_0 \in N_0).
$$% 
If $u = 0$ on $N_0 M_{1Q} \bar N_Q$, then $u = 0$ on $G.$ 
\end{thm}
\begin{rem} 
In the above, the space $\cD'(G)$ of distributions on $G$ is 
defined to be the continuous linear dual of the complete locally convex space $\cD(G):= C_c^\infty(G).$ The left and right regular actions are defined by 
$$ 
L_g u = u \after L_g^{-1},\; R_g u = u \after R_g^{-1}, \qquad (u \in \cD'(G), g \in G).
$$ 
\end{rem}

The following corollary is of immediate importance for our discussion.
We retain the assumption that $\gs$ is an irreducible unitary representation of $M_P.$ 
We will say that an element $\psi \in C^{-\infty}(G/\bar P:\gs:\nu)$ vanishes 
on an open subset $\cO$ of $G/\bar P$ if $\inp{f}{\psi} = 0$ 
for all $f \in C^\infty(G/\bar P:\gs: -\bar\nu)$ with $\supp f \subset \cO.$ 
\begin{cor}
\label{c: cor of main thm HC}
Let $\cchi$ be regular, $\nu \in \faPdc$ and $j \in C^{-\infty}(G/\bar P: \gs: \nu)_\cchi.$ 
If $j$ vanishes on the orbit $N_P\bar P$ in $G/\bar P$ then $j =0.$ 
\end{cor}

\proof 
Fix $v \in H_\gs\setminus \{0\}.$ For $\gf \in C_c^\infty(G)$ we define 
the function $T \gf: G \to H_\gs$ by 
$$ 
T\gf(x) = \int_{M_P} \int_{A_P} \int_{\bar N_P} a^{-\bar \nu - \rho_P} \gf(xma\bar n)  
[\gs(m)v] \; dm da d\bar n.
$$% 
Then it is readily verified that $T$ defines a continuous linear operator
$C_c^\infty(G) \to C^\infty(G/\bar P:\gs: -\bar \nu)$ which intertwines
the left regular actions of $G$ on these spaces. It follows that
$$ 
u: \gf \mapsto \inp{T\gf}{j}
$$% 
defines a distribution on $G.$ It is clear that $u$ is right $\bar N_P$-invariant.
By equivariance of $T$ we see that, for $\gf \in C^\infty_c(G)$ and  $n \in N_0,$ 
$$ 
    L_{n}u(\gf) =  u(L_n^{-1} \gf) = \inp{L_{n^{-1}}T\gf}j = \inp{T\gf}{L_n j}
=
\cchi(n)^{-1} u(\gf) 
$$% 
Now assume that $j =0$ on $N_P \bar P.$ 
If $\supp\, \gf \subset N_P \bar P$ then $\supp (T\gf)
\subset N_P \bar P$, from which it follows that $u = 0$ 
on $N_P \bar P.$ Since $\chi^{-1}$ is a regular character, it now follows from 
Theorem \ref{t: main thm HC} that
$u = 0.$ This implies that $j$ gives zero when applied to the space $T(C_c^\infty(G)).$ We will finish the proof by showing that the latter space is dense in $C^\infty(G/\bar P:\gs:-\bar \nu).$ 
In view of the natural decomposition $G \simeq K\times_{K_P} M_P \times A_P \times \bar N_P$ it suffices to show that the operator 
$$
S: C^\infty_c(K \times_{K_P} M_P) \to C^\infty(K/K_P: \gs)
$$%
defined by 
$$ 
S\psi(k) = \int_{M_P}  \psi(k,m) [\gs(m)v]\; dm
$$% 
has dense image. Let $f \in C^\infty(K/K_P: \gs)$ be $K$-finite from 
the left. Then it suffices to show that $f \in {\rm im}(S).$  There exists
a bi $K$-invariant finite dimensional subspace $F \subset C^\infty(K)$ such that 
$f$ belongs to $F \otimes (H_\gs)_{K_P}$ and is fixed under $R_{k_P}\otimes \gs(k_P)$ 
for all $k_P\in K_P.$  Thus, $f$ is a finite sum of terms
$f_j \otimes v_j$, with $f_j \in F$ and $v_j \in (H_\gs)_{K_P}.$ 
By irreducibility of $\gs$, there exist left $K_P$-finite $\psi_j \in C_c^\infty(M_P)$ 
such that $\gs(\psi_j) v = v_j.$ Put 
$$ 
\psi (k, m) = \sum_j \int_{K_P} f_j(k k_P) \psi_j(k_P^{-1} m)\; dk_P,
$$% 
where $dk_P$ denotes normalized Haar measure on $K_P.$ 
Then it is readily verified that $\psi$ defines an element of 
$C^\infty_c(K \times_{K_P} M_P)$ which has image $f.$ 
\qed
In the following it is not required that $\cchi$ is regular.  
We denote by $C(G, \bp H_{\gs, \nu}^{-\infty})^{\bar P}$ the space
of continuous functions $f: G \to H_\gs^{-\infty}$ such that 
$$ 
f(x ma\bar n) = a^{-\nu + \rho_P} \gs^{-\infty}(m^{-1}) f(x), 
$$ 
for $x \in G$ and $(m,a,\bar n) \in M_P \times A_P\times \bar N_P.$ 
Furthermore, we define
\begin{equation}
\label{e: defi C fs} 
C(G, \bp H_{\gs, \nu}^{-\infty})^{\bar P}_{\fds}
\end{equation}
to be the subspace consisting of 
$f\in C(G, \bp H_{\gs, \nu}^{-\infty})^{\bar P}
$
such that the support of $f$ is contained in $N_P\bar P$ and 
$$ 
\dim {\rm span}\, f(N_P) < \infty.
$$ 
 \begin{lemma} 
 \label{l: continuity f gf}
 If $f \in C(G, \bp H_{\gs, \nu}^{-\infty})^{\bar P}$, then 
 for every $\gf \in C^\infty(G/\bar P: \gs: - \bar \nu)$ the function $x \mapsto \inp{\gf(x)}{f(x)}_{\gs}$ 
 belongs to $C(G/\bar P : \gd_{\bar P}).$
 \end{lemma}
 
 \proof
 This is straightforward.
 \qed
 It follows that for $f$ as in the lemma, we may define the linear functional 
 $$
 \trw\!f_*: C^\infty(G/\bar P: \gs:-\bar \nu) \to \C, \;\; \gf \mapsto \int_{K/K_P} \inp{\gf(k)}{f(k)}_{\gs} \; dk .
 $$ 
 
 \begin{lemma}
 If $f \in C(G, \bp H_{\gs, \nu}^{-\infty})^{\bar P}_{\fds},$ then $\trw\!f_*$ is continuous
 hence is the image of a unique $f_* \in C^{-\infty}(G/P:\gs: \nu).$ Furthermore, the associated
 map $f \mapsto f_*$ is a linear injection
 $$ 
 C(G, \bp H_{\gs, \nu}^{-\infty})^{\bar P}_{\fds} \;\;\hookrightarrow \;\; C^{-\infty}(G/P:\gs : \nu).
 $$ 
 \end{lemma} 
\proof
For brevity, we write $\xi$ for the continuous representation $\gs \otimes -\bar \nu \otimes 1$
of $P$ in $H_\gs.$ There exists a compact subset $S\subset N_P$ such that $\supp f \subset S \bar P.$
From Lemmas  \ref{l: continuity f gf} and \ref{l: around gd P} it follows
that, for every $\gf \in C^\infty(G/\bar P: \xi),$ 
\begin{equation}
\label{e: int for f star gf}
\trw\! f_*(\gf) = \int_{N_P} \inp{\gf(n)}{f(n)}_\gs\; dn =  \int_S \inp{\gf(n)}{f(n)}_\gs\; dn.
\end{equation} 
In view of the hypothesis, the span $E$ of $f(N_P)$ is a finite dimensional subspace of 
$H^{-\infty}_\gs.$ The natural pairing $E \times H^\infty_\gs \to \C$ 
is continuous. It follows that there exist a norm $\|\dotvar\|_E$ on $E$ and a continuous 
seminorm $q$ on $H^\infty_\gs$ such that $|\inp{v}{\mu} |\leq \|\mu\|_E q (v)$ 
for all $\mu \in E$ and $v\in H_\gs^\infty.$ Using (\ref{e: int for f star gf}) we now infer that, 
for $\gf \in C^\infty(G/\bar P: \xi),$ 
$$ 
|\trw\!f_*(\gf)| \leq C \sup_{S} q (\gf(n))
$$ 
where $C:= \sup_{n \in S} \|f(n)\|_E.$ 
Let $\kappa: G \to K$ be the Iwasawa map associated with the decomposition 
$G = K A \bar N_0,$ and let $p_{\bar P}: G \to \bar P$ be defined by $p_{\bar P} (x) =
\kappa(x)^{-1} x.$ 
Then it follows that $p_{\bar P}(S)$ is a compact subset of $P.$ By uniform boundedness,
there exists a continuous seminorm $r$ on the Fr\'echet space $H_\gs^\infty$ such that 
$$
q(\bp\xi(p_{\bar P}(n))^{-1}v) \leq r(v), \quad (n \in S, v \in H_\gs^\infty).
$$ 
It now follows that, for 
$
\gf \in C^\infty(G/\bar P: \xi)
$
and $n \in S,$ 
\begin{equation}
\label{e: est q gf n}
q(\gf(n)) \leq r (\gf(\kappa(n))) \leq \sup_{k \in K} r (\gf(k)).
\end{equation}
Via the natural isomorphisms $C^\infty(K, H_\gs^\infty)^{K_P} \simeq C^\infty(K/K_P: \xi)
\simeq C^\infty(G/\bar P: \xi)$ we see that the expression on the right
of (\ref{e: est q gf n}) defines a continuous seminorm on $C^\infty(G/\bar P: \xi).$ 
The asserted continuity of $\trw\!f_*$ follows.

The element $f_*$ is  now defined by $\trw\! f_* = \inp{\dotvar}{f_*}.$ 
Since $f \mapsto \trw\!f_*$ is conjugate linear, 
linearity of $f \mapsto f_*$ is obvious. For injectivity, assume $f_* = 0.$ Let 
$\gf \in C_c^{\infty}(N_P)$ and $v \in H_\gs^\infty$ be arbitrary.
Then there exists a unique $\tilde \gf \in C^\infty(G/\bar P: \xi)$ with support in 
$N_P \bar P$ such that $\tilde \gf|_{N_P} = \gf\otimes v.$ From $\inp{\tilde \gf}{f_*} = 0$ 
it follows by using (\ref{e: int for f star gf}) that
$$ 
\int_{N_P} \gf(n) \inp{v}{f(n)}_\gs \; dn  = 0.
$$ 
Hence, for every $v \in H_\gs^\infty$ the continuous function 
$n \mapsto \inp{v}{f(n)}_\gs$ vanishes on $N_P$ and we  conclude that $f = 0.$ 
\qed

\begin{thm}
\label{t: evaluation of j}
Let $\nu \in \faPdc$ and $j \in C^{-\infty}(G/\bar P:\gs: \nu)_\cchi.$ There exists a unique
element $\eta \in H_\gs^{-\infty}$  such that the restriction of 
$j$ to $N_P\bar P$  equals the continuous function 
$j_\eta: N_P\bar P \to H_\gs^{-\infty}$ given 
by 
\begin{equation}
\label{e: formula j eta}
j_\eta(n_P m_P a_P \bar n_P) = a_P^{-\nu + \rho_P} \cchi(n_P)^{-1} \gs(m_P)^{-1} \eta. 
\end{equation}
Moreover, 
\begin{equation}
\label{e: behavior eta}
\eta \in (H_\gs^{-\infty})_{\cchi|_{N_0\cap M_P}}.
\end{equation} 
\end{thm}

\begin{rem} 
\label{r: meaning of equality vj j}
By the assertion that the restriction of $j$ to $N_P\bar P$ equals
the continuous function $j_\eta : N_P \bar P \to H^{-\infty}_\gs$ it is meant that
for every $\psi \in C^\infty(G/\bar P)$ with $\supp \,\psi \subset N_P \bar P$ we have
$(\psi j_\eta)_* = \psi j.$ 
\end{rem}

\proof
Let $i_{\gs, -\bar \nu}: C^\infty_c(N_P, H_\gs^\infty) \to 
C^\infty(G/\bar P\col \gs\col  -\bar \nu)$  be the continuous linear map 
defined as in (\ref{e: i nu})
with $\bar P$ in place of $P$ and with $\xi = \gs \otimes (-\bar\nu).$ 
Thus, if $f\in C^\infty_c(N_P, H_\gs^\infty)$ then $i_{\gs,-\bar \nu}(f) \in C^\infty(G:\gs:-\bar \nu)$
is uniquely determined by  
$\supp \, i_{\gs,-\bar \nu}(f) \subset N_P\bar P$ and by 
$$ 
i_{\gs,-\bar \nu}(f)|_{N_P} = f.
$$% 
Given $v \in H_\gs^\infty$ and $\gf \in C^\infty_c(N_P)$ 
we define 
$$ 
\jv(\gf): = \inp{i_{\gs, -\bar \nu}(\gf \otimes v)}{j}.
$$% 
Then $\jv$ defines a continuous linear functional on $C_c^\infty(N_P)$,
i.e., a distribution on $N_P$, which depends linearly on $v \in H_\gs^\infty.$  
\begin{lemma}
The map $(\gf, v) \mapsto {}_v j(\gf)$ is continuous bilinear
$C_c^\infty(N_P) \times H_\gs^\infty \to \C.$ 
\end{lemma}
\proof
Let us denote the above bilinear map by $b.$ 
For every compact set $S \subset N_P$ the bilinear map 
$(\gf, v)\mapsto \gf \otimes v$  evidently has a continuous restriction to 
$C^\infty_S(N_P) \times H_\gs^\infty.$ This implies that the restricted bilinear map 
$b: C^\infty_S(N_P) \times H_\gs^\infty \to \C$ is continuous. 

Since $L_n j = \chi(n)j,$ $(n \in N_P),$ it follows that 
$$
b(L_n \gf , v) = \chi(n) b(\gf, v), \qquad ((\gf, v)\in C_c^\infty(N_P) \times H_\xi^\infty, \; n \in N_P).
$$ 
Taking into account that the manifold $N_P$ is diffeomorphic to $\R^n$ for a certain $n,$ 
it now follows by application of Lemma \ref{l: b continuous bilinear} that $b$ is continuous.
\qed

\begin{lemma}
\label{l: b continuous bilinear}
Let $V$ be a (Hausdorff) locally convex space, and $b: C_c^\infty(\R^n) \times V \to \C$ 
a bilinear map. 
Suppose that there exists a compact neighborhood $\cK$ of $0$ in $\R^n$ such that the restriction 
of $b$ to $C_\cK^\infty(\R^n) \times V$ is continuous.
 Suppose in addition that  for every $x \in \R^n$ exists a 
diffeomorphism ${\ell}_x$ from an open neighborhood $\omega_x'$ of $0$ to an open neighborhood $\omega_x$ of $x$ and a constant $C_x > 0$ such that 
\begin{equation}
\label{e: condition on b}
|b( \gf , v)| \leq C_x |b(\ell_x^*\gf, v)|
\end{equation} 
for all $\gf \in C^\infty_c(\omega_x)$ and $v \in V.$ 
Then $b$ is (jointly) continuous. 
\end{lemma}

\proof
By hypothesis,
there exist continuous seminorms $p$ on $C_{\cK}^\infty(\R^n)$ and $q $ on $V$ such that
\begin{equation}
\label{e: first bil est b}
|b(\gf, v)| \leq p(\gf) q(v), \qquad (\gf \in C_{\cK}^\infty(\R^n), \,v \in V). 
\end{equation}
For $S \subset \R^n$ compact and $k \geq 0$ we define the seminorm $p_{S, k}$ on $C^\infty(\R^n)$ by 
$$
p_{S, k}(\gf) := \max_{|\ga|\leq k} \sup_S |\partial^\ga \gf|.
$$ 
There exist constants $k \in \N$ and $c> 0$ such that $p \leq c p_{\cK, k}$ 
on $C^\infty_{\cK}(\R^n).$ From now on we will keep $k$ fixed and write
$p_{S}$ for $p_{S, k}.$ 

Let $x \in \R^n$ be arbitrary and fix $\ell_x: \omega_x' \to \omega_x$ with the property
mentioned in the hypothesis. We select a compact neighborhood $S_x$ of $x$ 
such that $S_x \subset \ell_x(\omega_x\cap \cK).$  Then there exists a constant $D_x > 0$ such that
for all $\gf \in C^\infty_{S_x}(\R^n)$ 
$$ 
p_\cK(\ell^*_x\gf) \leq D_x p_{S_x}(\gf).
$$ 
Combining this with (\ref{e: condition on b}) and (\ref{e: first bil est b}) we find that, for all $\gf \in C^\infty_{S_x}(\R^n)$ and $v \in V,$ 
$$ 
|b(\gf, v)| \leq  c \,C_x \, p_\cK (\ell_x^*\gf) q(v) \leq c \,C_x D_x \,p_{S_x}(\gf) q(v).
$$
The sets ${\rm int}(S_x)$ cover $\R^n.$ By paracompactness of the latter space,
there exists a $C^\infty$ partition of unity $\{\psi_i\mid i \in I\} \subset C^\infty_c(\R^n)$ subordinate to the mentioned cover.
Thus, for each $i \in I,$ there exists $x_i \in \R^n$ such that $\gs_i:= \supp \,\psi_i \subset {\rm int} (S_{x_i}).$
Moreover, $\sum_{i\in I} \psi_i = 1,$ with locally finite sum. Write $S_i: = S_{x_i}.$ Then 
for each $i \in I$ it follows
by application of the Leibniz formula to $\partial^\ga(\psi_i \gf)$ that 
there exists
a constant $d_i > 0$ such that for all $\gf \in C^\infty(\R^n),$ 
$$ 
p_{S_i}(\psi_i \gf) \leq d_i p_{\gs_i}(\gf). 
$$ 
Noting that $\{\gs_i| i\in I\}$ is a locally finite collection of compact sets, we 
define the seminorm $\tilde p$ on $C_c^\infty (\R^n)$ by 
$$
\tilde p:= \sum_{i\in I}  c \,C_{x_i}D_{x_i} d_i \, p_{\gs_i}. 
$$ 
By a simple argument it follows that $\tilde p$ is continuous on $C^\infty_S(\R^n)$ 
for every compact subset $S \subset \R^n.$ Hence $\tilde p$ is continuous on $C_c^\infty(\R^n).$ 
Furthermore, for every $\gf \in C^\infty_c(R^n)$ and $v \in V$ 
we have
\begin{eqnarray*}
b(\gf, v) & = & \sum_{i\in I} b(\psi_i\gf, v) \\
& \leq &  \sum_{i\in I} c C_{x_i} D_{x_i} 
p_{S_i}(\psi_i \gf) q(v)\\
& \leq & \sum_{i \in I} c C_{x_i} D_{x_i} d_i p_{\gs_i}(\gf) q(v)
= \tilde p(\gf) q(v).
\end{eqnarray*}
This establishes the continuity of $b.$ 
\qed

{\em Completion of the proof of Theorem \ref{t: evaluation of j}.\ }
From Lemma \ref{l: b continuous bilinear} it follows that the map $v \mapsto {}_v j$ 
is continuous linear $H_\xi^\infty \to \cD'(N_P) =C_c^\infty(N_P)'.$ 
Let $L^\veec$ denote the contragredient of the left regular representation $L$ of $N_P$ 
in $C^\infty_c(N_P).$ 
Since the map $i_{\gs, -\bar \nu}$ 
is equivariant for the left regular actions of $N_P$, it follows that  
$$ 
L_{n}^\veec (\jv) =  \cchi(n)^{-1}  \jv,\quad (n \in N_P).
$$% 
We fix  $\psi \in C_c^\infty(N_P)$ such that
$$ 
\int_{N_P} \psi(n) \cchi(n)^{-1} \; dn = 1.
$$% 
Then it follows that 
$$ 
\jv = \int_{N_P} \psi(n) L_n^\veec (\jv ) \; dn = L^\veec(\psi)(\jv).
$$% 
In view of Lemma \ref{l: smoothing operator} below, it now follows that for every $v \in H_\gs^\infty$ 
there exists a unique function $J_v \in C^\infty(N_P)$ such that 
$$ 
\jv(\gf) = \int_{N_P} J_v(n) \gf(n)\; dn, \qquad (\gf \in C_c^\infty(N_P)).
$$% 
Moreover, by the same lemma, the map $v \mapsto J_v$ is continuous linear
$H_\gs^\infty \to C^\infty(N_P).$ 
By uniqueness, $J_v(n'n) = \cchi(n') J_v(n)$, for all $n, n' 
\in N_P,$ 
hence 
$J_v(n) = \cchi(n) J_v(e).$ 
Define
$ \gl : H_\gs^\infty \to \C$ by 
$$
\gl(v):= J_v(e),\qquad (v \in H_\gs^\infty).
$$%  
Since $\gd_e: C^\infty(N_P) \to \C$ is continuous, it follows that $\gl \in (H^\infty_\gs)'.$ 
We now apply the isomorphism (\ref{e: dual H1 infty to H2 minus infty}) in the setting  $\pi_1 = \pi_2 = \gs$ and with (\ref{e: sesquilinear pairing}) given by the inner product on $H_\gs$ which makes $\gs$ unitary. It follows that there 
exists a unique $\eta \in H_\gs^{-\infty}$
such that 
$$ 
\gl(v) = \inp{v}{\eta},\qquad (v \in H_\gs^\infty).
$$ 
Define 
$j_\eta: N_P \bar P \to H_\gs^{-\infty}$ as in (\ref{e: formula j eta}). Then it easily follows that $j_\eta$ is continuous. Furthermore, $j_\eta(N_P) \subset \C \eta.$ Let $\psi \in C^\infty(G/\bar P)$ 
have support contained in $N_P \bar P.$ Then it follows that 
$\psi j_\eta \in C(G, H_{\gs, \nu}^\infty)_{\fds}^{\bar P}$, see (\ref{e: defi C fs}). In view of Remark \ref{r: meaning of equality vj j} it suffices to show that $(\psi j_\eta)_* = \psi j.$ 
Let $\gf \in C_c^\infty(N_P)$ and $v \in H^\infty_\gs.$ Then 
\begin{eqnarray*}
\inp{i_{\gs, -\bar \nu}(\gf \otimes v)}{(\psi j_\eta)_*}&=& 
\inp{i_{\gs, -\bar \nu}(\gf \otimes v)}{\psi j_\eta} 
 =  \int_{N_P} \overline{\psi(n)}\inp{\gf(n)v }{j_\eta(n)}\; dn\\
& = &  \int_{N_P} \overline{\psi(n)}\gf(n)  \cchi(n) \inp{v}{\eta} \; dn =
\int_{N_P} \overline{\psi(n)}\gf(n)  \cchi(n) J_v(e) \; dn\\
&=& \int_{N_P} \overline{\psi(n)} \gf(n) J_v(n) \; dn 
= \jv(\bar \psi \gf) \\
& =&  \inp{i_{\gs, -\bar\nu}(\bar\psi \gf \otimes v)}{j}
= \inp{i_{\gs, -\bar \nu}( \gf \otimes v)}{\psi j} .
\end{eqnarray*}
Therefore,  the continuous linear functional 
$h : C^\infty(G/\bar P \col \gs \col -\bar \nu) \to \C$ defined by 
$h:= \inp{\dotvar}{\psi j - (\psi j_\eta)_*}$ 
vanishes on $i_{\gs,-\bar \nu}(C^\infty_c(N_P) \otimes H_\gs^\infty);$ 
here the algebraic tensor product has been taken.

Let $S_1\subset N_P$ be compact. We select $ S_2 \subset N_P$ compact such that $S_1 \subset {\rm int}\,(S_2);$ 
then $C^\infty_{S_1}(N_P, H_\xi^\infty)$ is contained in the closure of 
$C^\infty_{S_2}(N_P) \otimes H_\xi^\infty$ in $C^\infty_{S_2}(N_P, H_\xi^\infty).$
Since $i_{\gs, -\bar \nu}$ restricts to a topological isomorphism 
from $C^\infty_{S_j}(N_P, H_\xi^\infty)$ onto the space $C^\infty_{S_j}(G/\bar P: \gs :-\bar \nu),$ 
it follows that $h$ vanishes on $C^\infty_{S_1}(G/\bar P: \gs :-\bar \nu).$ As this is valid
with  $S_1$ an arbitrary compact subset of $N_P,$ we may assume that $S_1$ contains
an open neighborhood $\omega$ of $\supp \,\psi$ in $N_P.$ Then $\Omega_1:= \omega \bar P$ 
and $\Omega_2:= G \setminus \supp \,\psi \,\bar P$ form an open cover of $G/\bar P$ 
such that $\psi j - (\psi j_\eta)_*$ restricts to $0$ on $\Omega_j,$ $(j=1,2).$ By the restriction 
property mentioned in the text following (\ref{e: defi gen fun on Omega}) it follows that 
$\psi j = (\psi j_\eta)_*.$ 

This establishes the existence of $\eta$ such that (\ref{e: formula j eta}) is valid.
If $\eta'$ is a similar element, let $\psi\in C^\infty(G/\bar P)$ have support in $N_P \bar P$ 
and satisfy $\psi([e])\neq 0.$ Then $[\psi (j_\eta -j_{\eta'})]_* = 0.$ By injectivity of the map $f \mapsto f_*$ 
it follows that $\psi (j_\eta - j_{\eta'}) = 0.$ Evaluating this identity at $e$ 
we obtain $\psi([e]) (\eta - \eta') =0$ and conclude $\eta = \eta'.$ 
Uniqueness of $\eta$ follows.

It remains to show (\ref{e: behavior eta}). For this we note that,
for $n_0 \in N_0\cap M_P$, conjugation by $n_0$ 
leaves $N_P$ invariant and 
$$ 
L_{n_0} i_{\gs,\nu}(\gf \otimes v) = i_{\gs,\nu} (L_{n_0} R_{n_0} \gf \otimes \gs(n_0) v).
$$% 
From the definition of $\jv$ it now follows that
$$ 
\cchi(n_0) \;\jv (\gf) = {}_{\gs(n_0)v}j (L_{n_0} R_{n_0} \gf).
$$% 
This implies, in turn,
$$ 
 \cchi(n_0) J_{ v} (n) = J_{\gs(n_0)v} (n_0 n n_0^{-1}), \qquad (n \in N_P).
$$% 
Evaluation at $n =e$ gives
$$ 
\cchi(n_0)\inp{v}{\eta}_\gs = \inp{\gs(n_0) v}{\eta}_\gs,
$$% 
for all $v \in H_\gs^\infty$ and $n_0 \in N_0\cap M_P.$ 
Finally, this gives 
$$
\gs^{-\infty}(n_0)^{-1}\eta = \cchi(n_0)^{-1} \eta, \quad (n_0 \in N_0 \cap M_P),
$$% 
and (\ref{e: behavior eta}) follows. This finishes the proof of Theorem
\ref{t: evaluation of j}.
\qed

In the next lemma, we assume that $H$ is a Lie group, equipped
with a left Haar measure $dx.$ We equip $C^\infty_c(H)$ with the left
regular representation $L$ of $H$, and its dual $\cD'(H)$ with the 
contragredient representation $L^\veec.$ Accordingly, we define, for $\psi \in C_c^\infty(H)$, 
the continuous linear map $L^\veec(\psi): \cD'(H) \to \cD'(H)$ by
$$ 
L^\veec( \psi) (u) = \int_H \psi(x)\, u\after L_x^{-1}\; dx.
$$% 
For $v \in C^\infty(H)$ we define $i(v) \in \cD'(H)$ by 
$$ 
i(v)(f) = \int_H f(x) v(x)  dx, \qquad (f \in C^\infty_c(H)).
$$% 
Then $i: C^\infty(H) \to \cD'(H)$ is an equivariant injective continuous linear map with dense image.
In particular,  $i \after L(\psi) = L^\veec(\psi)\after i.$ 

\begin{lemma}
\label{l: smoothing operator}
 If $\psi \in C_c^\infty(H)$, then $L^\veec(\psi)$ is a smoothing 
operator in the sense that there exists a unique continuous
linear map $T_\psi: \cD'(H) \to C^\infty(H)$ such that
\begin{equation} 
\label{e: smoothing operator}
L^\veec(\psi) = i \after T_\psi.
\end{equation} 
\end{lemma}
 
\proof 
Uniqueness is obvious, since $i$ is injective. For $y \in H$ we define $R_y (\check\psi) \in C_c^\infty(H)$ 
by 
$$
R_{y^{-1}} (\check \psi)(x)  : = \check\psi(xy^{-1}) = \psi(yx^{-1}), \qquad (x \in H).
$$% 
Let $\gD: H \to ]0,\infty[$ be defined by $\gD(x) = |\det \Ad(x)|.$ Then the map 
$y  \mapsto \gD R_{y^{-1}}(\check\psi) $  is smooth $H \to C_c^\infty(H).$ 
Thus, if $u \in \cD'(H)$ then  $T_\psi (u): y\mapsto u(\gD R_{y^{-1}}(\check \psi))$ is a smooth function on $H.$ Moreover, the map $u \mapsto T_\psi(u)$ is continuous linear $\cD'(H) \to C^\infty(H).$ 
We will show that it satisfies (\ref{e: smoothing operator}). Since the expressions
on both sides of (\ref{e: smoothing operator}) are continuous linear endomorphisms of $\cD'(H),$ 
it suffices to establish the equality on the dense subspace of elements of the form $u = i(v),$ with $v \in C^\infty(H).$ 
Since $L^\veec(\psi  ) \after i = i \after L(\psi)$ and since $i$ is injective, it suffices to show that
$$ 
L(\psi)(v)  = T_\psi \after i(v).
$$% 
This identity of functions in $C^\infty(H)$ is established as follows. If $y \in H,$ then 
\begin{eqnarray*}
L(\psi)(v)(y) & = & \int_H \psi (x) v(x^{-1} y) \; dx  =   \int_H \psi(yx)  v(x^{-1})   dx\\
&=&  \int_H \psi(yx^{-1})  v(x)  \gD(x)  dx =   \int_H \gD(x) R_{y^{-1}}(\check \psi)(x) v(x) \; dx \\
 & = & i(v)(\gD R_{y^{-1}}(\check \psi)) = T_\psi(i(v))(y). 
\end{eqnarray*}
\qed
We return to the setting of Theorem \ref{t: evaluation of j}.
Given $j \in C^{-\infty}(G/\bar P:\gs: \nu)_\cchi $ we denote by by $\ev_e(j)$ the associated
element $\eta \in H_\gs^{-\infty}$ such that $j = j_\eta$ on $N_P\bar P.$ 
Then by uniqueness of $\eta$ combined with (\ref{e: behavior eta}),
we find that $\ev_e$ defines a linear map
\begin{equation}
\label{e: defi ev e}
\ev_e: C^{-\infty}(G/\bar P:\gs:  \nu)_\cchi \to (H_\gs^{-\infty})_{\cchi|_{N_0\cap M_P}}.
\end{equation}

\begin{cor}
\label{c: ev is injective}
If $\cchi$ is regular, then the map (\ref{e: defi ev e}) is injective.
\end{cor}

\proof
Let $j \in C^{-\infty}(G/\bar P:\gs: \nu)_\cchi$ and suppose that $\ev_e(j) = 0.$ 
Then it follows from Theorem \ref{t: evaluation of j} that 
$j = 0$ on $N_P \bar P.$ By Corollary \ref{c: cor of main thm HC} this implies 
that $j = 0.$ 
\qed

We retain the assumption that $\chi$ is regular. Then $\chi|_{N_0\cap M_P}$
is regular with respect to the roots of ${}^*\fa_P$ in $\fn_0 \cap \fm.$ 
Indeed, for each such root $\ga$, the associated root space satisfies $\fg_\ga \subset \fm_P\cap \fn_0. $ We agree to use the abbreviated notation
\begin{equation}
\label{e: abbreviation Whittaker data}
H_{\gs,\cchi_P}^{-\infty}:= (H_\gs^{-\infty})_{\cchi|_{N_0 \cap M_P}}.
\end{equation}
For $R \in \R$ we put
\begin{equation}
\label{e: defi faPR}
\faPdc(P,  R) :=  \{\nu \in \faPdc \mid \forall \ga \in \gS(\fn_P,\fa): \inp{\Re \nu}{\ga} > R\}.
\end{equation}
Given $\eta \in H_{\gs,\cchi_P}^{-\infty}$ and $\nu \in \faPdc(P,0),$ we define the function 
$
j_\nu  = j(P, \gs ,  \nu, \eta): G  \to  H_\gs^{-\infty}
$
by $j_\nu  = 0$ on $G \setminus N_P \bar P$ and by 
\begin{equation}
\label{e: defi j P gs nu eta}
j(P, \gs, \nu, \eta) (n m a \bar n ) =  a^{- \nu + \rho_P} \cchi(n)^{-1} \gs(m)^{-1} \eta,
\end{equation}
for $n \in N_P$, $  m \in M_P$,  $a \in A_P$ and $\bar n \in \bar N_P$, 
see also (\ref{e: formula j eta}). 

\begin{prop}
\label{p: integrability of j}
Suppose that $\gs \in \MPds$ and let $\eta \in H_{\gs,\cchi_P}^{-\infty}.$
\begin{enumerate}
\itema
If $\nu\in \faPdc(P,0)$ then the function $j_\nu = 
 j(P, \gs, \nu, \eta): G \to H_\gs^{-\infty}$  
satisfies 
 \begin{equation}
 \label{e: rule for j}
 j_\nu(n_0 x m a \bar n) =  \cchi(n_0)^{-1} \,a^{-\nu + \rho_P} \gs^{-\infty}(m)^{-1} j_\nu(x) , 
 \end{equation}
 for all $x \in G$,  $n_0 \in N_0$,  $m \in M_P$,  $ a \in A_P$, and $\bar n \in \bar N_P.$ 
 \itemb
The function $j_\nu$ is continuous $H_\gs^{-\infty}$-valued on $\cK:= K\cap N_P \bar P.$ There exists a continuous seminorm 
$s_\gs$ on $H_\gs^{\infty},$ and for every $R > 0$ a  Lebesgue integrable 
function $L_R: K \to [0,\infty[$ such that
$$ 
|\inp{v}{j(\bar P, \gs, \nu, \eta)(k)}_\gs | \leq L_R(k) s_\gs(v),
$$% 
for all $\nu \in \faPdc(P,R),$ $k \in K$ and $v \in H_\gs^\infty.$ 
\end{enumerate}
\end{prop}
 \proof 
Since $N_P \bar P$ is left $N_0$-invariant, 
so is the complement $G\setminus N_P \bar P$, 
where $j_\nu$ equals zero. Hence $j_\nu$ satisfies (\ref{e: rule for j}) for $x\in G\setminus N_P\bar P$ and all values of $n_0, m, a,\bar n.$ 

For $x \in N_P \bar P$, the rule (\ref{e: defi j P gs nu eta}) implies (\ref{e: rule for j}) 
for $m \in M_P, a \in A_P, \bar n \in \bar N_P$ and all $n_0 \in N_P.$ To obtain the rule for all $n_0 \in N_0$ 
we use that $N_0 = (N_0\cap M_P)N_P$, and note that for $n_0 \in N_0\cap M_P$ and 
$nma\bar n \in N_PM_PA_P\bar N_P$ we have, taking into account that
$n_0 n n_0^{-1} \in N_P$, 
\begin{eqnarray*} 
j_\nu(n_0 n man) & = & j_\nu(n_0 n n_0^{-1} (n_0 m)  an) \\
&=& = a^{-\nu + \rho_P}
\cchi(n_0 n n_0^{-1})^{-1} \gs^{-\infty} (n_0 m)^{-1} \eta \\
&=&  
a^{-\nu + \rho_P}  \cchi(n)^{-1} \gs^{-\infty}(m)^{-1} [\cchi(n_0)^{-1} \eta]  \\
 &=&  \cchi(n_0)^{-1} j_\nu(nman).
\end{eqnarray*}
Rule (\ref{e: rule for j}) now follows, and we turn to proving (b). 

Let $P_0'$ be the minimal parabolic subgroup containing $A$ with 
$N_{P_0'} = (N_0\cap M_P) \bar N_P.$ We consider the maps 
$\kappa': G \to K, \; h': G \to A$ and $n': G \to N_{P_0'}$ associated with the Iwasawa
decomposition 
\begin{equation}
\label{e: primed Iwasawa deco}
G = K A N_{P_0'},
\end{equation}
and put $H' := \log \after h': G \to \fa.$

The natural map $N_P \to G/\bar P$ is an embedding onto a dense open subset.
Via the natural diffeomorphism $K/K_P \simeq G/\bar P$ we have a corresponding open 
embedding $N_P \to K/K_P.$ Since $P_0' \subset \bar P,$ this open embedding is given 
by $n \mapsto \kappa'(n)K_P.$ The associated open embedding $N_P \times K_P \to K$, 
given by $(n, k_P) \mapsto \kappa'(n) k_P,$ has image $\cK.$ By transformation of variables
it is well known that a function $f: K \to \C$ is absolutely integrable if and only if 
\begin{equation}
\label{e: integral f over N P}
\int_{N_P\times K_P}  |f(\kappa'(n) k_P)| e^{ 2 \rho_P  H'(n)}\; dk_P dn < \infty 
\end{equation}
Furthermore, if (\ref{e: integral f over N P}) holds
then $\|f\|_{L^1(K)} $ equals the given integral 
(up to a positive scalar factor, depending on the normalization of measures).
 
In the following we put ${}^*\!A:= M_P\cap A$ and ${}^*\!N_0: = M_P \cap N_0,$ 
so that 
\begin{equation}
\label{e: Iwasawa deco MP}
M_P = K_P {}^*\!A {}^*\!N_0 
\end{equation}
 is the Iwasawa decomposition of $M_P$ 
associated with its minimal parabolic subgroup $P_0' \cap M_P = P_0 \cap M_P.$ 
Therefore, this Iwasawa decomposition for $M_P$ is compatible with the decomposition
(\ref{e: primed Iwasawa deco}). 
Let ${}^*\!\rho \in \fa^*$ and $\rho' \in \fa^*$ 
be defined by 
$$ 
{}^*\!\rho (H) = \half \tr [\ad(H)|_{{}^*\fn_0}],\quad \rho'(H) = \half \tr [\ad(H)|_{\fn_{P_0'}}],
$$% 
for  $H \in \fa.$ Then 
\begin{equation}
\label{e: deco rho prime}
\rho' = {}^* \rho -  \rho_P
\end{equation}
and this decomposition is compatible with $\fa = {}^* \fa_P \oplus \fa_P$ in the sense that
${}^*\rho = 0$ on $\fa_P$ and $\rho_P = 0$ on ${}^*\fa_P.$ 

We chose $N \in \N$ sufficiently large (the precise condition will appear later). 
By Corollary \ref{c: schwartz estimate for ds coeff} applied to $M_P$ in place of $G,$ 
and with $\gl = \inp{\dotvar}{\eta} \in \Wh_{\cchi_P}(H_\gs^\infty),$ 
 there exists a continuous seminorm $q$ on $H_\gs^{\infty}$ such that
\begin{equation}
\label{e: estimate coeff MP} 
|\inp{\gs(m)^{-1}v}{\eta}| \leq ( 1 + |{}^*H(m) |)^{-N} e^{-{}^*\!\rho({}^*H(m))} q(v)
\end{equation}
for all $v \in H_\gs^{\infty}$ and  $m \in M_P.$ Here ${}^*H$ denotes
the Iwasawa projection $M_P \to {}^*\fa$ associated with the decomposition (\ref{e: Iwasawa deco MP}).

In accordance with the decomposition $\fa = {}^* \fa \oplus \fa_{P},$  we decompose the Iwasawa
projection $H'(x)$ of an element $x \in G$ as 
\begin{equation} 
\label{e: deco H prime}
H'(x) = {}^*H'(x) +  H_{\bar P}'(x).
\end{equation}
Furthermore, we agree to write ${}^*h' = \exp \after {}^*H'$ and $h_{\bar P}' = \exp \after H_{\bar P}'.$ 
By compatibility of the decompositions (\ref{e: primed Iwasawa deco}) 
and (\ref{e: Iwasawa deco MP}) we note that  ${}^*H' |_{M_P} = {}^*H.$ 

We now observe that for $n \in N_P$ we have  
$\kappa'(n) = n n'(n)^{-1} h'(n)^{-1}.$
For $v \in H_\gs^\infty$ and $\nu \in \faPdc(P,0)$ the function 
$x \mapsto  \inp{v}{j_\nu(x)}_\gs$ is left $N_{0}$-equivariant and right $\bar N_P$-invariant. 
Hence, for $n \in N_P$ and $k_P\in K_P,$ 
\begin{eqnarray}
\nonumber
\inp{v}{j_\nu(\kappa'(n) k_P)}_\gs & = & 
\cchi(n)^{-1} \, \inp{\gs(k_P)v}{j_\nu(h'(n)^{-1})}_\gs\\
&=& \cchi(n)^{-1} \, e^{(\nu - \rho_P)H'_{\bar P}(n)} \,  \inp{\gs({}^*h'(n))^{-1}\gs(k_P)v}{\eta}_\gs
\label{e: formula j nu inp v}
\end{eqnarray}
Fom (\ref{e: deco H prime}) we see that $(\nu -\rho_P)H'_{\bar P}(n) = (\nu -\rho_P) H'(n)$ 
for all $\nu \in \faPdc.$ Furthermore, ${}^*\rho {}^*H'(n) = {}^* \rho H'(n).$ 
Applying the estimate (\ref{e: estimate coeff MP}) to (\ref{e: formula j nu inp v}) 
we now find that 
for all $\nu \in \faPdc(P,0)$ and $v \in H_\gs^\infty,$ all $n \in N_P$ and $k_P\in K_P,$ 
\begin{eqnarray} 
\nonumber
\lefteqn{| \inp{v}{j_\nu(\kappa'(n) k_P)}_\gs|  \,e^{2 \rho_P  H'(n)} }\\
\nonumber
&\leq& 
e^{(\Re \nu +  \rho_P)\,H'(n)}(1 + |{}^*H'(n)|)^{-N} e^{-{}^*\rho {}^*H'(n)}\,q(\gs(k_P)v) \\
&\leq& e^{(\Re \nu)\,H'(n)}(1 + |{}^*H'(n)|)^{-N} e^{-\rho' H'(n)} s_\gs(v).
\label{e: estimate f of kappa prime n}
\end{eqnarray}
For the last inequality we have applied (\ref{e: deco rho prime}).
Furthermore, $s_\gs$ is a continuous seminorm on 
$H_\gs^{\infty}$ such that $q(\gs(k_P)v) \leq s_\gs (v)$ for all $v \in H_\gs^\infty$ and $k_P \in K_P;$ 
it exists by compactness of $K_P.$ 

Since $\bar P$ contains $ P_0'$ it follows by 
Lemma \ref{l: image N Q under H} below, with $P_0'$ and $\bar P$ in place of $P_0$ and $Q$ 
respectively,
that $H'(N_P)$ equals the closed convex cone $\Gamma$ spanned by 
the root vectors $H_\ga$ for $\ga \in \gS(\bar \fn_P, \fa)$ (here $H_\ga \in \fa\cap (\ker \ga)^\perp$ 
and $\ga(H_\ga) = 2$). 
If $\ga \in \gS(\bar\fn_P, \fa )$, then $\rho_P(H_\ga) < 0$ 
and for all $\nu \in \faPdc(P, R)$ we have 
$$ 
\rho_P(H_\ga)^{-1} \Re \nu(H_\ga)= \inp{\rho_P}{- \ga}^{-1} \inp{\Re \nu}{ - \ga} >
\inp{\rho_P}{- \ga}^{-1} R \geq \ge R,
$$%
where $\ge > 0$ is the minimal value of $\inp{\rho_P}{\gb}^{-1}$,  for $\gb \in \gS(\fn_P, \fa).$
It follows that
$$
\Re\nu(H) \leq \ge R \rho_P(H)
$$% 
for all $H \in \Gamma.$ This implies the existence of a constant $C_{R,N}> 0$ such that 
for all $\nu \in \faPdc(P,R)$ and all $H \in\Gamma$ we have
$$ 
e^{\Re \nu(H)} \leq C_{R,N} ( 1 - \rho_P(H))^{-N}.
$$% 
Since $-\rho_P > 0$ on $\Gamma\setminus \{0\}$ there exists a constant $\gg > 0$ such that
$$
\gg |H| \leq -\rho_P(H), \quad (H\in \Gamma).
$$% 
If $H\in \fa,$ we write $H = {}^* H + H_{ P}$ according to the (orthogonal) decomposition
$\fa = {}^*\fa + \fa_P.$ Then $|H_P|\leq |H|. $ It follows that for all $\nu \in \faPdc(P,R)$ and all $H \in\Gamma$
we have
$$ 
e^{\Re \nu(H)} \leq C_{R,N} ( 1 + \gg |H_P|)^{-N} \leq \tilde C_{R,N}( 1+ |H_P|)^{-N};
$$% 
here $\tilde C_{R,N} = C_{R,N}\sup_{t\geq1} |(1 + t)^{-1}(1 + \gg t)|^{-N}.$
Finally, we infer that for $H \in \Gamma$ and $\nu \in \faPdc(P, R)$ we have the estimate 
$$ 
e^{\Re \nu (H)}(1 + |{}^* H|)^{-N} 
 \leq \tilde C_{R,N} (1 + |H_P|)^{-N} (1 + |{}^* H|)^{-N}  \leq \tilde C_{R,N} ( 1 + |H|)^{-N}.
$$% 
 Observing that in (\ref{e: estimate f of kappa prime n}) the element $H'(n)$ belongs to $\Gamma$, we infer that for all $\nu \in \faPdc(P, R)$ and all $v \in H_\gs^\infty$ we have 
\begin{equation}
\label{e: first estimate j}
| \inp{v}{j_\nu(\kappa'(n) k_P)}_\gs|  \,e^{2 \rho_P  H'(n)} 
\leq 
\tilde C_{R,N} e^{-\rho' H'(n)} \,(1 + |\log H'(n)|)^{-N} s_\gs (v),
\end{equation}
If $N$ is sufficiently large, then the integral of the latter function
over $N_P$ is absolutely integrable, see Lemma \ref{l: image N Q under H} (b).
For such a choice of $N$ 
the function $L: \cK \to [0,\infty[$ defined by 
\begin{equation}
\label{e: defi L R}
L_R(\kappa'(n) k_P) = \tilde C_{R,N} e^{- \rho' H'(n)} (1 + |\log H'(n)|)^{-N} e^{- 2 \rho_P  H'(n)}
\end{equation} 
satsfies the required conditions.
\qed

Given a root $\ga\in \gS$ we denote by $H_\ga$ the element of $\fa$ determined by 
$H_\ga \perp \ker \ga$ and $\ga(H_\ga) = 2.$
\begin{lemma}
\label{l: image N Q under H}
Let $Q$ be a standard parabolic subgroup. Then 
\begin{enumerate}
\itema
$H(\bar N_Q)$ 
equals the cone $\Gamma(\gS( \fn_Q, \fa))$ spanned by the elements $H_\ga$ for $\ga \in 
\gS(\fn_Q, \fa).$ 
\itemb
There exists a constant $m \in \N$ such that
$$ 
\int_{\bar N_Q} e^{-\rho H(\bar n)} (1 + |H(\bar n)|)^{-m} \; d\bar n < \infty.
$$% 
\end{enumerate}
\end{lemma}

\proof
We consider the minimal parabolic subgroup $R$ of $G$ determined by 
$
N_{R} = (M_Q\cap N_0) \bar N_Q .$ 
Then it is well-known, see e.g. \cite[Lemma 4.9]{BBconv}, that
$
H(N_R \cap \bar N_0) 
$ 
equals the cone spanned by the elements $H_\ga$ for $\ga \in \gS(\bar \fn_R \cap \fn_0).$ 
Now $N_R \cap \bar N_0 = \bar N_Q$ and $\bar \fn_R \cap \fn_0 = \fn_Q$ and (a) follows.

The validity of  (b) is due to Harish-Chandra, see e.g. \cite[\S 31]{HC_HA1}.
\qed

It follows from Proposition \ref{p: integrability of j}
that for every $\gf \in C^\infty(K/ K_P:\gsP)$ and $\nu \in \faPdc(P,0)$ 
the function
$$ 
k \mapsto \inp{\gf(k)}{j(\bar P, \gs, \nu, \eta )(k)}_\gs 
$$% 
is continuous on $\cK = K \cap N_P\bar P$ and dominated by a Lebesgue integrable function, hence integrable over $K.$ 
Accordingly, we define the linear functional $\trw\!j_*(\bar P, \gs, \nu, \eta): 
C^\infty(K/K_P:\gsP) \to \C$ by
\begin{equation}
\label{e: definition j star}
\trw\!j_*(\bar P,  \gs, \nu,\eta)(\gf):= \int_K \inp{\gf(k)}{j(\bar P, \gs, \nu, \eta )(k)}_\gs\; dk.
\end{equation}
It follows from the estimate (\ref{e: first estimate j}) that
\begin{equation}
\label{e: continuity estimate j}
|\trw\!j_*(\bar P,  \gs, \nu, \eta )(\gf)| \leq  I(L_R) \,  \sup_{k \in K} s_\gs(\gf(k)), \qquad (\nu \in \faPdc(P,R)),
\end{equation}
where $I(L_R):= \int_K L_R(k)\; dk.$
In particular, we see that $\trw\!j_*(\bar P, \gs, \nu,\eta) \in C^\infty(K/K_P:\gsP)'$ 
hence equals $\inp{\dotvar}{j_*(\bar P, \gs, \nu,\eta)}$ for a unique element
$j_*(\bar P, \gs, \nu,\eta) \in C^{-\infty}(K/K_P:\gsP), $ 
for 
$\nu \in \faPdc(P,0).$ From the text subsequent to (\ref{e:  intro norm f s}) we recall 
that $C^{-\infty}(K/K_P:\gsP)$ is the union of the Banach spaces $C^{-s}(K/K_P:\gsP),$ for 
$s\in \N.$ 

 \begin{prop}
\label{p: holomorphy j in dominant region}
If $R > 0$ there exists a constant $r > 0$ such that the following assertions hold for all 
$\eta \in H_{\gs, \cchi_P}^{-\infty}.$
\begin{enumerate} 
\itema
There exists a bounded subset of $C^{-r}(K/K_P:\gsP)$ to which 
$j_*(\bar P, \gs, \nu, \eta)$ belongs 
for every  $\nu \in \faPdc(P,R).$ 
\itemb
The map $\nu \mapsto j_*(\bar P, \gs, \nu, \eta)$ 
is holomorphic as a function on $\faPdc(P,R)$ with values in the Banach space
 $C^{-r}(K/K_P: \gsP)$. 
 \end{enumerate} 
\end{prop}

\proof
Let $\eta $ and $R$ be fixed. Let $s_\gs $ and $L_R$ be as in Proposition \ref{p: integrability of j}.
Then is it readily verified that $ \gf \mapsto \sup_{k\in K} s_\gs(\gf(k))$ is a continuous seminorm 
on $C^\infty(K/K_P:\gsP).$ As the topology on the latter space is generated by the seminorms
$\|\dotvar\|_{r},$ for $r \in \N,$ see 
the text accompanying (\ref{e:  intro norm f s}), 
there exists a constant $C >0$ such that 
\begin{equation}
\label{e: estimate sup s gs gf}
\sup_{k\in K} s_\gs(\gf(k)) \leq C \|\gf\|_r
\end{equation}
for all $\gf \in C^\infty(K/K_P:\gsP).$ Thus, from the estimate (\ref{e: continuity estimate j}) it follows that  there exists a constant $r \in \N$ such that $j_*(\bar P, \gs, \nu, \eta) \in 
C^{-r}(K/K_P:\gsP)$ for all $\nu \in \faPdc(P,R).$ By linearity in $\eta$ and finite dimensionality of  
$H_{\gs, \cchi_P}^{-\infty}$ the constant $r$ may be taken the same for all $\eta,$ 
and (a) follows, with the mentioned boundedness. 

For (b) we will first show that the map $\nu \mapsto j_*(\bar P, \gs, \nu, \eta)$ 
is continuous from $\faPdc(P,R)$ to $C^{-r}(K/K_P:\gsP).$ 
Let $\nu_0 \in \faPdc(P,R)$ be fixed. Then it suffices to show that
\begin{equation}
\label{e: limit j star for nu}
\| j_*(\bar P, \gs, \nu, \eta) - j_*(\bar P, \gs, \nu_0, \eta) \|_{-r} \to 0 \qquad (\nu \to \nu_0).
\end{equation}

In the notation of the proof of Proposition \ref{p: integrability of j} we define, for $\mu \in \faPdc$ the function 
$\ge_\mu: K \to \C$ by $\ge_\mu= 0$ outside $\cK$ and by 
$$ 
\ge_\mu(\kappa'(n) k_P)= e^{\mu H'(n)} , \qquad (n \in N_P, k_P \in K_P).
$$% 
If $\mu \in \faPdc(P,0)$ then $\Re \mu \leq 0$ on the cone $\Gamma = H'(N_P)$ so that
$|\ge_\mu(k)| \leq 1$ for all $k \in K.$ From (\ref{e: formula j nu inp v}) combined with 
the definitions in the proof of Proposition \ref{p: integrability of j} it is now readily checked that, for $0 < c < 1,$
$$ 
j_*(\bar P, \gs, \nu, \eta) - j_*(\bar P, \gs, \nu_0, \eta) = \left[ \ge_{\nu - c \nu_0} - \ge_{\nu_0 - c \nu_0}
\right] j_*(\bar P, \gs, c \nu_0, \eta).
$$% 
Fix $c$ sufficiently close to $1,$ so that $c \nu_0 \in \faPdc(P, R).$ 
Then combining (\ref{e: first estimate j}) and (\ref{e: defi L R}) it follows that,
for $v \in H_\gs,$ 
$$
|\inp{\left[ \ge_{\nu - c \nu_0} - \ge_{\nu_0 - c \nu_0}
\right] j_*(\bar P, \gs, c \nu_0, \eta)}{v} | \leq \left|\ge_{\nu - c \nu_0} - \ge_{\nu_0 - c \nu_0}
\right|  L_R \cdot s_{\gs}(v).
$$%
Substituting $v = \gf(k)$, integrating over $K$ and using the estimate (\ref{e: estimate sup s gs gf}) we find,
 for all $\gf \in C^\infty(K/K_P: \gs)$ and all $\nu \in \faPdc(P, R)$ that 
$$ 
|\inp{j_*(\bar P, \gs, \nu, \eta) - j_*(\bar P, \gs, \nu_0, \eta) }{\gf}|
\leq 
C I(|\ge_{\nu - c \nu_0} - \ge_{\nu_0 - c \nu_0}| \cdot L_R) \, \|\gf\|_{r},
$$% 
where $I$ denotes the integral over $K.$ It  
follows by application of the dominated convergence theorem that 
$$ 
I(|\ge_{\nu - c \nu_0} - \ge_{\nu_0 - c \nu_0}|\cdot L_R) \to 0 , \qquad (\nu \to \nu_0).
$$% 
The continuity (\ref{e: limit j  star for nu}) now follows. 

Now that the continuity has been established, it follows by a simple application 
of the Cauchy integral formula that it suffices to prove the holomorphy of  
\begin{equation}
\label{e: holomorphy of j map}
\nu \mapsto \inp{j_*(\bar P, \gs, \nu,\eta)}{\gf}, \quad \faPdc(P,R) \to \C
\end{equation}
for a fixed $\gf \in C^\infty(K/K_P:\gsP).$ 
According to (\ref{e: definition j star}) we have that 
\begin{equation}
\label{e: integral for j gf}
\inp{j_*(\bar P, \gs, \nu, \eta)}{\gf}  = \int_K  \inp{ j(\bar P,\gs, \nu, \eta)(k)}{\gf(k)}_\gs \; dk.
\end{equation}

By Proposition \ref{p: integrability of j} (b) the integrand is Lebesgue integrable in $k$ for every 
$\nu \in \faPdc(P,0)$, holomorphic in $\nu \in \faPdc(P,R)$ for every $k \in K$ and uniformly dominated by the  Lebesgue integrable function 
$[\sup_{k \in K} s_\gs(\gf(k))] \cdot L_R$ for $\nu \in \faPdc(P,R).$ This implies   
that the integral in (\ref{e: integral for j gf}) defines a holomorphic function of $\nu \in \faPdc(P,R).$ 
\qed

From now on we will omit the $*$ in the notation of the functional defined by (\ref{e: definition j star}),
thus identifying the function $j(\bar P, \gs, \nu, \eta): K \to H_\gs^{-\infty}$ with an element 
of $C^{-\infty}(K/K_P: \gsP)$  for every $\nu \in \faPdc(P,0)$. 
For such $\nu$ we will also write
$j(\bar P, \gs, \nu)$ for the element of 
\begin{equation}
\label{e: image space j}
(H_{\gs,\cchi_P}^{-\infty})^*\otimes C^{-\infty}(K/K_P: \gsP)
\simeq \Hom( H_{\gs,\cchi_P}^{-\infty}, C^{-\infty}(K/K_P: \gsP))
\end{equation}
defined by $\eta \mapsto j(\bar P, \gs ,\nu, \eta).$ 

\begin{prop}
\label{p: bijectivity j in dominant region}
Let $\nu \in \faPdc(P,0).$ Then the following assertions are valid.
\begin{enumerate}
\itema
If $\eta \in H_{\gs,\cchi_P}^{-\infty},$ then $ 
j(\bar P,\gs,\nu,\eta) \in C^{-\infty}(G/\bar P : \gs :\nu)_\cchi.
$
\itemb
The map 
\begin{equation}
\label{e: the map j}
j(\bar P, \gs, \nu): H_{\gs,\cchi_P}^{-\infty} \to C^{-\infty}(G/\bar P:\gs:\nu)_\cchi
\end{equation}
is bijective with inverse $\ev_e.$ 
\end{enumerate} 
\end{prop}

\proof 
Let $\eta \in H_{\gs,\cchi_P}^{-\infty}.$ Then for (a) it suffices to show that 
$$ 
\pi_{\bar P, \gs, \nu}^{-\infty}(n_0) \; j(\bar P , \gs, \nu,\eta) = \cchi(n_0)\, j(\bar P , \gs, \nu,\eta),
\qquad (n_0 \in N_0).
$$% 
Let $\gf \in C^\infty(K/K_P: \gsP)$ and let $\gf_{-\bar\nu}$ denote
the unique function in $C^\infty(G/\bar P:\gs: -\bar \nu)$ which restricts to $\gf$ on $K.$ 
For $n_0 \in N_0$ 
we may write $n_0 = n_1 n_2$ with $n_1 \in M_P\cap N_0$ and $n_2 \in N_P.$ 
Accordingly, in view of Lemma \ref{l: around gd P},
\begin{eqnarray*}
\inp{\gf}{\pi_{\bar P, \gs, \nu}^{-\infty}(n_0) \; j(\bar P , \gs, \nu, \eta)}
& = & 
\inp{\pi_{\bar P,\gs , -\bar \nu }(n_0)^{-1} \gf}{ j(\bar P , \gs, \nu, \eta)}\\
& =& 
\int_{N_P} \inp{\gf _{-\bar \nu}(n_0 n)}{ j(\bar P , \gs, \nu, \eta)(n)}_\gs \; dn\\
&=& 
\int_{N_P} \cchi(n)\; \inp{ \gf _{-\bar \nu}(n_1 n_2 n)}{\eta}_\gs \; dn\\
& = & 
\int_{N_P} \cchi(n_2^{-1} n)\; \inp{\gf_{-\bar \nu} (n_1 n)}{\eta}_\gs\; dn.
\end{eqnarray*}
Using that $n_1$ normalizes $N_P$ with Jacobian $1,$ whereas 
$\cchi(n_2^{-1} n) = \cchi(n_2^{-1} n_1^{-1} n n_1),$ we infer that 
\begin{eqnarray*} 
\inp{\gf}{\pi_{\bar P, \gs, \nu}^{-\infty}(n_0) \; j(\bar P , \gs, \nu, \eta)}
&=& 
\int_{N_P} \cchi(n_2^{-1} n) \;\inp{\gf_{-\bar \nu} (n)}{\gs(n_1)\eta}_\gs\; dn\\
&=& 
\cchi(n_0)^{-1} \inp{\gf}{ j(\bar P , \gs, \nu, \eta)}.
\end{eqnarray*} 
This establishes (a). For (b) we note that from (\ref{e: defi j P gs nu eta}) and 
the definition of $\ev_e$ in (4.5) it follows that $\ev_e \after j(\bar P,\gs, \nu, \eta) = \eta.$ 
This implies that the map (\ref{e: the map j}) is injective, with $\ev_e$ as left inverse.
Therefore, $\ev_e$ is surjective onto $H_{\gs,\cchi_P}^{-\infty}.$ In view of Cor. \ref{c: ev is injective}
it now follows that $\ev_e$ is bijective with two-sided inverse $j(\bar P, \gs, \nu).$ 
\qed

\section{The Whittaker integral}
\label{s: whit int}
We will now reformulate the results of the previous section in terms of 
what Wallach \cite[\S 15.4.1]{Wrrg2} calls the Jacquet integral. Given $f \in C^\infty(K/K_P:\gs)$ we write
$f_{\bar P,\nu}$ for the unique function in $ C^\infty(G/\bar P:\gs: \nu)$ whose restriction to $K$ equals $f.$ Thus,
\begin{equation}
\label{e: extension f}
f_{\bar P, \nu}(kman) = a^{-\nu + \rho_P} \gs(m)^{-1} f(k),
\end{equation}
for $k \in K$, $(m,a,n) \in M_P\times A_P\times \bar N_P.$ 
We recall the definition of the continuous linear isomorphism $\eta \mapsto \trw\! \eta$ 
from $H_{\gs, \chi_P}^{-\infty}$ onto $\Wh_{\chi_P}(H_\gs^\infty)$ as given in 
(\ref{e: from j to wj}) with $(M_P, N_0 \cap P, \chi_P)$ in place of $(G, N_0, \chi),$ 
and with $\pi_1 = \pi_2 = \gs.$

\begin{lemma}
\label{l: intro Jacquet}
Let $\cchi$ be regular, $\gs \in \widehat M_{P{\rm ds}}$ and $\nu \in \faPdc(P,0).$ Then for all $f \in C^\infty(K/K_P: \gs)$
we have
$$ 
\inp{f_{\bar P, -\bar\nu}} {j(\bar P, \gs, \nu, \eta)} 
=
\int_{N_P} \cchi(n)\,  \trw\! \eta (f_{\bar P, -\bar \nu}(n)) \; dn,
$$% 
with absolutely convergent integral.
\end{lemma}

\proof 
By the substitution of variables used in (\ref{e: integral f over N P}) it follows 
that, with absolutely convergent integrals, 
\begin{eqnarray*}
\inp{f_{\bar P, -\bar\nu}} {j(\bar P, \gs, \nu, \eta)} & = &  \int_{K/K_P} \inp{f(k)}{j(\bar P, \gs, \gl, \eta)(k)}_\gs\; dk\\
&=& 
\int_{N_P} \inp{f(\kappa'(n)} {j(\bar P, \gs, \nu, \eta)(\kappa'(n))}_\gs \;e^{2\rho_P H'(n)} \; dn\\
&=& 
\int_{N_P} \inp{f_{\bar P, - \bar \nu}(n )}{j(\bar P, \gs, \nu, \eta)(n)}_\gs\; dn\\
&=&
\int_{N_P} \inp{f_{\bar P, - \bar \nu}(n )}{\cchi(n)^{-1}\eta}_\gs\; dn\\
&=& 
\int_{N_P}  \cchi(n)\, \trw\! \eta( f_{\bar P, -\bar \nu}(n)) \; dn.
\end{eqnarray*}
\qed

If $\Re \nu$ is $\bar P$-dominant, and $\gl \in \Wh_{\cchiP}(H_\gs^\infty)$, 
this motivates the definition of the Jacquet integral 
\begin{equation}
\label{e: Jacquet integral}
J(\bar P, \gs, \nu, \gl )(f) := \int_{N_P} \cchi(n) \gl(f_{\bar P, \nu}(n)) \;  dn,
\end{equation} 
for $f \in C^\infty(K/K_P:\gs).$  

From Lemma \ref{l: intro Jacquet} we see that 
this integral is absolutely convergent for $\bar P$-dominant $\nu \in \faPdc$ 
and defines a Whittaker functional for $C^\infty(G/\bar P: \gs :\nu).$ 
In fact, the assertion of that lemma may be reformulated as 
\begin{equation}
\label{e: comparison j and J}
\trw\! j(\bar P, \gs, \nu, \eta) = J(\bar P,\gs, -\bar\nu, \trw\!\eta).
\end{equation} 
Here the expression on the left is viewed as an element
of $C^{\infty}(K/K_P:\gs)'$ according to the compact picture, see 
(\ref{e: compact picture smooth vectors}).
It follows from Wallach's work \cite[Thm~15.4.1]{Wrrg2} that the Jacquet 
integral has a weakly holomorphic extension to a map
$\faPdc \to C^\infty(K/K_P:\gs)'.$ Furthermore, for every $\nu \in \faPdc$ 
the extension gives a linear isomorphism
$$ 
 \Wh_{\cchiP}(H_\gs^\infty) \to \Wh_\cchi(C^\infty(G/\bar P, \gs, \nu)), \;\; \xi \mapsto J(P, \gs, \gl, \xi)
$$%  
At a later stage we will strengthen this result by deriving a functional equation
for $j(P,\gs,\nu)$ and applying it to show that for every $\eta \in H_{\gs, \cchi_P}^{-\infty}$ the function $j(P,\gs, \cdot, \eta)$ extends to a  holomorphic 
function $\faPdc \to C^{-\infty}(K/K_P,\gsP)$, where the image space is equipped
with the direct limit topology. By analytic continuation of the $N_0$-equivariance of 
$j(P,\gs, \cdot, \eta)$ combined with 
Corollary \ref{c: ev is injective}, 
it then follows that for every $\nu \in \faPdc$ 
the map 
$$
j(P, \gs, \nu, \dotvar): H_{\gs,\cchi_P}^{-\infty} \to C^{- \infty}(G/\bar P, \gs, \nu)_\cchi
$$% 
is a linear isomorphism.

Given a representation $\pi$ of the discrete series of $G$, we note that by 
Lemma \ref{l: finite dimensionality whit vectors} and Corollary \ref{c: schwartz estimate for ds coeff} the linear map $\mu_\pi: H_\pi^\infty \otimes \overline{H_{\pi, \chi}^{-\infty}} \to C^\infty(G/N_0: \chi)$ given by $\mu_\pi(v\otimes \eta)(x) = \inp{\pi(x)^{-1} v}{\eta}$ for $x \in G,$ 
is actually continuous linear into $\cC(G/N_0: \chi).$ Its image is denoted by 
$ 
\cC(G/N_0: \chi)_\pi.
$ 
We denote the closure of this subspace of $L^2(G/N_0:\chi_P)$ by 
$L^2(G/N_0: \chi_P)_\pi.$ Clearly the latter space is invariant under the left regular representation
$L$ of $G.$ The following result is contained in \cite[Thm.~15.3.4]{Wrrg2}.

\begin{lemma}
\label{l: extension mu pi}
The map $\mu_\pi$ has a unique extension to a topological linear isomorphism 
\begin{equation}
\label{e: extension mu pi} 
H_\pi \otimes \overline{H_{\pi, \chi}^{-\infty}} \;\;{ \buildrel \simeq \over \longrightarrow}\;\;L^2(G/N_0: \chi)_\pi.
\end{equation} 
This extension intertwines $\pi \otimes I$ with $L.$ 
\end{lemma}

\begin{cor} 
\label{c: ortho ds subs}
If $\pi_1, \pi_2 \in \widehat G_{\rm ds}$ and $\pi_1\not\sim \pi_2$ 
 then  $L^2(G/N_0: \chi)_{\pi_1}\perp L^2(G/N_0: \chi)_{\pi_2}$ 
\end{cor}

\proof 
Let $\iota_1$ denote the inclusion of the first of the spaces into $L^2(G/N_0:\chi)$ 
and let $p_2$ denote the orthogonal projection onto the second of these spaces.
Then $p_2 \after i_1$ is a $G$-equivariant operator 
$L^2(G/N_0:\chi)_{\pi_1} \to L^2(G/N_0:\chi)_{\pi_2}.$ From Lemma \ref{l: extension mu pi} and the inequivalence of $\pi_1$ and $\pi_2$ it readily follows that $p_2 \after i_1 = 0.$ 
\qed

\begin{lemma}
\label{l: induced inner product}
Let $(\pi, H)$ be an irreducible unitary representation of a Lie group $L,$ 
and let $V$ be a finite  dimensional 
linear space. Suppose that $H \otimes V$ is equipped with a Hermitian  inner product 
$\inp{\dotvar}{\dotvar}_{H\otimes V}$ for which $\pi \otimes 1_V$ is a unitary representation of $L.$ Then there exists a unique inner product $\inp{\dotvar}{\dotvar}_V$ on $V$ such that 
\begin{equation}
\label{e: tensor prod inner prod} 
\inp{x_1\otimes v_1}{x_2\otimes v_2}_{H\otimes V} = \inp{x_1}{x_2}_{H} \,
\inp{v_1}{v_2}_{V}, \qquad (x_1, x_2 \in H, \; v_1, v_2 \in V).
\end{equation} 
\end{lemma} 

\proof
We equip $\bar H'$ with the contragredient conjugate representation $\bar \pi^\veec$ 
of $L,$ and $V$ and $\bar V'$ with the trivial representation of $L.$ 
 
 By finite dimensionality of $V,$ it readily follows that the natural map
$(A, B) \mapsto A \otimes B$ induces a linear isomorphism 
$$
\Hom_L(H, \bar H') \otimes \Hom(V, \bar V') \simeq  \Hom_L(H \otimes V, \bar H' \otimes \bar V').
$$ 
By Schur's lemma, the space on the left equals
$$ 
\C i \otimes \Hom(V, \bar V') \simeq \Hom(V, \bar V'),
$$ 
where $i: H \to \bar H', x \mapsto \inp{x}{\dotvar}_H.$ 
We consider the $L$-equivariant linear map $h: H\otimes V \to \bar H' \otimes \bar V'$ determined by
$h(x \otimes v)(y \otimes w) = \inp{x\otimes v}{y \otimes w}.$ By the above isomorphism
there exists a unique linear map $j: V \to \bar V'$ such that $i\otimes j$ is mapped onto $h.$ 
It is now readily checked that for all $v, w\in V$ we have $j(v)(w) = \overline{j(w)(v)}$ 
and $j(v)(v) > 0.$ Thus, $\inp{\dotvar}{\dotvar}_V : V \times \bar V \to \C, (v,w) \mapsto j(v)(w)$ defines a Hermitian inner product which satisfies the requirement. Conversely, if such an inner product is given then $j: v\mapsto \inp{v}{\dotvar}$ is such that $i \otimes j$ is mapped onto $h,$ 
and uniqueness of $\inp{\dotvar}{\dotvar}_V$ follows.
\qed 

We retain the notation of  Lemma \ref{l: extension mu pi}.

\begin{cor}
\label{c: mu pi isometric}
There exists a unique Hermitian inner product on $\overline{H_{\pi, \chi}^{-\infty}}$ 
such that the isomorphism (\ref{e: extension mu pi}) becomes an isometry.
\end{cor} 

\proof 
Use Lemma \ref{l: finite dimensionality whit vectors} and apply Lemma \ref{l: induced inner product} with $L= G,$ $V = \overline{H^{-\infty}_{\pi, \chi}} \simeq \Wh_\chi(H_\pi^\infty)$ 
and with $H_\pi \otimes V$ equipped with the pull-back
of the $L^2$-inner product under $\mu_\pi.$
\qed

From now on, we assume that the finite dimensional spaces of Whittaker vectors
$H_{\gs, \chi_P}^{-\infty},$ for $\gs \in \MPds,$ are equipped with the Hermitian inner products satisfying the assertion of Corollary \ref{c: mu pi isometric}. 

We are now prepared to introduce Harish-Chandra's  Whittaker integral, which is an 
appropriate analogue of the Eisenstein integral for groups and symmetric spaces.
Let $\tau$ be a unitary representation of $K$ in a finite dimensional complex Hilbert space $V_\tau.$ 

We write $C^\infty(\tau: G/N_0 : \chi)$ for the space of smooth functions
$f: G\to V_\tau$ such that $f(k x n) = \chi(n)^{-1} \tau(k) f(x),$ for all $x \in G, k\in K, n\in N_0.$ 
Via the inverse of the natural isomorphism $C^\infty(G)\otimes V_\tau \to C^\infty(G, V_\tau)$ we have
$$ 
C^\infty(\tau: G/N_0 : \chi) \simeq (C^\infty(G/N_0:\chi) \otimes V_\tau)^K.
$$ 
Accordingly, we define the associated space of $\tau$-spherical Whittaker Schwartz functions by 
$$ 
\cC(\tau: G/N_0: \chi) : = C^\infty(\tau: G/N_0: \chi) \,\cap [\cC(G/N_0: \chi)\, \otimes \Vtau]
$$ 
Furthermore, we put 
$$ 
\cA_2(\tau: G/N_0:\chi)_\pi: = \cC(\tau: G/N_0: \chi) \cap [\cC(G/N_0: \chi)_\pi \otimes \Vtau].
$$ 
and   
\begin{equation}
\label{e: defi cA two G}
 \cA_2(\tau: G /N_0: \chi):= \oplus_{\pi\in \widehat G_{\rm ds}} \cA(\tau : G/N_0: \chi)_\pi.
\end{equation}
Note that this direct sum has only finitely many non-zero terms, since 
$(V_\tau \otimes H_\pi)^K \neq 0$ for only finitely many $\pi \in \widehat G_{\rm ds}.$ 
Moreover, each of the components is finite dimensional since $\pi$ is admissible and
$H_{\pi,\chi}^{-\infty}$ is finite dimensional. It follows that the space
(\ref{e: defi cA two G}) is finite dimensional. In particular, for each $\pi$ the corresponding 
summand is a closed
subspace of $L^2(G/N_0: \chi)_\pi \otimes \Vtau$. In view of Corollary \ref{c: ortho ds subs} 
it follows that (\ref{e: defi cA two G}) is a finite orthogonal direct sum of finite dimensional
subspaces of the Hilbert space $L^2(G/N_0: \chi)) \otimes \Vtau.$ 
Accordingly, we equip the space $\cA_2(\tau: G /N_0: \chi) $ with the restricted
Hilbert structure.

\begin{rem}
It is a result of both Harish-Chandra \cite{HCwhit} and Wallach \cite{Wrrg2} that 
the space (\ref{e: defi cA two G}) equals the space of $\fZ$-finite functions in $\cC(\tau: G/N_0: \chi).$
Equivalently, this means that the irreducible unitary representations which appear discretely in 
$L^2(G/N_0: \chi)$ belong to $\widehat G_{\rm ds}.$ However, we shall not 
need this in the present paper.
\end{rem}

Let 
$P = M_PA_PN_P$ be a standard parabolic subgroup of $G.$ 
We recall that $\chi_P:= \chi|_{M_P \cap N_0}$ is regular relative to $(M_P, M_P \cap N_0),$ 
put $\tau_P= \tau|_{K_P}$ 
and define the finite dimensional space 
\begin{equation}
\label{e: defi cAtwoP}
\cAtwoP:= \cA_2 (\tau_P: M_P/ M_P\cap N_0:\chi_P)
\end{equation}
as above. Then 
\begin{equation}
\label{e: deco cA P in ds types}
\cAtwoP = \oplus_{\gs \in \MPds} \cA_2 (\tau_P: M_P/M_P\cap N_0:\chi_P)_\gs.
\end{equation}
To keep notation manageable we will denote the summands by $\cAtwoPgs.$ 
Since $\tau$ will be kept fixed, this will not cause any ambiguity.

\begin{rem}
Note that for $P = P_0$ minimal we have $M_P \cap N_0 = \{e\}$ so that 
$$ 
\cA_{2, P_0} = C^\infty(\tau_M: M) \simeq \Vtau.
$$% 
At the other extreme, for $P =G$ we have $M_P = {}^\circ G,$ so that $M_P \cap N_0 = N_0$ 
and 
$$ 
\cA_{2, G} = \oplus_{\gs \in \widehat{{}^\circ G}_{\rm ds}} \cC(\tau \col {}^\circ G/N_0 \col \chi)_\gs.
$$
\end{rem}

We return to the setting of a parabolic subgroup $P$ containing $A.$ 
For $\gs\in \MPds$ we define $C^\infty(\tau_P\col K/K_P \col \gs)$ 
to be the space of smooth functions $\gf: K \to H_\gs^\infty \otimes \Vtau$ 
such that
\begin{equation}
\label{e: tau gs behavior} 
\gf(k_1 k m)= [ \tau(k_1) \otimes \gs(m)^{-1}] \gf(k), \qquad (k\in K, \; k_1,m \in K_P).
\end{equation}
We equip this space with the pre-Hilbert structure induced by the $L^2$-inner product on $L^2(K, H_\gs \otimes \Vtau),$ with respect to the Haar measure $dk$ on $K$ normalized
by $\int_K dk =1.$ 
 
For a finite subset $\vartheta \subset \widehat{K_P}$ we denote by 
$H_{\gs, \vartheta}$ the sum of the $K_P$-isotypical components of $H_\gs$ for the 
$K_P$-types in $\vartheta.$ 
We note that (\ref{e: tau gs behavior}) implies that 
$\gf(e) \in (H^\infty_\gs \otimes \Vtau)^{K_P}\subset H_{\gs,\vartheta} \otimes \Vtau,$ 
with 
$$ 
\vartheta = \{\gd\mid \Hom_{K_P}(\gd^\veec, \tau) \not= 0\}.
$$ 
 By sphericality 
this implies that $C^\infty(\tau_P\col K/K_P \col \gs)$ equals the space of 
smooth $\gf: K \to H_{\gs, \vartheta} \otimes \Vtau$ such that (\ref{e: tau gs behavior}).
In particular it is finite dimensional, hence Hilbert for the given pre-Hilbert structure.

We  define, for $\gs\in \MPds$ and $T = f \otimes \eta \in 
C^\infty(\tau: K/K_P: \gsP) \otimes \overline{ H_{\gs,\chi_P}^{-\infty}}$, the function 
$\psi_{T}: M_P \to V_\tau$ 
by 
\begin{equation}
\label{e: formula for psi T}
\psi_T(m) = \gg \after (\trw\!\eta \otimes I)\after (\gs(m)^{-1} \otimes I)(f(e)).
\end{equation}
Here $f$ is viewed as a function with values in $H_\gs^\infty \otimes V_\tau$,
$\trw\!\eta \in \Wh_\chi(H^\infty_\gs)$ is defined by $v\mapsto \inp{v}{\eta}$ 
and $\gg$ denotes the canonical linear map $\C \otimes \Vtau \to \Vtau.$  
It is readily verified that $\psi_T \in \cA_{2}(\tau_P: M_P/M_P\cap N_0: \cchiP)_\gs.$ 

\begin{lemma}
The linear map  $T \mapsto \psi_T$ is an isometric linear isomorphism
\begin{equation}
\label{e: iso by psi T}
C^\infty(\tau: K/K_P :\gs) \otimes \overline{ H_{\gs, \chi_P}^{-\infty}}\;\;{ \buildrel \simeq \over \longrightarrow}\;\; \cA_{2}(\tau_P: M_P/M_P\cap N_0:  \cchiP)_\gs.
\end{equation}
\end{lemma}
This is analogous to a result of Harish-Chandra in the case of the group, see
\cite[Lemmas 7.1, 9.1]{HCha3}. It is also analogous to \cite[Lemma 4.1]{Bps2} in the setting of 
symmetric spaces. 
\vspace{-8pt}
\medno
\proof
We equip $L^2(K, \Vtau)$ with the natural $L^2$-inner product corresponding to the 
fixed normalized Haar measure $dk.$ By restriction this induces an inner product 
on $C^\infty(\tau: K).$ Clearly, the map 
$C^\infty(\tau: K) \to V_\tau$, $f \mapsto f(e),$ is 
a linear isomorphism which is $K$-equivariant  for $R$ and $\tau.$ Furthermore,
for $f,g \in C^\infty(\tau: K)$ we have, by sphericality,
$$ 
\inp{f}{g} = \int_K \inp{f(e)}{g(e)}\;dk = \inp{f(e)}{g(e)}.
$$ 
Thus $f \mapsto f(e)$ defines an isometric linear isomorphism $C^\infty(\tau: K) \;{\buildrel \simeq
\over\to} \;V_\tau.$ 

It now follows that 
\begin{eqnarray*}
C^\infty(\tau: K/K_P :\gs) \otimes \overline{ H_{\gs,\cchiP}^{-\infty}}
& = & 
[C^\infty(\tau: K) \otimes H_\gs]^{K_P} \otimes \overline{ H_{\gs,\cchiP}^{-\infty}}\\
&\simeq& [V_\tau \otimes H_\gs]^{K_P}\otimes \overline{ H_{\gs,\cchiP}^{-\infty}}\\
&\simeq & [ V_\tau \otimes \cC(M_P/M_P\cap N_0)_\gs]^{K_P}\\
&= & \cA_{2, P}(\tau_P : M_P/M_P\cap N_0)_\gs.
\end{eqnarray*}

In the above array the identity signs indicate isometric isomorphisms via which spaces are naturally identified. 
The composition of the first two isomorphisms is given by $f \otimes \eta \mapsto f(e)\otimes \eta.$
By what we said in the above, this is an isometric isomorphism.
The application of the third isomorphism maps $f(e) \otimes \eta$ to the function $M_P \to \Vtau$ 
given by  $m \mapsto \trw \eta(\gs(m)^{-1}f(e)) = \inp{\gs(m)^{-1}f(e)}{\eta},$  which gives an isometric isomorphism 
in view of Corollary \ref{c: mu pi isometric}. From these descriptions it follows
that the composition of the isomorphisms in the array is isometric and gives $T \mapsto \psi_T.$ 
\qed

We now assume that $P$ is a standard parabolic subgroup of $G.$ 
For $\psi \in \cAtwoP,$ see (\ref{e: defi cAtwoP}),
 we define the associated Whittaker integral by 
\begin{equation}
\label{e: defi whittaker int}
\Wh(P, \psi, \nu) (x) = \int_{N_P} \cchi(n) \;\psi_{\bar P, -\nu}(xn)\; d \bn,
\end{equation}
where $\nu \in \faPdc$ and where $ \psi_{\bar P, -\nu} \in C^\infty(G/\bar P,  \gs,  -\nu) \otimes \Vtau$ is defined by 
$$
\psi_{\bar P, -\nu} (k m a \bar n) = a^{\nu + \rho_P} \tau(k) \psi(m),
$$%
for $k \in K, (m,a,\bar n) \in M_P\times A_P \times \bar N_P.$ 
This is precisely the definition given by Harish-Chandra, \cite[\S 1.7, p.147]{HCwhit}.
By rewriting this integral in terms of the Jacquet integral, we will see 
that it converges absolutely for $\nu \in \faPdc$ with $\inp{\Re \nu}{\ga} > 0$ for all $\nu \in \gS(P,\faP).$ 

The Whittaker integral can be related to matrix coefficients, hence to the Jacquet integral, 
as follows.
For $\gs$ a discrete series representation of $M_P$ and 
$T = f \otimes \eta \in C^\infty(\tau: K/K_P: \gsP) \otimes \overline{ H_{\gs,\chi_P}^{-\infty}}$, let 
$\psi_{T}: M_P \to V_\tau$ be defined
as in (\ref{e: formula for psi T}).

We note that $\gg \after (J(P,\gs,\nu,\gl) \otimes I_{\Vtau})$ defines a continuous linear map 
from $C^\infty(\tau: K/K_P: \gs)$ to $\Vtau$ which we shall denote by 
$$
J(P,\gs,\nu,\gl)_\tau: C^\infty(\tau: K/K_P: \gs) \to \Vtau.
$$ 
Accordingly, we have the following relation of the Whittaker integral 
with the Jacquet integral.

\begin{lemma} 
\label{l: whittaker integral as matrix coeff} 
Let $P =M_P A_P N_P$ be standard and $\gs \in \MPds.$
Let $f \in C^\infty(\tau:K/K_P:\gs)$ and $\eta \in H_{\gs, \chi_P}^{-\infty}.$ 
If $\nu \in \faPdc(P,0)$, then 
$$ 
\Wh(P, \psi_{f \otimes \eta} \,,\nu )(x) = J(\bar P, \gs,-\nu, \trw\!\eta)_\tau (\pi_{\bar P, \gs, -\nu}(x)^{-1} f),
$$% 
with absolutely convergent integral for the Whittaker integral on the left.
Here we have abused notation, by writing $\pi_{\bar P, \gs, -\nu}(x)$ 
for $ \pi_{\bar P, \gs, -\nu}(x) \otimes I_{V_\tau}.$ 
\end{lemma}
\proof 
We put $\psi:= \psi_{f\otimes \eta}$ 
and define $\psi_{\bar P, -\nu} :  G \to  \Vtau$ 
by 
$$ 
\psi_{\bar P, -\nu}(kman): =  a^{\nu + \rho_P}  \tau(k)  \psi(m).
$$% 
Furthermore, we define $f_{\bar P, -\nu}: G \to H_\gs \otimes V_\tau$ 
by 
$$ 
f_{\bar P, -\nu}(kman): = a^{\nu + \rho_P} \,  [\gs(m)^{-1} \otimes I]  f(k).
$$% 
Then 
\begin{eqnarray*}
\psi_{\bar P, -\nu}(kman) &=& a^{\nu + \rho_P}  \tau(k) \gg [\trw\!\eta \after \gs(m)^{-1}\otimes I]f(e) \\
&=& 
 \gg[\trw\! \eta \otimes I]  a^{\nu + \rho_P} [ \gs(m)^{-1} \otimes \tau(k)] f(e)\\
 &=& 
\gg [\trw\!\eta \otimes I]  a^{\nu + \rho_P} [ \gs(m)^{-1} \otimes I] f(k)\\
 &=& 
 \gg [\trw\!\eta \otimes I] f_{\bar P, -\nu}(kman).
\end{eqnarray*}
This in turn implies that 
\begin{equation}
\label{e: translates of psi and f} 
L_{x^{-1}} \psi_{\bar P, -\nu} = \gg [\trw\eta  \otimes I] 
([\pi_{\bar P, \gs , -\nu}(x^{-1} ) \otimes I]f)_{\bar P, -\nu}.
\end{equation}
The function on the right-hand side is integrable over $N_P$ with integral
$$
J(\bar P, \gs, -\nu, \trw\!\eta)_\tau ([\pi_{\bar P, \gs , -\nu}(x^{-1} ) \otimes I]f)
$$% 
(abusing notation). 
It follows that the function on the left-hand side of (\ref{e: translates of psi and f}) 
is also integrable over $N_P$, with integral being equal to $\Wh(P, \psi_{f\otimes \xi},\nu )(x),$
see (\ref{e: defi whittaker int}).
\qed

\begin{cor}
\label{c: Wh as matrix coeff with j}
Let the setting be as in Lemma \ref{l: whittaker integral as matrix coeff}. Then  
$$ 
\Wh(P, \psi_{f \otimes \eta}\, ,\nu )(x)  = 
\inp{\pi_{\bar P, \gs, -\nu}(x)^{-1} f}{j(\bar P, \gs , \bar \nu, \eta)}.
$$%
Here we have abused notation in the expression on 
the left, by suppressing 
trivial actions on  the tensor component $\Vtau$ 
and the role of the canonical isomorphism $\gg: \C \otimes \Vtau \to \Vtau.$  
\end{cor}

\proof 
This follows from Lemma \ref{l: whittaker integral as matrix coeff} combined with
(\ref{e: comparison j and J}).
\qed

\begin{cor}
\label{c: initial holomorphy Wh}
Let $P =M_PA_P N_P$ be standard.  
If  $\psi \in \cA_2(\tau_P: M_P/M_P\cap N_0: \cchiP)$, 
then for every $\nu \in \faPdc(P,0)$ we have  
$$ 
\Wh(P, \psi,\nu ) \in C^\infty(\tau: G/N_0\col\cchi).
$$ 
Furthermore, $\nu \mapsto \Wh(P ,\psi,\nu )$ is a holomorphic 
function on $\faPdc(P,0)$ with values in $C^\infty(\tau: G/N_0\col\cchi).$ 
\end{cor}

\proof 
By decomposition (\ref{e: deco cA P in ds types}) and linearity,
we may assume that $\psi$ belongs to the space $\cA_2(\tau_P: M_P/M_P\cap N_0: \cchiP)_\gs$ 
with $\gs$ a representation of the discrete series of $M_P.$ In view of the isomorphism
we may further assume that $\psi = \psi_{f \otimes \eta},$ with $f$ and $\eta$ as in Lemma
\ref{l: whittaker integral as matrix coeff}. The result now follows by application of 
Corollary \ref{c: Wh as matrix coeff with j} and Proposition \ref{p: holomorphy j in dominant region}
(recall that in the compact picture, $j_*$ is $j$ viewed as an element of $C^{-\infty}(K/K_P:\gsP)$). 
\qed

For $Z \in \fZ$ we note that the endomorphism $R_Z$ of $C^\infty(G)$ 
leaves the subspace $C^\infty(G/N_0 \col\cchi)$ invariant 
and induces a differential operator on that space, viewed as the space
of smooth section of the associated bundle $G\times_{N_0}\C_\cchi.$ This differential 
operator is denoted $R_Z$ as well.
The associated endomorphism $I \otimes R_Z$ of $\Vtau\otimes C^\infty(G/N_0\col \cchi)$ 
restricts to an endomorphism of $C^\infty(\tau: G/N_0\col \cchi).$ 
In a similar fashion we may equip the latter space with a left action $I \otimes L_Z.$ 
Since $I \otimes L_{Z^\veec} = I \otimes R_Z,$ we see that the latter operator 
preserves $\cA_{2}(\tau: G/N_0\col \cchi).$ 

Let $P \in \cP(A).$ We agree to equip $\cA_{2,P}$ with the structure 
of $\fZ(\fm_P)$-module induced by the right regular representation of
$\fm_P$ on $C^\infty(M_P),$ as in the preceding text with $G$ replaced by $M_{P}.$

Let 
$$ 
\mu_P: \fZ \to \fZ(\fm_{1P})
$$ 
be the canonical embedding. The decomposition $\fm_{1P} = \fm_P \oplus \fa_P$ 
induces the canonical isomorphisms
$$ 
\fZ(\fm_{1P}) \simeq  S(\fa_P) \otimes \fZ(\fm_P)  \simeq P(\faPd) \otimes \fZ(\fm_P).
$$ 
Thus, if $Z \in \fZ$ then $\mu_P(Z)$ may be viewed as a polynomial
function on $\faPdc$ with values in $\fZ(\fm_P).$ Accordingly, for 
$\nu \in \faPdc$ we put $\mu_P(Z, \nu): = \mu_P(Z)(\nu).$ 
We agree to denote by $\umu_P(Z,\nu)$ the endomorphism
by which $\mu_P(Z,\nu)$ acts on $\cA_{2,P}.$ Then $\umu_P(\dotvar ,\nu)$ may be viewed as 
an algebra homomorphism $\fZ \to \End(\cA_{2,P}),$ with polynomial dependence on $\nu.$ 

\begin{lemma}
Let $P \in \cP(A)$ be standard, $Z \in \fZ,$ $\psi \in \cAtwoP.$ Then for every $\nu \in \faPdc(P,0)$ we have 
\begin{equation}
\label{e: RZ on Wh}
R_Z \Wh(P,\psi ,\nu) = \Wh(P, \umu_P(Z,\nu)\psi, \nu).
\end{equation}
\end{lemma}

\proof
By linearity it suffices to fix a representation $\gs$ of the discrete series of $M_P$  
and to establish the identity for $\psi = \psi_{f \otimes \eta},$ with $f \in C^\infty(\tau_P:K/K_P:\gs)$ and $\eta \in H_{\gs,\chi_P}^{-\infty}.$ 
From Corollary \ref{c: Wh as matrix coeff with j} it follows that
\begin{eqnarray}\nonumber
R_Z \Wh(P, \psi, \nu)(x) &=& 
 \inp{\pi_{\bar P, \gs, -\nu}(x)^{-1} \bp\! f}{j(\bar P, \gs, \bar\nu, \eta)}\\
 &=& \Wh(P, \psi_{\bp\! f \otimes \eta}, \nu)(x), 
\label{e: RZ Wh and bp f}
\end{eqnarray}
with 
$$ 
\bp\! f (k) = \pi_{\bar P, \gs, -\nu}(Z^\veec) f (k) =  f_{\bar P, \gs, -\nu}(k; Z) =
\gs(\mu_{\bar P}(Z, \nu)^\veec) f(k) .
$$ 
From the definition of $\psi,$ see (\ref{e: formula for psi T}) it follows that
$$ 
\psi_{\bp\! f \otimes \xi} = R_{\mu_{\bar P}(Z, \nu)} \psi
= R_{\mu_P(Z, \nu)} \psi = \umu_P(Z, \nu) \psi.
$$ 
Substituting this in (\ref{e: RZ Wh and bp f}), we obtain (\ref{e: RZ on Wh}).
\qed

As mentioned in the introduction, the main purpose of the present paper is to show that the Whittaker integrals
extend holomorphically in the variable $\nu \in \faPdc$ and, for imaginary $\nu,$ satisfy estimates of a uniformly 
tempered type.

A first step into this direction is the following estimate, for $\Re \nu$ $P$-dominant.
The proof given below corresponds to the proof in \cite[Lemma 9.22.1]{HCwhit}.

\begin{lemma}
\label{l: first estimate Wh integral}
For every $\psi \in \cAtwoP$ there exists a constant $m > 0$ and for every $R > 0$ a constant
$C>0$ such that for all $a \in A$ 
and all $\nu \in \fadc(P,R),$ 
$$
\|\Wh(P, \psi, \nu)(a)\|_\tau \leq  C(1 + |\log a|)^m a^{\Re \nu -\rho}.
$$
\end{lemma}

\proof
Since $\psi \in \cAtwoP$ there exists for every $m >0$ a constant $C_m> 0$ such that
\begin{equation}
\label{e: estimate for psi}
\|\psi ({}^*\! a)\| \leq C_m (1 + |\log {}^*\! a|)^{-m} ({}^*\! a)^{-{}^*\!\rho_P},
\end{equation}
for ${}^*\!a \in {}^*A_P:= M_P\cap A.$ 

In the following we write $\psi_{-\nu} = \psi_{\bar P, -\nu}.$ 
Furthermore, we write $a = {}^*\!a  a_P$ according to the decomposition $A = {}^*\!A_P A_P.$  Then from (\ref{e: defi whittaker int}) with $x =a,$ we find,
by substituting $a n a^{-1}$ for $n$ that
$$ 
 \Wh(P, \psi, \nu)(a) =  a_P^{\nu - \rho_P}\int_{N_P}  \cchi(a^{-1}n a)) \psi_{-\nu}(n \,{}^*\!a ) dn.
$$ 
Let $P_0'$ be the minimal parabolic subgroup $MA (M_P\cap N_0) \bar N_P$ 
as in (\ref{e: primed Iwasawa deco}) and let $\kappa', h', n'$ be the
projection maps for the associated Iwasawa decomposition. Decomposing
$n = \kappa'(n) h'(n) n'(n)$ and $h'(n) = {}^*h'(n) h'_P(n)$ according to $A = {}^*A_P A_P,$ 
as well as $n'(n) = {}^*n'(n) \bar n_P'(n)$ according to $N_{P_0'} = (M_P \cap N_0) \bar N_P,$ 
we find
\begin{eqnarray*}
\psi_{-\nu}(n {}^* a) & =& \tau(\kappa'(n))\, \psi_{-\nu}(h'(n)n'(n){}^*\! a) \\
&=&  \tau(\kappa'(n))\, \psi_{-\nu}({}^*\!h'(n) {}^*\!n'(n) h_P'(n) {}^*\!a)\\
 &=& h_P'(n)^{\nu + \rho_P}\, \cchi({}^*\! a^{-1} {}^*\! n'(n){}^*\! a)^{-1}\,\tau(\kappa'(n))\, \psi({}^*\!h'(n) {}^*\! a).
\end{eqnarray*}
In view of the unitarity of $\tau$ and $\cchi$ this
leads to the estimate
\begin{equation}
\label{e: first estimate Wh pos nu}
\|\Wh(P, \psi, \nu)(a) \|\leq   
a_P^{\Re \nu -\rho_P} 
\int_{N_P} \|\psi({}^*\!h'(n) {}^*\! a)\|\; h_P'(n)^{\Re \nu + \rho_P} dn.
\end{equation}

From (\ref{e: estimate for psi}) it follows, taking account that
$a_P^\nu = a^\nu$ and $a_P^{\rho_P}({}^*\!a)^{{}^*\rho}  = a^\rho,$  that 
\begin{eqnarray*} 
\lefteqn{a_P^{\Re \nu -\rho_P} \|\psi({}^*\!h'(n) {}^*\! a)\|  \, h_P'(n)^{\Re \nu + \rho_P} }\\
& \leq & C_m \, a^{\Re \nu - \rho}(1 + |\log {}^*\!a|)^m (1 + |\log {}^*\!h'(n)|)^{-m} h'(n)^{- \rho'}h_P'(n)^{\Re \nu}.
\end{eqnarray*}
Applying Lemma \ref{l: image N Q under H} (a) with $P_0'$ in place of $P_0$ and 
$\bar P$ in place of $Q,$ we see that the image $\log \after h'_P(N_P)$ is contained in the cone
in $\fa_P$ spanned by the elements $- \pr_P H_\ga,$ for $\ga \in \gS(\fn_P, \fa);$ here
$\pr_P$ denote the orthogonal projection $\fa \to \fa_P.$  This implies
that there exists a constant $C_{m,R} > 0$ such that, for all $\nu \in \faPdc(P,R),$ 
$$ 
 h_P'(n)^{\Re \nu} \leq C_{m,R} ( 1+ | \log h'_P(n)|)^{-m}, \qquad (n \in N_P).
$$ 
Therefore,
\begin{eqnarray}
\lefteqn{a_P^{\Re \nu -\rho_P} \|\psi({}^*\!h'(n) {}^*\! a)\|  \, h_P'(n)^{\Re \nu + \rho_P}}\nonumber \\ 
& \leq &  C_m C_{m,R}\, a^{\Re \nu - \rho}(1 + |\log {}^*\!a|)^m 
(1 + |\log h'(n)|)^{-m} h'(n)^{- \rho'}.\label{e: estimate psi of star}
\end{eqnarray}
In view of Lemma \ref{l: image N Q under H} (b) we may fix $m>0$ such that 
\begin{equation}
\label{e: conv int m}
I_m:=  \int_{N_P} (1 + |\log h'(n)|)^{-m} h'(n)^{- \rho'}\; dn < \infty.
\end{equation} 
Combining (\ref{e: first estimate Wh pos nu}) with (\ref{e: estimate psi of star}) and (\ref{e: conv int m})
we find that, for $\nu \in \faPdc(P,R)$ and $a \in A,$ 
$$ 
\| \Wh(P, \psi, \nu)(a) \| \leq C_m C_{m,R} I_m \; a_P^{\Re \nu -\rho_P} (1 + |\log{}^*\!a |)^m.
$$ 
\qed

\section{Finite dimensional spherical representations}
\label{s: fin dim reps}
We assume that $\fh$  is a $\Cartan$-stable Cartan subalgebra of $\fg$ 
containing $\fa.$  Then $\fh = \ft \oplus \fa$ with $\ft$ a maximal torus in $\fm.$ 
Accordingly, we may naturally identify $\fadc$ with 
the space of $\gl \in \fhdc$ such that $\gl|_\ft = 0.$ 

The recall the definition of $B$ from (\ref{e: intro of B}), and denote its restriction to 
$\fh$ as well as the complex bilinear 
extension to $\fh_\iC$ by $\Rinp{\dotvar}{\dotvar}.$ The latter restricts to a positive definite inner product on $\fh_\iR:= i \ft \oplus \fa.$ Its complexified dual,
denoted by $\Rinp{\dotvar}{\dotvar}$ as well, is positive definite on $\fh_\iR^*:= i \ft^* + \fa^*.$ 
The restriction of $\Rinp{\dotvar}{\dotvar}$ to $\fa^*$ coincides with the dual of the restriction
of the inner product $\inp{\dotvar}{\dotvar}$ defined by (\ref{e: defi inp}).

We denote by $R(\fh)\subset \fhdc$ the root system of $\fh$ in $\fg_\iC$ 
and select a positive system $R^+(\fh)$ which is compatible with $\gS^+.$
The latter means that if $\ga \in R(\fh)$ and  $\ga|_\fa \in \gS^+$ 
then $\ga \in R^+(\fh).$ 

Let $\gL(\fh)$ denote the collection of weights of the pair $(\fg_\iC, \fh_\iC)$, 
i.e., the collection of $\gl \in \fhdc$ such that $2\Rinp{\gl}{\ga}/\Rinp{\ga}{\ga} \in \Z$ 
for all $\ga \in R(\fh).$ 
Let $\gL^+(\fh) \subset \fhdc$ be the associated
collection of dominant weights, i.e., the weights $\gl \in \gL(\fh)$ such that $\Rinp{\gl}{\ga} \geq 0 $ for all $\ga \in R^+(\fh).$ 

By the Cartan--Helgason classification 
\cite[Ch. 5, Thm. 4.1]{HelgasonGGA}, a finite dimensional irreducible 
representation $\pi$ of $G$ is spherical, i.e., has a $K$-fixed vector, if and only if $M$ acts trivially on its highest weight space. Furthermore, the latter condition implies that the highest weight of $\pi$ is an element of the set
$$ 
\gL^+(\fa) = \{\mu \in \fadc \mid \forall \ga\in \gS^+: \frac{\Rinp{\mu}{\ga}}{\Rinp{\ga}{\ga}} \in \N\}.
$$ 
Conversely, $\gL^+(\fa) \subset \gL^+(\fh)$ and if $\mu \in \gL^+(\fa)$, then up to equivalence there is a unique spherical representation of $G$ of highest weight $\mu.$

In \cite{HelgasonGGA} these results are proven for $G$ connected semisimple with 
finite center. The extension of this result to groups of the Harish-Chandra class is 
straightforward. Given an element $\mu \in \gL(\fa)$ we denote by 
$\pi_\mu$ the associated irreducible spherical representation of $G.$ 
  
The following result  is well-known.
\begin{lemma}
\label{l: double weight lattice}
$2 (\gL^+(\fh) \cap \fadc) \subset \gL^+(\fa).$
\end{lemma}
\proof
Let $\ga \in \gS$ and let $\tilde \ga \in R(\fh)$ be such that  $\ga = \tilde \ga|_{\fa}.$ Then 
$\Rinp{\tilde \ga}{\tilde \ga} = m \Rinp{\ga}{\ga}$, for a certain $m \in \{1,2,4\};$ if $m =4$ then $2\ga \in \gS$, see \cite[Ch. VII, Lemma 8.4]{HelgasonDS}.
Let $\gl \in\gL^+(\fh) \cap \fadc.$ Then $\Rinp{\gl}{\ga} = \Rinp{\gl}{\tilde \ga}$, so that 
$$ 
\frac{\Rinp{2 \gl}{\ga}}{\Rinp{\ga}{\ga}} = 2m \frac{\Rinp{\gl}{\tilde \ga}}{\Rinp{\tilde \ga}{\tilde \ga}} \in m\N \subset \N.
$$ 
\qed 

In the rest of this section we assume that $P$ is a standard parabolic subgroup of $G$. We write ${}^*\fh_P = \fh \cap \fm_P.$ Then ${}^*\fh_P$ 
is a real $\Cartan$-stable Cartan subspace of $\fm_P$, which decomposes as 
${}^*\fh_P = \ft \oplus (\fa \cap \fm_P)$ Note that $\fh = {}^*\fh_P \oplus \fa_P$, so that we may identify
${}^*\fh_{P\iC}^*$ and $\faPdc$ with subspaces of $\fhdc.$ We denote by $R({}^*\fh_P)$ the root system 
of ${}^*\fh_P$ in $\fm_{P\iC}.$ Then $R({}^*\fh_P)$ consists of the roots in $R(\fh)$ which vanish on $\fa_P.$ 
Furthermore, $R^{+}({}^*\fh_P) = R({}^*\fh_P) \cap R^+(\fh)$ is a positive system. The associated
weight lattice  is denoted by $\gL({}^*\fh_P)$ and the subset of dominant
ones by $\gL^+({}^*\fh_P).$ 

Via the decomposition $\fa = (\fa \cap \fm_P) \oplus \fa_P$ we view 
$\faPdc$ as the linear subspace of $\fadc$ consisting of all 
$\mu \in \fadc$ that  vanish on $(\fa \cap \fm_P).$  Accordingly, we define
$$ 
\gL^+(\fa_P):= \gL^+(\fa) \cap \faPdc. 
$$

\begin{lemma}
\label{l: star A on hw}
Let $\mu \in \gL^+(\fa)$ and let $\pi_\mu$ be the irreducible spherical representation of $G$ of highest weight $\mu.$ Then the following assertions are equivalent.
\begin{enumerate}
\itema $M_P$ acts trivially on the highest weight space of $\pi_\mu;$ 
\itemb $\mu \in \gL^+(\fa_P).$ 
\end{enumerate}
\end{lemma}

\proof
Let $F$ be a finite dimensional complex linear space on which $\pi_\mu$ is realized. Let $e_\mu \in F_\mu \setminus \{0\}$ be a non-zero highest weight vector.

Assume (a). Then $A \cap M_P$ acts trivially $e_\mu$, hence $\mu = 0 $ on $\fa \cap \fm_P$, which implies (b). 

For the converse, assume (b). Then $M_P = (M_P)_e M$, so that it suffices to show that $\fm_{P\iC}$ annihilates $e_\mu.$ Let $\fb$ be the Borel subalgebra of $\fm_{P\iC}$ determined by the positive system $R^+({}^*\fh_P).$ 
Then $\fb$ is contained in $ \fm +  (\fa \cap \fm_P) + \fn_0$ hence annihilates $e_\mu.$ 
This implies that $U(\fm_P) e_\mu$ is a finite dimensional cyclic highest weight $\fm_{P\iC}$-module of highest weight $0.$
Therefore, $U(\fm_P) e_\mu = \C e_\mu$ from which we obtain $\fm_P e_\mu = 0.$ 
\qed

We define 
\begin{equation}
\label{e: gL double plus}
\gL^{++}(\fa_P):= \{ \mu \in \gL^+(\fa_P)\mid \forall \ga \in \gS(\fn_P, \fa_P):\;\;  
\Rinp{\mu}{\ga} > 0\}.
\end{equation} 
The following lemma guarantees in particular that the set (\ref{e: gL double plus}) is 
non-empty.
\begin{lemma} 
\label{l: property 4 rho P}
The element $4 \rho_P$ belongs to $\gL^{++}(\fa_P).$
\end{lemma}
\proof 
Let $\Cartan_\iC$ denote the complex linear extension of the Cartan involution to 
$\fg_\iC.$ It restricts to a linear automorphism of  $\fh_\iC$ whose 
inverse transpose $\fh_\iC^* \to \fh_\iC^*$ is denoted by $\Cartan_\iC$ as well.
The latter map preserves both $R(\fh)$ and $R({}^*\fh_P).$  Since $R^+(\fh)$ 
is compatible with $\gS^+$ it follows that $-\Cartan_\C$ preserves 
the set $\widetilde \gS^+:= \{\ga \in R^+(\fh) \mid \ga|_\fa \neq 0\}.$  

We define  $\gd_P = \gd - \gd_{\fm_P}$ where $\gd$ and $\gd_{\fm_P}$ are half the sums 
of the positive roots from $R^+(\fh)$ and $R^+({}^*\fh_P)$, respectively. 
Then $\gd_P$ equals half the sum of the positive roots from $R^+(\fh)\setminus R^+({}^*\fh_P).$ 
The latter set equals $\widetilde \gS^+\setminus R({}^*\fh_P)$ hence is invariant under the map $-\Cartan_{\iC}.$ 
It follows that $- \Cartan_\iC \gd_P = \gd_P$, so that $\gd_P \in \fadc.$ Since clearly $\gd_P|_{\fa} = \rho_P$, 
we find that 
$$ 
\gd_P = \rho_P.
$$ 
In particular this implies that $2 \rho_P \in \gL(\fh).$ 
Let $\gb$ be a simple root from $ R^+(\fh).$ 
If it vanishes on $\fa_P$, then clearly, $\Rinp{\rho_P}{\gb} = 0.$ If $\gb$ does not vanish 
on $\fa_P$ then the simple roots $\gg$ from $R^+({}^*\fh_P)$ are simple  for $R^+(\fh)$ 
and not equal to $\gb$, hence satisfy $\Rinp{\gg}{\gb} \leq 0.$ For such a root $\gb$ we thus have 
$\Rinp{\gb}{\gd_{\fm_P}} \leq 0$ so that 
$$ 
\Rinp{\gb}{\gd_P} \geq \Rinp{\gb}{\gd}  = \Rinp{\gb}{\gb} > 0.
$$ 
We thus conclude that $2\rho_P \in \gL^+(\fh) \cap \fadc.$ 
By application of Lemma \ref{l: double weight lattice} it now follows that $4\rho_P \in \gL^+(\fa) \cap \faPdc = \gL^+(\fa_P).$ 

We finish the proof by establishing the inequalities of (\ref{e: gL double plus}). 
Let $\ga \in \gS(\fn_P, \fa).$ 
Then $\ga$ is the restriction to $\fa$ of a root $\hat\ga \in R^+(\fh)$ which does not 
vanish on $\fa_P.$ Now $\hat \ga$ can be written as a sum of simple roots $\gb \in R^+(\fh).$ 
For all these we have $\Rinp{\gd_P}{\gb} \geq 0,$ see above. For those not vanishing on $\fa_P$ 
we have $\Rinp{\gd_P}{\gb} >0.$ Therefore,
$$ 
\Rinp{\rho_P}{\ga} = \Rinp{\gd_P}{\hat \ga} > 0
$$ 
\qed

\section{Projection along infinitesimal characters}
\label{s: proj inf char}
Let  $V$ be an admissible $(\fg, K)$-module and suppose that 
$\fZ$, the center of $U(\fg),$ acts on $V$ in a finite way. By this we mean that
$V$ decomposes into a finite direct sum of generalized weight spaces for $\fZ.$ 
If $\xi$ belongs to the set $\widehat \fZ$ of characters of $\fZ,$ we denote the associated
generalized weight space by $V[\xi].$ Obviously, $V[\xi]$ is an admissible
$(\fg, K)$-submodule of $V.$ Let $X$ be the set of $\xi \in \widehat\fZ$ such that
$V[\xi] \neq 0;$ then $X$ is finite and $V$ is the direct sum of the weight spaces
$V[\xi]$ for $\xi \in X.$ For each $\xi \in \widehat \fZ$ the associated $\fZ$-equivariant
projection map $V \to V$ with image $V[\xi]$ is denoted $p^V_\xi = p_\xi.$ It is readily checked
that $p_\xi$ is $(\fg, K)$-equivariant. If $\fh$ is a Cartan subalgebra of $\fg$ and
$\gl \in \fhdc,$ then we agree to write $V[\gl]:= V[ \xi_\gl]$ and $p_\gl := p_{\xi_\gl};$ 
here $\xi_\gl: Z \mapsto \gg(Z, \gl)$ is the character of $\fZ$ defined via the canonical
isomorphism $\gg: \fZ \to P(\fad).$ 

Similar definitions can be given if $V$ is a complete locally convex space on 
which $G$ has a smooth admissible representation $\pi$ such that $\fZ$ acts finitely. 
This leads again to a finite decomposition into a sum of generalized weight spaces 
$$
V = \oplus_{\xi \in X} V[\xi]
$$
with $\fZ$-equivariant projection maps $p^V_\xi: V \to V$ with image $V[\xi].$ 
For a character 
$\xi \in \widehat \fZ$ the associated generalized weight space is the intersection
of the spaces $\ker (Z -\xi(Z))^p,$ for $Z \in \fZ,$ and $p\geq 1.$ As these spaces are all
$G$-invariant and closed, it follows that $V[\xi]$ is $G$-invariant and closed. This in turn implies that $p_\xi: V \to V$ is a $G$-equivariant continuous projection, for every $\xi \in \widehat \fZ.$ 
In view of admissibility we note that
$$ 
V[\xi] \cap V_K = V_K[\xi], \qquad V[\xi ]= \cl V_K[\xi]).
$$ 
Furthermore, $p_\xi|_{V_K}$ is the projection $p^{V_K}_\xi$ associated with $V_K$ 
and $p_\xi.$

We now assume  $(\rho, E)$ to be a smooth 
representation of $G$ in a complete locally convex space which is admissible and of finite length. 
The following result follows immediately from \cite[Thm.~5.1]{Ktensprod}.  
We assume that $(\pi, F)$ is a finite dimensional irreducible representation of $G$
of highest weight $\mu\in \fhdc$. 

\begin{lemma}
\label{l: action by Z on tensor prod}
Let $(\rho, E)$ be as above and have infinitesimal character $\gl\in \fhdc.$ 
Let $\{\mu_1 = \mu, \mu_2,\dots, \mu_m\} \subset \fhdc$ be 
the set of distinct weights of the finite dimensional representation $(\pi, F).$ 
Then for every $Z \in \fZ$, 
$$ 
\prod_{k =1}^m  \;(Z - \gg(Z, \gl + \mu_k)) \quad \mbox{\rm acts by zero on}\;\;\; E \otimes F.
$$ 
\end{lemma}
\proof 
The reference \cite[Thm.~5.1]{Ktensprod} gives this result for $E_K \times F,$ where
$E_K$ is the $U(\fg)$-module of $K$-finite vectors in $E.$ The required result now follows
by density of $E_K$ in $E$ and continuity of the action of $\fZ$ on $E.$ 
\qed

Let $Q$ a parabolic subgroup of  $G$ containing $A$ and let $\omega$ be a 
continuous representation of $Q$ in a Hilbert space $H_\omega.$ 
If $(\pi, F)$ is a finite dimensional representation of $G$, we have a natural $G$-equivariant topological linear isomorphism
\begin{equation}
\label{e: natural iso gf smooth}
\gf:  C^\infty(G/Q: \omega) \otimes F \;\;{\buildrel \simeq\over \longrightarrow} \;\;C^\infty(G/Q: \omega \otimes \pi|_Q)
\end{equation}
given by the formula 
$$ 
\gf(f \otimes v) (x) = f(x) \otimes \pi(x)^{-1} v, \qquad (x \in G),
$$ 
for $f \in C^\infty(G/Q: \omega)$ and $v \in F.$  The inverse to this isomorphism
is given by $\gf^{-1}(f)(x) = (1 \otimes \pi(x)) f(x),$ for $f \in C^\infty(G/Q: \omega \otimes \pi|_Q)$ 
and $x \in G.$ Clearly, all these assertions also hold with the bigger spaces of continuous functions that arise from replacing $C^\infty$ by $C$ everywhere.

\begin{lemma}
The isomorphism (\ref{e: natural iso gf smooth}) has a unique  extension 
to a continuous linear map 
\begin{equation}
\label{e: def gf mininfty tensprod}
 \gf^{-\infty}: C^{-\infty}(G/Q: \omega) \otimes F \to C^{-\infty}(G/Q: \omega \otimes \pi|_Q).
\end{equation}
This extension is a $G$-equivariant topological linear isomorphism.
\end{lemma}

\proof
Uniqueness is obvious, by density and continuity.
For existence, let ${}^*\gf$ denote the isomorphism (\ref{e: natural iso gf smooth}) for the conjugate representations 
$(\omega^*, H_\omega)$ and $(\pi^*, F)$ 
in place of  $(\omega, H_\omega)$ and $(\pi, F)$ (we assume $F$ to be equipped with a $K_Q$-invariant inner product). Then by taking the transpose of the isomorphism 
$$ 
({}^*\gf)^{-1}: C^{\infty}(G/Q: \omega^* \otimes \pi^*|_Q) \;{\buildrel \simeq \over\longrightarrow}\; C^{\infty}(G/Q:  \omega^*) \otimes F
$$ 
one obtains an extension of (\ref{e: natural iso gf smooth}) to a $G$-equivariant topological  linear isomorphism. 
\qed 

 At a  later stage we will use the notation $\gf_\nu$ for the map $\gf$ of 
(\ref{e: natural iso gf smooth}) in the case that $\omega = \gs \otimes \nu \otimes 1,$
with $\gs$ a unitary representation of $M_Q$ and $\nu \in \faQdc.$ 
\medbreak
Let  $\fh$ be a $\theta$-stable Cartan subalgebra of $\fg$ containing $\fa$ 
and retain the notation of the beginning of Section \ref{s: fin dim reps}.
In particular, $\ft = \fm \cap \fh$ and 
$$
\fh^*_\R := \{\xi \in \fh^*_\iC\mid \xi(i\ft + \fa) \subset \R\} = i\ft^* \oplus \fa^*.
$$ 
We denote by $W_Q(\fh)$ the centralizer of $\fa_Q$ in 
$W(\fh).$ 

Recall the definition of the complex bilinear form $\inp{\dotvar}{\dotvar}$ on $\fadc$ 
in the text following (\ref{e: defi inp}).
We denote by $\gS(\fn_Q, \fa_Q)$ the set of $\fa_Q$-weights in $\fn_Q.$ 

\begin{defi}
\label{d: affine root hyperplane}
By an affine $\gS(\fn_Q, \fa_Q)$-hyperplane in $\faQdc$ we mean a hyperplane
of the form 
\begin{equation}
\label{e: defi H ga c}
\hypp_{\ga, c} = \{\nu \in \faQdc \mid \inp{\nu}{\ga} = c\}, \qquad (\ga \in \gS(\fn_Q, \fa_Q), c \in \C).
\end{equation} 
The hyperplane is said to be real if $\hypp_{\ga , c} \cap \fa_Q^*\not= \emptyset,$ 
which is equivalent to $c \in \R.$ 
\end{defi}

\begin{lemma}
\label{l: generic W h conjugacy}
Let $P$ be a parabolic subgroup of $G$ containing $A$ and let $X \subset \fhdc$ be finite. 
Then there exists a finite collection $\cH= \cH(X)$ of affine $\gS(\fn_P, \fa_P)$-{\bfhyp} 
such that for each $\xi_1, \xi_2 \in X$,  $w \in W(\fh)$, and all $\nu\in \faPdc\setminus \cup \cH$,
$$ 
w(\xi_1+ \nu) = \xi_2  + \nu \implies w \in W_P(\fh).
$$ 
If $X \subset \fh_\R^*$, then $\cH$ may be chosen to consist of real $\gS(\fn_P, \fa_P)$-{\bfhyp}.
\end{lemma}

\proof 
It is easily verified that it suffices to prove the lemma for the case that $\fg$ is semisimple.
Assume this to be the case.

Let $\nu \in \faPdc$, $w \in W_P^c := W(\fh)\setminus W_P(\fh)$, and $\xi_1, \xi_2 \in X$, and assume that 
$w(\xi_1 + \nu) = \xi_2 + \nu.$  Put $X_w:= w(X) - X.$ Then $(I- w)(\nu) \in X_w$, 
hence $\nu \in (I-w)^{-1}(X_w) \cap \faPdc.$ It is sufficient to show that the latter set is contained in 
a finite collection $\cH_w$ of $\gS(\fn_P, \fa_P)$-{\bfhyp}. Then $\cH = \cup_{w \in W_P^c} \,\cH_w$ fulfils the requirements.

By our assumption on $w$, the linear space $\ker (I-w) \cap \faPd$ is properly contained in $\faPd, $
and therefore has a non-zero linear complement $T$ in $\faPd.$ We find that
\begin{equation}
\label{e: weyl conjugates} 
\nu \in  (T_\iC \cap (I- w)^{-1}(X_w)) +  (\ker (I-w) \cap \faPdc).
\end{equation}
We claim that the first of these sets is finite.
For this we note that the map $(I-w) \in \End(\fhdc)$  preserves $\fh^*_\R.$ 
Since $T\subset \faPd\subset \fh^*_\R$ and $\ker (I -w)\cap T = 0,$ it follows that 
$\ker (I-w) \cap T_\C = 0$ as well. This implies that the first set in (\ref{e: weyl conjugates})
has cardinality at  most $\# X_w$ and establishes the claim.

If $\xi \in \ker (I- w) \cap \faPdc$, then $w$ can be written 
as a product of reflections in $\fh$-roots vanishing on $\xi.$ 
At least one of these roots, say $\widehat \ga$, 
does not vanish on $\fa_P$, so that $\ga = \widehat\ga|_{\faPd} \in \gS(\fn_P, \fa_P).$ 
We see that $\xi \in \ker \ga.$ It follows that the set in (\ref{e: weyl conjugates}) is contained in $\cup \cH_w$, 
where $\cH_w$ is the collection of all {\bfhyp} of the form $\eta + \hypp_{\ga, 0}$, with $\eta \in (T_\iC \cap (I-w)^{-1}(X_w))$ 
and $\ga \in \gS(\fn_P, \fa_P).$ 

For the proof of the last assertion assume that $X\subset \fh^*_\R$. Then $X_w\subset \fh_\R^*.$ 
Let 
$$ 
\eta \in T_\iC \cap (I - w)^{-1}(X_w).
$$  
Then it suffices to show that $\eta \in T.$ Write $\eta = \eta_1 + i \eta_2.$ 
Then $(I-w)(\eta_j) \in \fh^*_\iR$ and  by considering real and imaginary parts
we conclude that $(I - w)(\eta_2) = 0.$ Now $T \cap \ker (I -w) = 0$ and we infer $\eta_2 = 0.$ 
 \qed

We now fix a standard parabolic subgroup $P$ and a unitary representation $(\gs, H_\gs)$ of $M_P$, 
which is admissible and of finite length, and such that $H_\gs^\infty$ is quasi-simple with infinitesimal 
character determined by $\gL \in (\fh \cap \fm_P)^*_\iC.$
 
Furthermore, we fix $\mu \in \gL^{++}(\fa_P)$, see
(\ref{e: gL double plus}),
and denote the associated irreducible finite dimensional representation of highest weight $\mu$ by 
$(\pi, F).$ Let $\{\mu_1,\ldots, \mu_m\}$ be the set of $\fh$-weights of $\pi$, 
ordered in such a way that $\mu_1 = \mu.$ 

Our goal is to describe the projection $p_{\gL + \nu + \mu}$ on the space
$C^{-\infty}(G/Q:\gs:\nu ) \otimes F$ of generalized vectors of the representation
$\Ind_Q^G(\gs \otimes \nu \otimes 1) \otimes \pi$, for $Q \subset G$ a parabolic subgroup with 
split component $A_Q = A_P.$ 
The (finite) set of all such parabolic subgroups is denoted by 
$\cP(A_P).$ 

\begin{lemma}
\label{l: second weyl conjugacy}
For every $j > 1$ the elements $\gL + \mu_j $ and $\gL + \mu$ are not conjugate under $W_P(\fh)$, 
the centralizer of $\fa_P$ in $W(\fh).$ 
\end{lemma}

\proof
Let $e_\mu$ be a (non-zero) highest-weight vector of $F.$ Then $M_P$ acts trivially on $e_\mu,$
see Lemma \ref{l: star A on hw},
so that $U(\bar\fn_P) e_\mu  = U(\fg) e_\mu $ is a subspace of $F$ which is invariant under
the action of $G_e$ and under the action of $M_P A_P$, hence under the action of 
$G.$ By irreducibility it follows that $U(\bar\fn_P) e_\mu = F.$ We thus see that each $\mu_j$, 
for $j > 1$, is of the form $\mu_j = \mu - \xi_j$, where $\xi_j \in \gS(\fn_P, \fh)\setminus \{0\}.$ 
The latter implies that $\xi_j$ does not vanish identically on $\fa_P.$ 

If $w \in W_P(\fh)$ and $j >1$, then $w(\gL + \mu) - (\gL + \mu_j) = w(\gL) - \gL + \xi_j.$ 
Now $w(\gL ) - \gL$ vanishes identically on $\faP$, and $\xi_j$ does not, so that 
$w(\gL + \mu) - (\gL + \mu_j) \neq 0.$ 
\qed 

\begin{cor}
\label{c: third Weyl conjugacy}
There exists a finite collection $\cH$ of affine $\gS(\fn_P,\faP)$-{\bfhyp} such that 
for all $\nu \in \faPdc\setminus \cup \cH$ and all $j > 1$ the element $ \gL + \nu +\mu_j$
does not belong to $W(\fh)(\gL + \nu + \mu).$ 

If $\gL \in \fh^*_\R,$  the assertion is 
valid with the additional  requirement that $\cH$ consists of real $\gS(\fn_P, \faP)$-{\bfhyp}.
\end{cor}

\proof 
Put $X = \{\gL+\mu_j\mid j \geq 1\}.$ Then $X \subset \gL + \fh^*_\R \subset \fhdc.$ 
Let $\cH$ be the finite collection of affine $\gS( \fn_P, \fa_P)$-{\bfhyp} 
satisfying the conclusions of Lemma \ref{l: generic W h conjugacy}.

If $\nu \in \faPdc\setminus \cup \cH$ and $\gL + \nu + \mu_j  \in W(\fh)(\gL + \nu + \mu)$
for $j >1,$ then it would follow that $\gL+ \mu_j \in W_P(\fh)(\gL + \mu),$ 
violating the assertion of Lemma  \ref{l: second weyl conjugacy}.
\qed

Let $\gg: \fZ \to P(\fhdc)^{W(\fh)}$ be Harish-Chandra's canonical isomorphism.
Following the notation of \cite{Bps2} we define, for $Z \in \fZ$,  the polynomial map $\Pi(Z) = \Pi_{\mu}(Z): \faPdc \to \fZ$ by 
\begin{equation}
\Pi(Z, \nu):= \prod_{j > 1} [Z   - \gg(Z, \gL + \nu + \mu_j)] 
\end{equation}
\begin{lemma}
\label{l: center zero on ker p}
Let $\cH$ be a finite collection of affine $\gS(\fn_P,\faP)$-{\bfhyp} as in 
Corollary \ref{c: third Weyl conjugacy}.  Let $Q \in \cP(A_P).$ 
If $\nu \in \faPdc\setminus \cup \cH$, then 
\begin{equation}
\label{e: vanishing Pi on ker p}
[\Ind_Q^G(\gs \otimes \nu \otimes 1) \otimes \pi](\Pi(Z, \nu))  = 0  
\quad {\rm on }\;\;\; \ker (p_{\gL + \nu + \mu}).
\end{equation}
\end{lemma}

\proof
Fix $\nu \in \faPdc\setminus \cup \cH$ and put  
$\rho = \Ind_Q^G(\gs \otimes \nu \otimes 1) \otimes \pi.$  Write $\Pi(Z): = \Pi(Z, \nu)$, for $Z \in \fZ.$ 
For each $Z \in \fZ$ the operator $\rho(Z)$ 
commutes with $p_{\gL + \nu + \mu}$ hence leaves the kernel
of the latter invariant. Let $\gd \in \widehat K.$ Since the operator $\rho(Z)$ is $G$-equivariant, 
it restricts to an endomorphism $\rho(Z)_\gd$ of the finite dimensional 
isotypical component $\cK_\gd:= \ker (p_{\gL + \nu + \mu})_\gd.$ For each 
character $\cchidiff$ of $\fZ$ let $\cK_{\gd, \cchidiff}$ denote the associated 
union of the spaces $\ker (\rho(Z)_\gd - \cchidiff(Z))^k$, for $k \geq 1.$ Let $\fX$ be the  set 
of characters $\cchi$ for which $\cK_{\gd, \cchidiff}$ is non-zero. 
Then $\cK_{\gd}$ 
is the finite direct sum of the generalized weight spaces $\cK_{\gd,\cchidiff}$, for $\cchidiff \in \fX.$ 

Fix $\cchidiff \in \fX.$ For each $Z \in \fZ$, the endomorphism $\rho(Z)$ restricts to an endomorphism 
$\rho(Z)_{\gd,\cchidiff}$ of $\cK_{\gd, \cchidiff}.$ It suffices to show that $\rho(\Pi(Z))_{\gd,\cchidiff} = 0$ 
for all $Z \in \fZ.$ Let $\cE$ be a finite dimensional linear subspace of $\fZ$ 
which generates the algebra $\fZ.$ Since $\fZ$ is the union of such subspaces, it suffices to 
show that $\rho(\Pi(Z))_{\gd,\cchidiff} = 0$ for all $Z \in \cE.$ 

For $Z \in \fZ$ we define $l(Z) = Z - \gg(Z, \gL +\nu + \mu).$ Then it follows from 
Lemma \ref{l: action by Z on tensor prod} with $\gL + \nu$ in place of $\gl$  that 
\begin{equation}
\label{e: prod l and D is zero}
\rho(l(Z))_{\gd,\cchidiff} \rho(\Pi(Z))_{\gd, \cchidiff} = 0.
\end{equation}
For $Z \in \fZ$ the endomorphism $\rho(Z)_{\gd, \cchidiff}$ has the single eigenvalue
$\cchidiff(Z)$, so that the endomorphism $\rho(l(Z))_{\gd,\cchidiff} = \rho(Z)_{\gd,\cchidiff} - \gg(Z, \gL + \nu +\mu)$ 
has the single eigenvalue $\cchidiff(Z) - \gg(Z, \gL + \nu + \mu).$ 

By definition of $p_{\gL + \nu + \mu}$ each character from $\fX$ is different from 
$Z \mapsto \gg(Z, \gL + \nu + \mu).$ Since $\cE$ generates the algebra 
$\fZ$ there must be an element $Z \in \cE$ such that $\gg(Z, \gL + \nu + \mu) \neq \cchidiff(Z).$ 
It follows that the subspace $\cE_0$ of $Z \in \cE$ such that $\gg(Z, \gL + \nu + \mu) =  \cchidiff(Z)$ 
is a proper hyperplane in $\cE.$ For $Z \in \cE\setminus \cE_0$ the endomorphism
$\rho(l(Z))_{\gd, \cchidiff}$ has a single non-zero eigenvalue, hence is invertible.
Taking (\ref{e: prod l and D is zero}) into account we infer that $\rho(\Pi(Z))_{\gd,\cchidiff} = 0$ 
for all $Z \in \cE \setminus \cE_0.$ By density this extends to all $Z \in \cE$. 
\qed

From now on we assume that $\gs$ is a representation from 
the discrete series of $M_P.$ In particular, its infinitesimal character $\gL$ belongs
to $\fh_\R^*$ and is regular. 
\begin{lemma}
\label{l: iso image pr}
Let $Q \in \cP(\fa_P)$ and $\mu \in \gL^{++}(\fa_P).$ 
Then there exists a locally finite collection $\cH = \cH(Q, \gs, \mu)$ of  
affine $\gS(\fn_P, \fa_P)$-{\bfhyp} in $\faPdc$ such that, 
for every $\nu \in \faPdc\setminus \cup \cH,$ 
$$ 
p_{\gL + \nu + \mu} [C(G/Q: \gs: \nu)_K \otimes F] \simeq C(G/Q:\gs: \nu + \mu)_K
$$ 
as $(\fg, K)$-modules.  
\end{lemma}

\proof 
We first assume that $Q = P$ and denote by $\gf_\nu$ the isomorphism
(\ref{e: natural iso gf smooth}) for $\omega = \xi_\nu:= \gs \otimes \nu \otimes 1.$  Then 
$\gf_\nu$ restricts to an equivariant isomorphism
$$ 
C(G/P: \gs: \nu)_K \otimes F \;{\buildrel\simeq \over\longrightarrow} \; C(G/P: \xi_\nu \otimes \pi|_P)_K
$$ 
and therefore $p_{\gL + \nu + \mu} \after \gf_\nu = \gf_\nu \after p_{\gL + \nu +\mu}.$ 
Since $\mu$ is $P$-dominant,  the highest weight space $F_\mu$ is a $P$-submodule 
of $F.$ Let $\cH_0$ be the finite collection of affine $\gS(\fn_P, \fa_P)$-{\bfhyp} of Cor. \ref{c: third Weyl conjugacy}. We claim that
for $\nu \in \faPdc\setminus \cup \cH_0$ we have
$$ 
p_{\gL + \nu +\mu}(C(G/P: H_{\gs, \nu} \otimes F)_K) = C(G/P: H_{\gs ,\nu} \otimes F_\mu)_K.
$$ 
As the latter space is isomorphic to $C(G/P : \gs : \nu + \mu)_K$, the result for $Q = P$ will follow from the claim.
To prove the claim, we use exactness of the induction functor $\omega \mapsto \Ind_P^G(\omega)$ 
from the category of admissible $(\fm_{1P}\oplus \fn_P, K_P)$-modules 
to the category of $(\fg, K)$-modules. 
Let $F_1 \subset F_2$ be a sequence
of $P$-submodules of $F$ containing $F_\mu$, 
such that $F_2/F_1$ is irreducible. Then by the mentioned exactness it suffices
to show that $p_{\gL + \nu +\mu} = 0$ on $\Ind_P^G(H_{\gs, \nu} \otimes (F_2/F_1)).$

It follows from the irreducibility that $N_P$ acts trivially on $F_2/F_1$ and that 
the $M_{1P}$-action is irreducible, with a set of $\fh$-weights of the form $\{\mu_j \mid j \in J\}$ 
with $J \subset \{2, \ldots, m\}.$ Note that $\mu_j|_{\fa_P}$ is independent of $j \in J.$ 
 It follows from Lemma \ref{l: action by Z on tensor prod} that the infinitesimal characters of $\fZ(\fm_{1P})$ in $H_\gs \otimes (F_2/F_1)$ 
are all of the form $\gL + \mu_j$, with $j \in J.$
We conclude that the infinitesimal characters appearing in 
$\Ind_P^G(H_{\gs,\nu} \otimes (F_2/F_1))$ are all of the form $\gg(\dotvar, \gL + \mu_j + \nu)$, 
with $j > 1.$ By our choice of $\cH_0$ these characters are different
from $\gg(\dotvar, \gL +\mu + \nu)$ for $\nu \in \faPdc\setminus \cup \cH_0.$ 
Hence $p_{\gL + \nu +\mu} $ vanishes on $\Ind_P^G(H_{\gs, \nu} \otimes (F_2/F_1)).$ 
This establishes the result for $Q = P.$ 

Let now $Q \in \cP(A_P)$ be arbitrary. Then by the rank one product formula  
there exists a locally finite collection
$\cH_Q$ of affine $\gS(\fn_P, \faP)$-{\bfhyp} such that the standard intertwining operator $A(Q,P,\gs,\nu)$ 
is regular and invertible for $\nu \in \faPdc\setminus\cup \cH_Q.$ 

Then 
$$
\Ind_Q^G(\gs \otimes \nu \otimes 1) \otimes \pi \simeq \Ind_P^G(\gs \otimes \nu \otimes 1) \otimes \pi 
$$ 
and 
$$ 
\Ind_Q^G(\gs \otimes (\nu + \mu) \otimes 1)\simeq \Ind_P^G(\gs \otimes (\nu + \mu) \otimes 1)
$$ 
for $\nu \notin \cup \cH_Q \cup (-\mu + \cup \cH_Q).$ The required result now follows
with the hyperplane collection 
$$ 
\cH(Q,\gs, \mu) = \cH_Q \cup (-\mu + \cH_Q) \cup \cH_0;
$$ 
here we have written $- \mu + \cH_Q = \{-\mu + H \mid H\in \cH_Q\}$.
\qed

\begin{rem}
In this paper we shall not need the deep result that the collection 
$\cH_Q$ in the preceding proof 
may be chosen to consist of real affine $\gS(\fn_P,\faP)$-hyperplanes. 
Indeed, the singular sets of the standard intertwining 
operators $\nu \mapsto A(Q,P, \gs, \nu)$ are locally finite unions of 
real affine $\gS(\fn_P, \fa_P)$-hyperplanes. For this it is only required that 
$\gs$ is irreducible unitary with real infinitesimal character, see \cite[Thm.~6.6]{KS2}. 
The mentioned deep result  requires in addition that the zero set of $\nu \mapsto \eta(Q,P,\gs,\nu)$ 
be contained in a locally finite union or real affine $\gS(\fn_P,\faP)$-hyperplanes.
If $\gs$ belongs to the discrete series of $M_P,$ then Harish-Chandra's
explicit determination of the Plancherel measure, see 
\cite[Lemma 35.3]{HCha3}, guarantees this.
\end{rem}

For $Z \in \fZ$ we define the polynomial function
$b(Z) : \faPdc \to \C$ by 
$$ 
b(Z, \nu):= \prod_{j > 1} \, [\gg(Z, \gL + \nu +\mu) - \gg(Z, \gL + \nu + \mu_j)].
$$ 
\begin{cor}
\label{c: vanishing of b}
Let $\cH$ be a finite set of affine $\gS(\fn_P,\fa_P)$-{\bfhyp}
in $\faPdc$ as in Corollary \ref{c: third Weyl conjugacy}. If  $\nu \in \faPdc$ is such that 
the function $b(\dotvar, \nu)$  vanishes identically on $\fZ,$ then $\nu \in \cup \cH.$ 
\end{cor}

\proof
Assume that $\nu$ satisfies the hypothesis, then it follows that $b(Z, \nu)$ vanishes for all $Z \in \fZ.$ 
Since $\fZ$ is finitely generated, we may fix a finite dimensional 
linear subspace $\fZ_0$  of $\fZ$ which generates $\fZ.$ Since $Z \mapsto b(Z,\nu), \fZ_0 \to \C$ 
is a polynomial function on $\fZ_0$ which vanishes identically on $\fZ_0$ it follows that there exists $j > 1$ such that 
the factor
$$ 
Z \mapsto \gg(Z, \gL +\nu + \mu) - \gg(Z, \gL +\nu +\mu_j)
$$ 
vanishes identically on $\fZ_0.$ Since $\gg$ is an algebra homomorphism,
whereas $\fZ_0$ generates $\fZ$, it follows that the above factor 
vanishes identically on $\fZ.$ In turn, this implies that 
$\gL +\nu + \mu$ and $\gL +\nu +\mu_j$ are $W(\fh)$ conjugate. By application of Corollary \ref{c: third Weyl conjugacy} it now follows that $\nu \in \cup \cH.$ 
\qed

\begin{lemma}
\label{l: p and b and D}
There exists a locally finite collection $\cH = \cH_{\gs,\mu}$ of affine $\gS(\fn_P, \fa_P)$-hyper\-planes in $\faPdc$ such that for all $\nu \in \faPdc \setminus \cup \cH$ and $Z \in \fZ$, 
$$ 
b(Z, \nu)\; p_{\gL + \nu + \mu} = [\Ind_Q^G(\gs \otimes \nu \otimes 1) \otimes \pi](\Pi(Z, \nu))
$$ 
on 
$C^{-\infty}(G/Q: \gs :\nu) \otimes F.$
\end{lemma}

\proof 
The two mentioned maps are continuous linear, hence by density and continuity it
suffices to prove the identity on the level of $K$-finite vectors.
Let $\cH_1$ be a finite collection of affine $\gS(\fn_P,\faP)$-{\bfhyp} as
in Lemma \ref{l: center zero on ker p}.
Let $\cH_2$ be a locally finite collection of such {\bfhyp} as in Lemma \ref{l: iso image pr}. 
We will prove the result with $\cH = \cH_1 \cup \cH_2.$ 
Let $\nu \in \faPdc\setminus \cup \cH.$ Then by Lemma \ref{l: center zero on ker p} the required identity is valid on $\ker p_{\gL + \nu + \mu}.$ By Lemma \ref{l: iso image pr} the image of $p_{\gL + \nu +\mu}$ 
is isomorphic to $C(G/Q:\gs: \nu + \mu)_K$ on which $Z \in \fZ$ acts by the 
scalar  $\gg(Z, \gL + \nu + \mu).$ Therefore, the identity is also valid on the image
of $p_{\gL + \nu +\mu}.$ Since $p_{\gL + \nu + \mu}$ is a projection (on the level of $K$-finite vectors), the result follows.
\qed

The induced representation $\Ind_Q^G(\gs\otimes \nu \otimes 1)$
has infinitesimal character $\gL + \nu$, hence it follows from 
Lemma \ref{l: action by Z on tensor prod} that 
$$ 
(Z - \gg(\gL + \nu + \mu)) \Pi(Z, \nu) = 0 \quad{\rm on} \;\; C^{-\infty}(G/Q: \gs :\nu) \otimes F.
$$

If $Q \in \cP(A_P)$ we define the algebra homomorphism 
$$ 
I_{Q,\nu}: \fZ \to \End(C(K/K_P :\gs)_K\otimes F), \;\;Z \mapsto  
[\Ind_Q^G(\gs\otimes \nu \otimes 1) \otimes \pi](Z).
$$%
Note that $I_{Q,\nu}(Z)$ depends polynomially on $\nu$, for fixed $Z \in \fZ.$ 

\begin{lemma}
\label{l: intro Z mu}
There exists a polynomial map $\ulrmZ: \faPdc \to \fZ$ and a polynomial 
function $q \in P(\fa_P^*)$ which is a finite product of factors of the form 
$\inp{\ga}{\dotvar} - c$, with $\ga \in \gS(\fn_P, \fa_P)$ and $c \in \C$,  such that
for every $Q \in \cP(A_P)$,
$$ 
b(Z,\nu) \ulrmZ(\nu)  -  q(\nu) \Pi(Z, \nu) \in \ker(I_{Q,\nu}) 
$$ 
for all $Z \in\fZ$ and  $\nu \in \faPdc.$ 
\end{lemma}

\proof
We first assume that $Q = P$ and follow the ideas of \cite[proof of Prop. 8.3]{Bps2}. 
It follows from Lemma \ref{l: p and b and D} that for all $Z_1, Z_2 \in \fZ$ 
$$ 
b(Z_2, \nu) \Pi(Z_1, \nu) - b(Z_1, \nu) \Pi(Z_2, \nu) \in \ker(I_{P, \nu})
$$% 
for $\nu$ in an open dense subset of $\faPdc.$ By continuity the above identity 
actually holds for all $\nu \in \faPdc.$ 

Let $\cI$ be the ideal in the ring $S(\fa_P)\simeq P(\fa_P^*)$ generated by the polynomials $b(Z)$, for $Z \in \fZ.$ Let $V_\cI$ be the associated common zero set in $\faPdc.$ 
Let $\cH$ be a finite collection of affine $\gS(\fn_P, \fa_P)$-{\bfhyp} in $\faPdc$ 
as in Corollary \ref{c: vanishing of b}. Then it follows from the mentioned corollary that
$$ 
V_\cI \subset \cup \cH.
$$

We select $\tilde q \in P(\faPd)$ a product of linear factors of the form 
$\inp{\ga}{\dotvar} - c,$ with $\ga \in \gS(\fn_P, \faP)$ and $c \in \C,$ 
such that $\tilde q$ vanishes on $\cup \cH.$ 
Then it follows that $\tilde q$ vanishes on $V_\cI.$ By Hilbert's Nulstellen Satz,  
there exists a positive integer $N$ such that $q: = \tilde q^N$ belongs to $\cI.$ By the Noetherian
property, the ideal $\cI$ is already generated by finitely many of its elements,
say $b(Z_1), \ldots , b(Z_l).$ It follows that there exist $a_j \in S(\fa_P)$ such that,
for all $\nu \in \faPdc,$  
$$ 
q(\nu) = \sum_{j=1}^l a_j(\nu) b(Z_j , \nu). 
$$ 
We define
$$ 
\ulrmZ(\nu) := \sum_{j=1}^l a_j(\nu)\, \Pi(Z_j , \nu).
$$ 
Then for all $Z \in \fZ$ and $\nu \in \fadc$ we have that, modulo $\ker (I_{P,\nu})$, 
\begin{eqnarray*}
q(\nu)\,\Pi(Z, \nu) & = & \sum_{j=1}^l a_j(\nu) b(Z_j , \nu)\, \Pi(Z,\nu)\\
& \equiv & \sum_{j=1}^l a_j(\nu) b(Z , \nu) \,\Pi(Z_j,\nu)\\
& =  & b(Z,\nu) \ulrmZ(\nu).
\end{eqnarray*}
This establishes the result for $Q= P.$ The general result follows from an easy application
of the standard intertwining operator $A(\nu) = A(Q,P,\gs,\nu),$ by noting that 
$I_{Q,\nu} = A(\nu) \after I_{P, \nu}$ for generic $\nu,$ combined with a density argument.
\qed

\begin{cor}
\label{c: q p is Z nu}
Let $Q \in \cP(\fa_P).$ Then there exists a locally finite collection $\cH$ 
 of affine $\gS(\fn_P, \fa_P)$-{\bfhyp} in $\faPdc$ such that, for all $\nu \in \faPdc\setminus \cup \cH$,
 $$ 
q(\nu) p_{\gL + \nu + \mu} = [\Ind_Q^G(\gs \otimes \nu \otimes 1) \otimes \pi](\ulrmZ(\nu)) 
$$ 
on $C^{-\infty}(K/K_Q: \gs) \otimes F$. 
\end{cor}

\proof
By a simple density argument we see that it suffices to establish the identity on the 
$K$-finite level.
In that case, let $\cH_{\gs,\mu}$ be the collection of affine $\gS(\fn_P, \faP)$-{\bfhyp} of Lemma \ref{l: p and b and D}.
Then it follows from combining that lemma with Lemma \ref{l: intro Z mu} that
$$
q(\nu) b(Z,\nu) p_{\gL + \nu +\mu} = b(Z, \nu) I_{Q,\nu}( \ulrmZ(\nu))
$$ 
on $C(K/K_Q:\gs ) \otimes F,$ for all $Z \in \fZ$ and all $\nu \in \faPdc\setminus \cup \cH_{\gs, \nu}.$ 
Let $\cH_2$ be a finite collection of affine $\gS(\fn_P, \faP)$-{\bfhyp} as in Corollary \ref{c: vanishing of b}.
Then for $\nu \in \faPdc\setminus \cup \cH_2$ there exists $Z \in \fZ$ such that
$b(Z, \nu) \neq 0.$ The required result now follows with $\cH = \cH_{\gs, \mu} \cup \cH_2.$
\qed

\begin{rem}
\label{r: uniqueness ulrmZ}
 It follows from Corollary \ref{c: q p is Z nu} that if $\ulrmZ': \faPdc \to \fZ$ is a second
polynomial map as in Lemma \ref{l: intro Z mu}, then for all $Q \in \cP(\fa_P)$ and all $\nu \in \faPdc$ 
we have $\ulrmZ(\nu) - \ulrmZ'(\nu) \in \ker({ I}_{Q,\nu}).$
\end{rem}

\section{The functional equation}
\label{s: functional equation}
 We retain the assumption of the previous section that  $P$ is a standard parabolic
 subgroup of $G,$ that $\gs$ is a representation of the discrete series of $M_P,$ 
 that $\mu \in \gL^{++}(\fa_P),$ and that 
 $(\pi, F)$ is the irreducible finite dimensional spherical representation of  
 highest weight $\mu.$ 

Let 
$e_K \in F$ be a nonzero $K$-fixed vector, and
$e^{-\mu} \in F^*$ a non-zero lowest weight  vector. 
We define the matrix coefficient map $i_\mu: F \to C^\infty(G)$ 
by
\begin{equation}
\label{e: intro e min mu}
i_\mu(v)(x) = \inp{v}{\pi(x) e^{-\mu}} ,\qquad (x \in G).
\end{equation}  
Then $i_\mu$ intertwines $\pi$ with the left regular representation 
and is readily seen to define a $G$-equivariant embedding
\begin{equation}
\label{e: embedding i mu}
i_\mu: F \embeds C^\infty(G/\bar P : 1 : \mu + \rho_P).
\end{equation} 
In particular, we note that $i_\mu(e_K)(kan) = a^{-\mu}\inp{e_K}{e^{-\mu}}$, 
for $(k , a, n) \in K \times A \times \bar N_0$, and see that $i_\mu(e_K)$ 
is a nowhere vanishing function on $G.$  By renormalizing $e_K$ we may arrange 
that $\inp{e_K}{e^{-\mu}} = i_\mu(e_K)(1) = 1.$ 
Accordingly, we define the map 
\begin{equation}
\label{e: defi M mu} 
\cM_\mu : C^{-\infty}(G/\bar P: \gs: \nu + \mu) \to C^{-\infty}(G/\bar P : \gs: \nu) \otimes F\end{equation}
by $ f \mapsto i_\mu(e_K)^{-1} f \otimes e_K.$ 
In the compact picture of the induced representations,
$\cM_\mu$ corresponds to the map 
$$ 
C^{-\infty}(K/K_P:\gs) \to 
C^{-\infty}(K/K_P:\gs)\otimes F, \;\; f \mapsto f \otimes e_K.
$$%
We fix a non-zero highest weight vector $e_\mu \in F.$ In order to emphasize
the feature that it is $N_0$-fixed, we agree to also write 
$e_{N_0} = e_{\mu}.$
 Likewise, we
denote by $e^{N_0}$ a fixed choice of non-zero highest weight vector in 
the complex linear dual space
$F^*.$ 

 We now define
the map 
\begin{equation}
\label{e: defi ge N zero}
\ge^{N_0} := m \after (I \otimes e^{N_0}): C^{-\infty}(G/\bar P:\gs: \nu) \otimes F \to 
C^{-\infty}(G/\bar P:\gs: \nu),
\end{equation}  
where $m$ denotes the natural linear isomorphism  
from 
$C^{-\infty}(G/\bar P:\gs: \nu) \otimes \C$ onto $C^{-\infty}(G/\bar P:\gs: \nu)$,
induced by multiplication. In the compact picture, this becomes a map 
$
\ge^{N_0}: C^{-\infty}(K/K_P:\gs) \otimes F \to  C^{-\infty}(K/K_P: \gs),
$ 
constant in the variable $\nu \in \faPdc.$ 

For $\nu \in \faPdc$ we define the endomorphism $\ulZ_{\bar P, \mu}(\nu)$ of $C^{-\infty}(G/\bar P : \gs: \nu )\otimes F$ by 
\begin{equation}
\label{e: defi ulZ bP mu nu} 
\ulZ_{\bar P,\mu}(\nu): = [\Ind_{\bar P}^G(\gs\otimes \nu \otimes 1) \otimes \pi](\ulrmZ(\nu)),
\end{equation} 
with $\ulrmZ: \faPdc \to \fZ$ a polynomial map as in Lemma \ref{l: intro Z mu}. Note that
the endomorphism (\ref{e: defi ulZ bP mu nu}) is independent of the particular choice
of $\ulrmZ$, in view of Remark \ref{r: uniqueness ulrmZ}.

Finally, we define the operator 
$$ 
D_\mu(\gs,\nu): C^{-\infty}(G/\bar P:\gs: \nu + \mu) \to C^{-\infty}(G/\bar P:\gs:\nu)
$$ 
by 
\begin{equation}
\label{e: defi D mu}
D_\mu(\gs, \nu)= \ge^{N_0} \after \ulZ_{\bar P, \mu}(\nu)  \after \cM_\mu,
\end{equation}
 
\begin{prop}
\label{p: estimate of generalized operator norm D mu}
The operator $D_\mu(\gs,\nu)$, viewed as an endomorphism of  the  space 
$C^{-\infty}(K/K_P:\gs)$ is continuous and depends polynomially
on $\nu \in \faPdc.$ There exists a constant $d \in \N$ such that the following is valid.

There exists an $r \in \N$ and for every $s \in \N$ a constant $C > 0$ such that for all $\nu \in \faPdc$ the endomorphism $D_\mu(\gs,\nu)$
maps the Banach space $C^{-s}(K/K_P:\gs)$ continuously to the Banach space $C^{-s-r}(K/K_P:\gs)$ 
with operator norm satisfying the estimate
$$ 
\| D_\mu(\gs, \nu) \|_{\rm op} \leq C(1 + |\nu |)^d.
$$ 
\end{prop}

For the proof we need the following lemma. 

\begin{lemma}
\label{l: estimate infinitesimal action}
Let $u \in U(\fg)$ be an element of degree at most $d.$ Then  for every $\nu \in \faPdc$ the endomorphism
$\pi_{\bar P, \gs, \nu}^{\infty}(u)$ of $C^{\infty}(K/K_P :\gsP)$ is continuous and support preserving. Furthermore, the following assertions are valid.
\begin{enumerate}
\itema The function $\nu \mapsto \pi_{\bar P, \gs, \nu}^{\infty}(u)$ is polynomial 
$\End(C^\infty(K/K_P:\gsP))$-valued of degree at most $d.$ 
\itemb There exists a constant $t \in \N$ and for every  $s \in \N$ a constant $C > 0$ such 
that for all $f \in C^\infty(K/K_P: \gsP),$ 
$$ 
\| \pi_{\bar P, \gs, \nu}^{\infty}(u) f\|_{s} \leq C (1 + |\nu|)^d \|f\|_{s +t}.
$$ 
\end{enumerate}
\end{lemma}

\proof 
See \cite[Lemma 2.1]{BScfunc}.
\qed

\begin{cor}
\label{c: estimate dual infinitesimal action}
Let $u \in U(\fg)$ be an element of degree at most $d.$ Then  for every $\nu \in \faPdc$ the endomorphism
$\pi_{\bar P, \gs, \nu}^{- \infty}(u)$ of $C^{-\infty}(K/K_P :\gsP)$ is continuous and support preserving. Furthermore,
\begin{enumerate}
\itema The map $\nu \mapsto \pi_{\bar P, \gs, \nu}^{- \infty}(u)$ is polynomial 
$\End(C^{- \infty}(K/K_P:\gsP))$-valued of degree at most $d.$ 
\itemb There exists a constant $t \in \N$ and for every  $s \in \N$ a constant $C > 0$ such that 
for all $\nu \in \faPdc$ the endomorphism $\pi_{\bar P, \gs, \nu}^{- \infty}(u)$ maps $C^{-s}(K/K_P:\gsP)$ 
continuous linearly into $C^{-s-t}(K/K_P:\gsP)$  with operator norm
$$ 
\| \pi_{\bar P, \gs, \nu}^{- \infty}(u)\|_{\rm op} \leq C (1 + |\nu|)^d.
$$ 
\end{enumerate}
\end{cor}

\proof 
This follows from Lemma \ref{l: estimate infinitesimal action} by taking adjoints.
\qed
\medno
{\em Proof of Proposition \ref{p: estimate of generalized operator norm D mu}. }
We start by observing that $\ulrmZ: \faPdc \to \fZ$ is polynomial in the variable 
$\nu.$ For $Z \in \fZ$ we define 
$$ 
D(Z)(\nu): = \ge^{N _0}\after [\pi^{-\infty}_{\bar P,\gs, \nu} \otimes \pi_\mu](Z) \after \cM_\mu.
$$ 
Then it clearly suffices to prove the assertions of the proposition for $D(Z)$ in place of $D_\mu(\nu).$ 
\qed

In terms of the above maps we can now present the functional 
equation for the Whittaker vectors. We will stay close to the 
notation of \cite[Thm. 9.3]{Bps2} in order to emphasize the strong 
analogy. Recall the definition of $D_\mu(\gs, \nu)$ in (\ref{e: defi D mu}).

\begin{thm}{\rm (Functional equation)\;}
\label{t: the functional equation}
Let $\mu \in \gL^{++}(\fa_P).$ Then there exists
a rational $\End(H_{\gs,\cchi_P}^{-\infty})$-valued function 
$\nu \mapsto R_\mu(\gs,\nu)$ on $\faPdc$ such that 
\begin{equation}
\label{e: functional equation}
j(\bar P,\gs , \nu ) = D_\mu(\gs, \nu) \after j(\bar P,\gs, \nu + \mu) \after R_\mu(\gs, \nu).
\end{equation}
\end{thm}

\begin{rem} In the next section we will show that $\nu \mapsto p(\nu) R_\mu(\gs, \nu)$ 
is polynomial for a suitable polynomial function $p: \faPdc \to \C$ 
which can be written as a product of linear factors of the form $\nu \mapsto \inp{\nu}{\ga} + c,$ 
$\ga \in \gS(\fn_P, \faP)$ and $c \in \C.$ 
\end{rem}

We will prove Theorem \ref{t: the functional equation} in a sequence of lemmas occupying the rest of this section
and the next. A key ingredient in our proof is the map 
$$ 
\Phi_\mu(\nu): C^{-\infty}(G/\bar P: \gs: \nu)\otimes F \to C^{-\infty}(G/\bar P: \gs: \nu+\mu), \;\;f\otimes v \mapsto i_\mu(v) f,
$$ 
with $i_\mu$ as in (\ref{e: embedding i mu}),
which is readily verified to be $G$-equivariant.

\begin{defi}
We will say that an assertion depending on a parameter $\nu \in \faPdc$ holds for 
$\gS(\fn_P,\fa_P)$-generic $\nu$, if there exists a locally finite collection
$\cH$ of $\gS(\fn_P, \faP)$-{\bfhyp} in $\faPdc$ such that the assertion is valid
for all $\nu \in \faPdc\setminus\cup \cH$. 
\end{defi}

\begin{lemma}
\label{l: Phi after p is Phi}
For $\gS(\fn_P, \faP)$-generic $\nu \in \faPdc$, 
\begin{equation}
\label{e: Phi mu after p}
\Phi_\mu(\nu) \after p_{\gL + \nu +\mu} = \Phi_\mu(\nu) .
\end{equation}
\end{lemma}
\proof
By equivariance of $\Phi_\mu$, we have 
$$
\Phi_\mu(\nu) \after p_{\gL + \nu +\mu} =  p_{\gL + \nu +\mu} \after \Phi_\mu(\nu).
$$ 
The map on the right of this equation equals $\Phi_\mu(\nu)$, 
since $\Ind_{\bar P}^G(\gs \otimes (\nu + \mu) \otimes 1)$ has infinitesimal character $\gL +\nu +\mu.$ 
\qed
The following identity, for $\nu \in \fadc$,  is a straightforward consequence of the definitions,
\begin{equation}
\label{e: identity Phi and M}
\Phi_\mu(\nu)\after \cM_\mu = I \qquad {\rm on} \quad C^{-\infty}(G/\bar P: \gs: \nu).
\end{equation}
The map $\cM_\mu$ is not equivariant. However, for the map 
\begin{equation}
\label{e: defi Psi mu}
\Psi_\mu(\nu):= q(\nu)^{-1} \ulZ_{\bar P, \mu}(\nu) \after \cM_\mu
\end{equation}
we have the following result.
\begin{lemma}
\label{l: inversion Phi mu nu}
The map $\nu \mapsto q(\nu) \Psi_\mu(\nu)$ is polynomial as a map 
with values in the space of equivariant continuous linear operators from  
$(C^{-\infty}(K/K_P:\gs), \pi_{\bar P,\gs,\nu})$  to 
%the representation 
$(C^{-\infty}(K/K_P:\gs)\otimes F, \pi_{\bar P,\gs,\nu} \otimes \pi)$. 
Furthermore, for all $\nu \in \faPdc\setminus q^{-1}(0)$, 
\begin{eqnarray}
\Phi_\mu(\nu) \after \Psi_\mu(\nu) &=&  I,\label{e: first identity Phi Psi}\\
 \Psi_\mu(\nu) \after\Phi_\mu(\nu) & =&  q(\nu)^{-1} \ulZ_{\bar P, \mu}(\nu). \label{e: second identity Phi Psi}
\end{eqnarray}
\end{lemma} 

\proof 
We first observe that from the definitions it follows that $\nu \mapsto q(\gl) \Psi_\mu(\nu)$ is polynomial as a map into the space of continuous linear operators $C^{-\infty}(K/K_P:\gs) \to C^{-\infty}(K/K_P:\gs)\otimes F.$ 
The equivariance of the operators in the image of that polynomial map will be addressed in a moment. 

In view of Corollary \ref{c: q p is Z nu} we have 
\begin{equation}
\label{e: Z equals p}
q(\nu)^{-1} \,\ulZ_{\bar P,\mu}(\nu) = p_{\gL + \nu + \mu} \quad {\rm on} \quad C^{-\infty}(G/\bar P :\gs :\nu) \otimes F,
\end{equation}
for $\gS(\fn_P, \faP)$-generic $\nu \in \fadc.$ From (\ref{e: defi Psi mu}) it now follows that
$$ 
\Phi_\mu(\nu) \after \Psi_\mu(\nu) = \Phi_\mu(\nu) \after \pr_{\gL + \nu + \mu} \after \cM_\mu.
$$ 
Taking (\ref{e: Phi mu after p}) into account we infer the validity of (\ref{e: first identity Phi Psi}), for $\gS(\fn_P, \faP)$-generic $\nu \in \faPdc.$ 
By analytic 
continuation the identity (\ref{e: first identity Phi Psi}) follows for all $\nu \in \faPdc\setminus q^{-1}(0).$

It follows from (\ref{e: identity Phi and M}) that $\Phi_\mu(\nu)$ is surjective for all $\nu.$ 
Furthermore, for $\gS(\fn_P, \faP)$-generic $\nu,$ the image  
$
{\rm im}(p_{\gL +  \nu+ \mu}) \subset C^{-\infty}(G/\bar P:\gs:\nu) \otimes F
$  is a closed $G$-invariant subspace, which  satisfies 
\begin{equation}
\label{e: im p on K finite}
{\rm im}(p_{\gL +  \nu+ \mu})_K \simeq C^{-\infty}(G/\bar P: \gs: \nu + \mu)_K
\end{equation}
as $(\fg, K)$-modules, in view of Lemma \ref{l: iso image pr}. 
From (\ref{e: Phi mu after p}) we see, 
still for $\gS(\fn_P, \faP)$-generic $\nu, $ that 
$\Phi_\mu(\nu)$ is a $G$-equivariant surjective continuous linear map from the space ${\rm im}(p_{\gL + \gs + \nu})$ onto $C^{-\infty}(G/\bar P: \gs: \nu + \mu),$ which by (\ref{e: im p on K finite}) is injective,
hence bijective.

It follows from (\ref{e: Z equals p}) and (\ref{e: defi Psi mu}) that 
$\Psi_\mu(\nu)$ maps into ${\rm im}(p_{\gL + \nu + \mu}),$ for generic $\nu,$ 
so the equivariance of $\Psi_\mu(\nu)$ follows from (\ref{e: first identity Phi Psi}) and 
the equivariance of $\Phi_\mu(\nu),$ for generic $\nu \in \faPdc.$ By analytic continuation 
the equivariance of $q(\nu) \Psi_\mu(\nu)$ follows for all $\nu \in \faPdc.$ 

We finally turn to proving the identity (\ref{e: second identity Phi Psi}). By analytic continuation, it suffices to establish that identity for $\gS(\fn_P, \faP)$-generic $\nu \in \faPdc.$ Since the maps on both sides of (\ref{e: second identity Phi Psi}) map into the space ${\rm im} (p_{\gL + \nu +\mu})$,  on which 
$\Phi_\mu(\nu)$ restricts to an injective map, for generic $\nu,$ it suffices to check that 
$$ 
\Phi_\mu(\nu)\after \Psi_\mu(\nu) \after \Phi_\mu(\nu)  =  \Phi_\mu(\nu) \after q(\nu)^{-1} \ulZ_{\bar P, \mu}(\nu).
$$%
The expression on the left simplifies to $\Phi_\mu(\nu),$ in view of (\ref{e: first identity Phi Psi}). The expression on 
the right equals $\Phi_\mu(\nu) \after p_{\gL + \nu + \mu}$ by (\ref{e: Z equals p}), and
$\Phi_\mu(\nu)$ by  (\ref{e: Phi mu after p}).
\qed

Recall from (\ref{e: intro e min mu}) that $e^{-\mu}\in F^*$ is a
non-zero lowest weight vector (of $\fa$-weight $-\mu$)
and put 
$$
m_\mu := i_\mu(e_{N_0})(1) = \inp{e_{N_0}}{e^{-\mu}}.
$$
 From $m_\mu = 0$ it would follow
that $i_\mu(e_{N_0})$ vanishes on $N_0 \bar P$ hence on $G,$ contradicting the injectivity of 
$i_\mu.$ Therefore, $m_\mu$ is a nonzero complex number. 

\begin{lemma}
For every $\eta\in H_\gs^{-\infty,\cchi}$ and all $\nu\in \fadc$ with 
$\inp {\Re \nu}{\ga} > 0$, $(\ga \in \gS(\fn_P, \fa))$,
\begin{equation}
\label{e: Phi on tensor j}
\Phi_\mu(\nu)[ j(\bar P , \gs, \nu, \eta) \otimes e_{N_0}] = j(\bar P, \gs,\nu + \mu) (m_\mu  \eta).
\end{equation}
\end{lemma}

\proof
It is readily verified that the expression on the left belongs
to $C^{-\infty}(G/\bar P :\gs: \nu + \mu)^{\cchi}$ hence is of the form 
$j(\bar P, \gs,\nu + \mu)(\eta')$ for some $\eta' \in H_\gs^{-\infty,\cchi}.$ 
On the other hand, on $N_P \bar P$ the expression on the 
left hand side is the continuous $H_\gs^{-\infty}$-valued function 
whose value at the unit element $e=1$ is 
$$ 
\ev_e j(\bar P , \gs, \nu, \eta)\, i_\mu(e_{N_0})(1) = m_\mu \eta.
$$
It follows that $\eta' = m_\mu \eta.$ 
\qed

{\em Completion of the proof of Theorem \ref{t: the functional equation}. } 
Applying the operator $q(\nu) \ge^{N_0} \after \Psi_\mu(\nu)$  to  
(\ref{e: Phi on tensor j}) and taking into account Lemma \ref{l: inversion Phi mu nu}
we obtain
$$ 
 \ge^{N_0}\after \ulZ_{\bar P, \mu}(\nu) [ j(\bar P , \gs, \nu) \eta \otimes e_{N_0}] = 
q(\nu) \ge^{N_0} \after \Psi_\mu(\nu) j(\bar P, \gs,\nu + \mu) (m_\mu  \eta).
$$ 
From (\ref{e: defi D mu}) and (\ref{e: defi Psi mu}) we see that
$
q(\nu) \ge^{N_0} \after \Psi_\mu(\nu) = D_\mu(\gs, \nu).
$
Hence, 
$$ 
\ge^{N_0} \after \ulZ_{\bar P, \mu}(\nu) [ j(\bar P , \gs, \nu)\eta \otimes e_{N_0}] = 
D_{ \mu}(\gs, \nu) j(\bar P, \gs,\nu + \mu)( m_\mu  \eta),
$$ 
for all $\eta \in H_{\gs,\cchi_P}^{-\infty}$ and all $\nu \in \faPdc$ with $\Re \nu$ $P$-dominant. 
The functional equation now follows with $R_{\mu}(\gs, \nu) = M_{\mu}(\gs, \nu)^{-1} m_\mu$ 
by application of Proposition \ref{p: intro M mu nu} below.
\qed
\begin{prop}
\label{p: intro M mu nu}
There exists a unique polynomial function 
$M_\mu = M_\mu(\gs, \dotvar): \faPdc \to \End(H_{\gs,\cchi_P}^{-\infty})$ 
such that for all $\nu \in \faPdc(P,0)$ we have
$$ 
\ge^{N_0} \after \ulZ_{\bar P, \mu}(\nu) [ j(\bar P, \gs, \nu) \eta \otimes e_{N_0}]
=  j(\bar P, \gs, \nu)( M_\mu (\nu) \eta).
$$ 
The polynomial function $\det \after M_\mu: \faPdc \to \C$ is not identically zero.
\end{prop}

To prepare for the proof, we introduce the space 
 $C^\infty_{c,N_P}(G/\bar P: \gs: \nu)$ of functions $f \in C^\infty(G/\bar P: \gs: \nu)$ 
whose support $\supp f$ has compact intersection with $N_P.$ Restriction to $N_P$ induces a linear
isomorphism $r_\nu$ from $C^\infty_{c,N_P}(G/\bar P: \gs: \nu)$ onto $C^\infty_c(N_P, H_\gs^\infty).$ 
We note that $C^\infty_{c, N_P}(G/\bar P: \gs: \nu)$ is invariant under the left regular action 
by $U(\fg)$ and denote by ${}^n \pi_{\bar P, \gs, \nu}$ the unique representation of $U(\fg)$ in $C^\infty_c(N_P, H_\gs^\infty)$ such that  $r_\nu$ intertwines the left regular representation of $U(\fg)$ with ${}^n \pi_{\bar P, \gs, \nu}.$ 

Our first step in the proof of Proposition \ref{p: intro M mu nu} is the following observation.
\begin{lemma}
\label{l: pol dep for Ccinfty}
Let $\gf \in C_c^\infty(N_P, H_\gs^\infty).$ Then for  each element $u \in U(\fg)$ the function $\faPdc \to C^\infty_c(N_P, H_\gs^\infty),$ 
$\nu \mapsto  {}^n \pi_{\bar P, \gs, \nu}(u) \gf$ is polynomial, i.e., it belongs to the space
$P_k(\faPdc) \otimes C^\infty_c(N_P, H_\gs^\infty)$ for some $k \in \N.$ 
\end{lemma}
\proof
It suffices to prove this assertion for $u = X  \in \fg,$ with $k=1.$ 
Let $\Omega $ be a bounded open neighborhood of $\supp\, \gf$ in $N_P.$ 
Let $X \in {\rm Lie}(\bar P).$ Then there exists an open interval $I \ni 0$ in $\R$ 
such that for every $t \in I$ and $n \in \Omega$ we have $ \exp (-tX)n  \in N_P \bar P.$ 
Consequently, there exist smooth functions $U: I\times \Omega \to N_P$ and 
$V: I\times \Omega \to \bar P$ 
such that
$$ 
\exp(-tX) n =  n U(t,n) V(t,n), \qquad (n \in \Omega,\; t \in I).
$$ 
We note that $U(0,n) = e = V(0,n)$ for all $n \in N_P.$ 
Let $\gf_\nu \in  C^\infty(G/\bar P, \gs, \nu)$ be defined by $\supp \gf_\nu \subset \Omega \bar P$ 
and 
$$ 
\gf_\nu(n m a\bar n) := a^{-\nu + \rho_P} \gs(m)^{-1} \gf(n),\qquad ((n, m, a, \bar n ) \in N_P \times M_P \times A_P \times\bar N_P).
$$
Then $\gf_\nu|_{N_P} = \gf,$ so that ${}^n \pi_{\bar P, \gs, \nu}(X) \gf(n) = L_X \gf_\nu(n).$ 
For $n \in N_P$ and $t \in I$ we have
\begin{eqnarray*} 
\gf_\nu(\exp (- tX) n) & = & \gf_\nu( n U(t,n) V(t,n)) \\
&=&  [\gs \otimes (-\nu + \rho_P) \otimes 1](V(t,n))^{-1} \gf(n U(t,n )).
\end{eqnarray*} 
Differentiating this expression in $t$ at $t =0$ we find, for $n \in \Omega,$
\begin{eqnarray*}
{}^n \pi_{\bar P, \gs, \nu}(X) \gf(n) &= &
L_X(\gf_\nu)(n) \\
& = & - [\gs \otimes (-\nu + \rho_P) \otimes 1](\partial_t V(0,n))\gf(n)
+ R_{\partial_t U(0,n)}\gf (n).
\end{eqnarray*} 
We thus see that $\nu \mapsto {}^n \pi_{\bar P, \gs, \nu}(X) \gf$ belongs to $P_1(\faPdc) \otimes C^\infty_c(\Omega, H_\gs^\infty).$ 
\qed

The next step in the proof of Proposition \ref{p: intro M mu nu} is formulated 
in the following lemma.

\begin{lemma}
\label{l: intro m Z nu}
For $Z \in \fZ$ and $\nu \in \faPdc(P,0)$ there exists 
a unique endomorphism $m(Z ,\nu) \in \End(H_{\gs,\cchi_P}^{-\infty})$ such that,
for all $\eta \in H_{\gs,\cchi_P}^{-\infty },$ 
\begin{equation}
\label{e: m Z mu}
\ge^{N_0} [\pi_{\bar P, \gs, \nu}  \otimes \pi_\mu](Z) [j(\bar P,\gs,\nu, \eta) \otimes e_{N_0}]  = 
j(\bar P,\gs,\nu ,m(Z, \nu)(\eta)).
\end{equation}
The map $\nu \mapsto m(Z, \nu)$ is polynomial on $\faPdc$ with values in $\End(H_{\gs,\cchi_P}^{-\infty}).$
\end{lemma}

\proof
For $\nu$ as stated, the expression on the left-hand side belongs
to the space $C^{-\infty}(G/\bar P: \gs:\nu)_\cchi$, hence by Proposition 
\ref{p: bijectivity j in dominant region} (b) can be written as
the expression on the right-hand side, with a uniquely determined 
$m(Z, \nu)(\eta) \in H_{\gs,\cchi_P}^{-\infty}.$ By uniqueness, $m(Z, \nu)(\eta)$ 
depends linearly on $\eta$, hence, $m(Z, \nu) \in \End(H_{\gs,\cchi_P}^{-\infty}).$ 

It remains to be shown that $m(Z,\dotvar )$ is polynomial. Let $\eta \in H_{\gs,\cchi_P}^{-\infty}$ and let 
$v \in H_\gs^\infty$ be an arbitrary element. Then it suffices to show that
$\nu \mapsto \inp{m(Z,\nu)\eta}{v}_\gs$ is polynomial in the indicated range.
See (\ref{e: reversed order pairing}) for the definition of the pairing used.
For this we recall that for $\nu$ in that range,  $j(\bar P, \gs, \nu)(\eta)$ 
restricted to $N_P$ is the continuous function $N_P\to H_{\gs, \cchi_P}^{-\infty}$ 
given by $n \mapsto \cchi(n)^{-1}  \eta.$ Fix a function $\psi \in C^{\infty}_c(N_0)$ such that
$$ 
\int_{N_P} {\cchi(n)} \psi(n) \; dn = 1.
$$ 
We define $f  \in C_c^\infty(N_P, H_\gs^\infty) $ 
by 
$f(n) = \psi(n) v$ for $n \in N_P.$ Furthermore, we denote by $f_{-\bar \nu}$ the extension 
of $f$ to an element of $C^\infty(G/\bar P: \gs: -\bar \nu)$ with support contained in $N_P \bar P.$ 
Then, for $\eta' \in H_{\gs,\cchi_P}^{-\infty},$ 
$$ 
\inp{j(\bar P, \gs ,\nu)(\eta' )}{f_{-\bar \nu}} = 
\int_{N_P}  \cchi(n)^{-1} \inp{\eta'}{\psi(n) v}_\gs \; dn = \inp{\eta'}{v}_\gs.
$$ 
Substituting $\eta' = m(Z, \nu) \eta$ and combining the result with (\ref{e: m Z mu})
we now find that
$$ 
\inp{v}{m(Z, \nu)\eta}_\gs= 
\inp{\ge^{N_0} [\pi_{\bar P, \gs, \nu}  \otimes \pi_\mu](Z)
 [j(\bar P, \gs,\nu, \eta) \otimes e_{N_0}]}{f_{-\bar \nu}}.
 $$%
By the Leibniz rule for tensors, the expression on the right-hand side is a sum of terms of the form
$$ 
\inp{\ge^{N_0} \pi_{\bar P, \gs, \nu}(U)  j(\bar P, \gs, \nu)( \eta) \otimes \pi_\mu (V)e_{N_0} } {f_{-\bar \nu}},
$$%
with $(U,V)$ ranging over a subset of $U(\fg)\times U(\fg)$ independent of $\nu.$ The above term 
may be rewritten as
$$ 
\ge^{N_0}(\pi_\mu(V) e_{N_0}) \cdot \inp{\pi_{\bar P, \gs, \nu}(U) j(\bar P, \gs, \nu)( \eta)}{f_{-\bar \nu}}.
$$%
Thus, it suffices to show that the latter expression is polynomial in $\nu.$ 
We now observe that 
\begin{eqnarray*}
\inp{\pi_{\bar P, \gs, \nu}(U) j(\bar P,\gs,\nu)( \eta)}{f_{-\bar\nu}}&=& 
\inp{ j(\bar P, \gs, \nu)( \eta)}{\pi_{\bar P, \gs, -\bar \nu}(U^*)f_{-\bar\nu} }\\
&=& 
\inp{ j(\bar P,\gs, \nu)( \eta)} {({}^n \pi_{\bar P, \gs, -\bar \nu}(U^*)f)_{-\bar\nu} }\\
&=& 
\int_{N_P} \cchi(n)  \inp{\eta}{{}^n \pi_{\bar P, \gs, -\bar \nu}(U^*)f (n)}_\gs \; dn.
\end{eqnarray*}
By virtue of Lemma \ref{l: pol dep for Ccinfty} 
 the latter integral depends polynomially on $\nu.$ 
\qed

{\em Proof of Proposition \ref{p: intro M mu nu}.}
We note that the map $\faPdc \to \fZ,$ $\nu \mapsto \ulZ_\mu(\nu)$ is polynomial. Moreover, 
$$ 
\ulZ_{\bar P, \mu}(\nu) = [\pi_{\bar P, \gs, \nu}  \otimes \pi_\mu](\ulZ_\mu(\nu)).
$$%
By application of Lemma \ref{l: intro m Z nu} the first assertion now follows with 
$$ 
M_\mu(\nu) = m(\ulZ_{\mu}(\nu), \nu).
$$% 
For completing the proof of Proposition \ref{p: intro M mu nu} it thus remains to establish the lemma below.
\qed

\begin{lemma}
\label{l: det M mu not zero}
The polynomial  function $\faPdc \to \C,$ $ \nu \mapsto \det M_\mu(\nu)$ is not identically zero.
\end{lemma}

\proof
It suffices to show that $\nu \mapsto \det M_\mu(\nu)$ is non-zero for a suitable $\nu.$ 
For this it suffices to show that there exists a $\nu \in \faPdc(P,0)$ 
such that $M_\mu(\nu)$ is an injective endomorphism of $H_{\gs, \cchi_P}^{-\infty}.$  
Taking the characterization of $M_\mu$ in Proposition \ref{p: intro M mu nu} into account 
and using that 
\begin{equation}
\label{e: Z q p}
\ulZ_{\bar P, \mu}(\nu) = q(\nu) p_{\gL+ \nu + \mu} \quad {\rm on} \quad C^{-\infty}(G/\bar P: \gs: \nu) \otimes F
\end{equation}
we infer that it suffices to show that for generic $\nu \in \faPdc$ the map 
$j \mapsto \ge^{N_0} p_{\gL + \nu +\mu} (j \otimes e_{N_0})$
is an injective endomorphism from $C^{-\infty}(G/\bar P : \gs: \nu)_\cchi$ to itself. 
The latter statement follows from Lemma  \ref{l: generic injectivity of map with j} below, which will be proven in the next section. 
\qed
\begin{lemma}
\label{l: generic injectivity of map with j}
Let $Q \in \{ P, \bar P\}.$ Then for $\gS(\fn_P,\fa_P)$-generic $\nu \in \faPdc$ the map 
\begin{equation}
\label{e: j to e after proj j}
j \mapsto \ge^{N_0} \after p_{\gL + \nu + \mu} ( j \otimes e_{N_0})
\end{equation} 
is an injective linear endomorphism of $C^{-\infty}(G/Q: \gs:\nu)_\cchi.$ 
\end{lemma}

\section{Proof of Lemma \ref{l: generic injectivity of map with j}}
We retain the notation of the previous section. 
The following lemma serves as a first step in the proof of Lemma \ref{l: generic injectivity of map with j}.
 
\begin{lemma}
\label{l: injectivity p after j e}
Let $Q \in \{P, \bar P\}.$ For $\gS(\fn_P,\fa_P)$-generic $\nu \in \faPdc$, the map 
$$
j \mapsto  p_{\gL + \nu + \mu}[j \otimes e_{N_0}] : \;C^{-\infty}(G/Q:\gs:\nu)_\cchi \to C^{-\infty}(G/Q:\gs:\nu) \otimes F
$$  
is injective with values in $( C^{-\infty}(G/Q:\gs:\nu) \otimes F)_\cchi.$ 
\end{lemma}

\proof 
That the given map attains values in $(C^{-\infty}(G/Q:\gs:\nu) \otimes F)_\cchi$ 
is a straightforward consequence of the $G$-equivariance of $p_{\gl + \nu + \mu}.$ 
We therefore focus on the asserted injectivity.

There exists a locally finite union $\cH$ of affine $\gS(\fn_P, \faP)$-{\bfhyp}
such that for $\nu \in \faPdc\setminus \cup \cH$ 
the standard intertwining operator
$$
A(\nu) = A(\bar Q , Q, \gs, \nu): C^{-\infty}(G/Q :\gs:\nu) \to 
C^{-\infty}(G/\bar Q:\gs:\nu)
$$ 
is bijective and    
maps the subspace $C^{-\infty}(G/Q:\gs:\nu)_\cchi$ bijectively onto 
$C^{-\infty}(G/ \bar Q: \gs: \nu)_\cchi.$ 
Furthermore, by the intertwining property of $A(\nu)$ we have, for such $\nu,$ that  
\begin{equation*} 
(A(\nu) \otimes I)\after p_{\gL + \nu + \mu}[j \otimes e_{N_0}] =  
p_{\gL + \nu + \mu} [A(\nu) j \otimes e_{N_0}].
\end{equation*} 
We thus see  that it suffices to establish the assertion for $Q = \bar P.$ 
Then by Lemma \ref{l: Phi after p is Phi} we have 
$$ 
\Phi_\mu(\nu) p_{\gL + \nu + \mu} (j \otimes e_{N_0}) = \Phi_\mu(\nu) ( j \otimes e_{N_0}) = i_\nu(e_{N_0}) j,
$$ 
for $\gS(\fn_P,\fa_P)$-generic $\nu\in \faPdc$ and $j \in C^{-\infty}(G/\bar P: \gs :\nu)_\cchi.$
Since  the function $i_\nu(e_{N_0})$ is non-zero on $N_0 \bar P$, it follows that 
the expression on the right of the above equality is zero if and only if $j|_{N_0 \bar P} = 0.$ 
By Corollary \ref{c: cor of main thm HC} this in turn is  equivalent to $j =0.$ The asserted injectivity follows.
\qed

To prepare for the proof of Lemma \ref{l: generic injectivity of map with j}, we need 
to introduce certain particular subspaces of generalized vectors of induced representations.
Let $Q$ be any parabolic subgroup of $G$ containing $A.$ 
We consider a continuous Hilbert representation $(\xi, H_\xi)$ of $Q.$

We denote by $W_Q(\fa)$ the centralizer of $\fa_Q$ in $W(\fa).$
From the Bruhat decompositions for $G$ and $M_{1Q}$ it follows that the map $v \mapsto N_0 v Q, \; N_K(\fa) \to N_0\backslash G / Q$ induces a bijection
from $W(\fa)/W_Q(\fa) $ onto the double coset space $N_0\backslash G / Q.$ 
Precisely one of these cosets is open in $G;$ it will be denoted by $\cO_Q.$ 
In fact, for $v\in W(\fa)$ we have
$$ 
N_0 v Q = \cO_Q\iff v \bar Q v^{-1} \supset P_0.
$$

\begin{defi}
\label{d: defi sub N zero}
\begin{equation}
\label{e: smooth subspace of induced, without nu}
C^{-\infty}(G/Q:\xi)_{N_0} \subset C^{-\infty}(G/Q:\xi)
\end{equation}
is defined to be the subspace of elements $u \in C^{-\infty}(G/Q:\xi)$ such that 
\begin{enumerate}
\itema $u$ is $N_0$-finite from the left;
\itemb $u|_{\cO_Q}$ is a continuous function
$\cO_Q \to H_\xi^{-\infty}$
\end{enumerate}
\end{defi}

Assertion (b) means that there exists a continuous function 
$\tilde u: \cO_Q \to H_\xi^{-\infty}$ such that $\tilde u(nvq) = \bp\xi^{-\infty}(q)^{-1}  \tilde u(nv)$ 
$(q\in Q, n \in N_0)$ and such that for all $\psi \in C^\infty(G/Q: \xi^*)$ 
with $\supp\, \psi \subset \cO_Q,$ 
the following identity is valid:
\begin{equation}
\label{e: partial K pairing}
\inp{u}{\psi} = \inp{\tilde u}{\psi}_K : = \int_{K/K_Q} \inp{\tilde u}{\psi}_\xi(k)  \; d\dot{k}.
\end{equation} 
Note that the integrand is a continuous complex valued function with support
contained in $K \cap \cO_Q.$ Note also that $\tilde u$ is uniquely determined.

If $v \in \cO_Q,$  the evaluation map 
$$ 
\underline\ev_v: u \mapsto \tilde u(v), \qquad C^{-\infty}(G/Q:\xi)_{N_0} \to \cH_\xi^{-\infty}
$$
is well defined and linear. 

In the special setting $\xi = \gs \otimes \nu $ with $\gs$ a unitary representation 
of $M_Q$ and $\nu \in \faQdc,$ the space on the left in (\ref{e: smooth subspace of induced, without nu})
is also denoted by $C^\infty(G/Q: \gs:\nu)_{N_0}.$ 

The following observation will allow us to connect to the case that $\bar Q$ is standard.
Let $v \in N_K(\fa)$ be such that $v \bar Q v^{-1}$ is standard. Then $\cO_Q= N_0 v Q.$ 
We write $v \xi$ for the representation of $vQv^{-1}$ in $H_\xi$ given by $v\xi(q) = \xi(v^{-1} q v).$ 
Then the right regular action by $v$  
induces  a $G$-equivariant topological linear isomorphism 
\begin{equation}
\label{e: Rv on gen vectors, xi} 
R_v: C^{-\infty}(G/Q: \xi) \;\;{\buildrel \simeq \over \longrightarrow}\; \; C^{-\infty}(G/vQv^{-1}: v\xi).
\end{equation} 
If $f$ belongs to the subspace $C^{\infty}(G/Q: \xi)$ then $R_v f$ is given by $x \mapsto f(xv)$ and belongs to $ C^{\infty}(G/vQv^{-1}: v\xi)$.

\begin{lemma}
\label{l: Nzero space invariant under Rv}
The isomorphism (\ref{e: Rv on gen vectors, xi}) restricts to a linear isomorphism
$$ 
\underline R_v:  C^{-\infty}(G/Q: \xi)_{N_0} \;\; {\buildrel \simeq \over \longrightarrow} \;\; C^{-\infty}(G/vQv^{-1}: v\xi)_{N_0}.
$$ 
Furthermore, $ \underline \ev_e \after R_v = \underline \ev_v .$
\end{lemma}

\proof 
This is a straightforward consequence of the definitions.
\qed

\begin{lemma} 
\label{l: ev and ev}
Let $Q \in\cP(A),$ $\gs$ a unitary representation of $M_Q$ and $\nu \in \faQdc.$ 
Then
\begin{enumerate}
\itema  
$
C^{-\infty}(G/ Q:\gs: \nu)_\cchi \subset  C^{-\infty}(G/Q: \gs: \nu)_{N_0}.
$ 
\itemb
If $\bar Q$ is standard, then the evaluation map $\ev_e$ defined on the space on the left is the restriction of $\underline\ev_e$ defined on the space on the right.
\itemc
Let $v\in N_K(\fa)$ be such that $\cO_Q= N_0 v Q .$ Then the evaluation map $\underline \ev_v$ 
restricts to a linear map
\begin{equation}
\label{e: map ev v on j}
\ev_v: C^{-\infty}(G/ Q:\gs: \nu)_\cchi  \to (H_{v\gs}^{- \infty})_{\cchi|_{N_0 \cap v M_Q v^{-1}}}.
 \end{equation}
 If $\cchi$ is regular, the map (\ref{e: map ev v on j})
is injective.
\end{enumerate}
\end{lemma}

\proof
We will first prove (a) - (c) under the assumption that  $Q =\bar P$ with 
$P$ standard. Then assertions (a) and (b) follow by application of 
Theorem \ref{t: evaluation of j}, Equation (\ref{e: smooth subspace of induced, without nu})
with $\xi = \gs \otimes \nu$
and the definitions of $\ev_e $ and $\underline \ev_e.$ Assertion (c) follows
from Corollary \ref{c: ev is injective}.

Now assume that 
$Q$ is general and fix $v \in N_K(\fa)$ as in (c). Then $vQv^{-1} = \bar P,$ 
with $P$ standard. We consider the isomorphism $R_v$ of (\ref{e: Rv on gen vectors, xi})
for $\xi:= \gs\otimes \nu.$ Note that $v\xi = v\gs \otimes v\nu.$ 
By Lemma \ref{l: Nzero space invariant under Rv}, 
 $R_v$ maps $C^{-\infty}(G/Q: \gs : \nu)_{N_0}$ 
onto the space $C^{-\infty}(G/vQv^{-1}: v\gs : v \nu)_{N_0}$ and the evaluation maps 
in $v$ and $e$ respectively are related by 
$\underline ev_e \after R_{v}= \underline ev_v.$

By $G$-equivariance it also follows that $R_v$ restricts to a linear isomorphism
$$ 
R_v: C^{-\infty}(G/Q: \gs : \nu)_\cchi \;\; {\buildrel \simeq \over \longrightarrow}\;\; C^{-\infty}(G/vQv^{-1}: v\gs : v \nu)_\cchi.
$$ 
As the space on the right is contained in $C^{-\infty}(G/vQv^{-1}: v\gs : v\nu)_{N_0}$ 
by the first part of the proof, it follows that the space on the left is
contained in $C^{-\infty}(G/Q: \gs : \nu)_{N_0}$ and we have the following 
commutative diagram with evaluation maps:
$$ 
\begin{array}{ccc}
C^{-\infty}(G/Q: \gs : \nu)_\cchi 
& {\buildrel R_v \over \longrightarrow} 
& C^{-\infty}(G/vQv^{-1}: v\gs : v \nu)_\cchi\\
{\!\!\!\!\!\!\!\!\!\!\!\!\!\!\!\scriptstyle \ev_v} \; \downarrow\;\;\;\;\; & &\!\!\!\! \downarrow\, {\scriptstyle \ev_e}\\
(H_{v\gs}^{-\infty})_{\cchi|_{N_0 \cap vM_Qv^{-1}}} 
& {\buildrel I \over \longrightarrow}  & 
(H_{v\gs}^{-\infty})_{\cchi|_{N_0 \cap vM_Qv^{-1}}}
\end{array}
$$ 
If $\cchi$ is regular, then $\ev_e$ is injective, and the injectivity of $\ev_v$ follows.
\qed

We return to the setting that $Q \in \cP(A)$ and that $(\xi, H_\xi)$ is a continuous Hilbert
representation of $Q.$  The representation $\xi^*$ of $Q$ in $H_\xi$ is defined as in Remark \ref{r: dual of pi}.

\begin{prop}
\label{p: nilpotent space invariant under g}
The subspace $C^{-\infty}(G/Q:\xi)_{N_0}$ of 
$C^{-\infty}(G/Q:\xi) $ is invariant under the left action by $\fg.$
\end{prop}
\proof
We fix $v \in  N_K(\fa) $ such that $v Q v^{-1} = \bar P,$ with $P$-standard.
By Lemma \ref{l: Nzero space invariant under Rv} and since $R_v$ 
of (\ref{e: Rv on gen vectors, xi}) is $\fg$-equivariant, it suffices to 
establish the assertion of the proposition with $vQv^{-1}$ in place of $Q.$
In other words, without loss of generality we may and will assume from the start
that $Q = \bar P,$ with $P$ standard.

We recall that the left action by an element $X \in \fg$ on $u \in \cF:= C^{-\infty}(G/\bar P: \xi)$ 
is defined by 
$$ 
\inp{L_X u}{\psi} = \inp{u}{L_{X^\veec} \psi},
$$ 
for all $\psi \in C^\infty(G/Q:  \xi^*).$

It is sufficient to show that for $u \in \cF_0:= \cF_{N_0}$ and $X \in \fg$ the element $ L_X u \in \cF$ restricts 
to a continuous $H_\xi^{-\infty}$-valued function on $N_0\bar P.$ Indeed the $N_0$-finiteness is obvious, since $L_{n} L_X = L_{\Ad(n) X} L_n$ for all $n \in N_0.$ For the first statement we need a suitable interpretation of $L_X$ on the space $\cF_{N_0}.$

Let $u\in \cF_0.$ The span  $E$
of the left $N_0$-translates of $u$ is a finite dimensional subspace of $\cF_0.$ 
The restriction of $L|_{N_0}$ to $E$ is a finite dimensional
representation of $N_0$, which we denote by $\omega.$ As $\omega$ is the restriction
of the continuous representation $L^{-\infty}|_{N_0}$ to $E$, it is continuous. By finite 
dimensionality it follows that $(\omega, E)$ is smooth.
We claim that for $X \in \fn_0$, we have 
$L_X = \omega(X)$ on $E.$ To see this, note that for $\psi \in C^\infty(G/\bar P: \xi^*),$
\begin{eqnarray*}
\inp{L_X u}{\psi}
&=& \inp{u}{- L_X \psi} 
=
\left. \frac{d}{dt}\right|_{t=0} \inp{u}{L_{\exp tX}^{-1} \psi }\\
&=& 
\left. \frac{d}{dt}\right|_{t=0} \inp{L_{\exp tX}u}{\psi} 
 =   \left. \frac{d}{dt}\right|_{t=0} \inp{\omega(\exp tX)u}{\psi}\\
&=& \inp{\omega(X) u}{\psi}.
\end{eqnarray*} 
It follows from this that $L_X u \in \cF_0,$ hence $L_X u$ is continuous on $N_0\bar P,$ 
for all $u \in \cF_0$ and $X \in \fn_0.$ 

We now fix a general element $X \in \fg$ and will establish the continuity asssertion
for $L_X u.$ Since  $N_0 \bar P = N_P \bar P,$ 
the assertion of continuity of $u \in \cF$ on $N_0 P$  means that there exists
a unique continuous function $\tilde u : N_P\bar P \to H_\xi^{-\infty}$ 
such that $\tilde u(x\bar p) = \bp\xi(\bar p)^{-1} \tilde u(x)$ $(x \in N_P\bar P, \bar p \in \bar P)$
and such that
(\ref{e: partial K pairing}) is valid. By the usual transformation of variables corresponding
to the open embedding $N_P \to G/\bar P\simeq K/K_P$ that equation is equivalent
to 
$$ 
\inp{u}{\psi}= \inp{\tilde u}{\psi}_{N_P}:= \int_{N_P} \inp{\tilde u}{\psi}_\xi (n) \; dn.
$$ 
Since $\fg = \fn_P \oplus {\rm Lie}(\bar P)$ it follows  
that for every $n\in N_P$ we may 
write 
$$
\Ad(n)^{-1} X = \Ad(n)^{-1}  Y(n) + Z(n)
$$  with $Y: N_P \to \fn_P$ and $Z: N_P \to {\rm Lie}(\bar P)$ smooth functions.
Let $Y_1, \ldots, Y_k$ be a basis for $\fn_P.$ 
Then we see that $Y(n) = \sum_i y^i(n) Y_i$  with $y^i: N_P \to \R$ smooth functions.

Let now $\psi \in C^\infty(G/\bar P: \xi^*)$ have support contained in 
$N_P \bar P.$ 
Then 
$$
\inp{L_X u}{\psi} = - \inp{u}{L_X \psi} = - \inp{\tilde u} {L_X\psi}.
$$ 
Furthermore, for $n \in N_P,$ 
\begin{eqnarray*}
L_X \psi(n) & = & [ L_{Y(n)}\psi] (n) +[ L_{\Ad(n) Z(n)} \psi](n)\\ 
& = & \sum_i y^i(n) L_{Y_i} \psi (n)-  [R_{Z(n)}\psi](n) \\
&=&  \sum_i L_{Y_i} [y^i \psi](n) - \sum_i L_{Y_i}(y^i)(n)\psi(n) + \bp \xi(Z(n))\psi(n).
\end{eqnarray*}
By what we established above, $L_{Y_i} u = \omega(Y_i) u$ is given by 
a continuous function $\tilde u_i$ on $N_0 \bar P.$ 
Let $\hat y^i$ denote the unique element of $C^\infty (N_P\bar P)$ given  
by $\hat y^i (n \bar p) = y^i (n).$
Then $\hat y^i \psi \in C^\infty(G/\bar P : \xi^*)$ (extension by zero outside $N_P \bar P$), and we see that, for each $1 \leq i \leq k,$
$$
- \inp{\tilde u} {L_{Y_i} [y^i \psi]}_{N_P} = \inp{L_{Y_i}u }{\hat y^i\psi} =  \inp{\tilde u_i}{\hat y^i \psi}_{N_0} =
\inp{y^i \tilde u_i }{\psi}_{N_0}.
$$ 

This leads to
$$ 
- \inp{L_X u}{\psi} = \inp{\sum_i  y^i \tilde u_i  - L_{Y_i}(y^i) \tilde u - [\bp \xi\after Z]\tilde u}{\psi}_{N_0}.
$$ 
As $\psi|_{N_P}$ ranges over all functions of $C_c^\infty(N_P, H_\xi^\infty),$ 
it follows that on $N_0\bar P = N_P \bar P,$ the generalized function 
$- L_X u$ is represented by 
$$ 
\sum_i [y^i  \tilde u^i  - L_{Y_i}(y^i) \tilde u] - [\bp\xi\after Z]\tilde u.
$$ 
The latter function is obviously continuous $N_0\bar P \to H_\xi^{-\infty}.$ 
\qed
We now assume that $Q \in \cP(A),$ that $\gs$ is a unitary representation of $M_Q,$ 
and that  $\nu \in \faQdc,$ and define the representation $\xi_\nu$ of $Q= M_Q A_Q N_Q$ in $H_\gs$ by 
$$
\xi_\nu = \gs \otimes \nu \otimes 1.
$$
Furthermore, we assume that $(\pi, F)$ is a continuous finite dimensional representation of $G.$ 
From (\ref{e: def gf mininfty tensprod}) we recall the existence
of a unique $G$-equivariant topological linear isomorphism 
\begin{equation}
\label{e: gf gl} 
\gf_\nu^{-\infty}: C^{-\infty}(G/Q:\xi_\nu) \otimes F \;\;{\buildrel \simeq \over \longrightarrow}\;\; C^{-\infty}(G/Q:\xi_\nu \otimes \pi|_Q) 
\end{equation} 
determined by 
$$ 
\gf_\nu^{-\infty}(u \otimes e)  = (1 \otimes \pi(\dotvar)^{-1}) ( u  \otimes  e)
$$ 
on the subspaces with $C^\infty$ in place of $C^{-\infty}.$ 
The inverse is given by $w \mapsto ( I \otimes \pi(\dotvar) ) w$ 
on the mentioned subspaces.

\begin{cor}
\label{c: proj preserves Nzero space}
For $\gS(\fn_Q,\fa_Q)$-generic $\nu \in \faQdc$, the endomorphism $p_{\gL + \mu + \nu}$ of 
the space 
 $C^{-\infty}(G/Q: \xi_\nu \otimes \pi|_Q)$ preserves the subspace 
 $C^{-\infty}(G/Q:  \xi_\nu \otimes \pi|_Q)_{N_0}.$ 
\end{cor}

\proof
This follows by combining Proposition \ref{p: nilpotent space invariant under g} applied
to $\xi = \xi_\nu$ and the characterization of $p_{\gL +\nu +\mu}$ in Corollary \ref{c: q p is Z nu}.
\qed

\begin{lemma}
\label{l: gf nu and ev v}
Let $Q \in \cP(A)$ and let $v \in N_K(\fa)$ be such that $v \bar Q v^{-1}$ is standard.
The map $\gf_\nu^{-\infty}$ of (\ref{e: gf gl}) restricts to a $(\fg, N_0)$-equivariant linear isomorphism 
\begin{equation}
\label{e: underline gf nu}
 \underline\gf_\nu:  C^{-\infty}(G/Q:\gs: \nu)_{N_0} \otimes F\;\;{ \buildrel \simeq \over \longrightarrow}\;\; C^{-\infty}(G/Q: \xi_\nu \otimes \pi|_Q)_{N_0}
 \end{equation}
 which satisfies
\begin{equation}
\label{e: underline gf eval}
\underline{\ev}_v \after \underline \gf_\nu = \underline \ev_v \otimes \pi(v)^{-1}.
\end{equation}
 \end{lemma}
\medbreak
\proof
We note that $(\xi_\nu)^* = \gs \otimes -\bar \nu =\xi_{-\bar \nu}.$ We assume that
$F$ is equipped with a Hermitian inner product, and $F_*$ denote
the finite dimensional Hilbert space $F$, equipped with the conjugate representation $\pi^*.$ 
Let ${}^*\gf_{\nu}$ be the equivariant isomorphism from $C^\infty(G/Q: \xi_{-\bar \nu}) \otimes F_*$ onto 
$C^\infty(G/Q: \xi_{-\bar \nu} \otimes \pi^*|_Q)$ 
as defined in (\ref{e: natural iso gf smooth}).

Let $f \in C^{-\infty}(G/Q:\gs:\nu)_{N_0} \otimes F.$ Then the restriction of $f$ to $\cO_Q =N_0 v Q$ 
is continuous $\cO_Q \to H_{\gs}^{-\infty} \otimes F.$ It follows from the definition 
of $\gf_\nu^{-\infty}$ as in (\ref{e: def gf mininfty tensprod}) that 
$$ 
\inp{\gf_\nu^{-\infty}(f) }{g} = \inp{f}{({}^*\gf_\nu)^{-1}(g)},
$$ 
for all $g \in C^\infty(G/Q: \xi_{-\bar\nu}) \otimes F_*.$
In particular this is true for all such $g$ with support contained in $N_0 v Q.$ 
In that case  the above equality tells us that
\begin{eqnarray*}
\inp{\gf_\nu^{-\infty}(f) }{g} &= &\int_{K/K_Q} \inp{f(k)}{({}^*\gf_\nu)^{-1}(g)(k)}_\gs \, d\dot{k}\\
&=&\int_{K/K_Q} \inp{f(k)}{(I \otimes \pi^*(k))g(k)} _\gs \,d\dot{k}\\
&=& \int_{K/K_Q} \inp{(I \otimes \pi(k)^{-1})f(k)}{g(k)}_\gs\,  d \dot{k}.
\end{eqnarray*}
It follows from this that, for $x \in K \cap \cO_Q,$ 
\begin{equation}
\label{e: gf nu eval}
\gf_\nu^{-\infty}(f)(x) = (1 \otimes \pi(x))^{-1} (\ev_x \otimes I)(f).
\end{equation}
By $Q$-equivariance this equality is true for all $x \in \cO_Q.$ 
From this it is immediately clear that $\gf_\nu^{-\infty}(f)$ belongs to the space on the right in 
(\ref{e: underline gf nu}). Moreover, by substituting $x = v$ in (\ref{e: gf nu eval}) 
we obtain (\ref{e: underline gf eval}) when applied to $f.$ 

By using a similar argument involving the maps $[\gf_\nu^{-\infty}]^{-1}$ and 
${}^*\gf_\nu$ one sees that the map $\underline \gf_\nu$ is a linear isomorphism as asserted.

Finally, the $(\fg, N_0)$-equivariance of $\underline \gf_\nu$ follows from the similar 
equivariance of $\gf_\nu$ combined with Proposition \ref{p: nilpotent space invariant under g}.
\qed

\begin{cor}
\label{c: proj preserves Nzero space two} 
For $\gS(\fn_Q,\fa_Q)$-generic $\nu \in \faQdc$, the endomorphism $p_{\gL + \mu + \nu}$ of 
the space 
 $C^{-\infty}(G/Q: \gs: \nu) \otimes F$ preserves the subspace 
 $C^{-\infty}(G/Q:  \gs: \nu)_{N_0}\otimes F.$ 
\end{cor}

\proof 
Since the map $\gf_\nu^{-\infty}$ of (\ref{e: gf gl}) is a $G$-equivariant isomorphism 
from the space $C^{-\infty}(G/Q: \gs : \nu) \otimes F$ onto $C^{-\infty}(G/Q : (\gs \otimes \nu)\otimes \pi|_Q),$
we have  $\gf_\nu^{-\infty} \after p_{\gL + \mu + \nu} = \gf_\nu^{-\infty} \after p_{\gL + \mu + \nu}.$ 
The result now follows by  combining Corollary \ref{c: proj preserves Nzero space}
with  Lemma \ref{l: gf nu and ev v}.
\qed

Finally, we are prepared to complete the proof announced in the title of this section.
\medno
{\em Proof of Lemma \ref{l: generic injectivity of map with j}.\ }
By application of standard intertwining operators as in the proof of Lemma \ref{l: injectivity p after j e},
 we may reduce to the case that $Q = P$ (recall that $P$ is standard). 
 In this case, we argue as follows. For $\nu \in \faPdc$ we 
 denote by  $\bp H_{\gs, \nu}$ the space $H_\gs$ on which $P=M_P A_PN_P$ acts by 
 $man \mapsto a^{\nu + \rho_P}\gs(m).$ The space $\C e_{N_0}$ is a $P$-invariant subspace of $F,$ 
 on which $P $ acts by $man \mapsto a^\mu.$ We write $\pi|_{P,\C e_{N_0}}$ 
 for the restriction of $\pi|_P$ to this subspace. The inclusion map $\C e_{N_0} \to F$ 
 is $P$-equivariant, hence 
 induces a $G$-equivariant continuous linear map 
$$
\iota_\nu:  C^{-\infty}(G /P: \xi_\nu \otimes \pi|_{P,\C e_{N_0}}) \embeds C^{-\infty}(G /P: 
\xi_\nu \otimes \pi|_P),
 $$
 see Lemma \ref{l: functoriality C min infty} for details.
 
In the sequel we shall briefly write $p_\nu$ for $p_{\gL + \mu +\nu}.$ 
Since $H_{\gs, \nu} \otimes \C e_{N_0}$ is naturally isomorphic to $ H_{\gs, \nu + \mu}$ as
a $P$-module, it follows that for ($\gS(\fn_P,\fa_P)$-)generic $\nu\in \faPdc$ the projection 
\begin{equation}
\label{e: proj p nu in End}
p_\nu \in \End(C^{-\infty}(G /P: \xi_\nu \otimes \pi|_P))
\end{equation}
equals the identity on the image of 
$\iota_\nu.$  On the other hand, it follows from Lemma \ref{l: iso image pr}
combined with the isomorphism 
$\Ind_P^G (\xi_\nu ) \otimes \pi \simeq \Ind_P^G (\xi_\nu  \otimes \pi|_P)$
that for generic $\nu \in \faPdc,$ 
$$ 
{\rm im}(p_{\nu})_K \simeq C(G/P: \gs: \mu + \nu)_K \simeq {\rm im}(\iota_\nu)_K,
$$
as $(\fg, K)$-modules. Since (\ref{e: proj p nu in End}) is a continuous projection 
 its image ${\rm im} (p_\nu)$ is a closed subspace of $C^{-\infty}(G /P: \xi_\nu \otimes \pi|_P).$
We now infer that for generic $\nu \in \faPdc$ we have
$$ 
{\rm im}(p_{\nu}) = \cl({\rm im}(\iota_\nu)_K),
$$ 
where ${\rm cl}$ indicates that the closure in $C^{-\infty}(G/P:\xi_\nu \otimes \pi|_P)$ is taken.
To characterize this closure in a useful way, we fix a Hermitian inner product on $F$
and define the continuous representation $\pi^*$ of $G$ on it by $\pi^*(p) = \pi(p^{-1})^{*}.$ 
The restriction of $\pi^*|_P$ to the $\pi^*(P)$-invariant subspace $E = (\C e_{N_0})^\perp$ of $F$ is denoted by $\pi^*_E.$ Clearly, $\C e_{N_0} = E^\perp.$ We view $C^\infty(G/P: \xi_{-\bar \nu} \otimes \pi^*_E)$ as an invariant subspace of $C^\infty(G/P: \xi_{-\bar \nu} \otimes \pi^*_E)$.
Via the sesquilinear pairing
\begin{equation}
\label{e: ann pairing}
C^{-\infty}(G/P: \xi_\nu \otimes  \pi) \times C^\infty(G/P: \xi_{-\bar \nu} \otimes \pi^*)\to \C
\end{equation}
we accordingly define ${\rm Ann}_\nu$ to be the annihilator of $C^\infty(G/P: \xi_{-\bar \nu} \otimes \pi^*_E)$
in the first of the spaces in (\ref{e: ann pairing}). This annihilator is closed and on the $K$-finite level it is readily seen that
$({\rm Ann}_\nu)_K = {\rm im} (\iota_\nu)_K.$ It follows that the annihilator equals the closure of
${\rm im} (\iota_\nu)_K.$ Hence,
$$ 
 {\rm im}(p_{\nu}) = {\rm Ann}_\nu .
 $$ 
 We now select $v \in N_K(\fa)$ such that $\cO_P = N_0 v P$ is open in $G.$ 
Note that $\cO_P = v \bar N_P P.$ 
 \begin{lemma}
 \label{l: char annihilator on N zero}
 Let $\nu \in \faPdc,$ write $\xi_\nu = \gs \otimes \nu \otimes 1$ and let
 $ 
 u \in  C^{-\infty}(G/P:\xi_\nu \otimes \pi|_P)_{N_0}.
 $
 If $u \in {\rm Ann}_\nu $ then 
 $
 u|_{\cO_P} \in C(\cO_P, H_\gs^{-\infty} \otimes \C e_{N_0})^P.
 $
 \end{lemma}
 
 \proof
Let $u$ fulfill the hypothesis.
Then it follows from Definition \ref{d: defi sub N zero}
that the restriction of $u$ to $\cO_P$ 
is continuous with values in $H_\gs^{-\infty} \otimes F.$ 
This means that there exists a continuous function $\tilde u \in N_0 v P \to H_\gs^{-\infty} \otimes F$ 
such that $\tilde u(n v p) = [\bp \xi^{-\infty}(p) \otimes \pi(p) ]^{-1} \tilde u(n v)$  $(p \in P, n \in N_0),$ 
and such that for all $\psi \in C^\infty(G/P: \xi_{-\bar \nu} \otimes \pi^*)$  with support contained
in $\cO_P$ we have (\ref{e: partial K pairing}), which by the substitution of variables $k = 
v \kappa_P(\bar n)K_P$ may be rewritten 
as
$$ 
\inp{u}{\psi} = \int_{\bar N_P} \inp{\tilde u ( v \bar n) }{\psi(v \bar n)}_{\gs \otimes \pi} \,d \bar n. 
$$ 
 Let now $h \in H_\gs^\infty,$ $f \in E= (\C e_{N_0})^\perp$ and  $\phi \in C^\infty_c(\bar N_P),$ and define 
$ \psi \in C^\infty(G/P:\xi_{-\bar \nu} \otimes \pi^*)$ by the requirements
$$ 
\psi(v \bar n ) = \phi(\bar n) ( h \otimes f) , \qquad \psi = 0 \;\;{\rm on} \;\;G\setminus v\bar N_P P.
$$ 
Then it follows from $u \in {\rm Ann}_\nu$ that
$$ 
\int_{\bar N_P} \phi(n) \inp{\tilde u (v \bar n)}{h \otimes f} \; d\bar n = 0.
$$ 
As this is valid for all $\phi$ as above, the continuous function $\bar n \mapsto \inp {\tilde u(v \bar n)}{ h \otimes f}$ 
is zero. This implies that for each $\bar n \in \bar N_P$ the element $\tilde u(v\bar n) \in H_\gs^{-\infty} \otimes F$
satisfies 
$$ 
\inp{\tilde u(v\bar n)}{h \otimes f)} = 0,\qquad (h \in H_\gs^\infty, f \in E).
$$ This in turn implies that $u(v \bar n) \in 
H^{-\infty}_\gs \otimes E^\perp = H^{-\infty}_\gs \otimes \C e_{N_0}.$ 
Since $v \bar N_P P = N_0 v P $ this finishes the proof.
\qed

We proceed with the completion of the proof of Lemma \ref{l: generic injectivity of map with j}. 
From Lemma \ref{l: injectivity p after j e} it follows that for generic $\nu \in \faPdc$ 
the map
\begin{equation}
\label{e: j to gf nu after p} 
j \mapsto   \gf_\nu^{-\infty} \after p_{\nu} (j \otimes e_{N_0})
\end{equation}
is injective from $C^{-\infty}(G/P:\gs:\nu)_\cchi$ to 
$C^{-\infty}(G /P: (\gs\otimes\nu) \otimes \pi|_P)_\cchi.$
Furthermore, since $\gf_\nu$ is an intertwining isomorphism,
it follows that  
$$ 
\gf_\nu^{-\infty} \after p_{\nu } = p_{\nu} \after \gf_\nu^{-\infty}.
$$ 
We thus see that the map (\ref{e: j to gf nu after p}) maps $C^{-\infty}(G/P:\gs:\nu)_\cchi$
injectively 
to 
$$ 
{\rm im}(p_\nu) \cap C^{-\infty}(G/P: (\gs \otimes \nu) \otimes\pi|_P)_\cchi
\subset {\rm Ann}_\nu \cap C^{-\infty}(G/P: (\gs \otimes \nu) \otimes\pi|_P)_{N_0}.
$$ 
Applying Lemmas \ref{l: ev and ev} and \ref{l: char annihilator on N zero}
we now
see that the map 
$$ 
j \mapsto   \ev_v [ \gf_\nu^{-\infty} \after p_{\gL +\nu + \mu} (j \otimes e_{N_0})]
$$
is injective from   $C^{-\infty}(G/P:\gs:\nu)_\cchi$ to $H_{\gs}^{-\infty}\otimes \C e_{N_0}.$
Put 
\begin{equation}
e^v:= e^{N_0} \after \pi(v) \in F^*;
\end{equation}
then $e^v(e_{N_0}) \neq 0$, since otherwise
the matrix coefficient $x \mapsto e^v(\pi(x) e_{N_0})$ would be zero
on $N_0 vP$ hence on $G$, contradicting
the irreducibility of $\pi.$ It thus follows that 
\begin{equation}
\label{e: j mapsto vs 2}
j \mapsto   (I \otimes e^v) [ \ev_v  [ \gf_\nu^{-\infty} \after p_{\nu} (j \otimes e_{N_0}) ]]
\end{equation}
is injective  $C^{-\infty}(G/P:\gs:\nu)_\cchi \to H_{\gs}^{-\infty}\otimes \C.$
On the other hand, $j \mapsto j \otimes e_{N_0}$ maps
$C^{-\infty}(G/P:\gs:\nu)_\cchi$ into $C^{-\infty}(G/P:\gs : \nu)_{N_0} \otimes F.$ 
By Lemmas \ref{l: ev and ev}, \ref{l: gf nu and ev v} and Corollary \ref{c: proj preserves Nzero space two} 
it now follows that (\ref{e: j mapsto vs 2}) equals
\begin{equation}
\label{e: j mapsto vs 3}
j \mapsto (I \otimes e^v)\after \underline{\ev}_v \after \underline\gf_\nu \after  \underline p_\nu (j \otimes e_{N_0})
\end{equation}
which is therefore  an injective map 
$C^{-\infty}(G/P:\gs:\nu)_\cchi \to H_{\gs, \nu}^{-\infty}.$ 
By Lemma \ref{l: gf nu and ev v} the above map (\ref{e: j mapsto vs 3}) equals
\begin{eqnarray*}
j   & \mapsto  &  (I \otimes e^v) \after (\underline{\ev}_v  \otimes \pi(v)^{-1}) \after \underline p_{\nu} (j \otimes e_{N_0})\\
&=& (\underline\ev_v \otimes I) \after (I \otimes e^{N_0}) \after \underline p_{\nu} (j \otimes e_{N_0}).
\end{eqnarray*}
The injectivity of the latter map implies the injectivity of 
$$ 
j \mapsto (I \otimes e^{N_0}) \after \underline p_{\nu }(j\otimes e_{N_0})
=
(I \otimes e^{N_0}) \after p_{\nu }(j\otimes e_{N_0}).
$$ 
as a map from $C^{-\infty}(G/P:\gs : \nu)_\cchi $ to $C^{-\infty} (G/P: \gs:\nu) \otimes \C$.
Since $\ge^{N_0} = m \after (I \otimes e^{N_0}),$ with $m$ injective, 
see (\ref{e: defi ge N zero}), the required injectivity of the map (\ref{e: j to e after proj j}) 
follows.
\qed

\section{Holomorphy and uniformly moderate estimates}
\label{s: holo and ume}
In this section, we assume that
$P$ is a standard parabolic subgroup and $(\gs, H_\gs)$ a discrete series representation of $M_P.$ 
We will first prove the following result which is inspired by Wallach \cite[Thm. 15.4]{Wrrg2}.
Let $\gd > 0.$

\begin{thm}
\label{t: estimate by functional equation}
For every $R \leq \gd$ the function $\nu \mapsto j(\bar P,\gs, \nu),$ originally defined for $\nu \in \faPdc(P,\gd)$ 
allows a meromorphic extension to $\faPdc(P,R)$ as a function with values in the space
$(H_{\gs,\cchi_P}^{-\infty})^* \otimes C^{-\infty}(K/K_P:\gsP).$ 

Furthermore, there exists a non-trivial polynomial function $p_R \in P(\faPd)$ and constants $s, N \in \N$
and $C > 0$ such that for all $\eta \in H_{\gs,\cchi_P}^{-\infty}$ the extended function $\nu \mapsto 
p_R(\nu) j(\bar P, \gs, \nu,\eta)$ is holomorphic $C^{-s}(K/K_P: \gsP)$-valued on $\faPdc(P,R)$ 
and satisfies the estimate
$$ 
\|p_R(\nu) j(\bar P,\gs, \nu, \eta)\|_{-s} \leq C (1 + |\nu|)^N \|\eta\|, \qquad (\nu \in \faPdc(P,R)).
$$ 
\end{thm}

\proof 
First of all, by Proposition \ref{p: holomorphy j in dominant region} 
the above result is true for $R = \gd,$ with $s = N=0$   Let $\mu \in \gL^{++}(\faP).$ Then $\inp{\mu}{\ga} > 0$ for all $\ga \in \gS(\fn_P, \faP).$  Let $m>0$ be fixed and strictly smaller than the minimum of the numbers
$\inp{\mu}{\ga},$ for $\ga \in \gS(\fn_P, \faP).$ Then 
$$ 
\faPdc(P, R-m) + \mu \subset \faPdc(P,R). 
$$ 
We will show that if the assertions of the theorem are valid for $R \leq \gd,$ 
they are also valid with $R$ replaced by $R - m.$ The result then follows by 
induction.

Assume the result to be proven for a given $R,$ with constants $s, N_R$ and $C_R$ 
in place of $s, N, C.$ 
Let $D_\mu(\gs, \nu)$ and $R_\mu(\gs, \nu)$ be as in Theorem \ref{t: the functional equation}.
By holomorphic continuation, this functional equation is still valid on $\faPdc(P,R) $
for the extended function $j(\bar P, \gs, \dotvar).$ 
Let $q: \faPdc \to \C$ be a  non-trivial polynomial function such that $q R_\mu(\gs, \dotvar)$ is polynomial
of degree $d',$ with values in $\End (H_{\gs,\cchi_P}^{-\infty}).$ 

For  $\nu \in \faPdc(R - r)$ we define
$$
j_e(\bar P:\gs: \nu)(\eta) = D_\mu(\gs, \nu) j(P, \gs ,  \nu + \mu) R_{\mu}(\gs, \nu) \eta.
$$ 
Let $r,d, C \in \N$ be the constants of Proposition \ref{p: estimate of generalized operator norm D mu}. Then it follows by application
of the mentioned proposition that $\nu \mapsto  q(\nu) p_R(\nu +\mu ) j_e(\bar P:\gs: \nu)(\eta)$ 
is holomorphic on $\faPdc(P, R - m)$ with values in $C^{-s-r}(K/K_P:\gsP).$ Furthermore, 
\begin{eqnarray*}
\lefteqn{\!\!\!\!\!\!\!\!\|q(\nu) p_R(\nu + \mu) j_e(\bar P:\gs: \nu)(\eta)\|_{-s-r} } \\
& \leq &  C( 1 + |\nu|)^{d}  \| p_R(\mu + \nu) j(P,\gs, \nu + \mu, q(\nu)R_\mu(\gs, \nu)\eta)\|_{-s}\\
&\leq & C C_R (1 + |\nu|)^d ( 1 + |\nu + \mu|)^{N_R} \|q(\nu) R_\mu(\gs, \nu)\|_{\rm op} \|\eta\| \\
&\leq & C_{R-m} (1 +|\nu|)^{d + N_R + d'}\|\eta\|,
\end{eqnarray*}
with $C_{R-m} > 0$ a constant which is uniform for $\nu \in \faPdc(P, R-m).$ 
By the functional equation of Theorem \ref{t: the functional equation} it follows that
$$ 
q(\nu) p_R(\nu + \mu) j_e(\bar P,\gs, \nu)(\eta) = q(\nu) p_R(\nu + \mu) j(\bar P, \gs, \nu)(\eta)
$$ 
for all $\nu \in \faPdc(P,R).$ 
This shows that the $(H_{\gs \cchi_P}^{-\infty})^* \otimes C^{-\infty}(K/K_P:\gsP)$-valued function $\nu \mapsto j_e(\bar P, \gs, \nu)$ is the meromorphic extension of the original 
$\cH_{\gs,\cchi_P}^{-\infty}$-valued function $j(\bar P, \gs, \dotvar)$ defined on $\faPdc(P, \gd).$ 
The proof is complete.
\qed

It follows from the above result that as a  $(H_{\gs, \cchi_P}^{-\infty})^*\otimes C^{-\infty}(K/K_P: \gsP)$-valued 
function, the function $\nu \mapsto j(\bar P, \gs, \nu)$ has a meromorphic extension to all of 
$\faPdc.$ This meromorphic extension will be denoted by the same symbol.

\begin{lemma}
For $\eta \in H_{\gs,\cchi_P}^{-\infty}$ and a regular point $\nu \in \faPdc$
the element $j(\bar P, \gs, \nu, \eta) \in C^{-\infty}(K/K_P: \gsP)$ satisfies the transformation 
rule 
$$ 
\pi_{\bar P, \gs, \nu}^{-\infty} (n) j(\bar P, \gs, \nu, \eta) = \cchi(n)  j(\bar P, \gs, \nu, \eta), \qquad  (n \in N_0).
$$ 
\end{lemma}

\proof 
This follows by analytic continuation.
\qed

\begin{lemma} 
\label{l: restricted j holomorphic}
Let $\gf \in C^{\infty}(K/K_P:\gsP)$ have compact support contained in 
the set $K \cap N_0\bar P.$  Then for every $\eta \in H_{\gs,\cchi_P}^{-\infty},$ the 
meromorphic function $\nu \mapsto \inp{j(\bar P, \gs, \nu)(\eta)}{\gf}$ is holomorphic.
\end{lemma}

\proof
It follows from the definition
of $j(\bar P,\gs, \nu, \eta )$ for $\nu \in \faPdc(P,0)$ that the restriction $j(\bar P,\gs, \nu, \eta)|_{K\cap N_P \bar P}$ is given by the 
continuous function $K\cap N_0 \bar P  \to \cH_\gs^{-\infty}$ described by the formula
\begin{equation}
\label{e: defi j eta nu} 
j_\eta(\nu): k \mapsto a(k)^{-\nu + \rho_P}  \cchi(n_P(k)) \gs(\mu_P(k))^{-1}\eta.
\end{equation}
By this we mean that for a function $\gf \in C^{\infty}(K/K_P\col \gsP)$ 
with compact support contained in $K \cap N_0\bar P$ we have 
$$ 
\inp{j(\bar P, \gs, \nu)(\eta)}{\gf} = \int_{K/K_P} \inp{ j_\eta(\nu)(k)}{\gf(k)}_\gs dk
$$ 
From
$$ 
j_\eta(\gl) = a(\dotvar)^{(\nu-\gl)} j_\eta(\nu),\qquad (\gl \in \faPdc),
$$ 
we see that $j_\eta$ extends to a holomorphic function from $\faPdc$ to $C(K\cap N_0 \bar P, H_\gs^{-\infty}).$ 
Accordingly, it follows that for $\gf \in C^{-\infty}(K/K_P:\gsP)$ with compact support in 
$K \cap N_P\bar P$ the $\C$-valued function 
$$ 
\nu \mapsto \inp{j(\bar P, \gs, \nu)(\eta)}{\gf}
$$ 
is holomorphic on $\faPdc.$ 
\qed

\begin{thm}
\label{t: j is holomorphic}
The map $\nu \mapsto j(\bar P,\gs, \nu)$ is holomorphic as a function on $\faPdc$ with values in 
the complete locally convex space $(H_{\gs,\cchi_P}^{-\infty})^*\otimes C^{-\infty}(K/K_P: \gsP)$.
Here $C^{-\infty}(K/K_P:\gsP)$ is understood to be equipped with the direct limit topology, see
the text below (\ref{e: limit topology}).

\end{thm}

\proof 
Let $R \leq 0$ and let $p = p_R$ be as in Theorem \ref{t: estimate by functional equation}. 
Let $\Omega:= \faPdc(P, R)$ and let $X$ be the zero set of $p$ in $\Omega.$ Fix 
$\eta \in H_{\gs,\cchi_P}^{-\infty}.$ It follows from the theorem that 
$j: \nu \mapsto j(\bar P, \gs, \nu)$ is holomorphic as a function from $\Omega \setminus X$ 
to the complete locally convex space $V: = C^{-\infty}(K/K_P : \gsP).$

By Theorem \ref{t: Hartog type continuation} (Appendix) it suffices to show that
$j$ admits an extension to a holomorphic function $\Omega \setminus X_s \to V,$ 
with $X_s = X \setminus X_r,$ where $X_r$ is the set of points $\nu_0 \in X$ at which
$X$ is a complex differentiable submanifold of co-dimension $1.$  It is readily verified
that $X_r $ is open in $X;$ therefore, $X_s$ is closed in $X$ hence in $\Omega.$ 

Let $\nu_0 \in X_r.$ Then it suffices to show that there exists an open neighborhood  $\Omega_0$
of $\nu_0$ in $\Omega$ such that $j|_{\Omega_0 \setminus X}$ admits an extension to a holomorphic
function $\Omega _0 \to V.$ See also Lemma \ref{l: observation about extension} (Appendix).

By definition of $X_r$ there exists an open neighborhood $\omega$ of $\nu_0$ in $\Omega$ 
such that $X_0:= X\cap \omega$ is a connected complex differentiable submanifold of codimension $1.$ 
Let $\xi \in \faPdc$ be such that 
$$
T_{\nu_0} X_0 \,\oplus\,  \C \xi = \faPdc.
$$ 
Then it follows that the 
map $\gf: (\nu , z)\mapsto  \nu + z \xi$ is a local holomorphic diffeomorphism at $(\nu_0, 0).$ Replacing $\omega$ by a smaller neighborhood if necessary, and taking $r >0$ sufficiently small,  we arrive at the situation that $\gf: X_0 \times D(0,r) \to \faPdc$ is a holomorphic diffeomorphism onto an open neighborhood $\Omega_0$ of $\nu_0$ in $\faPdc$ and that $\gf(X_0 \times \{0\}) = X \cap \Omega_0.$ 

Put $D = D(0, r).$ Then $j^* = j\after \gf|_{X_0 \times D\setminus \{0\}}$ is a holomorphic function $X_0 \times D\setminus \{0\} \to V$ and it suffices to show that this function extends to a 
holomorphic function $X_0 \times D \to V.$ 

Since $j^*: (\nu, z) \mapsto j(\nu + z\xi)$ is holomorphic on $X_0 \times (D\setminus \{0\})$ 
it has a Laurent series expansion in $z$ of the form
$$ 
j(\nu + z\xi) = \sum_{k\in \Z} c_k(\nu) z^k,
$$ 
with $c_k: X_0 \to V$ holomorphic, for all $k \in \Z.$ 
 Thus, it suffices to show that $c_k = 0 $ for $k < 0.$ 

The zero set of $p\after \gf$ equals $\gf^{-1}(X) = X_0 \times \{0\}. $ Hence, 
there exists a constant $d\geq 1$ such that  
$p(\nu + z\xi) = z^d q(\nu, z),$ with $q: X_0 \times D \to \C$ a holomorphic function 
that is not identically zero on $X_0 \times \{0\}.$ Let $X_0'$ be the open dense subset of 
$\nu \in X_0$ such that $q(\nu, 0) \neq 0.$ Fix $\nu \in X_0',$ an open neighborhood
$X_1$ of $\nu$ whose closure is contained in $D_0'$ and a disk $D' \subset D$ centered at $0$ such that $q(\nu, z) \neq 0$ for all $\nu \in X_1$ and $z \in D'.$  
Then for every $\nu \in X_1$ the function 
$z \mapsto z^d j(\nu + z \xi)$ extends to a holomorphic function $D' \to V. $ It follows that 
 $c_k(\nu) = 0$ for $k < - d$ and $\nu \in X_1.$ By analytic continuation it now follows that $c_k = 0$ on $X_0$ for $k < - d.$ 
 
Let $m \in \Z$ be the maximal number such that $c_{-m} \neq 0.$ Arguing by contradiction
we will show that $m \leq 0,$ thereby completing the proof. Thus, suppose $m > 0.$ Then there
exists $\nu_1 \in X_0$ such that $c_{-m}(\nu_1) \in V\setminus \{ 0 \}.$ We claim 
that for $n \in N_0$ we have 
\begin{equation}
\label{e: transformation c min N}
\pi_{\nu_1}^{-\infty}(n) c_{-m}(\nu_1) = \cchi(n) c_{-m} (\nu_1). 
\end{equation} 
Indeed, fix $n \in N_0$ and $\psi \in C^\infty(K/K_P:\gs).$ Then
it suffices to prove the identity evaluated at $\psi.$ Writing $\pi_\nu = \pi_{\bar P,\gs, \nu},$ 
we start with the known identity expressing that $j(\nu)$ is a Whittaker vector in the induced representation:
$$ 
\inp {j(\nu) }{\pi_{-\bar\nu}^\infty (n^{-1})\psi} = \cchi(n) \inp{j(\nu)}{\psi},
$$ 
for $\nu \in \Omega\setminus X.$ Substituting $\nu = \nu_1 + z \xi$ 
for $z \in D\setminus \{0\},$ we obtain the identity 
\begin{equation}
\label{e: Whittaker behavior parametrized by z}
\inp{z^m j(\nu_1 + z \xi)}{\pi_{-\bar\nu_1 -\bar  z\bar \xi}^{\infty}(n^{-1})\psi} = 
\cchi(n)\inp{z^m j(\nu_1 + z \xi)}{\psi}
\end{equation}
of holomorphic functions on $D(0,r)\setminus \{0\}.$ We observe that we may write 
$$ 
z^m j(\nu_1 + z\xi) = c_{-m}(\nu_1) + z R(z),
$$ 
as an identity of holomorphic $V$-valued functions in $z \in D(0,r),$ with $R: D(0,r) \to V$ 
holomorphic. We may also write 
$$ 
\pi_{\nu_1 - z\bar \xi}^{\infty}(n^{-1})\psi =  \pi_{\nu_1}^\infty(n^{-1})\psi+ z \Psi(z), \qquad (z \in D(0,r)),
$$ 
with $\Psi$ a holomorphic function $D(0,R) \to C^\infty(K/K_P:\gsP).$ Substituting these 
expressions in (\ref{e: Whittaker behavior parametrized by z}) we obtain
\begin{equation}
\label{e: Whittaker and rest in z}
\inp{c_{-m}(\nu_1)}{\pi_{\nu_1}^{\infty}(n^{-1})\psi} =\cchi(n)  \inp{c_{-m}(\nu_1)}{\psi} + z F(z),
\end{equation}
where $F: D(0, r) \to \C$ is given by 
$$ 
F(z) =  \inp{R(z)}
{\cchi(n)\psi - \pi_{\nu_1}(n^{-1}) \psi - \bar z \Psi(\bar z)} -  \inp{c_{-m}(\nu_1)}{ \Psi(\bar z) }.
$$ 
If $z$ is restricted to a compact neighborhood of $0$ in $D(0,r),$ then $\Psi(z)$ stays in a bounded 
subset of $C^\infty(K/K_P:\gsP)$ and $R(z)$ stays in a bounded subset of $C^{-\infty}(K/K_P:\gsP).$ 
This implies that $F(z)$ remains bounded, so that $\lim_{z \to 0} zF(z) = 0.$ By taking  the limit 
of (\ref{e: Whittaker and rest in z}) we find that (\ref{e: transformation c min N}) is valid after pairing
both sides with $\psi \in C^\infty(K/K_P:\gsP).$ Since $n$ and $\psi$ were arbitrary, (\ref{e: transformation c min N}) follows.
Thus, in the induced picture we have 
$$ 
c_{-m}(\nu_1) \in C^{-\infty}(\bar P:\gs: \nu_1)_{\cchi}.
$$ 
Furthermore, by application of Lemma \ref{l: restricted j holomorphic}, the generalized function 
$c_{-m}(\nu_1)$ vanishes on the open orbit $N_0 \bar P.$ 
In view of Corollary \ref{c: cor of main thm HC} applied to $c_{-m}(\nu_1)$ in place of 
$j$ it finally follows that $c_{-m}(\nu_1) = 0,$ contradicting the condition involved in 
 the choice of $\nu_1.$ 
\qed

\begin{cor} 
\label{c: j bar P bijection}
For every $\nu \in \faPdc$ the map
$$ 
j(\bar P,\gs,\nu): \;\; H_{\gs,\cchi_P}^{-\infty} \to C^{-\infty}(G/\bar P: \gs :\nu)_\cchi
$$ 
is a linear isomorphism with inverse equal to the map $\ev_e$ defined in (\ref{e: defi ev e}).
\end{cor}

\proof 
It follows from the definition of $j(\bar P , \gs, \nu)(\eta)$ for $\nu \in \faPdc(P,0)$ that
the restriction of $j(\bar P, \gs, \nu)(\eta)$ to $K\cap N_P\bar P$ is equal to the continuous 
function $j_\eta(\nu): K\cap N_P \bar P \to H^{-\infty}_\gs $ defined in (\ref{e: defi j eta nu}). In the proof
of Lemma \ref{l: restricted j holomorphic} it is shown that $j_\eta(\nu)$ extends to a holomorphic function
of $\nu \in \faPdc$ with values in $C(K\cap N_P \bar P, H_\gs^{-\infty}).$ By analytic continuation it
follows that $j(\bar P , \gs, \nu, \eta)|_{K\cap N_P \bar P}$ is given by the function
$j_\eta(\nu).$ It now follows that $\ev_e j (\bar P, \gs, \nu ,\eta) = \eta.$ 
Hence, $\ev_e$ is a left inverse to $j(\bar P, \gs, \nu)$ and we see that $ev_e$ is surjective
$C^{-\infty}(G/\bar P: \gs :\nu)_\cchi \to H_{\gs,\cchi_P}^{-\infty}.$ If we combine this with the 
injectivity of $\ev_e,$ asserted in Corollary \ref{c: ev is injective}, the required result follows.
\qed 

\begin{lemma}
\label{l: det M mu product}
Let $M_\mu = M_{\mu}(\gs, \dotvar): \faPdc \to \End(H_{\gs,\cchi_P}^{-\infty})$ be the polynomial function
$\faPdc \to \End (H_{\gs,\cchi_P}^{-\infty})$ introduced in Proposition \ref{p: intro M mu nu}. 
Then the function $\det M_\mu$ is a non-zero constant times a finite product of first order polynomial 
functions on $\faPdc$  of the form $\nu \mapsto \inp{\nu}{\ga}  + c$ with $\ga \in \gS(\fn_P, \fa_P)$ 
and $c \in \C.$
\end{lemma}

\proof 
We observe that $\det M_\mu$ is a polynomial function $\faPdc \to \C.$ By Proposition \ref{p: factorization in first order} it suffices to show that $\det M_\mu$ is non-zero on the complement of a locally finite collection
$\cH$ of affine $\gS(\fn_P, \faP)$-hyperplanes. Increasing $\cH$ we may assume that 
$\cup \cH$ contains the zero set of the polynomial function $q: \faPdc \to 
\C$ introduced in Lemma \ref{l: intro Z mu}.
From Lemma \ref{l: generic injectivity of map with j}
we know that there exists such a collection of {\bfhyp} such for $\nu \in \faPdc\setminus \cup \cH$ 
the map 
\begin{equation}
\label{e: endo on j} 
\psi_\nu: j \mapsto \ge^{N_0} \after p_{\gL + \nu + \mu} (j \otimes e_{N_0} )
\end{equation}
is a bijective endomorphism of $C^{-\infty}(G/\bar P:\gs: \nu)_\cchi.$ 
By Corollary \ref{c: j bar P bijection} this implies the existence of a unique linear automorphism  
$b_\nu$ of $ H_{\gs,\cchi_P}^{-\infty}$ such that
$$ 
\psi_\mu(j(\bar P, \gs, \nu)) = j(\bar P, \gs,\nu) \after b_\nu.
$$ 
Since $q(\nu) \neq 0$ it follows from (\ref{e: Z q p}) 
that $p_{\gL + \mu +\nu} = q(\nu)^{-1} Z_{\bar P, \mu(\nu)}$ on $C^{-\infty}(G/\bar P: \gs :\nu)\otimes F.$ 
If we combine this with (\ref{e: endo on j}) and Proposition \ref{p: intro M mu nu} we find that 
$$ 
\psi_\mu(j(\bar P, \gs, \nu))  =  j(\bar P, \gs,\nu) \after [q(\nu)^{-1} M_\mu(\nu)],
$$ 
for $\gS(\fn_P,\fa_P)$-generic $\nu.$ By uniqueness, this implies $q(\nu)^{-1}M_\mu(\nu) = b_\nu $ 
for $\gS(\fn_P,\fa_P)$-generic $\nu.$ We conclude that $\det M_\mu(\nu) \neq 0$ for $\nu$ outside $\cup \cH.$ 
\qed

\begin{cor} There exists a polynomial function $p: \faPdc\to\C$ which is a finite product of linear factors of the form $\inp{\ga}{\dotvar} - c,$ with $\ga \in \gS(\fn_P, \faP)$ and $c \in \C$ 
such that $p R_\mu(\gs,\dotvar)$ is a polynomial $\End(H_{\gs,\cchi_P}^{-\infty})$-valued map.
\end{cor} 

\proof 
Let $p = \det M_\mu,$ then by Lemma \ref{l: det M mu product} the function $p$ has the required form.
By application of Cramer's rule the result now follows from the formula $R_\mu(\gs, \nu) = M_\mu(\nu)^{-1} m_\mu,$ 
given at the end of the proof of Theorem \ref{t: the functional equation} just before Proposition 
\ref{p: intro M mu nu}.
\qed

\begin{thm} 
\label{t: uniform estimate of j}
For every $R \leq 0$ there exist constants $C> 0, N \in \N$ and $r \in \N$ such that
for all $\eta \in H_{\gs,\cchi_P}^{-\infty}$ and $\nu \in \faPdc(P, R)$ we have  
$j(\bar P,\gs,\nu)\eta \in C^{-r}(K/K_P:\gsP)$ and 
$$ 
\|j(\bar P,\gs,\nu)\eta\|_{-r} \leq C (1 + |\nu|)^N \|\eta\|.
$$ 
\end{thm}

\proof 
We agree to write $j(\nu, \eta) = j(\bar P, \gs, \nu)\eta$. 
Following the induction in the proof of Theorem \ref{t: estimate by functional equation} one
sees that its assertion is valid with $p = p_R$ a polynomial function $\faPdc \to \C$ which is 
a finite product of linear factors of the form $l_{\ga, c}: \nu \mapsto \inp{\nu}{\ga} - c$ 
with $\ga \in \gS(\fn_P, \faP)$ and $c \in \C.$ From the mentioned theorem we know that there
exist constants $C' > 0$ and $s, N \in \N$ such that for all $\eta\in H_{\gs,\cchi_P}^{-\infty}$ 
and all $\nu \in \faPdc(P,R)$ we have that $\nu \mapsto p(\nu) j(\nu, \eta)$ 
is holomorphic on $\faPdc(P,0)$ with values in $C^{-s}(k/K_P: \gsP)$ and satisfies the estimate
\begin{equation}
\label{e: estimate p j}
\|p(\nu) j(\nu, \eta)\|_{-s} \leq  C' (1 +|\nu|)^{N} \|\eta\|.
\end{equation}
Let $l: \faPdc \to \C,$ $\nu \mapsto \inp{\nu}{\ga} - c$ be a linear polynomial dividing $p.$ Then 
 it suffices to prove the assertion and estimate of the above
type with $\bp p  = p/l$ in place of $p$ and $2^{N+1} C'$ in place of $C'.$ 

Put $\xi  =\ga/|\ga|$ and let $H = l^{-1}(0).$ Then $\gf: H\times \C \to \faPdc,$ $(\nu_0, z) \mapsto 
\nu_0 + z \xi$ is an affine isomorphism such that $l\after \gf(\nu_0, z) = z.$ 
Then clearly,  
$(\nu_0, z) \mapsto \bp p(\nu_0 + z\xi)  j(\nu_0 + z \xi, \eta)$ 
is holomorphic on $H \times \C\setminus\{0\}.$ Furthermore, for $\nu$ in the complement 
of $\gf(H \times D(0,\frac12))$ we have the estimate (\ref{e: estimate p j}) with $\bp p$ in place of $p$ and 
$2 C'$ in place of $C'.$ 

Let now $(\nu_0 , z) \in H \times D(0,\frac12).$ Then by the Cauchy integral formula  we have
\begin{eqnarray*}
\bp p (\nu_0 + z \xi) j(\nu_0 + z \xi , \eta) & = & \frac{1}{2\pi i} \int_{|w| = 1} 
\frac{\bp p (\nu_0 + w \xi) j(\nu_0 + w \xi , \eta)}{(w - z)}\; dw\\
&=& 
\frac{1}{2\pi i} \int_{|w| = 1} 
\frac{p (\nu_0 + w \xi) j(\nu_0 + w \xi , \eta)}{w(w - z)}\; dw.
\end{eqnarray*}
The formula holds a priori as an integral formula of $C^{-\infty}(K/K_P:\gsP)$-valued functions.
However, as the integrand has values in $C^{-s}(K/K_P:\gsP),$ it readily follows that
$\bp p j(\dotvar ,\eta)$ is $C^{-s}(K/K_P:\gsP)$-valued. Furthermore, by a straightforward  estimation
we obtain:
\begin{eqnarray*}
\|\bp p (\nu_0 + z \xi) j(\nu_0 + z \xi , \eta)\|_{-s}
& \leq & 
\sup_{|w| =1} \frac{C' (1 + |\nu_0 + w \xi|)^N}{|w|- |z|} \leq  2 C' (2 + |\nu_0|)^N \\
& \leq & 
2^{N+1} C' ( 1+ |\nu_0 + z\xi|)^N .
\end{eqnarray*}
\qed

The above estimates give rise to the following uniformly moderate estimates for matrix coefficients
of Whittaker vectors with smooth vectors. 

\begin{thm}
\label{t: uniformly moderate estimate}
Let $R \in \R.$ Then there exist constants $N \in \N, r, s > 0$ and a
continuous seminorm $n$ on $C^\infty(K/K_P:\gsP)$ such that 
for all $f \in C^\infty(K/K_P:\gsP),$  all $\eta \in H_\gs^{\infty, \cchi}$, $\nu \in \faPdc(P,R)$
and $x \in G$, 
$$ 
|\inp{f}{\pi_{\bar P, \gs, \nu}(x) j(\bar P, \gs, \nu,\eta)}| \leq  
(1 + |\nu|)^N e^{s|\Re \nu||H(x)|} e^{r|H(x)|}  n(f).
$$ 
\end{thm}

We prepare for the proof with a few lemmas.

\begin{lemma}
There exists a constant $s > 0$ such that for all $g \in C^\infty(K/K_P:\gsP),$  $\nu \in \faPdc$ and $a \in A$,
$$ 
\|\pi_{\bar P, \gs, \nu}(a^{-1}) g\|_0 \leq e^{s (|\Re \nu| + |\rho_P|)|\log a|}\, \|g\|_0.
$$ 
\end{lemma}

\proof 
The constant $s$ is built in to make the assertion independent of the choice of norms on $\fa$ and $\fad.$ 
Here we will need that $|\nu(H)|\leq s |\nu| |H|$ for $\nu \in \faPd$ and $H \in \fa.$ 
Define 
$g_\nu: G \to H_\gs$ by 
$$ 
g_\nu( k m a_P \bar n ) =  a_P^{- \nu + \rho_P} \gs(m^{-1}) g(k)
$$  
for $(k, m, a_P, \bar n) \in K \times M_P \times A_P \times \bar N_{P}.$ 
Then for $k \in K$ we have $[\pi_{\bar P, \gs, \nu}(a^{-1})g](k) = g_{\nu}(ak)$ from which it follows that
\begin{eqnarray*}
\|[\pi_{\bar P, \gs, \nu}(a^{-1})g](k)\|_\gs  & \leq & e^{(- \Re \nu +\rho_P)(H_{\bar P}(ak))}
 \|g(\kappa_{\bar P}(a k))\|_\gs \\
&\leq& e^{s |\Re \nu - \rho_P||H_{\bar P}(ak)|} \|g\|_0.
\end{eqnarray*}
Since $|H_{\bar P}(ak)| \leq |\log a|$ the required estimate follows.
\qed

\begin{lemma}
\label{l: estimate u a f}
Let $u \in U(\fg).$ Then there exist constants $N \in \N$, $r > 0$ and $t \in \N$ such that 
for all $f \in C^\infty(K/K_P:\gsP),$ $\nu \in \faPdc$ and $a \in A,$
$$ 
\|\pi_{\bar P, \gs, \nu}(u)\pi_{\bar P, \gs, \nu}(a^{-1}) f\|_0  \leq ( 1+ |\nu|)^N e^{s|\Re \nu||\log a|}e^{r|\log a|} \|f\|_{t}.
$$ 
\end{lemma}
\medbreak

\proof 
In view of the PBW theorem we may assume that $\fa$ acts via the adjoint action by a weight $\xi$ on $u.$ 
Accordingly we have, for $a \in A,$ 
\begin{eqnarray*}
\|\pi_{\bar P, \gs, \nu}(u)\pi_{\bar P, \gs, \nu}(a^{-1}) f\|_0 
& \leq & a^{\xi} |\pi_{\bar P, \gs, \nu}(a^{-1}) \pi_{\bar P, \gs, \nu}(u) f\|_0\\
& \leq & e^{s(|\Re \nu| + |\rho_P|)|\log a|} e^{\xi(\log a)} \|\pi_{\bar P, \gs, \nu}(u) f\|_0
\end{eqnarray*}
by application of the previous lemma. The proof is finished by application of Lemma
\ref{l: estimate infinitesimal action}.
\qed
\medno
{\em Completion of the proof of Theorem \ref{t: uniformly moderate estimate}.\ }
It is clear that it suffices to prove the estimate  for $x =a \in A.$
By Theorem \ref{t: uniform estimate of j} there exists a constant $N > 0$ and a finite collection
$F \subset U(\fg)$ such that for all $f \in C^\infty(K/K_P:\gsP)$, all $\nu \in \faPdc(P,R)$ 
and all $a \in A,$ 
$$ 
|\inp{f}{\pi_{\bar P, \gs, \nu}(a) j(\bar P, \gs, \nu,\eta)}| \leq  
(1 + |\nu|)^N \|\eta\| \max_{u \in F} \| \pi_{\bar P, \gs, \nu}(u) \pi_{\bar P, \gs, \nu}(a^{-1})f \|_0.
$$ 
The proof is now readily completed by application of Lemma \ref{l: estimate u a f}.
\qed

\section{Uniformly tempered estimates}
\label{s: uni temp estimates}
The purpose of this section is to obtain uniformly tempered estimates for holomorphic families 
of Whittaker functions satisfying requirements of moderate growth.
 
Let $P$ be a standard parabolic subgroup of $G$ and  $(\gs, H_\gs)$ a representation
of the discrete series of $M_P.$ For $\ge > 0$ we put 
$$
\faPdc(\ge) = \{\nu \in \faPdc\mid |\Re (\nu)| < \ge\}.
$$
\begin{defi}
\label{d: whittaker family}
By a holomorphic family of Whittaker maps associated
with $(P,\gs)$ and $\ge_0 > 0$ we mean a family of maps
\begin{equation}
\label{e: intro w nu}
\wh_\nu: C^\infty(K/K_P:\gs) \to C^\infty(G/N_0\col \cchi) , \qquad (\nu \in \faPdc(\ge_0)),
\end{equation}
given by the matrix coefficient formula
\begin{equation}
\label{e: w nu as matrix coeff}
\wh_\nu(f)(x) = \inp{\pi_{\bar P,\gs, -\nu}(x)^{-1} f}{j_{\bar \nu}}
\end{equation}  
with $ j_{\bar \nu} \in C^{-\infty}(\bar P: \gs: \bar \nu)_\cchi$, $(\nu \in \faPdc(\ge_0)),$ 
such that $\nu \mapsto j_{\nu}$ is holomorphic as a $C^{-\infty}(K/K_P:\gs)$-valued 
function.
\end{defi}

\begin{rem}
\label{r: properties Whit map}
Let $(\wh_\nu)_{\nu \in \faPdc(\ge_0)}$ be a holomorphic family of Whittaker maps as above. 
We note that the following assertions
are valid.
\begin{enumerate}
\itema
The map $(\nu, f) \mapsto \wh_\nu(f)$,  $\faPdc(\ge_0) \times C^\infty(K/K_P:\gs) \to C^\infty(G/N_0\col \cchi)$ is continuous and holomorphic in the variable $\nu.$ 
\itemb
For each $\nu \in \faPdc(\ge_0)$ the map (\ref{e: intro w nu}) intertwines the generalized
principal series $\pi_{\bar P, \gs, -\nu}$
with the left regular representation $L.$ 
\end{enumerate}
\end{rem}

\begin{rem}
\label{r: example Whittaker family}
It follows from  Theorem \ref{t: j is holomorphic}
that for any $\ge_0 > 0$ and all $\xxi \in H_{\gs,\cchi_P}^{-\infty}$ 
the family $(\wh_\nu)_{\nu \in \faPdc(\ge_0)}$ of maps $H_\gs^\infty \to C^\infty(G/N_0\col  \cchi)$ defined by
$$ 
\wh_{\nu}(f)(x) = \inp{\pi_{\bar P, \gs, -\nu}(x)^{-1} f}{j(\bar P, \gs, \bar \nu, \xxi)}
$$ 
is a holomorphic family of Whittaker maps associated with $(P,\gs)$. 
Furthermore, 
by Theorem \ref{t: uniformly moderate estimate} it satisfies
the condition of uniform moderate growth mentioned below.
\end{rem}

To keep notation manageable, we will write
\begin{equation}
\label{e: notation IPgs}
\IPgs:= C^\infty(K/K_P:\gsP).
\end{equation}
Furthermore, we will use
the notation 
$$ 
|(\nu, a)|: =  (1 + |\nu|)(1+ |\log a|), 
$$%
for $a \in A$ and $\nu \in \faPdc.$ 

\begin{defi}
A Whittaker family $(\wh_\nu)$ as in Definition \ref{d: whittaker family} 
is said to have {\em uniform moderate growth} if there exist constants 
$r, s, N > 0$ and  
a continuous seminorm $\snorm$ on $\IPgs$  such that 
\begin{equation}
\label{e: estimate of uniform moderate growth}
|\wh_\nu(f)(a)| \leq  |(\nu,a)|^N  e^{s |\Re \nu| |\log a| } e^{r|\log a|} \snorm(f)
\end{equation} 
for all $f \in \IPgs$, $\nu \in \faPdc(\ge_0)$ and $a \in A.$ 
\end{defi}
If we combine the estimate (\ref{e: estimate of uniform moderate growth})
with Lemma \ref{l: basic comparison of estimates}, then
we see that for any linear functional $\xxi \in \fad$ with $\xxi \geq r|\dotvar|$ 
on $\fa^+$
we may adapt the continuous seminorm $\snorm$ so that 
for all 
$f\in \IPgs$ we have
\begin{equation}
\naam{e: domination whittaker family}
|\wh_\nu(f)(a)| \leq  |(\nu,a)|^N e^{s |\Re \nu| |\log a| } a^\xxi \snorm(f)
\end{equation}
for all $a \in A$ and $\nu \in \faPdc(\ge_0).$ 
We fix such a choice of $\xxi$ and observe that in particular $\xxi \geq -\rho $ 
on $\fa^+.$ 

The above estimate can be improved to a much sharper estimate 
of {\em uniform tempered growth} for $\ge_0 > 0$ taken sufficiently small.
More precisely, we have the following result.

\begin{thm}
\label{t: uniformly tempered estimates}
Let $G = {}^\circ G,$ and let $(\wh_\nu)_{\nu \in \faPdc(\ge_0)}$ be a holomorphic family of 
Whittaker maps as in Definition \ref{d: whittaker family}. Assume the family satisfies the condition 
of uniform moderate growth formulated in (\ref{e: estimate of uniform moderate growth}).
Then for $\ge > 0$ sufficiently small there exist constants $s > 0, N >0$ and a continuous seminorm $\snorm$ on $\IPgs$ such that for all $f \in \IPgs$, all $\nu \in \faPdc(\ge)$ and all $a \in A,$ 
\begin{equation}
\label{e: first uniform temp estimate}
|\wh_\nu(f)(a)| \leq  |(\nu,a)|^N e^{s |\Re \nu| |\log a| } a^{-\rho} \snorm(f).
\end{equation} 
\end{thm}

Before turning to its proof,  we first formulate a useful  consequence of the above result.

\begin{cor}
\label{c: uniformly tempered estimates} 
Let hypotheses be as in Theorem \ref{t: uniformly tempered estimates} and let $\ge >0$
be such that the conclusions of the theorem are valid.

Let $u \in U(\fg).$ Then there exist $s > 0, N >0$ and a continuous seminorm 
$\snorm$ on $\IPgs$ such that 
for all $f \in \IPgs$~, all  $\nu \in \faPdc(\ge)$ and all $x \in G$,
\begin{equation}
\label{e: third uniform temp estimate}
|L_u[\wh_\nu(f)](x)| \leq  (1 + |\nu|)^N (1 + | H(x)|)^N e^{s |\Re \nu| |H(x)| } e^{-\rho H(x)} \snorm(f).
\end{equation}
\end{cor}

\proof 
In view of Remark \ref{r: properties Whit map} (b) 
and since 
$\pi_{\bar P, \gs, -\nu}^{\infty}(u)$ acts continuously on $\IPgs,$ 
with polynomial dependence on $\nu,$  the estimate (\ref{e: first uniform temp estimate}) gives rise 
to an estimate
\begin{equation}
\label{e: second uniform temp estimate}
|L_u[\wh_\nu(f)](a)| \leq  |(\nu,a)|^N e^{s |\Re \nu| |\log a| } a^{-\rho} \snorm(f),
\end{equation}
for all $f \in \IPgs$,  all $\nu \in \faPdc(\ge)$ and all $a \in A,$   provided $N$ and $\snorm$ are suitably enlarged.

For a given finite subset $S\subset U(\fg)$ we may enlarge $N$ and $n$ further to arrange that
the estimate (\ref{e: second uniform temp estimate}) is valid for all $u \in S$ and all 
$f,\nu, a$ as before.

For  $k\in K, a\in A$ and $n \in N_0$ we have
$$ 
|L_u(\wh_\nu(f))(kan)| = |L_{k^{-1}} L_u (\wh_\nu(f)) (a)| = | L_{\Ad(k^{-1})u } \wh_\nu(L_{k^{-1}}f )(a)|.
$$ 
We may write $\Ad(k^{-1}) u = \sum_i c_i(k) u_i$ with $S=\{u_i\}$ a finite subset of $U(\fg)$ 
and such that $c_i: K \to \C$ are functions with sup-norm bounded by $C>0.$ 

We may enlarge $N$ and $n$ so that the estimate (\ref{e: second uniform temp estimate}) is valid with 
$u$ replaced by any element of $S$ and for all $f,\nu, a.$ 
There exist a continuous seminorm $\snorm'$ on $\IPgs$ such
that $\snorm  \after L_k \leq \snorm'$ for all $k \in K.$ From the last estimate mentioned above we 
now readily infer that for all $f \in \IPgs,$ all $\nu \in \faPdc(\ge)$ and all $(k,a,n) \in K\times A\times N_0,$
$$ 
|L_u(\wh_\nu(f))(kan)| \leq |(\nu, a)|^{N} e^{s |\Re \nu||\log a|} e^{- \rho \log a} \;C|S|\snorm'(f).
$$  
Enlarging $N$ and $n$ once more, we obtain the required estimate (\ref{e: third uniform temp estimate}).
\qed

In the proof of Theorem \ref{t: uniformly tempered estimates} the following terminology will
be useful.

\begin{defi}\label{d: domination}
We will say that a functional $\xxi \in \fad$ dominates the
given Whittaker family $ (\wh_\nu)$ if there exist $\ge >0, s> 0, N>0$ and 
$\snorm$ as above such that the estimate (\ref{e: domination whittaker family}) is valid
for all $f \in \IPgs$,  $\nu \in \faPdc(\ge)$ and $a \in A.$
\end{defi}

{\em Proof of Theorem \ref{t: uniformly tempered estimates}. \ }
Clearly it is sufficient 
to prove that $-\rho$ dominates the Whittaker family $(\wh_\nu).$  
We will achieve this by improving $\xxi$ in a finite number
of steps, each step corresponding to a simple root $\ga \in \gD$, by using the asymptotic behavior of $\wh_\nu(f)$ along standard maximal parabolic subgroups.

Since ${}^\circ G = G,$ the collection $\gD$ of simple roots in $\gS^+$ is a basis of $\fad.$ 
Let $(h_\ga)_{\ga \in \gD}$ be the dual basis in $\fa.$ 
We will now establish an improvement step for each $\ga \in \gD.$ Put 
$\Phi = \gD\setminus \{\ga\}.$ Then $P_\Phi$ is a maximal parabolic subgroup with
 split component $A_\Phi = \exp(\R h_\ga).$

 We define the partial ordering $\preceq $ on $\fad$ by 
 \begin{equation}
 \label{e: defi majorisation} 
 \gl \preceq \mu \iff \;\gl(H) \leq \mu(H) \;\;\;\mbox{\rm for all}\;\; H \in \fa^+.
 \end{equation}
 The condition on the right is equivalent to 
 $\gl(h_\ga) \leq \mu(h_\ga)$ for all $\ga \in \gD.$ 
  
 Given $\xxi \in \fad$ we write $i_\ga(\xxi)$ for the element in $\xxi + \R \ga$ 
satisfying  $i_\ga(\xxi)(h_\ga)=  -\rho(h_\ga).$  Equivalently, $i_\ga(\xxi) \in \fad$ is determined by
$$
i_\ga(\xxi)(h_\gb) = \left\{ 
\begin{array}{lll}
\xxi(h_\gb) & {\rm for} & \gb \in \gD\setminus \{\ga\};\\
-\rho(h_\ga) & {\rm for }&  \gb = \ga.
\end{array}
\right.
$$ 

\begin{lemma}
\label{l: i ga eta majorizes min rho}
If $\xxi \succeq -\rho$, then for every simple root $\ga \in \gD$ it holds that 
$i_\ga(\xxi) \succeq -\rho.$ 
\end{lemma}

\proof 
This is straightforward.
\qed

\begin{lemma}
\label{l: improvement step}
{\bf (Improvement step)\ } 
Suppose that $\xxi \in \fad$ dominates the Whittaker family $(\wh_\nu)$ 
and satisfies  $\xxi \succeq -\rho.$ Let $\ga \in \gD.$ 
\begin{enumerate}
\itema If $\xxi(h_\ga) - 1 \geq - \rho(h_\ga)$ then for every $c \in [0,1)$, 
the functional $\xxi':= \xxi - c\ga$ dominates $(\wh_\nu)$ and satisfies   $\xxi' \succeq -\rho.$ 
\itemb If $\xxi(h_\ga) -1  <  - \rho(h_\ga)$, then $i_\ga(\xxi)$ dominates $(\wh_\nu)$ 
and satisfies $i_\ga(\xxi) \succeq - \rho.$ 
\end{enumerate} 
\end{lemma}
\medno
The rest of this section will be dedicated to establishing this lemma. Before turning to the proof of 
the lemma we will show
how Theorem  \ref{t: uniformly tempered estimates} can be deduced from it.
\medbreak

{\em Completion of the proof of Theorem  \ref{t: uniformly tempered estimates}.\ }
Let $\ga \in \gD$ and assume that $\xxi\in \fad$ dominates $(\wh_\nu)$ and satisfies 
$\xxi \succeq -\rho.$ Then $\xxi(h_\ga) \geq  - \rho(h_\ga).$ Let $k$ be the smallest natural number 
such that $\xxi(h_\ga) - k < -\rho(h_\ga).$ Then there exists a $c \in [0,1)$ such that 
$\xxi(h_\ga) - kc =  -\rho(h_\ga).$ By applying (a) of the above lemma $k$-times
successively, we find that $\xxi'': = \xxi - kc\ga $ dominates $(\wh_\nu),$ while $\xxi'' \succeq -\rho.$  Since $\xxi''(h_\ga) - 1 < - \rho(h_\ga),$ we may apply (b) of the above lemma
to conclude that $i_\ga(\xxi'')$ dominates $(\wh_\nu)$ and satisfies $i_\ga(\xxi'') \succeq -\rho$ 
Since $i_\ga(\xxi'') = i_\gs(\xxi),$ we conclude that $i_\ga(\xxi)$ 
dominates $(\wh_\nu)$ and satisfies $i_\ga(\xxi) \succeq -\rho.$

Let now $\ga_1, \ldots, \ga_r$ be a numbering of the simple roots from $\gD.$ 
Then by the above reasoning it follows that 
$\xxi''':= i_{\ga_r} \after \cdots \after i_{\ga_1} (\xxi)$ dominates $(\wh_\nu)$ while 
$\xxi''' \succeq -\rho.$ Since $(h_\ga)_{\ga \in \gD}$ is the basis of $\fa$ dual to $\gD,$ 
it is readily verified that  $\xxi''' = -\rho.$ 
\qed

{\em Start of proof of Lemma \ref{l: improvement step}.\ }
We assume that  $(\wh_\nu)_{\nu \in \faPdc(\ge)}$ is a Whittaker family associated
with $(P,\gs)$ which is  dominated by $\xxi \in \fad. $ Moreover we assume that $\xxi \succeq -\rho.$ We fix a simple root $\ga \in \gD$ and put $\Phi:= \gD\setminus \{\ga\}.$ 
Our goal is to establish the two assertions (a) and (b) of Lemma \ref{l: improvement step}.
For this we will need a proper exploitation 
of the differential equations satisfied by the given Whittaker family. 

For $X \in \fg_\iC$ 
we denote by $\bar X$ the complex conjugate of $X$ relative to the real form $\fg.$ 
Let $U(\fg)$ denote the universal enveloping algebra of $\fg_{\iC}.$ 
The map $X \mapsto \bar X$ has a unique extension to a conjugate linear algebra 
isomorphism  $U(\fg) \to U(\fg),$ which is denoted by $u\mapsto \bar u.$ 
In particular, this means that $\overline{uv} = \bar u\bar v$ for $u,v \in U(\fg).$ 

For $X \in \fg$ we have 
\begin{equation}
\label{e: R X on wh nu} 
R_X [\wh_\nu (f)](x) = \inp{\pi(x^{-1})f}{\pi_{\bar P, \gs, \bar\nu}( \bar X) j_{\bar \nu}},
\end{equation}
for $f \in \IPgs,$ $\nu \in \faPdc(\ge)$ and $x \in G.$ 
By complex linearity (\ref{e: R X on wh nu}) is valid for all $X \in \fg_\iC$, which leads to the similar 
formula with $X$ replaced by a general element of  $U(\fg).$ 
Let $\ft$ be a maximal 
torus in $\fm_1;$ then 
$$ 
\fh := \ft \oplus \fa
$$ 
is a Cartan subalgebra of $\fg.$ 
We put $\fh^*_\iR:= i \ft^* \oplus \fa^*.$ Then $\fh^*_\iR$ is the real 
span of the roots of $\fh_\iC$ in $\fg_\iC.$  Let $\gg: \fZ \to P(\fh)^{W(\fh)}$ 
be the Harish-Chandra isomorphism for  $(G,\fh)$ and let 
$\gg_{M_{1P}}: \fZ(M_{1P}) \to S(\fh)^{W_P(\fh)}$ be the similar isomorphism
for $(M_{1P}, \fh).$ 
Let $\gL_\gs \in \fh_\iC^*  \cap \fm_{P\iC}^*$ be an infinitesimal character for 
the representation $\gs$ of the discrete series of $M_P,$ by which we mean that
$\gg_{1P}(\dotvar, \gL_\gs)$ is the character of $\fZ(M_P)$ by which it acts
on $H_\gs^\infty.$ We note that $\gL_\gs $ belongs to the real span of the
roots of $(\fm_{1P}, \fh),$ hence to $i\ft^* \oplus {}^*\fa_P^*.$ 
Applying formula (\ref{e: R X on wh nu}) with $Z \in \fZ$ in place of $X,$ we find
that 
\begin{equation}
\label{e: whittaker system of differential operators}
R_Z [\wh_\nu(f)] = \overline{\gg(\bar Z, \gL_\gs + \bar \nu)}\, \wh_\nu(f) =  \gg(Z, -\gL_\gs + \nu)\, \wh_\nu(f),
\end{equation} 
for all $Z \in \fZ.$ 
Following an idea similar to the one in Section \ref{s: sharp estimates of whit coeff}, but with dependence on parameters, we will exploit this system to establish the improvement 
step of Lemma \ref{l: improvement step}. 

The standard parabolic subgroup $P_\Phi$ has split component 
$A_\Phi:=  \exp \fa_\Phi = \exp \R h_\ga. $ We agree to write $\bp \fa_\Phi$ 
for the real linear span of the elements $h_\gb$ with $\gb\in \Phi$ and, accordingly,
$\bp A_\Phi= \exp(\bp \fa_\Phi).$  We note that $\fa = \bp \fa_\Phi  \oplus \fa_\Phi$ 
and $A = \bp A_\Phi  A_\Phi,$ see also (\ref{e: bp deco of fa}).

We denote by $\fZ_{1\Phi} = \fZ_{\fm_{1\Phi}}$ the center of $U(\fm_{1\Phi}).$
In view of the PBW theorem we have $U(\fg) = U(\fm_{1\Phi}) \oplus (\bar \fn_\Phi U(\fg) + U(\fg) \fn_\Phi).$ 
The associated projection $U(\fg) \to U(\fm_{1\Phi}),$ restricted to $\fZ,$ defines an 
algebra homomorphism
$$ 
p: \fZ \to \fZ_{1\Phi}.
$$ 
It is well known that $p$ is injective and that $\fZ_{1\Phi}$ is a free $p(\fZ)$-module of finite rank. Let $u_1, \ldots, u_\ell$ be a free basis of this module, and let $E_\Phi$ the complex linear span of this basis.
Then it follows that 
$$ 
\fZ_{1\Phi} = E_\Phi p(\fZ).
$$ 
Moreover, the map $(u, Z)\mapsto  u p(Z)$ induces a linear isomorphism
$
E_\Phi\otimes \fZ \simeq \fZ_{1\Phi}.
$
In the formulation of the following lemma, $\fh$ is a $\Cartan$-stable Cartan subalgebra of 
$\fg$ containing $\fa.$ Thus, $\fh = \ft \oplus \fa$ with $\ft$ a maximal torus in $\fm.$ 

\begin{lemma}
\label{l: weights in quotient center}
If $I$ is a cofinite ideal of $\fZ,$ then $I_\Phi:= E_\Phi p(I)$ is a cofinite ideal in $\fZ_{1\Phi}.$ 
Furthermore, if $\gl_\Phi \in \fhdc$ is an infinitesimal character for $\fZ_{1\Phi}$ appearing in the 
quotient module $\fZ_{1\Phi}/I_\Phi,$ then there exists an infinitesimal character $\gl \in \fhdc$ 
for the $\fZ$-module $\fZ/I$ such that 
\begin{equation}
\gl_\Phi \in W(\fh)\gl - \rho_\Phi.
\end{equation}
\end{lemma}

\proof 
Let $u \in E_\Phi$ and $W \in \fZ_{1\Phi}$ then $W u = \sum_{j} u_j p(Z_j)$ 
with $Z_j \in \fZ.$ It follows that 
$W up(I) \subset \sum_{j} u_j p(Z_j) p(I) = \sum_{j} u_j p(Z_j I) \subset E_\Phi p(I).$ This shows
that $E_\Phi p(I)$ is an ideal. In view of the linear isomorphism $E_\Phi \otimes \fZ \to \fZ_{1\Phi} $ 
we have $\fZ_{1\Phi}/ E_\Phi p(I) = E_\Phi p(Z)/p(I)$ as complex vector spaces.
Since $p(Z)/p(I)$ is finite dimensional, the cofiniteness of $I_\Phi$ follows.

Let $\xi \in \widehat \fZ_{1\Phi}$ be a character appearing in $\fZ_{1\Phi}/I_{\Phi}.$ 
Then there exists a $k \in \N$ and an element $v\in \fZ_{1\Phi}\setminus I_\Phi$ 
such that $(W - \xi(W))^{k}v \in I_\Phi$ for all $W \in \fZ_{1\Phi}.$ In particular
the latter is valid for $W = p(Z),$ $Z \in \fZ.$ Decompose $v = \sum u_i p(Z_i)$ 
with $Z_i \in \fZ.$ Then it follows that $[p(Z) - \xi(p(Z))]^k p(Z_i) \in p(I)$ for all $Z \in \fZ$
and all $1\leq i \leq \ell.$ By injectivity of $p$ this implies that 
\begin{equation}
\label{e: ann in Z mod I}
[Z - \xi(p(Z))]^k Z_i \in I,
\end{equation}
for all $Z \in \fZ$ and all $i.$  
On the other hand,  $v \notin I_{\Phi}$ implies that $Z_i \notin I$ for at least one $i.$
Combining this with (\ref{e: ann in Z mod I})  we infer that  the character $\xi\after p \in \widehat \fZ$ 
appears in $\fZ / I.$ 

Let $\gg: \fZ \to P(\fhdc)^{W(\fh)}$ and $\gg_{1\Phi}: \fZ_{1\phi} \to P(\fhdc)^{W_{1\Phi}(\fh)}$ 
denote the canonical isomorphisms. Then it is well known that 
$$ 
\gg_{1\Phi} \after p = T_{\rho_{\Phi}} \after \gg,
$$ 
where $T_{\rho_\Phi}\in {\rm Aut}(P(\fhdc))$ is the translation $p\mapsto p(\dotvar + \rho_\Phi).$ 

Let $\gl_\Phi$ be as stated. Then $\xi = \gg_{1\Phi} (\dotvar, \gl_{\Phi})$ 
is a character of $\fZ_{1\Phi}$ which appears in $\fZ_{1\Phi} /I_\Phi.$ 
It follows that $\xi \after p$ is a character of $\fZ$ appearing in $\fZ / I$ 
hence of the form $\gg(\dotvar, \gl),$ with $\gl \in \fhdc.$ 
We now conclude that for all $Z \in \fZ$ we have
$$ 
\gg(Z, \gl_\Phi + \rho_\Phi) = \gg(Z, \gl).
$$ 
This in turn implies that $\gl_\Phi + \rho_{\Phi} \in W(\fh) \gl.$
\qed

If $I \ideal \fZ$ is cofinite, we denote by $\spec(\fZ/I)$ the (finite) collection
of infinitesimal characters appearing in $\fZ/I.$ Since $U(\fa_\Phi)$ is a submodule
of $\fZ_{1\Phi}$, the $\fZ_{1\Phi}$-module $\fZ_{1\Phi}/ E_\Phi p(I)$ 
is a $U(\fa_\Phi)$-module as well. 

\begin{cor}
Let $I \ideal \fZ$ be cofinite. Then the set of $\fa_\Phi$-weights in 
$\fZ_{1\Phi}/ E_\Phi p(I) $ is contained in 
$$ 
\cup_{\gl \in \spec (\fZ/I)} W(\fh) \gl |_{\fa_\Phi} - \rho_\Phi.
$$ 
\end{cor}
\proof 
Apply Lemma \ref{l: weights in quotient center}.
\qed

We now specialize to the ideal $I_{\nu} = \ker \gg(\dotvar, - \gL_\gs + \nu)$ for 
$\nu \in \faPdc.$ Then 
$$
\spec(\fZ / I_\nu) = W(\fh)(- \gL_\gs + \nu). $$

\begin{cor} The inclusion map induces a linear isomorphism 
\begin{equation}
\label{e: iso E Phi to quotient}
E_\Phi \;{\buildrel  \simeq \over \longrightarrow}\;\;\fZ_{1\Phi}/\fZ_{1\Phi} p(I_{\nu}).
\end{equation}
The set $\wt(\nu)$ of (generalized) $\fa_\Phi$-weights in the displayed quotient equals 
$$ 
\wt(\nu) = W(\fh)(-\gL_\gs + \nu)|_{\fa_\Phi }- \rho_\Phi.
$$ 
\end{cor}

We write $E_{\Phi,\nu}$ for the space $E_\Phi$ equipped with the $\fa_\Phi$ action 
for which the map (\ref{e: iso E Phi to quotient}) becomes an isomorphism 
of $\fa_\Phi$ modules. We agree to write $B_1(\nu)$ for the linear map  
by which $h_\ga$ acts on $E_{\Phi, \nu}.$ 
There exist unique $Z^k_j \in \fZ$ such that 
\begin{equation}
\label{e: Z k j}
h_\ga u_j = \sum_{k = 1}^\ell u_k p(Z^k_j), \qquad (1 \leq j \leq \ell).
\end{equation}
Since $Z^k_j - \gg(Z^k_j, -\gL_\gs + \nu)\in I_{\nu}$ it follows that, for all $\nu \in \faPdc,$ 
$$ 
B_1(\nu) u_j  \in  \sum_k \gg(Z^k_j, - \gL_\gs + \nu)u_k + \fZ_{1\Phi} p(I_\nu).
$$ 
Therefore, the matrix of $B_1(\nu)$ relative to the basis $u_1, \ldots, u_\ell$ of $E_\Phi$ 
is given by 
\begin{equation}
\label{e: B and Z k j}
B_1(\nu )_j^k = \gg(Z^k_j, -\gL_\gs + \nu).
\end{equation} 
In particular, it follows that $\nu \mapsto B_1(\nu)$ is polynomial $\faPdc \to \End(E_\Phi).$ 

Let $I_\cchi$ denote the left ideal in $U(\fg)$ generated by the elements $Y - \cchi_*(Y)$ 
for $Y \in \fn_0.$ Let $I_{\gs, \nu}$ denote the left ideal in $U(\fg)$ generated by the 
ideal $I_\nu$. 

\begin{lemma}
There exist  elements $v_j \in \bar \fn_\Phi U(\bar \fn_0 + \fm_{1})$ such that
$$ 
h_\ga u_j-  B_1(\nu) u_j  - v_j \in I_\cchi + I_{\gs, \nu},\qquad (1 \leq j \leq \ell).
$$ 
\end{lemma}

\proof
Let $Z_j^k \in \fZ$ be as in (\ref{e: Z k j}). 
Now $p(Z_j^k) - Z_j^k \in \bar \fn_\Phi U(\fg),$ 
and since $u_j \in U(\fm_{1\Phi})$ it follows that also 
$$
u_k p(Z_j^k) - u_k Z_j^k \in \bar \fn_\Phi U(\fg) = \bar \fn_\Phi U(\bar \fn_0 + \fm_{1})U(\fn_0).
$$ 
We note that any $W \in U(\fn_0)$ equals $\cchi_*(W)$ modulo $I_\cchi,$ hence 
$$ 
u_k p(Z_j^k) - u_k Z_j^k  \in  v_j^k + I_\cchi,
$$ 
with $v_j^k \in \bar \fn_\Phi U(\bar \fn_0 + \fm_{1}).$ It follows that for all 
$\nu \in \faPdc,$ 
$$ 
u_k p(Z_j^k) - u_k \gg(Z_j^k, - \gL_\gs + \nu) - v_j^ k \in I_\cchi + I_{\gs, \nu}.
$$ 
Summing the above over $k$, putting $v_j = \sum_k v^k_j$ and using 
(\ref{e: B and Z k j}) we obtain the desired assertion.
\qed

We proceed with the proof of the improvement step of Lemma \ref{l: improvement step}.
For $f \in \IPgs$ and $\nu \in \faPdc(\ge)$ we define
the function $F(f, \nu) : A \to \C^\ell$ by 
$$ 
F(f, \nu, a)_j := \inp{ \pi_{-\nu}(a)^{-1} f}{ \pi_{\bar\nu}(u_j) j_{\bar\nu}},
$$
where we briefly wrote $\pi_\nu$ for $\pi_{\bar P,\gs, \nu}.$  We define the function
$R(f, \nu): A \to \End(\C^\ell)$ by 
$$ 
R(f, \nu, a)_j := \inp{\pi_{-\nu}(a)^{-1}f}{\pi_{\bar \nu}(v_j)j_{\bar \nu}}.
$$ 
Furthermore, let $B(\nu) \in \End(\C^\ell)$ be the endomorphism with matrix equal to 
$\overline{B_1(\bar \nu)}^{\rm T}.$ 

For the following result we recall the decomposition $A = \bp A_\Phi A_\Phi.$ It allows us to
decompose any element $a \in A$ in a unique way as 
$$ 
a = \bp a a_t , 
$$ 
with $\bp a \in \bp A_\Phi$ and with $a_t = \exp t h_\ga$, $(t \in \R).$

\begin{lemma}
\label{l: diff eq for F}
The function $F$ introduced above satisfies the equation
$$ 
\frac{d}{dt} F(f,\nu, \bp a a_t) = B(\nu) F(f, \nu, \bp aa_t) + R(f, \nu, \bp a a_t),
$$ 
for every $f \in \IPgs$ and all $\nu \in \faPdc(\ge_0),$ $\bp a \in A_\Phi$ and  $t \in \R.$ 

Here $B$ is a polynomial function $\faPdc \to \End(\C^\ell).$ For every $\nu \in \faPdc$ the spectrum 
of $B(\nu)$ satsfies 
$$ 
\spec (B(\nu)) \subset  \{w(-\gL_\gs + \nu)(h_\ga) - \rho(h_\ga)\mid w \in W(\fh)\},\qquad (\nu \in \faPdc).
$$ 
\end{lemma}

\proof
Noting that $I_{\gs, \bar \nu}$ and $I_\cchi$ vanish on $j_{\bar \nu},$ 
we obtain that the functions $F_j$ introduced above satisfy the equations
\begin{eqnarray*} 
\frac{d}{dt} F_j (f, \nu, \bp aa_t) & = & \inp{\pi_{-\nu} (a)^{-1} f}{\pi_{\bar \nu}(h_\ga u_j) j _{\bar \nu}}
\\
&=&  \inp{\pi_{-\nu} (a)^{-1} f}{ \pi_{\bar \nu}(B_1(\bar \nu) u_j + v_j) j_{\bar \nu}} ,
\end{eqnarray*}
for all $f \in \IPgs,\; \nu \in \faPdc(\ge_0), \; \bp a\in \bp A_\Phi$ and $t \in \R.$ 
This gives the required equation, with $B(\nu)$ as asserted. The spectrum 
of $B(\nu)$ equals that of $\overline{B_1(\bar \nu)}$ hence consists of the 
elements 
\begin{equation}
\label{e: conj els of spec}
\overline{w(-\gL_\gs + \bar \nu)(h_\ga)} - \rho_\Phi(h_\ga),
\end{equation} 
for $w \in W(\fh).$ 
Since $\gL_\gs \in \fh_\R^* = i \ft^* \oplus \fa^*$ whereas $W(\fh)$ leaves $\fh_\R$ invariant,
it follows that for each $w \in W(\fh)$ the value $w(-\gL_\gs)(h_\ga)= -\gL_\gs(w^{-1}h_\ga)$ is real. Likewise, it follows that $\overline{w\bar \nu(h_\ga)} = w \nu(h_\ga).$ Finally, since $\rho(h_\ga) = \rho_\Phi(h_\ga),$ it follows that the element (\ref{e: conj els of spec}) equals 
$w(-\gL_\gs +\nu) (h_\ga) -\rho(h_\ga).$ 
\qed

By integration the 
equation of Lemma \ref{l: diff eq for F} leads to the equality
\begin{equation}
\label{e: integral equation with gl}
F(f,\nu, \bp aa_t) = e^{tB(\nu)} F(\bp a) + e^{tB(\nu)} \int_0^t e^{-\tau B(\nu)} R(f, \nu, a)\; d\tau,
\end{equation}
for $f \in \IPgs,$ $\nu \in \faPdc(\ge_0),$ $\bp a \in A_\Phi$ and $t \in \R.$
We will first derive estimates for $F$ and $R$, following from the information that
$\xxi \in \fad$  dominates $\wh = (\wh_\nu).$ This means that there exist $\ge > 0, s >0, N >0$ 
and a continuous seminorm $\snorm$ on $\IPgs$ such that (\ref{e: domination whittaker family}) is valid 
for all $f \in \IPgs$, $\nu \in \faPdc(\ge)$ and $a \in A.$   

\begin{lemma}
\label{l: estimate F}
There exist $\ge > 0, s>0, N > 0$ and a continuous seminorm $n$ on $\IPgs$ such that 
\begin{equation}
\label{e: estimate F}
\|F(f, \nu, a)\| \leq |(\nu,a)|^N e^{s |\Re \nu| |\log a| } a^\xxi \snorm(f)
\end{equation}
for all $f \in I_\gs^\infty,$  $\nu \in \fad(\ge)$ and $a\in A.$  
\end{lemma}

\proof 
We take $\ge >0$ sufficiently small such that the estimate (\ref{e: domination whittaker family})
is valid. 
Since 
$$ 
F(f, \nu, a)_j = \inp{   \pi_{-\nu}(a)^{-1} \pi_{-\nu}(\bar u_j^\veec)f     }{  j_{\bar\nu}}
=
\wh_\nu( \pi_{-\nu}(\bar u_j^\veec)f )(a)
$$ 
the estimate (\ref{e: estimate F}) follows from (\ref{e: domination whittaker family}) with 
the same $s$ and possibly enlarged constant $N> 0$ and enlarged  
seminorm $n.$
\qed

\begin{rem}{\rm (structure of proof).\ }
Throughout the proof of Lemma \ref{l: improvement step} we will prove assertions of the form
\begin{equation}
\label{e: assertion A}
\exists (\ge,s, N, n): A(\ge, s, N, n)
\end{equation} where $A(\ge, s , N, n)$ is an assertion (usually containing an estimate) depending on positive constants $\ge, s, N$ and a continuous seminorm $n$ on $\IPgs.$ 
Moreover, the assertion has the property that for any $(\ge', s', N', n')$ with 
$\ge' \leq \ge, s' \geq s, N' \geq N$ and $n' \geq n,$ 
\begin{equation}
\label{e: property of assertion}
A(\ge, s, N, n) \implies A(\ge', s', N', n').
\end{equation}
A typical assertion of this type is the assertion 
$$ 
\forall (f \in \IPgs, \nu \in \faPdc, a \in A): 
\|F(f, \nu, a)\| \leq |(\nu,a)|^N e^{s |\Re \nu| |\log a| } a^\xxi \snorm(f)
$$ 
of Lemma \ref{l: estimate F}.
The proofs of assertions of this type will make use of finitely many valid similar assertions
$\exists (\ge, s, N, n) : A_i(\ge, s, N, n),$ for $i \in I,$ with $I$ a finite index set. If all $A_i$ have the property (\ref{e: property of assertion}) then it follows that also 
$\exists  (\ge, s, N, n) : \wedge_{i\in I} A_i(\ge, s, N, n) $ is valid. Indeed, if $A_i(\ge_i , s_i, N_i, n_i)$ 
is true for every $i \in I,$ then $A_i(\ge, s, N, n),$ for $i \in I,$ are simultaneously valid 
as soon as $\ge \leq \min_i \ge_i, \; s \geq \max_i s_i , N \geq \max_i N_i$ and $n \geq \max_i n_i.$ 
In the proof we shall indicate this informally by saying 
that $A_i(\ge, s, N, n)$ are valid for sufficiently small 
$\ge$ and sufficiently large $s,N$ and $n.$ A logical reasoning will then give
the validity of $A(\ge', s', N', n')$ for suitably chosen $\ge', s', N', n',$ which finally allows 
the conclusion that (\ref{e: assertion A}) is valid.
\end{rem} 

\begin{lemma}
\label{l: estimate R}
There exist $\ge >0, s>0, N > 0$ and a continuous seminorm $n$ on $\IPgs$ such that 
\begin{equation}
\label{e: estimate R}
\|R(f,\nu, a)\| \leq |(\nu,a)|^{N} e^{s |\Re \nu| |\log a| } a^{\xxi-\ga} \snorm(f)
\end{equation}
for all $f \in \IPgs$, $\nu \in \faPdc(\ge)$ and 
$a \in A$.
\end{lemma} 
\proof
We recall that $R(f,\nu)$ is the $\C^\ell$-valued function defined by 
$$ 
R(f, \nu, a)_j  := \inp{  \pi_{-\nu}(a)^{-1} f  }{  \pi_{\bar \nu}(v_j) j_{\bar \nu}}.
$$ 
We will now derive the estimate for $R$ with the required properties
of uniformity.
For this we note that $\pi_{\bar \nu}(v_j)  j_{\bar\nu}$ may be written as
a finite sum of terms of the form $\pi_{\bar \nu}(U )j_{\bar\nu}$, 
with $U \in U(\bfn_\Phi)U(\bar \fn_0 + \fm_1)$ of $\fa$-weight $-\mu \in - \sum_{\gb \in \gD} \N \gb$ such that  $\mu(h_\ga) \geq 1.$ Each corresponding term $r(f,\nu,a)$  in $R(f,\nu,a)_j$ 
may be rewritten as 
\begin{eqnarray*} 
r(f,\nu, a) &= & \inp{ \pi_{-\nu}(U^*) \pi_{-\nu}(\bp aa_t)^{-1} f}{j_{\bar \nu}}\\
&=& (\bp a a_t)^{-\mu} \inp{ \pi_{-\nu}(\bp aa_t)^{-1} \pi_{-\nu}(U^*) f}{j_{\bar \nu}}. 
\end{eqnarray*}
We note that the restriction of $\mu$ to ${\bp \fa_\Phi}$ equals the restriction
of $\mu_\Phi := \sum_{\gb\in \Phi} \mu_\gb \gb$ to this space. For each $\gb \in \Phi$ we may 
select a simple root vector $X_\gb \in \fg_{\gb}$ such that $\cchi_*(X_\gb) = 1.$ 
The product $X:= \prod_{\gb\in \Phi} X_\gb^{\mu_\gb}$ belongs to $U(\fn_0),$ 
satisfies $\cchi_*(X) = 1$ and has $\fa$-weight $\mu_\Phi.$ Therefore, 
$$ 
\Ad(\bp a a_t) X = (\bp a)^\mu X
$$ 
and it follows that
\begin{eqnarray*}
r(f,\nu, a)   
&=&  (\bp a a_t)^{-\mu} \inp{ \pi_{-\nu}(X^*) \pi_{-\nu}(\bp aa_t)^{-1} \pi_{-\nu}( U^*) f}{j_{\bar \nu}}\\
&=&   (a_t)^{-\mu}\inp{ \pi_{-\nu}(\bp aa_t)^{-1} \pi_{-\nu}(X^*U^*) f}{j_{\bar \nu}}\\
&=&  (a_t)^{-\mu} \wh_\nu( \pi_{-\nu}(X^*U^*) f)(\bp a a_t).
\end{eqnarray*}
We now select $\ge, s, N >0$ and $n$ a continuous seminorm 
on $\IPgs$ which make (\ref{e: domination whittaker family}) valid for  $f, \nu,a$ 
in the indicated sets.
Then it follows that there exists a constant $C_\mu>0,$ only depending on 
$\mu,$ such that 
\begin{eqnarray*}
|r(f,\nu, a)|
& \leq &
C_\mu (a_t)^\mu (\bp aa_t)^\xxi |(\nu, a)|^N e^{s|\Re \nu||\log (\bp aa_t)|} \snorm(\pi_{-\nu}(X^* U^*) f)\\
 & \leq & 
(\bp a a_t)^{\xxi -\ga} |(\nu,a)|^{N} e^{s|\Re \nu||\log (\bp aa_t)|} \snorm'(f) 
\end{eqnarray*} 
with $\snorm'$ a seminorm on $\IPgs$, independent of $f.$ Combining the above estimates and enlarging $\snorm, N$ if necessary, we find that 
$$ 
\|R(f,\nu, a)\| \leq |(\nu,a)|^{N} e^{s |\Re \nu| |\log a| } a^{\xxi-\ga} \snorm(f),
$$
for all $f \in \IPgs$, $\nu \in \faPdc(\ge)$ and $a \in A$. 
\qed

Our next goal is to show how this stronger estimate on the remainder term leads to an improved estimate
of $F$, hence of $\wh_\nu(f)= F_1(f,\nu).$ For this we need to decompose the formula (\ref{e: integral equation with gl}) into parts corresponding to the spectrum of $B(\nu).$ 
In the course of the argument we will impose finitely many conditions on the constant $\ge >0$, ensuring it is sufficiently small. 

Let us first analyze the spectrum of $B(\nu),$ in particular its dependence on $\nu.$ 
For $w \in W(\fh)$ we  
write $\xx_w = - w(\gL_\gs)(h_\ga) - \rho(h_\ga)$ which is a real number as shown in the proof
of Lemma \ref{l: diff eq for F}. For every $\nu \in \faPdc$ the spectrum of $B(\nu)$ 
consists of complex numbers of the form 
$$ 
\xx_w + w(\nu)(h_\ga),
$$ 
for $w \in W(\fh).$ 
Put 
$$ 
X: = \{ \xx_w \mid w \in W(\fh)\}.
$$ 
Let $\gg >0 $ be a positive real number such that all distinct elements $\xx_1$ and $\xx_2$ 
of $X$  satisfy $|\xx_1 - \xx_2| > 2 \gg.$  Fix $\ge_1 > 0$ such that 
$$
\ge_1 < \frac12 \gg |h_\ga|^{-1}.
$$
In the course of this section we will always assume that 
$0< \ge \leq \ge_1.$

We note that $W(\fh)$ preserves the subspace $\fh_\R = i\ft + \fa$ of $\fh_\iC.$ 
Let $p_\fa$ denote the projection $\fh_\R \to \fa$ along $i\ft.$ 
Then for all $\nu \in \faPdc(\ge)$ and  $w \in W(\fh)$
we have
\begin{eqnarray*} 
| \Re(w\nu(h_\ga))| & = &  |\Re [\nu(w^{-1}(h_\ga))]| \\
&=& |\Re [\nu(p_\fa (w^{-1}(h_\ga)))]|\\
&\leq &  |\Re \nu| |p_\fa(w^{-1}(h_\ga))| \leq |\Re \nu| |h_\ga|.
\end{eqnarray*}
Likewise, for all $\nu \in \faPdc$ we have
$$ 
|\Im (w\nu(h_\ga))| \leq |\Im \nu||h_\ga|.
$$ 
Given $\xx \in X$ and $w \in W(\fh)$ such that $\xx_w= \xx$ we thus see that for all
$\nu \in \faPdc(\ge)$ we have
\begin{equation}
\label{e: estimate real part spec}
|\Re ( \xx_w  + w\nu(h_\ga)) - \xx| \leq |\Re \nu ||h_\ga| < \gg/2,
\end{equation}
and 
\begin{equation}
\label{e: estimate ima part spec} 
|\Im (\xx_w + w\nu(h_\ga))| \leq |\Im(\nu)||h_\ga|.
\end{equation} 
For $t\in \R$ and $\nu \in \faPdc$ we define 
$$ 
C_{t,\nu} := \frac{\gg}{2} ( 1+ |t|)^{-1}(1 + |\nu|)^{-1} .
$$ 
and the rectangle $R(t,\nu) \subset \C$ 
to be the set of points $z \in \C$ such that 
$$ 
|\Re z| \leq |\Re \nu| |h_\ga)| + C_{t,\nu}, \quad |\Im z| \leq |\Im \nu||h_\ga| + C_{t,\nu}.
$$ 
Then it is readily seen that, for $t \in \R$ and $\nu \in \faPdc(\ge),$ 
$R(t, \nu) \subset [-\gg, \gg] + i\R,$ from which it readily follows that the translated 
rectangles $\xx + R(t,\nu)$ are mutually disjoint. For $\nu \in \faPdc(\ge)$ 
we define $S(\nu)$ to be the rectangle of points with 
$$ 
|\Re z| \leq |\Re \nu| |h_\ga)|, \quad |\Im z| \leq |\Im \nu||h_\ga|.
$$ 
Then it is clear that $S(\nu)$ is contained in the interior of $R(t,\nu)$ 
so that the translated rectangles $\xx + S(\nu),$ $(\xx \in X),$ are mutually disjoint as well.
Furthermore, it follows from the estimates (\ref{e: estimate real part spec}) and 
(\ref{e: estimate ima part spec}) that 
\begin{equation}
\label{e: deco spec Bnu}
\spec\, B(\nu) = \bigcup_{\xx \in X}\;\; \spec\, B(\nu) \cap [\xx + S(\nu)].
\end{equation} 
For $\xx \in X$ and $\nu \in \faPdc(\ge)$ we denote by $P_{\xx}(\nu)$  the 
spectral projection of $B(\nu)$ onto the sum of the generalized eigenspaces 
corresponding to the eigenvalues from $\spec\, B(\nu) \cap [\xx + S(\nu)].$ 
Then $P_\xx$ is a holomorphic function on $\faPdc(\ge)$ with values in $\End(\C^\ell).$ We note that 
$$ 
I = \sum_{\xx \in X} P_\xx(\nu).
$$

\begin{lemma}
There exists a $C > 0$ such that for every $\xx \in X$, all $t \in \R$ and $\nu \in \faPdc(\ge_1),$ 
\begin{equation}
\label{e: estimate projection}
\| e^{t B(\nu)} P_\xx(\nu) \| \leq C (1 + |t|)^p (1 + |\nu|)^p e^{\xx t  + |h_\ga||\Re \nu | |t|}.
\end{equation}
\end{lemma}

\proof  
Since for all $t \in \R,$ $ \nu \in \faPdc(\ge_1)$ we have that 
\begin{equation}
\label{e: spec Bnu in Snu}
\spec \,B(\nu) \cap (\xx + R(t,\nu)) \subset \xx + S(\nu) \subset {\rm int} (x + R(t,\nu)),
\end{equation}
it follows that 
$$ 
P_\xx (\nu)e^{tB(\nu)}= \frac1{2\pi i} \int_{\xx + \partial R(t,\nu)} e^{tz} (zI - B(\nu))^{-1}\; dz.
$$ 
We will complete the proof by estimation of the integral. First of all, the length of
the boundary of $\xx + R(t, \nu )$ is estimated by
\begin{equation}
\label{e: length partial R}
{\rm length}(\partial R(t, \nu)) \leq 2 \ge |h_\ga| + 2|\nu||h_\ga| + 4 C_{t,\nu} \leq 
3\gg\ge^{-1}(1 + |\nu|).
\end{equation}
For $z \in \xx + \partial R_{t,\nu}$ we have  
$$ 
|e^{-t\xx} e^{tz}| \leq e^{|t||\Re z -\xx|}  \leq e^{|t| |\Re \nu| |h_\ga| + |t| C_{t,\nu}},
$$ 
so that
\begin{equation}
\label{e: estimate power e t z}
|e^{tz}|\leq e^{t\xx} e^{|t| |\Re \nu| |h_\ga|}e^{\gg/2}.
\end{equation}
If $\xx, \xx' \in X$ then the distance from $\xx +\partial R(t,\nu)$ to $\xx' + S(\nu)$ 
is at least $C_{t,\nu}.$ 
From (\ref{e: deco spec Bnu}) and (\ref{e: spec Bnu in Snu}) we now see that for $z \in \xx + \partial R_{t,\nu}$ the distance of $z$ to the spectrum of $B(\nu)$ is at least $C_{t,\nu}$ so that
$$ 
|\det (zI - B(\nu))|^{-1} \leq C_{t,\nu}^{- \ell} \leq (2/\gg)^\ell (1 + |t|)^\ell ( 1+ |\nu|)^\ell.
$$ 
In view of Cramer's rule, there exists a constant $C_\ell > 0$ such that,
for all $A \in {\rm GL}(\ell, \C),$
$$ 
\|A^{-1}\| \leq C_\ell |\det A|^{-1}( 1+ \|A\|)^{\ell -1}.
$$ 
Applying this with $A = (z I - B(\nu))$ we see that for $z \in \xx + \partial R(t,\nu)$,
 $$ 
\|(zI - B(\nu))^{-1}\| \leq C_\ell (2/\gg)^\ell (1 + |t|)^\ell ( 1+ |\nu|)^\ell (|z| + \|B(\nu)\|)^{\ell -1}.
$$
As $\nu \mapsto B(\nu)$ is polynomial in $\nu,$ there exist
constants $N \in \N$ and $C' >0$ such that for all $t\in \R$, $\nu\in \faPdc$
and $z \in \xx + \partial R(t,\nu),$ 
\begin{equation}
\label{e: estimate resolvent}
\|(zI - B(\nu))^{-1}\|   \leq C' (1 + |t|)^\ell (1 +|\nu|)^N.
\end{equation} 
Combining the estimates (\ref{e: length partial R}), (\ref{e: estimate power e t z})
and (\ref{e: estimate resolvent}) we infer the existence of a constant $C > 0$ 
such that for all $t \in \R$ and $\nu \in \faPdc(\ge_1)$ the estimate 
(\ref{e: estimate projection}) is valid with $ p = \max (\ell , N+1).$ 
\qed
We now decompose $F(f, \nu, a)$ and $R(f,\nu, a)$ into components
$$ 
F_\xx (f, \nu, a) := P_\xx(\nu) F(f, \nu, a) , \qquad  R_\xx (f, \nu, a) = P_\xx(\nu) R(f, \nu, a).
$$ 
Then using (\ref{e: estimate projection}) with $t = 0$ we obtain, after decreasing $\ge >0$ 
and increasing $N$ and $\snorm$ suitably, 
$$ 
\|F_\xx(f, \nu, a) \| \leq  |(\nu,a)|^N e^{s|\Re  \nu||\log a|}a^\xxi \snorm(f)
$$ 
and
$$ 
\|R_\xx(f, \nu, a) \| \leq  |(\nu,a)|^N e^{s|\Re \nu||\log a|}a^{\xxi- \ga} \snorm(f)
$$ 
for all $f \in \IPgs,$ $\nu\in \faPdc(\ge),$ $a \in A$. 

We will now obtain sharper estimates for $F_\xx$, for each $\xx \in X. $ 
Our main tool will be the following identity which follows from  (\ref{e: integral equation with gl}) by application
of $P_\xx(\nu);$ 
\begin{equation}
\label{e: integral equation with gl and xi}
F_\xx(f,\nu, \bp aa_t) = e^{tB(\nu)} F_\xx(f,\nu, \bp a) + e^{tB(\nu)} \int_0^t e^{-\tau B(\nu)} R_\xx(f, \nu, a)\; d\tau.
\end{equation}
In the course of this proof we will need to distinguish two cases, depending 
on which of the following sets $\xx$ belongs to:
\begin{equation}
\label{e: defi X pm} 
X_+ : = X \cap \; ]\xxi(h_\ga) - 1, \infty[, \quad{\rm and} \quad X_- :=X\setminus X_+.
\end{equation}
There exists a constant $c_0 \in [0,1[$ such that $(\xxi(h_\ga) - [c_0, 1[\,)  \cap X = \emptyset.$ 
Furthermore, we fix an arbitrary $c\in \, [c_0 ,1[.$  Then for all $\xx \in X$ the following
assertion is valid:
\begin{equation}
\label{e: condition on c}
\xx > \xxi(h_\ga) - 1  \iff   \xx > \xxi(h_\ga) - c.
\end{equation}
We put $\xxi':= \xxi - c\ga;$ then by (\ref{e: condition on c})
we have,
for $\xx \in X,$  
$$ 
\xx \in X_\pm  \iff \pm(\xx -  \xxi'(h_\ga)) > 0.
$$ 
We fix $\ge_2> 0$ such that 
$$
\xxi'(h_\ga) + [-\ge_2|h_\ga|, \ge_2|h_\ga|] \cap X = \emptyset.
$$ 
In the course of this section we will always assume that  $0< \ge \leq \ge_2.$ Then for 
every $\xx \in X$ and all $\nu \in \faPdc(\ge)$ the real part of the spectrum  $\;\spec [B(\nu)|_{{\rm im}\,(P_\xx(\nu) )}]$ 
is contained in $[\, \xx-\ge|h_\ga|, \xx + \ge|h_\ga|],$ which is contained
in the interval $] \xxi'(h_\ga) , \infty[$ if $\xx \in X_+$ and in 
 the interval 
$]-\infty, \xxi'(h_\ga)[$ if $\xx \in X_-$.

\begin{lemma}
\label{l: estimate integral of R xi}
Assume that $\xxi \in \fad$ dominates $(\wh_\nu)$ and let $\ga > 0$ be simple.
Then there exist $\ge>0$, $N > 0$, $  s>0 $ and a continuous seminorm $\snorm$ on $I_\gs^\infty$ such that for every $f \in I_\gs^\infty$, $\nu \in \faPdc(\ge)$, $\bp a \in \bp A_\Phi$,
and $t \geq 0$ the following estimates are valid.
\begin{enumerate}
\itema
If $\xx \in X_-$, then 
$$ 
\| \int_0^t e^{(t -\tau) B(\nu)} R_\xx(f, \nu, \bp aa_\tau)\; d\tau \;\| \leq 
|(\nu, \bp aa_t)|^{N} (\bp aa_t)^{\xxi'}  e^{s|\Re \nu||\log (\bp a a_t) | }\snorm(f).
$$ 
\itemb
If $\xx \in X_+,$ then 
$$ 
\| \int_t^\infty e^{(t -\tau)B(\nu)} R_\xx(f, \nu, \bp aa_\tau)\; d\tau\;  \| \leq 
|(\nu, \bp aa_t)|^{N} (\bp aa_t)^{\xxi'}  e^{s|\Re \nu||\log (\bp a a_t) | }\snorm(f),
$$
\end{enumerate}
with absolutely converging integral.
\end{lemma}

\proof
After decreasing $\ge >0$ and increasing $s,N$ and $\snorm$ if necessary, we obtain, for all $f\in 
\IPgs,$ $\nu \in \faPdc(\ge)$ $\bp a \in \bp A_\Phi$ and $t, \tau \in \R $, 
\begin{equation} 
\label{e: estimate R xi}
\|e^{(t-\tau) B(\nu)} R_\xx(f,\nu, \bp aa_\tau) \| \leq   C_N(\nu, \bp a)\, D_N(\nu, t, \tau)\, \snorm(f),
\end{equation} 
where
$$ 
C_N(\nu, \bp a) = |(\nu, \bp a)|^N (\bp a)^\xxi e^{s|\Re \nu| |\log \bp a|}
$$
and  
\begin{eqnarray*}
D_N(\nu,t, \tau) &  =  & 
|(t,\tau)|^{N} e^{(t - \tau)\xx} e^{ (s+ 1) |\Re \nu| |h_\ga||\tau|}
(a_\tau)^{\xxi - \ga}\\
&=& |(t,\tau)|^{N} e^{t\xx}e^{ (s+1) |\Re \nu| |h_\ga||\tau|}  e^{[ -\xx + \xxi'(h_\ga)]\tau} e^{(c-1)\tau} .
\end{eqnarray*}
Here we have used the notation
$
|(t,\tau)|^{N} = ( 1 + |\tau|)^N (1 +|t|)^N
$ 
for $\tau, t \in \R.$ 

In order to prove (a), assume that $\xx < \xxi'(h_\ga).$ 
Then we may fix $\ge_s >0$ so that 
$$ 
\ge_s (s+1) |h_\ga|  - \xx +  \xxi'(h_\ga)  > 0. 
$$ 
By decreasing $\ge$ further if necessary, we may asssume that $0 < \ge \leq \ge_s.$ 
Then for $t \geq 0$ and $\nu \in \fadc(\ge)$ we have that the function 
$$ 
\tau \mapsto e^{t\xx} e^{ (s+1) |\Re \nu| |h_\ga||\tau|}  e^{[ -\xx   + \xxi'(h_\ga)]\tau}
$$ 
is increasing on $[0,t]$, hence dominated by its value at $t$, so 
that  
\begin{eqnarray*} 
\int_0^t D_N(\nu, \tau,t) d\tau 
& \leq & e^{ [ (s+1)|\Re \nu| |h_\ga|  + \xxi'(h_\ga) ] t   } 
\int_0^t (1+ |(t,\tau)|)^N e^{(c-1)\tau} \; d \tau \\
& \leq & C'    e^{(s+1) |\Re \nu||\log a_t|}   (a_t)^{\xxi'} (1+ |\log a_t|)^N , 
\end{eqnarray*} 
where 
\begin{equation}
\label{e: int for C p}
C' = \int_0^\infty (1 + \tau)^N e^{(c-1)\tau} \; d\tau.
\end{equation}
We observe that 
$(\bp a)^\xxi = (\bp a)^{\xxi'}.$ Accordingly, it now follows
that 
\begin{equation*}
C_N(\nu, \bp a) \int_0^t D_N(\nu, \tau,t) d\tau 
  \leq   
C' |(\nu, \bp aa_t)|^{2N} 
e^{2(s+1)|\Re \nu||\log (\bp a a_t) | }  (\bp aa_t)^{\xxi'}.
\end{equation*}
In view of (\ref{e: estimate R xi}) we finally obtain the desired estimate
of (a), with $s, N$ and $\snorm$ chosen large enough.

We now turn to (b), and assume $\xx > \xxi'(h_\ga).$ Then there exists 
$\ge_s' > 0$ such that 
$$ 
\ge_s' (s+1) |h_\ga|  -\xx + \xxi'(h_\ga)  < 0.
$$ 
Decreasing $\ge$ further if necessary, we may assume that 
$0<\ge < \ge_s'.$ Then for $t \geq 0$ and $\nu \in \fadc(\ge)$ 
the function
$$ 
\tau \mapsto e^{t\xx} e^{ (s+1) |\Re \nu| |h_\ga||\tau|}  e^{[ -\xx   + \xxi'(h_\ga)]\tau}
$$ 
is decreasing on $[t, \infty)$ hence dominated by its value  at $t$, so that 
\begin{eqnarray*} 
\int_t^\infty D_N(\nu, \tau,t) d\tau 
& \leq & e^{ [ (s+1)|\Re \nu| |h_\ga|  + \xxi'(h_\ga) ] t   } 
\int_t^\infty (1+ |(t,\tau)|)^N e^{(c-1)\tau} \; d \tau \\
& \leq & C'   (1+ |\log a_t|)^N  e^{(s+1) |\Re \nu||\log a_t|}   (a_t)^{\xxi'}, 
\end{eqnarray*} 
with $C'$ given by (\ref{e: int for C p}).
As in the first part of the proof, we now infer that 
$$ 
C_N(\nu, \bp a) \int_t^\infty D_N(\nu, \tau,t) d\tau \leq 
C' |(\nu, \bp aa_t)|^{2N} 
e^{2(s+1)|\Re \nu||\log (\bp a a_t) | }  (\bp aa_t)^{\xxi'}.
$$
Using (\ref{e: estimate R xi}) and further enlarging $s, N$ and $\snorm$ if necessary, we 
obtain the desired estimate of (b).
\qed

\begin{prop}
\label{p: xi and domination by eta}
Let the Whittaker family $(\wh_\nu)$ be dominated by $\xxi \in \fad.$ 
Assume that $\xxi\succeq -\rho.$ 
Let $\xx \in X.$ 
\begin{enumerate}
\itema If $\xx \leq \xxi(h_\ga) -1$, then $F_\xx$ is dominated by 
$\xxi' = \xxi - c \ga$ for all $c \in [0,1).$ 
\itemb If $\xx > \xxi(h_\ga) - 1 \geq  -\rho(h_\ga)$, 
then $F_\xx$ is dominated by 
$\xxi' = \xxi - c \ga$ for all $c \in [0,1).$ 
\itemc
If $\xx > \xxi(h_\ga) - 1 $ and $\xxi(h_\ga) - 1 < -\rho(h_\ga)$,
then $F_\xx$ is dominated 
by $i_\ga(\xxi).$ 
\end{enumerate}
\end{prop}

Before we continue with the proof of Proposition \ref{p: xi and  domination by eta}, we will first argue that the proposition is sufficient for the proof of the improvement step asserted in
Lemma \ref{l: improvement step}. We first observe that $F$ is dominated 
by a functional $\theta \in \fad$ if and only if every component $F_\xx$, for 
$\xx \in X$ is dominated by $\theta$, where the obvious extension  of the notion 
of domination is assumed. 
\medno
{\em Proof of Lemma \ref{l: improvement step}. }
We begin by observing that the hypothesis on $\xxi$ in Lemma \ref{l: improvement step} (a) 
guarantees that $\xxi' \succeq -\rho.$ In (b) of the lemma, $i_\ga(\xxi) \succeq -\rho$ 
by virtue of Lemma \ref{l: i ga eta majorizes min rho}. Thus it is sufficient to establish the asserted dominations.

To establish Lemma \ref{l: improvement step} (a), assume that $\xxi(h_\ga) - 1 \geq  -\rho(h_\ga).$  Let $\xx \in X$ and assume that 
$\xx > \xxi(h_\ga) -1.$ Then by (b) of the above proposition, it follows 
$F_\xx$ is dominated by $\xxi':= \xxi - c\ga$ for each $c \in [0,1).$ 
By (a) of the above proposition, the same is true 
for all remaining $\xx \in X$.  
This establishes Lemma \ref{l: improvement step} (a).

To prove (b) of Lemma \ref{l: improvement step}, assume that  $\xxi(h_\ga) - 1 <  -\rho(h_\ga).$ 
Since $\xxi(h_\ga) \geq -\rho(h_\ga)$, it follows that 
$i_\ga(\xxi) = \xxi - d\ga $ for a unique $d \in [0,1).$  

If $\xx \in X$ satisfies $\xx \leq \xxi(h_\ga) - 1$, then according to (a) of the above proposition,
$F_\xx$ is dominated  by $\xxi' = \xxi - d \ga = i_\ga(\xxi).$ 

On the other hand, if $\xx > \xxi(h_\ga) - 1$, then it follows from  
(c) of the above proposition that $F_\xx$ is dominated by  $i_\ga(\xxi).$  

We conclude that every $F_\xx$ is dominated by $i_\ga(\xxi)$, hence so is $F.$ 
This establishes Lemma \ref{l: improvement step} (b).
\qed

{\em Proof of Proposition \ref{p: xi and domination by eta}.\ }
Before we proceed with the proof, we note that in the cases (a) and (b) it suffices to 
show that $F_\xx$ is dominated by $\xxi' = \xxi -  c\ga$ for all $c \in [c_0,1 [.$ Indeed 
assume this to be the case and let $c' \in [0, c_0]$ and put $\xxi''= \xxi - c_\ga.$ 
Then $\xxi'' \leq \xxi'$ on $\bp\! \fa_\Phi + \R_+ h_\ga$ whereas 
$\xxi'' \leq \xxi $ on $\bp\! \fa_\Phi + \R_- h_\ga.$ It therefore follows from the domination
of $F_\xx$ by both $\xxi$ and $\xxi'$ that $F_\xx$ is dominated by $\xxi''.$

For the actual proof, let $\xx \in X.$ 
To establish (a) assume $\xx \leq \xxi(h_\ga) -1.$ Let $c \in [c_0, 1[.$ 
Then obviously $\xx < \xxi'(h_\ga)$, where $\xxi' = \xxi - c\ga.$ 
We now have an estimate of the type of Lemma \ref{l: estimate integral of R xi} (a).
On the other hand, we also have the identity (\ref{e: integral equation with gl and xi}). From the domination assumption on $F$ and the estimate of Lemma \ref{l: estimate integral of R xi} (a) it follows that there exists 
$\ge > 0$ such that for $N, s$ and the continuous seminorm $\snorm$ 
all chosen large enough, we have, for all $f \in \IPgs$, $\nu \in \faPdc(\ge)$, 
$\bp a \in \bp A_\Phi$, and $t \geq 0,$ 
 $$ 
e^{t B(\nu)} F_\xx(f, \nu, \bp a) \leq (1 + |t|)^N e^{t(\xx + |\Re \nu|) } |(\nu, \bp a)|^N a^{s|\Re \nu||\log \bp a|} 
(\bp a)^\xxi \snorm(f).
$$  
Using that $(1 + |t|)^N = (1 + |\log a_t|)^N \leq (1 + | \log a|)^N$, that 
$(\bp a)^\xxi = (\bp a)^{\xxi'}$, and that 
$$
e^{t(\xx + |\Re \nu|) } \leq e^{t\xxi'(h_\ga)}e^{t |\Re \nu||h_\ga|} =  e^{|\Re \nu ||\log a_t|}(a_t)^{\xxi'}
$$ 
we obtain
the estimate
$$ 
e^{t B(\nu)} F_\xx(f,\nu, \bp a) \leq |(\nu,a)|^{2N} a^{(s+1)|\Re \nu||\log a|} a^{\xxi' } \snorm(f),
$$
where $a = \bp a a_t$, $t \geq 0.$ 
Combining this with the identity (\ref{e: integral equation with gl and xi}) and  Lemma \ref{l: estimate integral of R xi} (a) we see that we may enlarge
$N, s$ and the continuous seminorm $\snorm$ such that for all $f \in \IPgs$, 
$\nu \in \faPdc(\ge)$ all $\bp a \in \bp A_\Phi$ and all $t \geq 0$ we have the estimate
\begin{equation}
\label{e: improved estimate F xi}
F_\xx(f, \nu, \bp a a_t) \leq |(\nu, a)|^N e^{s|\Re \nu||\log a|} a^{\xxi'} \snorm(f).
\end{equation}
After decreasing $\ge,$ and increasing $s, N, \snorm$, the same estimate becomes valid 
for all $f \in \IPgs$, $\nu \in \faPdc(\ge)$, $\bp a\in \bp A_\Phi$ and $t \leq 0$, 
provided $\xxi'$ is replaced by $\xxi.$ Since $\xxi' \leq \xxi$ 
on $\bp \fa_\Phi + (-\infty, 0] h_\ga$, we see that the estimate 
(\ref{e: improved estimate F xi}) is in fact valid for all $t \in \R.$ 
It follows that $\xxi'$ dominates $F_\xx.$ This establishes (a).

We turn to (b) and (c) and assume that $\xx > \xxi(h_\ga) -1.$ 
Fix $c\in [c_0, 1[$.
Then by (\ref{e: condition on c}) we have 
$\xx > \xxi'(h_\ga)$, where $\xxi' = \xxi - c\ga.$ 
From Lemma \ref{l: estimate integral of R xi} (b) with $t =0$ we now obtain that 
the integral
$$ 
I_\xx(f, \nu, \bp a) := \int_0^\infty  e^{-\tau B(\nu)} R_\xx(f,\nu, \bp aa_\tau) \; d\tau
$$ 
converges for $\nu \in \faPdc(\ge)$ and satisfies the estimate
\begin{equation}
\label{e: estimate for I xi}
\| I_\xx(f, \nu, \bp a)\| \leq |(\nu, \bp a)|^N  e^{s|\Re \nu| |\log \bp a|} (\bp a)^{\xxi'} \snorm(f),
\end{equation}
provided that $\ge>0$ is taken sufficiently small, and $s, N, n$ suffiently large.
Put 
\begin{equation}
\label{e: defi F infty}
F^\infty_\xx(f,\nu, \bp a): =  F_\xx(f, \nu, \bp a) + I_\xx(f,\nu, \bp a).
\end{equation}
From (\ref{e: integral equation with gl}) we see that 
\begin{equation}
\label{e: expression F xi with int J}
F_\xx(f, \nu, \bp a a_t) = e^{t B(\nu)} F^\infty_\xx(f, \nu, \bp a)  + R^\infty_\xx(f,\nu, \bp aa_t),
\end{equation}
where the last term is  given by the convergent integral
\begin{equation}
\label{e: defi R infty 2}
R_\xx^\infty(f, \nu, \bp aa_t) = - \int_t^\infty e^{(t -\tau)B(\nu)} R_\xx(f, \nu, \bp aa_\tau)\; d\tau.
\end{equation} 
From Lemma \ref{l: estimate integral of R xi} (b) it follows that
$\ge>0$ can be decreased, and $s, N$ and $n$ increased such that
\begin{equation}
\label{e: estimate R infty} 
\| R_\xx^\infty(f, \nu, \bp aa_t)\| \leq 
|(\nu, \bp aa_t)|^{N} (\bp aa_t)^{\xxi'}  e^{s|\Re \nu||\log (\bp a a_t) | }\snorm(f),
\end{equation}
for all $f \in \IPgs,$ $\nu \in \faPdc(\ge),$ $\bp\! a\in \bp\!A_\Phi$ and $t \geq 0.$ 

From the domination assumption on $F$ we obtain the estimate 
(\ref{e: estimate for I xi}) for $F_\xx(f, \nu, \bp \!a)$ in place of $I_\xx(f, \nu, \bp\!a)$, provided 
we shrink $\ge$ and enlarge $N,s, \snorm$ if necessary. Here we need that
$(\bp a)^{\xxi} = (\bp a)^{\xxi'}$ for  $\bp a \in \bp A_\Phi.$ This observation leads to the estimate
$$
\|F^\infty_\xx(f, \nu, \bp a)\| \leq |(\nu, \bp a)|^N a^{s|\Re \nu||\log \bp a|} 
(\bp a)^{\xxi'} \snorm(f),
$$ 
for all $f \in \IPgs $, $\nu \in \faPdc(\ge)$, and $\bp a \in \bp A_\Phi$. 
Combining this estimate with (\ref{e: estimate projection}) we see that we may increase $N, \snorm$ further to arrange that
\begin{equation}
\label{e: estimate F xi infty}
\|e^{tB(\nu)} F^\infty_\xx (f, \nu, \bp a) \|  \leq  
|(\nu, \bp a)|^N(\bp a)^\xxi e^{s|\Re\nu||\log \bp a|}e^{\xx t + |\Re \nu||\log a_t|}
\snorm(f)
\end{equation}
for all $f \in \IPgs$,  $\nu \in \faPdc(\ge)$,   $\bp a \in \bp A_\Phi$, and  $t \geq 0$. 

We will first consider the case that $F^\infty_\xx= 0.$ 
Then it follows from 
(\ref{e: expression F xi with int J}) and 
(\ref{e: estimate R infty}) that 
 $$ 
\|F_\xx(f, \nu, \bp aa_t)\| \leq |(\nu,\bp aa_t)|^N (\bp aa_t)^{\xxi'} e^{s|\Re \nu||\log \bp aa_t)|}\, n(f).
$$ 
for $f \in \IPgs$, $\nu \in \faPdc(\ge)$, $\bp a \in \bp A_\Phi$, $t \geq 0.$ 
Here we used that $\xxi'  = \xxi$ on $\bp \fa_\Phi.$ 
In view of the assumed domination of $(\wh_\nu)$ by $\xxi$, we have the similar estimate
for $t \leq 0$, with $\xxi$ in place of $\xxi'.$ Since $\xxi' \geq \xxi$ on $\bp \fa_\Phi + \fa_\Phi^-$, it follows that $F_\xx$ is dominated by $\xxi'$. In case (b) we still have that  $\xxi' \succeq -\rho.$ 
If the hypothesis of (c) is fulfilled this need not be the case. However, we may chose $c\in [c_0, 1[$ 
such that $\xxi' \leq i_\ga(\xxi) \leq \xxi$ on $\bp A_\Phi A_\Phi^+$ and by an argument similar to 
the previous argument, if follows that $F_\xx$ is dominated
by $i_\ga(\xxi).$ This establishes both (b) and (c) under the assumption $F^\infty_\xx = 0.$

To finish the proof, we assume that $F_\xx^\infty \neq 0.$ Then it follows from the proposition 
below that $\xx \leq -\rho(h_\ga)$ 
so that: $\xxi(h_\ga) - 1 < \xx \leq -\rho(h_\ga).$ In particular, we are in (a subcase of) case (c). 
From (\ref{e: condition on c})
it follows that 
$$
\xxi'(h_\ga) < \xx  \leq  - \rho(h_\ga) = i_\ga(\xxi)(h_\ga). 
$$ 
It now follows from 
(\ref{e: expression F xi with int J}), (\ref{e: estimate R infty}) and (\ref{e: estimate F xi infty})
that, for a suitably enlarged seminorm $n,$ 
$$ 
\|F_\xx(f, \nu, \bp aa_t)\| \leq |(\nu,\bp aa_t)|^N (\bp aa_t)^{i_\ga(\xxi)} e^{s|\Re \nu||\log \bp aa_t)|}
$$ 
for $f \in \IPgs$, $\nu \in \faPdc(\ge)$, $\bp a \in \bp A_\Phi$, $t \geq 0.$ 
Here we used that $i_\ga(\xxi) = \xxi$ on $\bp \fa_\Phi.$ 
In view of the assumed domination of $(\wh_\nu)$ by $\xxi$, we have the similar estimate
for $t \leq 0$, with $\xxi$ in place of $i_\ga(\xxi).$ Since $i_\ga(\xxi) \geq \xxi$ 
on $\bp \fa_\Phi + \fa_\Phi^-$, it follows that $F_\xx$ is dominated by $i_\ga(\xxi).$ 
This completes the proof of Proposition \ref{p: xi and domination by eta}.
\qed

\begin{prop}
\label{p: vanishing F xi infty}
Let $\xxi \in \fad$ satisfy $\xxi \succeq -\rho$ and dominate the Whittaker family $(\wh_\nu).$ 
If  $\xx \in X$ is such that $\xx > \max(\xxi(h_\ga) - 1, -\rho(h_\ga)),$  then, for $\ge>0$ sufficiently small,
$$ 
F^\infty_\xx(f, \nu , a) = 0
$$ 
for all $f \in \IPgs$, $  \nu \in \faPdc(\ge)$, and $a \in A.$  
\end{prop}
To prove the proposition, we need some preparation.
\begin{lemma}
For every left $K$-finite $f \in C^\infty(K/K_P: \gsP)$ and every $\nu \in \faPdc(P,0)$ 
there exist constants $m, C >0$ such that
\begin{equation}
\label{e: estimate w nu f}
|\wh_\nu(f)(a) | \leq C (1 + |\log {}^*\!a|)^m  a^{\Re \nu - \rho}.
\end{equation}  
\end{lemma}
\proof
By $K$-finiteness, there exists a unitary representation $(\tau, V_\tau)$  of $K$, 
a function $g \in C^\infty(\tau_P: K/K_P: \gsP)$ and a linear functional $\mu \in V_\tau^*$ 
such that $f = (I \otimes \mu)\after g.$ Define the function 
$\wh_\nu(g): G \to V_\tau$ by 
$$
\wh_\nu(x) = \inp{g}{\pi_{\bar P, \gs, \bar \nu} (x) j_{\bar \nu}}.
$$%
Then clearly, $\wh_\nu(f )= \mu \after \wh_\nu(g).$ We now note that there exists
a $\psi \in \cAtwoP$ such that $\wh_\nu(g) = \Wh(P,\psi, \nu).$ From Lemma \ref{l: first estimate Wh integral} it follows
that there exist $m, C' >0$ such that
$$  
\|\wh_\nu(g)(a) \| \leq C' (1 +| \log {}^*\!a|)^m a^{\Re \nu - \rho}.
$$% 
\noindent%
This implies (\ref{e: estimate w nu f}).
\qed

\begin{cor}
\label{c: estimate F for K-finite f}
For every left $K$-finite $f \in C^\infty(K/K_P: \gsP)$ and every $\nu \in \faPdc(P,0)$ 
there exist constants $m, C >0$ such that
\begin{equation}
\label{e: estimate F nu f}
\|F(f, \nu, a)\| \leq C (1 + |\log a|)^m  a^{\Re \nu - \rho}.
\end{equation}  
\end{cor}

\proof 
The $j$-th component of $F$ is given by 
$
F_j(f, \nu , a) = \wh_\nu( \pi_{\bar P, \gs, -\nu}(u_j) f),
$
with $u_j \in U(\fg).$ The function $\pi_{\bar P, \gs,  -\nu}(u_j) f$ 
is $K$-finite in $C^\infty(K/K_P, \gsP)$, hence $F_j$ satisfies an estimate
of the form (\ref{e: estimate F nu f}). The proof is completed by the observation
that $|\log {}^*\!a|\leq |\log a|,$ for all $a \in A,$ by orthogonality of the sum
$\fa = {}^*\!\fa_\Phi + \fa_\Phi.$ 
\qed

{\em Proof of Proposition \ref{p: vanishing F xi infty}.\ }
Put $\Omega: = \faPdc(P,0).$ 
From (\ref{e: expression F xi with int J}) it follows that 
$$ 
e^{tB(\nu)} F^\infty_\xx(f, \nu, \bp a) = F_\xx(f,\nu, \bp aa_t)  -  R^\infty_\xx(f, \nu, \bp a a_t)
$$% 
We first assume that $f\in \IPgs$ is $K$-finite. 
Let $\ge' <\ge.$ Then for 
$\nu \in \faPdc(\ge') \cap \Omega$ and $\bp a \in \bp A_\Phi$
it follows from Corollary \ref{c: estimate F for K-finite f} that there exist constants 
$m >0$ and $C > 0$ 
such that for all $t \geq 0,$ 
\begin{equation}
\label{e: estimate F xi infty 2}
\|F_\xx(f, \nu \bp a a_t)\|\leq C ( 1 + |t|)^m e^{t(\ge'|h_\ga| - \rho(h_\ga))}.
\end{equation}
From (\ref{e: defi R infty 2}) and Lemma \ref{l: estimate integral of R xi} (b)
it follows, possible after adapting $m$ and $C,$ that also 
\begin{equation}
\label{e: estimate R infty 2}
\|R^\infty_\xx(f, \nu, \bp a a_t)\| \leq C (1 + |t|)^m e^{t ( \ge'|h_\ga| + \xxi'(h_\ga))}),
\end{equation}
again for $t \geq 0.$ 
From the hypothesis on $\xx$ combined with (\ref{e: condition on c}) we see that there exists $\gd > 0$ such that
$\xx - \gd > \max(-\rho(h_\ga), \xxi'(h_\ga)).$ Keeping this in mind when combining 
the estimates (\ref{e: estimate F xi infty 2}) and (\ref{e: estimate R infty 2}) we obtain
\begin{equation}
\|e^{tB(\nu)} F^\infty_\xx(f, \nu, \bp a) \| \leq 2C(1 +|t|)^m e^{t(\ge' |h_\ga| + \xx - \gd)}
\end{equation}
as $t \geq 0.$ 
On the other hand, the expression inside the norm on the left-hand side is exponential polynomial in $t$  with exponents whose real part is at least $\xx - \ge'|h_\ga|.$ 
For $\ge'>0$ sufficiently small we have $\xx -\ge'|h_\ga| > \ge' |h_\ga| + \xx - \gd$ 
so that by uniqueness of asymptotics for $t \to \infty$ we find that 
$$ 
e^{tB(\nu)} F^\infty_\xx (f, \nu, \bp a)  = 0
$$% 
provided $\nu \in \faPdc(\ge')\cap \Omega.$ 
Since the expression on the left is holomorphic in $\nu \in \faPdc(\ge)$, 
it follows by analytic continuation that the assertion of the lemma holds for $K$-finite $f.$ 

For every $\nu \in \faPdc(\ge)$ and $\bp a \in \bp A_\Phi$ it follows from the definitions given 
that $f \mapsto F^{\infty}_\xx(f, \nu, \bp a)$ is a  linear
map $C^\infty(K/K_P:\gsP) \to \C^\ell.$ This linear map is continuous in view of 
(\ref{e: defi F infty}) and 
the estimates (\ref{e: improved estimate F xi}) and (\ref{e: estimate for I xi}) .
As it vanishes on the dense subspace of $K$-finite functions, it follows
that the given map is zero on the entire space $C^\infty(K/K_P:\gsP).$ 
This finishes the proof.
\qed

\section{Uniform temperedness of the Whittaker integral}
\label{s: uniform temp W int}
Let $P$ be a standard parabolic subgroup. We recall from (\ref{e: defi cA two G}) and (\ref{e: defi cAtwoP}) the definition of the space $\cA_2(\tau_P:M_P/M_P \cap N_0:\cchi_P).$ For $\psi $ in this space and 
for $\nu \in \faPdc(P,0) $ the Whittaker integral
$\Wh(P,\psi,\nu),$ defined by (\ref{e: defi whittaker int}), is a function in $C^\infty(\tau: G /N_0: \cchi).$

\begin{prop}
Let $\psi$ be as above. Then $\nu \mapsto \Wh(P,\psi, \nu)$ extends to a holomorphic function
$\faPdc \to C^\infty(\tau: G/N_0 : \cchi).$ 
\end{prop}

\proof As in the proof of Corollary \ref{c: initial holomorphy Wh} it suffices to prove this for $\psi = \psi_{f\otimes \xi}$ 
with $f \in C^\infty(\tau: K/K_P: \gs)$ and $\xi \in \Wh_{\cchi_P}(H_\gs^\infty).$ 
In that case the result follows from Corollary \ref{c: Wh as matrix coeff with j} combined with 
Theorem \ref{t: j is holomorphic}. 
 \qed

\begin{thm}{\bf (uniformly tempered estimate)}{\ }
\label{t: uniformly tempered for Wh}
Let $P$ be a standard parabolic subgroup. Then there exists an $\ge > 0,$ and for each $u\in U(\fg)$ 
constants $s > 0$ and $C, N > 0$ such 
that 
$$ 
\| L_u[\Wh (P, \psi, \nu)] (x) \| \leq C \|\psi\| (1 + |\nu|)^N (1 + |H(x)|)^N e^{-\rho H(x) + s|\Re \nu|} ,
$$% 
for all $\nu \in \faPdc(\ge),$ all $x \in G$ and all $\psi \in \cA_2(\tau_P:M_P/M_P\cap N_0:\cchi_P).$ 
\end{thm}

\proof
We first assume that $G = \circG.$ By finite dimensionality of $\cAtwoP,$ it suffices to prove the result
for a fixed $\psi$ of unit length. 
By linearity of the Whittaker integral in $\psi,$ and using the decomposition
(\ref{e: deco cA P in ds types}) and the isomorphism (\ref{e: iso by psi T}) we may as well assume 
in addition that $\psi = \psi_{f \otimes \xi},$ with $\xi \in \Wh_{\cchi_P}(H_\gs^\infty)$ and 
$$ 
f \in C^\infty(\tau: K/K_P: \cchi_P) \subset V_\tau \otimes  \IPgs.
$$% 
By analytic continuation we then have that 
$$ 
\Wh (P, \psi, \nu) (x) = \inp{f}{\pi_{\bar P, \gs, \bar \nu}(x) j(\bar P, \gs, \bar\nu, \bar\xi)} 
$$% 
for all $\nu \in \faPdc$ and $x \in G.$  In view of Remark \ref{r: example Whittaker family}
the map $\wh_\nu: C^\infty(K/K_P:\gs) \to C^\infty(G/N_0\col\cchi)$ defined by 
$$
 \wh_\nu(\bp f )(x)  = \inp{\bp f}{\pi_{\bar P,\gs, \bar\nu}(x) j(\bar P, \bar\nu, \bar\xi)}
$$% 
defines a holomorphic family of Whittaker functions of moderate growth.  By Cor. \ref{c: uniformly tempered estimates}, there exists an $\ge>0,$ and for each $u \in U(\fg)$ constants $N,s >0$ and a continuous seminorm $\bp\snorm$ on $C^\infty(K/K_P:\gs)$ such that 
$$ 
|L_u[\wh_\nu(\bp f)](x)|\leq (1 +|\nu|)^N (1 + |H(x)|)^N e^{s|\Re \nu ||H(x)|} e^{-\rho H(x)} \bp\snorm(\bp f)
$$% 
for all $\bp f \in \IPgs,$ $x \in G$ and $\nu \in \faPdc(\ge).$ Applying this with $\bp f = (\zeta \otimes I)f$
for $\zeta \in V_\tau^*$ we find 
$$
|\zeta \after L_u[ \Wh (P, \psi, \nu)] (x) | \leq (1 +|\nu|)^N (1 + |H(x)|)^N e^{s|\Re \nu ||H(x)|} e^{-\rho H(x)} \bp\snorm( 
(I \otimes \zeta) f). 
$$% 
Now there exists a constant $C > 0$ such that $|\bp\snorm( I \otimes \zeta) f)| \leq C |\zeta|$ 
for all $\zeta \in V_\tau^*.$ We conclude that 
$$ 
\|L_u[\Wh (P, \psi, \nu)] (x) \| \leq  C (1 +|\nu|)^N (1 + |H(x)|)^N e^{s|\Re \nu ||H(x)|} e^{-\rho H(x)},
$$% 
for all $x \in G$ and $\nu \in \faPdc(\ge)$ 
as required. This finishes the proof for the case $G = {}^\circ G$ 

In general, the group decomposes as  $G = {}^\circ G \times A_\Delta,$ 
where $A_\Delta = \exp \fa_\gD,$ with $\fa_\gD := \cap_{\ga \in \gD} \ker \ga$ central in $\fg.$ From the definitions it then readily follows that
the spaces $\cAtwoP$ for $G$ and ${}^\circ G$ coincide and that 
$$ 
\Wh(G, P,\psi, \nu)(x a_{\scriptscriptstyle \gD}) = a_{\scriptscriptstyle \gD}^\nu \, \Wh({}^\circ G, {}^\circ G\cap P, \psi, \bp \nu)(x),
$$% 
for $\psi\in \cAtwoP,$ $x \in {}^\circ G,$ $a_{\scriptscriptstyle \gD}\in A_\gD$ and $\nu \in \faPdc$
with $\bp \nu$ the restriction of $\nu$ to ${}^\circ\fg \cap \fa.$ All assertions now readily generalize
from ${}^\circ G$ to $G.$ 
\qed

\begin{cor}
\label{c: estimate diff Wh}
Let $P \in \cPst.$ Then for all $u \in U(\fg)$ and $v \in S(\faPd)$ 
there exist constants $C > 0$ and $N>0$ such that
$$ 
\|L_u [\Wh(P, \psi, \nu; v)](x)\| \leq C\|\psi\| (1+ |\nu|)^N (1 + |H(x)|)^N e^{-\rho H(x)} 
$$% 
for all $\psi \in \cAtwoP,$ $\nu \in i\faPd$ and $x \in G.$ 
\end{cor}

\begin{rem}
In the displayed equation, we have used Harish-Chandra's convention to denote the action
of a differential operator by putting it next to the variable relative to which it is applied, 
separated from the variable by a semi-colon.
 In the present context, if $v\in S(\faPd),$ then $v$ 
is viewed as a constant coefficient complex differential operator on $\faPdc,$ and 
if $\gf: \faPdc \to V$ is a holomorphic function with values in a locally convex space then 
$\gf(\nu ; v)$ stands for $v \gf$ at the point $\nu \in \faPdc.$ 
\end{rem} 

\proof 
This follows from the estimates of Theorem \ref{t: uniformly tempered for Wh} by using the Cauchy integral 
formula in the variable $\nu,$ with polydiscs of polyradius $\ge (\dim\faP)^{-1} (1 + |H(x)|)^{-1}.$ 
\qed

We may now define a Fourier transform $\cF_P: \cC(G/N_0\col \cchi) \to C^0(i\faPd, \cAtwoP)$ 
by the formula 
$$ 
\inp{\cF_P(f) (\nu)}{\psi} = \inp{f}{\Wh(P,\psi, \nu)}_2 := \int_{G/N_0} f(x) \overline{\Wh(P, \psi, \nu)(x)}\; d \dot{x},
$$%
for $f \in \cC(G/N_0\col\cchi),$ $\nu \in i\faPd$ and $\psi \in \cAtwoP.$ 
Indeed, let $\ell = \dim A + 1$ then by Lemma \ref{l: convergence integral over G/N_0}
the function
$x \mapsto (1+ |H(x)|)^{-\ell} e^{-2\rho H(x)}$ is absolutely integrable over $G/N_0.$ 
By application of Theorem \ref{t: uniformly tempered for Wh}, we infer the existence of a constant $N > 0$ 
and a continuous seminorm 
$\snorm$ on $\cC(\tau: G/N_0\col \cchi)$ such that for all $f \in \cC(\tau: G/N_0\col\cchi) $ we have
\begin{equation}
\label{e: first estimate Fourier}
|f(x) \overline{\Wh(P, \psi, \nu)(x)}| \leq  (1 +|\nu|)^N (1+ |H(x)|)^{-\ell} e^{-2\rho H(x)} \snorm(f),\qquad (x \in G).
\end{equation} 
It follows from this that the Fourier transform is defined by an absolutely converging integral, and defines a continuous linear operator $\cC(G/N_0\col\cchi) \to C^0(i\faPd, \cAtwoP).$ By application of Cor. \ref{c: estimate diff Wh} it follows that differentiation
under the integral is allowed, and that $\cF_P$ is continuous linear
$\cC(G/N_0\col\cchi) \to C^\infty(i \faPd, \cAtwoP).$ 

\begin{lemma}
\label{l: Fourier Z f}
Let $P \in \cPst.$ Then for all $f \in \cC(\tau : G/N_0\col \cchi)$ 
and every $Z \in\fZ$ we have
$$ 
\cF_P (R_Z f)(\nu)  = \umu_P(Z, \nu ) \cF_Pf(\nu), \qquad (\nu \in i\faPd).
$$% 
\end{lemma}

\proof Let $\psi \in \cAtwoP.$ 
Since $\umu(Z, \dotvar)$ is polynomial with values in $\End(\cAtwoP)$ 
it follows from () by analytic continuation that 
$$
R_Z \Wh(P, \psi, \nu) = \Wh(P, \umu_P(Z, \nu) \psi, \nu)
$$% 
for all $\nu \in \faPdc.$  Hence, by differentiation under the integral sign,
\begin{eqnarray*}
\inp{\cF_P( R_Z f)(\nu)}{ \psi} &=& \inp{L_{Z^\veec} f}{\Wh(P, \psi, \nu)}_2\\
&=& 
\inp{f}{R_{\bar Z^\veec} \Wh(P, \psi, \nu) }_2\\
&=& 
\inp{f}{\Wh(P, \umu_P(\bar Z^\veec, \nu)  \psi, \nu)}_2\\
&=& 
\inp{\cF_P f}{ \umu_P(\bar Z^\veec, \nu)  \psi}.
\end{eqnarray*}
As this holds for arbitrary $\psi \in \cAtwoP$ we conclude that 
$$ 
\cF_P (R_Z f)(\nu) = \umu_P(\bar Z^\veec, \nu)^*\; \cF_P f(\nu) ,\quad (\nu \in i\faPd).
$$% 
where the star indicates that the adjoint is taken with  respect to the $L^2$-Hilbert structure on $\cAtwoP.$ By a straightforward calculation it follows that
$$ 
\umu_P(\bar Z^\veec, \nu)^* = \overline{\umu_P(\bar Z^\veec, \nu)^\veec} =
\umu_P(Z, -\bar \nu ) = \umu_P(Z, \nu )
$$% 
for all $\nu \in i \faPd.$ 
\qed

Let $\cS(i\faPd)$ denote the usual space of Schwartz functions on the real vector space $i\faPd.$ 
Then the following result is valid.

\begin{thm}
\label{t: Fourier continuous into Schwartz}
$\cF_P$ maps  $\cC(\tau \!:\!G/N_0\col \cchi)$ continuous linearly to 
$\cS(i\faPd)\otimes \cAtwoP.$ 
\end{thm}

The proof follows the usual strategy of applying partial integration,
involving minus the Casimir operator associated with the invariant bilinear form 
$B,$ see (\ref{e: intro of B}).
The following lemma prepares for this.
\begin{lemma}
\label{l: limit umu L}
 Let $L \in \fZ$ be minus the Casimir operator. Then 
$$ 
\lim_{\topping{\nu \in i\faPd}{ |\nu| \to \infty}} (1 + |\nu|)^{-2} \, \umu_P(L, \nu) = I
$$% 
in $\End(\cAtwoP).$ 
\end{lemma}

\proof 
By finite dimensionality of $\cAtwoP$ it suffices to prove the
identity for the restriction of the endomorphisms to the subspace 
$\cA_{2,P, \gs} :=\cA_\gs(\tau_P\!:\! K/K_P\col \cchi_P),$ 
with $\gs$ a representation of  the discrete series of $M_P.$ 
Let $\gL_\gs \in \fhdc$ 
be the infinitesimal character of $\gs$ and let $\gd$ be half of the sum of a choice
of positive roots for the root system of $\fh$ in $\fg_\iC.$ 

The restriction of $\umu_P(L, \mu)$ to $\cA_{2, P, \gs}$ 
equals the restriction of $R_{\umu_P(L, \mu)} = L_{\umu_P(L, \mu)^\veec},$
which is given by multiplication by the scalar 
\begin{eqnarray*}
\gg_{M_P} (\umu_P(L, \nu)^\veec, \Lambda_\gs) & = & 
\gg_{M_P} (\umu_P(L, \nu), - \Lambda_\gs) \\
&= &\gg(L, \nu - \gL_\gs)=  |\nu|^2 + C_\gs
\end{eqnarray*} 
with $C_\gs = - B^*(\gL_\gs , \gL_\gs) + B^*(\gd, \gd) \in \R,$ where $B^*$ is the dual 
of $B.$ 
Accordingly, the restriction of the limit equals
$$ 
\lim_{|\nu| \to \infty} (1 + |\nu|)^{-2} ( |\nu|^2 +  C_\gs) I = I.
$$%
\qed
\medbreak
{\em Proof of Theorem \ref{t: Fourier continuous into Schwartz}.\ }
In the above we already showed that $\cF_P$ maps $\cC(\tau: G/N_0\col \cchi)$ 
continuously to $C^\infty(i\faPd) \otimes \cA_{2,P}.$ By Lemma \ref{l: limit umu L}
there exists
a constant $R > 0$ such that for all $\nu \in i\faPd$ with $|\nu| \geq R$ 
the endomorphism $\umu_P(L, \nu)$ of $\cAtwoP$ is invertible,
whereas the operator norm of its inverse satisfies 
$$ 
\| \umu_P(L, \nu)^{-1} \| \leq 2(1 + |\nu|)^{-2}. 
$$%
We will finish the proof by showing that for every $u \in U(\faPd)$ and all
$k \in \N$ there exists a continuous seminorm $\bf n$ on $\cC(\tau: G/N_0\col \cchi)$ 
such that for all $f \in \cC(\tau: G/N_0\col \cchi)$ we have 
$$
\|\cF_Pf(\nu; u)\| \leq (1 + |\nu|)^{-k}  {\bf n}(f), \qquad (|\nu| \geq R).
$$% 
For this we proceed by induction on the order of $u.$ Clearly, the result is true
for $u$ of order $-1$ since then $u= 0.$ Thus, assume that $u$ 
has order $m \geq 0$ and assume the result has been established for
$u$ of order strictly smaller than $m.$ 

It follows from (\ref{e: first estimate Fourier})
that for all $f \in \cC(\tau: G/N_0\col \cchi)$ and $\nu \in i\faPd$ we have 
\begin{equation}
\label{e: initial estimate cFP f}
\|\cF_Pf(\nu;u )\| \leq (1 + |\nu|)^{N}  {\bf n}(f).
\end{equation} 
Here $N \in \N$ and ${\bf n}$ is a continuous seminorm on $\cC(\tau: G/N_0\col \cchi).$ 

Fix $\ell\in N$ such that $N - 2\ell \leq -  k.$ Then it follows by application of Lemma \ref{l: Fourier Z f} and the Leibniz rule that there exist a finite collection of polynomial functions $q_j \in P(\faPd)$ and 
a finite collection of elements $u_j \in S(\faPd)$ of order strictly smaller than 
the order of $u$ such that 
$$
\umu_P(L, \nu)^\ell \cF_P f(\nu; u) = \cF_P( L^\ell f)(\nu; u) 
+\sum_{j=1}^r q_j(\nu) \cF_P f(\nu; u_j) ,
$$% 
for all $f \in \cC(\tau: G/N_0\col \cchi)$ and all $\nu \in i\faPd.$ 
By application of the initial estimate (\ref{e: initial estimate cFP f} ) and the inductive hypothesis there exists a continuous seminorm $\bf n'$ on $\cC(\tau: G/N_0\col \cchi)$ 
such that for all $f$ and all $|\nu|\geq R,$ we have
$$ 
\|\umu_P(L,\nu)^\ell \cF_Pf(\nu; u)\| \leq (1 + |\nu|)^N {\bf n'}(f).
$$% 
This implies that
\begin{eqnarray*}
\|\cF_Pf(\nu; u)\| & \leq & \| \umu_P(L, \nu)^{-1}\|^\ell \|\umu_P(L,\nu)^\ell \cF_Pf(\nu; u)\| 
\\
&\leq & 2^\ell (1 + |\nu|)^{N-2\ell} {\bf n'} ( f)\\
&\leq & (1 + |\nu|)^{-k} {\bf n'} (2^\ell  f).
\end{eqnarray*} 
completing the induction.
\qed

\section{Appendix: factorization of polynomial functions}
In this section we will prove the following result, which is needed in Section \ref{s: functional equation}.
\begin{prop}
\label{p: factorization in first order}
Let $f: \C^n \to \C$ be a polynomial function of degree $d \geq 1$ and assume that $\cH$ is a locally finite
collection of affine {\bfhyp} in $\C^n$ such that $f^{-1}(0) \subset \cup \cH.$ Then $f$ can be expressed 
as a finite product $f =\ell_1 \cdots \ell_d$  with  $\ell_j: \C^n \to \C$ a
linear polynomial function whose zero set $\ell_j^{-1}(0)$ belongs to $\cH,$ for $1 \leq j \leq d.$ 
\end{prop}

The following lemma is a first step in the proof.

\begin{lemma}
\label{l: f is zero on hyperplane}
Let $f , \ell : \C^n \to \C$ be non-zero polynomial functions,  with $\deg \ell = 1.$ 
If $f$ vanishes on the hyperplane $\ell^{-1}(0)$ then $f/\ell$ is polynomial.
\end{lemma}

\proof 
By application of a suitable 
affine coordinate transformation we may reduce to the case $\ell (z) = z_n.$ In view of  the hypothesis, $f(z' , 0) = 0$ for all $z' \in \C^{n-1}.$ This implies that all partial derivatives
$\partial^\ga f(0),$ with $\ga \in \N^{n}$ and $\ga_n =0$ are zero. Hence, 
$$ 
f(z) = \sum_{\ga\in F} c_\ga z^\ga
$$% 
with $F \subset \N^{n-1} \times \N_+.$  The result now follows.
\qed

In the following we denote by $\cO = \cO_{\C^n}$ the sheaf of holomorphic functions
on $\C^n.$ For a point $a\in \C^n$ we denote by $\cO_a$ the stalk at $a,$ i.e., the ring of germs
at $a$ of locallly defined holomorphic functions. If no confusion is possible, we will switch 
between elements of $\cO_a$ and local representatives for them without explicitly mentioning this.

\begin{lemma}
\label{l: linear factor lemma}
Let $f \in \cO_0$ and let $\Xi$ be a finite collection of non-zero linear functionals $\xi: \C^n \to \C.$ 
Suppose that $f(0) = 0$ and that $f^{-1}(0) \subset \cup_{\xi \in \Xi}\ker \xi$ (in the sense of germs).
Then there exists a functional $\xi \in \Xi$ such that $f = 0$ on $\ker \xi$ (in the sense of germs). 
\end{lemma} 

\proof 
By a suitable linear change of coordinates we may reduce to the case that 
$f$ and each $\xi \in \Phi$ does not vanish identically on the coordinate axis $\C e_n.$ 
Then by the Weierstrass preparation theorem, there exists an invertible element $q \in \cO_0$ 
such that $F= q^{-1}f$ is a Weierstrass polynomial given by 
\begin{equation}
\label{e: f in z n}
F(z', z_n) = z_n^N + \sum_{k < N} c_k(z') z_n^k , 
\end{equation} 
for $z = (z', z_n) \in \C^{n-1}\times \C$ sufficiently close to $(0,0),$ 
with $c_k \in \cO_0' = \cO_0(\C^{n-1})$ such that $c_k(0) = 0$ for $0\leq k < N.$ 
Furthermore, for every $\xi \in \Phi$ there exists a linear functional
$\eta_\xi: \C^{n-1} \to \C$ such that $\ker \xi \subset \C^{n-1}\times \C$ equals 
the graph of $\eta_\xi.$

Let $D\subset \C^n$ be a polydisk centered at $0$ such that $f, q$ and $F$ admit representatives
in $\cO(D).$ We decompose $D = D' \times D_n$ according to the decomposition $\C^n = \C^{n-1} 
\times \C.$ Then $f^{-1}(0) \cap D = F^{-1}(0)\cap D.$ 
Let $\Omega$ be the open dense subset of $\C^{n-1}$ consisting of $z' \in \C^{n-1}$ 
such that $\eta_{\xi_1}(z') \neq \eta_{\xi_2}(z')$ for all distinct $\xi_1, \xi_2 \in \Xi.$ 

By continuous dependence of the roots of (\ref{e: f in z n}) on the coefficients $c_k(z')$ we may shrink $D'$ sufficiently so that for all $z' \in D'$ there exists an element 
$z_n \in D_n$ such that $F(z', z_n) = 0.$ Accordingly, we may select $\ga' \in \Omega\cap  D'$ and 
$\ga_n \in D_n$ such that $F(\ga', \ga_n) = 0.$ Then $f(\ga', \ga_n) = 0$ so there exists a $\xi_0 \in \Xi$ such that
$\xi_0(\ga', \ga_n) = 0 ,$ or, equivalently, 
$$ 
\ga_n = \eta_{\xi_0}(\ga'). 
$$% 
Since $\ga' \in \Omega,$ there exists $\ge >0$ such that for $z' = \ga'$ we have 
\begin{equation}
\label{e: estimate eta H and ga n}
|\eta_\xi(z') - \ga_n| > \ge \qquad (\forall \xi \in \Xi\setminus\{\xi_0\}). 
\end{equation}
By continuity there exists an open neighborhood $U$ of $\ga' $ in $D' \cap \Omega $ 
such that the estimates (\ref{e: estimate eta H and ga n})  are still valid for all $z' \in U.$ 

By the continuity of roots mentioned above, we may shrink
$U$ so that in addition for every $z' \in U$ there exists a $z_n \in D(\ga_n, \ge)$ 
such that $f(z', z_n) = 0.$ The latter implies that $z_n  =  \eta_\xi(z')$ for a $\xi \in\Xi.$ 
Now $|\eta_\xi(z') - \ga_n| < \ge $ implies that $\xi = \xi_0$ and we see that 
\begin{equation} 
\label{e: vanishing f z pr eta}
F(z', \eta_{\xi_0} (z')) = 0
\end{equation} 
for all $z' \in U.$ By analytic continuation, (\ref{e: vanishing f z pr eta}) is valid for all 
$z \in D'.$ Therefore, $f = 0$ on $\{z \in D \mid z_n = \eta_{\xi_0}(z') \} = \ker \xi_0.$ 
\qed

\begin{cor}
\label{c: linear factor cor}
Let $f: \C^n \to \C$ be a polynomial function of positive degree 
and assume that $\cH$ is a locally finite collection of affine {\bfhyp} in $\C^n$ such that $f^{-1}(0) \subset \cup \cH.$ Then there exists a hyperplane $H \in \cH$ such that $f$ vanishes on $H.$ 
\end{cor} 

\proof 
Since $f$ has positive degree, $f^{-1}(0) \neq 
\emptyset.$ By application of a suitable translation, we may reduce
to the case that $f(0) = 0.$ Let $\cH_0$ be the finite collection of $H \in \cH$ with $H \ni 0.$ 
For each $H \in \cH_0$ we fix $\xi_H \in (\C^n)^*$ such that $H = \ker \xi_H.$ Put 
$\Xi =\{\xi_H \mid H \in \cH_0\}.$

By application of Lemma \ref{l: linear factor lemma} it follows that there exists a $\xi_0\in \Xi$ 
and a polydisk $D \subset \C^n$ centered at $0$ such that $f = 0$ on $\ker (\xi_0)\cap D.$ 
 By analytic continuation of $f|_{\ker \xi_0}$ this implies that $f$ vanishes
on $\ker\xi_0.$ 
\qed
\medno
{\em Proof of Proposition \ref{p: factorization in first order}.\ }   In view of Corollary \ref{c: linear factor cor}
and Lemma \ref{l: f is zero on hyperplane} the proof follows by a straightforward induction 
on $d.$
\qed

\section{Appendix: a Hartog type continuation result}
\label{s: Hartog thm}
In this paper we will need the following continuation result for holomorphic functions
on a domain in $\C^n$ with values in a quasi-complete locally convex space. 

Let $\Omega \subset \C^n$ be a connected
open subset, $p: \Omega \to \C$ a non-zero holomorphic function, and $X \subset \Omega$ 
its zero locus $p^{-1}(0).$  We denote by $X_r$ the set of points $z \in X$ such 
that $X$ is a smooth complex hypersurface at $z.$ By this we mean that there should exist an open neighborhood $U$ of $z$ in $\Omega$ such that $X\cap U$ is a complex submanifold of dimension $n-1.$ The complement of $X_r$ in $X$ is denoted by $X_s.$ Clearly, $X_s$ is closed in $X$ hence in 
$\Omega.$ Since $\Omega$ is connected, the set $X,$ hence also $X_s,$ has empty interior.  

\begin{thm}
\label{t: Hartog type continuation}
Let $V$ be a quasi-complete locally convex (Hausdorff) space. Then every holomorphic function $f: \Omega \setminus X_s \to V$ admits a unique extension to a holomorphic function $\Omega \to V.$ 
\end{thm}

\begin{rem}
For $V = \C$ the result is well known and can be obtained as a consequence
of \cite[Thm.~6.12]{FGhol}, which asserts that the above result is valid for $V =\C$ and 
with $X_s$ replaced by any analytic subset $Y$ of $\Omega$ which is 
everywhere locally of codimension at least two. That result actually extends to $V$-valued holomorphic functions, but we have not found a decent reference to the literature for this.  
We have therefore chosen to present a self-contained proof of Theorem \ref{t: Hartog type continuation}, following a strategy suggested in Exercise 4.26 in the set of lecture notes
{\em Several Complex Variables}, by Jaap Korevaar and Jan Wiegerinck, version August 23, 2017.
This turned out to be possible since Cauchy's integral formula is valid in the setting of $V$-valued holomorphic functions.
\end{rem}

To prepare for the proof of Theorem \ref{t: Hartog type continuation} we make the following general observation. We assume that $V$ is a quasi-complete locally convex space. 

\begin{lemma}
\label{l: observation about extension}
Let $\Omega \subset \C^n$ be an open subset and let $Y$ be a closed subset of $\Omega$ 
which has empty interior. Let $f: \Omega\setminus Y \to V$ be holomorphic. Then the following
assertions are equivalent.
\begin{enumerate}
\itema $f$ extends to a holomorphic function $\Omega \to V;$
\itemb for every $y \in Y$ there exists an open neighborhood $\omega \ni y$ in 
$\Omega$ such that $f|_{\omega \setminus Y}$
has a holomorphic extension to $\omega.$ 
\end{enumerate}
\end{lemma}

\proof 
That (a) implies (b) is obvious. Assume (b). Then one 
may cover $\Omega$ with open subsets 
$\omega_j$ for $j$ in an index set $I,$ such that for each $j\in I$ the function 
$f|_{\omega_j \setminus Y}$ has a holomorphic extension $f_j : \omega_j \to V.$ Clearly, if 
$i,j \in I$ and $\omega_i \cap \omega_j \neq \emptyset,$
then $f_i = f_j$ on $(\omega_i \cap \omega_j)\setminus Y.$ By density this 
implies that $f_i = f_j $ on $\omega_i \cap \omega_j.$ From this (a) readily follows.
\qed

\medno
{\em Proof of Theorem \ref{t: Hartog type continuation}.\ } 
By Lemma  \ref{l: observation about extension} applied with $Y = X_s$ it suffices to show that 
 for every $z \in X_s$ there exists an open neighborhood $\omega \ni z$ in $\Omega$ 
 such that $f|_{\omega \setminus X_s}$ extends to a holomorphic function $\omega \to V.$ 

Let $z^0 \in X_s.$ Then 
we may apply an affine coordinate transformation
to arrange that $z^0 = 0$ and that locally at $0$ the function $p$ is $z_n$-regular,
see \cite[p.~109]{FGhol}. By the Weierstrass preparation theorem, see \cite[Thm.~III.2.7]{FGhol}, locally at $0$ the function $p$ 
factors as a product of holomorphic functions $p_0 \cdot W,$ with the germ at $0$ of $p_0$ being a unit in the
ring $\cO_0 = \cO_0(\C^n)$  of germs of holomorphic functions defined locally at $0$ in $\C^n$ 
and with $W \in \cO_0'[z_n]$ 
a Weierstrass polynomial of order $d$ over the ring $\cO_0' = \cO_0(\C^{n-1}).$ Then for 
$D$ a sufficiently small polydisk in $\C^n$ centered at $0$ the germ $W$ has a representative which is holomorphic
on $D$ and such that $D \cap X$  is contained in $W^{-1}(0).$ Note that $0\in D\cap X.$ 

The ring $\cO_0$ is a unique factorization domain, see \cite[Thm.~III.3.3]{FGhol}.
Let $W=  f_1^{m_1} \cdots f_r^{m_r}$ be a decomposition into 
irreducibles of $\cO_0$ with the $f_j$ mutually prime. 
Since $W$ is $z_n$-regular, each $f_j$ is $z_n$-regular as well. By the Weierstrass
preparation theorem we may write $f_j = \ge_j W_j$ with $\ge_j$ a unit in $\cO_0$ 
and $W_j$ a Weierstrass polynomial in $\cO_0'[z_n].$ 
Then $W = \ge W_1^{m_1}\cdots W_r^{m_r} $ 
with $\ge = \ge_1^{m_1}\cdots \ge_r^{m_r}$ a unit in $\cO_0$. Clearly the product
$ W_1^{m_1}\cdots W_r^{m_r} $ is a Weierstrass polynomial. By the uniqueness
statement of the Weierstrass preparation theorem, it now follows that
$W = W_1^{m_1} \cdots W_r^{m_r}.$ The $W_j$ are mutually distinct and irreducible 
in $\cO_0$ hence in $\cO_0'[z_n],$ see \cite[Lemma II.5]{GRcomp}.

Let $\bp W = W_1 \cdots W_r;$ then after sufficiently shrinking
$D,$ keeping it centered at $0,$ we obtain that $\bp W$ has a holomorphic representative on $D$ and that $D \cap X$ equals the zero locus of $\bp W$ in $D.$ 
We agree to write $D = D' \times D_n$ in accordance with the decomposition
$\C^n = \C^{n-1} \times \C.$ 

If $z \in D$ satisfies $\bp W(z) = 0 $ and $\partial_n \bp W (z) \neq 0,$ then it follows by application
of the submersion theorem that  $\bp W^{-1}(0)$ is a complex differentiable submanifold of 
codimension $1$ locally at $z.$ This implies that $z \notin X_s.$ It follows that  $X_s \cap D$ is contained in the zero locus of both $\bp W$ and $\partial_n \bp W.$ 

Let now $z \in X_s\cap D.$ Then it follows that the polynomial functions $\bp W(z', \dotvar)$ and $\partial_n W(z', \dotvar)$ have $z_n\in D_n$ as a common zero. Hence, $z_n$ 
is a root of higher multiplicity of $\bp W(z', \dotvar)$ and it follows that 
the discriminant $\gD(z')$ of $\bp W(z', \dotvar)$ is zero. Since $\gD(z')$ is 
a polynomial expression in the coefficients of $\bp W(z', \dotvar),$ it follows that 
$\gD \in \cO(D').$ We view $\gD$ as a polynomial function on $\C$ 
with coefficients in $\cO(D')$ and conclude that 

\begin{equation}
\label{e: X as zero set of W gD}
D \cap X_s \subset D \cap \bp  W^{-1}(0)  \cap \gD^{-1}(0).
\end{equation}

We will now establish the claim that $\gD$ does not vanish identically on $D.$
Let $Q$ be the quotient field of $\cO'_0 := \cO_0(\C^{n-1}).$ Then
it follows from \cite[Cor.~3.2 (2)]{FGhol} with $I = \cO_0'$ that $W_1,\ldots , W_r$ are irreducible in 
  the polynomial ring $Q[z_n].$ Furthermore, since $Q\setminus\{0\}$ is the set of 
units in $Q[z_n]$, no distinct $W_i$ and $W_j$ are related by a unit factor in $Q[z_n].$ 
Since $Q[z_n]$ is a unique factorization domain, it follows that the factors $W_j$ are prime.
Furthermore, the elements $\bp W$ and $\partial_n \bp W$ have a greatest common divisor $\gg \in Q[z_n],$
which, up to a unit factor, may be written as the product of those factors $W_j$ that 
divide $\partial_n \bp W.$ By application of Leibniz's rule for differentiation one sees that 
such a factor $W_j$ must divide $W_1 \cdots W_{j-1} \cdot \partial_n W_j \cdot W_{j+1} \cdots W_r,$ 
hence $W_j$ must divide $\partial_n W_j$ which is impossible since the latter has lower degree
than the former. We conclude that $\gg$ is a unit hence belongs to $Q\setminus\{0\}.$ 

By the Euclidean division algorithm
there exist $\gl_1, \gl_2 \in Q[z_n]$ such that
$ 
1 = \gl_1\bp W + \gl_2 \,\partial_n \bp W.
$ 
Let $a \in \cO_0'$ be a non-zero element such that ${}^\circ \!\gl_j:= a \gl_j$ belong to
$\cO_0'[z_n],$  for $j =1,2.$ Then
\begin{equation}
\label{e: W and a}
a = {}^\circ\! \gl_1 \bp W + {}^\circ\! \gl_2 \,\partial_n \bp W.
\end{equation}
Shrinking the polydisk $D$ sufficiently we may arrange that this equation is valid 
for all $z \in D.$ Shrinking the polydisk $D'$ sufficiently, we may also arrange that for every $z' \in D'$ the polynomial functions $\bp W(z', \dotvar)$ and $\partial_n \bp W(z', \dotvar)$ 
have all their roots contained in $D_n.$ 
Suppose now that $z' \in D$ and $\gD(z') = 0.$ 
Then the polynomial functions $\bp W(z', \dotvar)$ and 
$\bp \partial_n \bp W(z' , \dotvar)$ have a common root $\gz$, which must be contained in $D_n.$
Evaluating (\ref{e: W and a}) in $(z', \gz)$ we find that $a(z') = 0.$ We infer that $D' \cap \gD^{-1}(0)
\subset D' \cap \gD^{-1}(0).$ Since $a$ is non-zero, we infer
that $\gD$ is not vanishing on all of $D'$; the validity of the claim follows.
  
Since the germ of $\gD$ at $0$ is non-zero we may apply 
a linear transformation in the first $n-1$ coordinates of $\C^n$ to
 arrange that $\gD = \gd W'$ in $\cO'_0 := \cO_0(\C^{n-1}),$ with $\gd$ a unit in $\cO'_0$ and with $W' \in \cO''_0[z_{n-1}]$ a Weierstrass polynomial of degree $d'$; here $\cO''_0 =\cO_0(\C^{n-2}).$ Shrinking $D'$ to a polydisk centered at $0$ with 
respect to the new coordinates, we may arrange to be in the situation that 
(\ref{e: X as zero set of W gD}) is still true,
and such that $W'$ can be represented by a holomorphic function on $D'\,$ whereas
$\gd$ can be represented by a nowhere vanishing function on $D'.$ Accordingly,
$D \cap \gD^{-1}(0) = D \cap (W')^{-1}(0).$ 

We fix a sufficiently small $r >0$ such that $\partial D(0, r) \subset D_n.$ 
Furthermore, we shrink the polydisk $D',$ keeping it centered at $0$, such that $\bp W$ is nowhere zero on an open neighborhood $D' \times A$ of $D' \times \partial D(0, r)$ in $D.$ 
It follows that $(D'\times A)  \cap X_s = \emptyset,$ so that $f$ is well-defined and holomorphic on 
$D'\times A.$ We write $D'= D'' \times D_{n-1}$ according to the decomposition
$\C^{n-1} = \C^{n-2} \times \C$ and fix $r'> 0$ such that $\partial D(0,r') \subset D_{n-1}.$ 
Then  it follows by application of Cauchy's integral formula to the $(n-1)$-th coordinate 
that for $(z'', z_{n-1}, z_n) \in D'' \times D(0, r') \times A,$
\begin{equation}
\label{e: first application of Cauchy}
f(z'', z_{n-1}, z_n) = \frac{1}{2\pi i} \int_{\partial D(0, r') } \frac{f(z'', \zeta_{n-1}, z_n)}{\zeta_{n-1} - z_{n-1}}\; d\zeta_{n-1}.
\end{equation}
There exists a sufficiently small polydisk $D_0'' \subset D''$, centered at $0,$ such that 
$W'$ does not vanish on  an open neighborhood $D_0'' \times B$ of 
$D_0'' \times \partial D(0, r').$  Since $W'$ is constant as a function of the 
$n$-th variable, it follows that $X_s$ has empty intersection with $D_0'' \times B \times D_n,$ 
so that $f$ is holomorphic on the latter set as well. In particular, it follows that the integrand
of (\ref{e: first application of Cauchy}) is a holomorphic function of $z_n \in D_n$ as long 
as $z'' \in D_0''$ and $\zeta_{n-1} \in \partial D(0, r').$ Applying Cauchy's integral formula
to this holomorphic function of $z_n,$ we obtain
$$ 
f(z'', z_{n-1}, z_n) = \left(\frac{1}{2\pi i}\right)^2\int_{\partial D(0,r') \times \partial D(0,r)}
 \frac{f(z'', \zeta_{n-1}, \zeta_n)}{(\zeta_{n-1} - z_{n-1})(\zeta_n - z_n)}\; d\zeta_n d\zeta_{n-1},
$$% 
for $z \in D_0'' \times D(0,r')  \times A.$
On the other hand, since $f$ is holomorphic on $D_0''\times B \times A,$ 
the integral on the right defines a $V$-valued holomorphic  function $F(z)$ 
of $z \in \bp D : = D_0'' \times D(0, r') \times D(0,r).$ 
It follows from the last displayed equality that $F = f$ on the non-empty open subset 
$\cO:= D_0'' \times D(0, r') \times (A\cap D(0; r))$ of $\bp D.$ 
By analytic continuation it follows that $F = f$ on the connected open set 
$\bp D\setminus X_s.$ Thus, $\bp D$ 
is an open neigborhood of $z^0$ such that 
$f|_{\bp D\setminus X_s}$ extends to the holomorphic function $F: \bp D \to V.$  
\qed

%Bibliography ------------------------------------------------------------
%\bibliographystyle{plain}
%\bibliography{unitemp.bib}

%Names -------------------------------------------------------------------
\def\adritem#1{\hbox{\small #1}}
\def\distance{\hbox{\hspace{3.5cm}}}
\def\apetail{@}
\def\addVdBan{\vbox{
\adritem{E.~P.~van den Ban}
\adritem{Mathematical Institute}
\adritem{Utrecht University}
\adritem{PO Box 80 010}
\adritem{3508 TA Utrecht}
\adritem{The Netherlands}
\adritem{E-mail: E.P.vandenBan{\apetail}uu.nl}
}

}
\mbox{}
\vfill
\hbox{\vbox{\addVdBan}\vbox{\distance}}
\end{document}